\documentclass[reqno, 11pt]{amsart}  

\usepackage{mathtools} 
\usepackage{amssymb}
\usepackage[nonewpage]{imakeidx}
\makeindex[columns=3, title=Index]

\makeatletter
\def\ckech{\mathaccent"\accentclass@014}
\def\tah{\mathaccent"\accentclass@05E}
\makeatother

  \renewcommand\check{\bm\ckech} 
\renewcommand\hat{\bm\tah}

\newcommand\mysection[1]{% 
\section{#1}\setcounter{equation}{0}}

\usepackage{bm}

\newtheorem{theorem}{Theorem}[section]
\newtheorem{lemma}[theorem]{Lemma}
\newtheorem{corollary}[theorem]{Corollary}

\theoremstyle{definition}

\newtheorem{definition}[theorem]{Definition}

\theoremstyle{remark}
\newtheorem{example}[theorem]{Example}
\newtheorem{remark}[theorem]{Remark}
 
\newcommand\loc{\textnormal{loc}}
\newcommand{\shharp}{=\!\!\!\!\|}  
\newcommand{\vsharp}{\asymp\kern -.5em\|}

 \makeatletter 
 \def\dashint{\operatorname{\,\,\,\mathclap{\!\int}\! 
 \!\text{\bf--}\!\!}}  
 \makeatother

\def\sfG{{\sf G}}

\newcommand\sfp{{\sf p}}

\newcommand\bB{\mathbb{B}}
\newcommand\bE{\mathbb{E}}
\newcommand\bC{\mathbb{C}}
 
\newcommand\bM{\mathbb{M}}
\newcommand\bR{\mathbb{R}}

\newcommand\bS{\mathbb{S}}

\newcommand\bZ{\mathbb{Z}}

\newcommand\frB{\mathfrak{B}}

\newcommand\frQ{\mathfrak{Q}}

\newcommand\cB{\mathcal{B}}

\newcommand\cE{\mathcal{E}}
\newcommand\cF{\mathcal{F}}

\newcommand\cL{\mathcal{L}}
\newcommand\cM{\mathcal{M}}

\newcommand\cO{\mathcal{O}}

\newcommand\cQ{\mathcal{Q}}

\newcommand\E{{\sf{E}}}
\renewcommand\L{{\sf{L}}}  {{}}
\newcommand\W{{\sf{W}}}

\renewcommand\){{\rm)}}

\def\+){\tmspace+\thinmuskip{.05em}\)}

\def\dashnorm{\,\,\text{\bf--}\kern-.5em\|}

\newcommand{\nlimsup}{\operatornamewithlimits{\overline{lim}}}
\newcommand{\nliminf}{\operatornamewithlimits{\underline{lim}}}

\newcommand{\osc}{\operatornamewithlimits{osc}}

\renewcommand{\eqref}[1]{\text{(\ref{#1})}}

\begin{document}

\title[]{Essentials of Real Analysis and Morrey-Sobolev spaces for
second-order elliptic and parabolic PDEs with singular first-order
coefficients}
\author[]{N.V. Krylov}
\address{School of Mathematics, University of Minnesota, Minneapolis, MN, 55455}
\email{nkrylov@umn.edu}

\renewcommand{\subjclassname}{\textup{2010} 
Mathematics Subject Classification}
 


\begin{abstract}
In recent years we witness growing interest
in using   Real Analysis methods and results
in the theory of nondivergence form partial differential equations (PDEs) 
and the goal of this article is to give a brief
and concise introduction into  the
applications of several results in Real Analysis to the theory of
elliptic and parabolic equations in Sobolev and Morrey-Sobolev spaces. 
In particular, we concentrate on such results as Hardy-Littlewood maximal function
theorem, Fefferman-Stein theorem,
theory of Muckenhoupt weights, and Rubio de Francia
extrapolation theorem and their role in Sobolev or Morrey-Sobolev space
theory of parabolic equations with mixed norms.

In our exposition we do not try to give the strongest
known results for particular equations in particular spaces.
We only show how the Real Analysis results, we present
with all proofs, can be used in model cases such as the Laplace
and the heat equations with singular first order terms. The only exception is the last
section where we present    new results.
 
\end{abstract}

\maketitle

In recent years we witness growing interest
in using   Real Analysis methods and results
in the theory of  partial differential equations (PDEs) and the goal of this article is to give a brief
and concise introduction into  the
applications of several results in Real Analysis to the theory of nondivergence form 
elliptic and parabolic equations in Sobolev and Morrey-Sobolev spaces. 
The origin of some of these results
is rather old, such as of  the  Hardy-Littlewood maximal function
theorem,  the  Fefferman-Stein theorem,  the 
theory of Muckenhoupt weights, and  the  Rubio de Francia
extrapolation theorem. Some of them later were 
extended from Lebesgue spaces to Morrey spaces
with mixed norms in the case of parabolic equations
with successful applications to PDEs.

Since the famous results of Cald\'eron and Zygmund (1952)
until   now  many results on solvability of nondivergence form elliptic
and parabolic equations and systems are based on the theory
of singular integrals developed by them.
An impressive impact of Real Analysis on PDEs.
The theory of singular integrals, basically, requires explicit 
integral representation for  the parametrix of the problem and allows
treating differential operators with continuous main coefficients.
In 1991 Chiarenza, Frasca, and Longo by using
a commutator theorem from the theory  
of singular integrals made another long lasting impact
allowing higher order coefficients of elliptic equations to
be in VMO (vanishing mean oscillation, or having small
local BMO norm). Their method is based again on explicit
formulas for parametrix  and is still in 
recent use  
(see, for instance \cite{BSS_24}).

Since 2007 another way of using Real Analysis results emerged.
This approach   based on the Fefferman-Stein theorem
is developed in a series of works
by Krylov, Dong, and Kim (see \cite{Kr_08}, \cite{DK_09}, \cite{DK_11},
\cite{DK_11_1}, \cite{DK_18}) in some of which the Muckenhoupt weights
(one more tool from Real Analysis)
play major role. Using the Adams (1986) or the Chiarenza-Frasca
(1990) theorem (from Real Analysis) allowed treating singular lower-order terms
in Sobolev space theory, much more singular than the Sobolev
embedding theorems would allow (see, for instance, \cite{TTV_95}, \cite{CLV_96}, \cite{MT_13}).

Finally, Morrey-Sobolev spaces showed themselves as a powerful
tool of investigating equations with {\em singular\/} lower order terms.
They were extensively used in many papers on {\em general\/} issues
in the theory of elliptic equations, with \cite{DR_93},  \cite{TTV_95} , \cite{FLY_98}
as ones of the first ones. 

 A recurrent gap in many such papers (including \cite{Li_03}, \cite{GS_15}) is that
the authors  estimate the Morrey norm of only the second derivatives
of the unknown function without paying any attention to
the first derivatives, which automatically
appear, for instance, after flattening the boundary of a domain. They, probably, tacitly assumed that the
latter  can be estimated as usual by interpolation inequality.
This is indeed possible although not completely trivial
and proved rather late in \cite{Kr_22}. The first attempt to
prove it, known to the author, was made also rather late in 
\cite{FHS_17} with an error in the proof.

It is somewhat unfortunate that there are no single
publication where the reader could
 find a thorough
exposition of Hardy-Littlewood maximal function
theorem, Fefferman-Stein theorem,
theory of Muckenhoupt weights,  Rubio de Francia
extrapolation theorem and the elements of Morrey-Sobolev spaces
relevant to PDEs, all in one place.
Of course, there are books on Morrey spaces \cite{Ad_15}, \cite{SDI_20},
the latter citation consisting of two books partly aimed at applications to PDEs
with total volume over 800 pages, and none of them
treats even the gap mentioned above.

This was one of motivations to write the present article:
to have all relevant tools in one place.

We do not touch another set of tools 
used to investigate
elliptic and parabolic equations based on such notions as
$H^{\infty}$-calculus, $R$-boundedness, maximal regularity,
$\gamma$-radonifying operators, UMD Banach spaces. The reader can get
acquainted with this approach starting with \cite{GV_17}, \cite{HNVW_23},
and the references therein.

In our exposition we do not try to give the strongest
known results for particular equations in particular spaces.
We only show how the Real Analysis results, we present
with all proofs, can be used in model cases arising
in PDEs. Nevertheless the last section of the article
contains   new results.

The paper consists of seven sections: 1. Some tools from
Real Analysis, 2. Morrey-Sobolev spaces, 3. Muckenhoupt
weights, 4. Parabolic equations, 5. Rubio de Francia
extrapolation theorem and mixed norm inequalities,
6. Parabolic equations in mixed norm Sobolev and
Morrey-Sobolev spaces, 7. Parabolic equations with
variable main part. In Section 1 we prove the Fefferman-Stein
and the Hardy-Littlewood theorems in $\bR^{d}$. Various versions
of these theorems in parabolic setting, with Muckenhoupt
weights, and in Morrey spaces are given in later section.
Here we also prove the solvability of elliptic equations
with BMO coefficients. In Section 2 we investigate
interpolation and embeddings for the Morrey-Sobolev spaces
and prove the solvability in these spaces
 of the Laplace equation
with first-order coefficients in Morrey spaces. 
Section 3 is devoted to the theory of Muckenhoupt weights
and its application to the solvability of
the Laplace equation
with Morrey first-order coefficients  
in weighted Sobolev spaces. In Section 4
we consider applications of the previous theory
to parabolic equations with Morrey first-order coefficients
 in parabolic Morrey-Sobolev
spaces. To be able to treat them in mixed norm
spaces we need the 4-page long Section 5, where, in particular,
we present Dong-Kim mixed norm theorem. Section 6
deals with parabolic equations, whose   main part
  is the heat operator and the first-order
coefficients are in mixed norm Morrey space, in mixed norm 
Sobolev and Morrey-Sobolev spaces. The final Section 7
is devoted to the same issues as Section 6 but
for operators with BMO second-order coefficients
and in the mixed norm spaces, in which the norm
is defined with inside integral taken with respect
to $t$ and not $x$ as usual.

 We finish the introduction  with some basic notation. 
 All other notations can be tracked down by using 
 the index at the end of the article.

By $\bR^{d}$ we denote the Euclidean space
of points $x=(x^{1},...,x^{d})$,
\index{$B$@Sets!$\bR^{d}$}% 
 $\bR^{d+1}=\{(t,x):t\in\bR,x\in\bR^{d}\}$.
\index{$B$@Sets!$\bR^{d+1}$}% 
In
the notation of function spaces over
$\bR^{d}$ this symbol is dropped, so that, for instance,  $C^{\infty}_{0}=C^{\infty}_{0}(\bR^{d})$ 
\index{$A$@Sets of functions!$C^{\infty}_{0}$}%
is the set
of infinitely differentiable functions
on $\bR^{d}$ with compact support.
On the other hand, we use $C^{\infty}_{0}$
to denote also $C^{\infty}_{0}(\bR^{d+1})$
with the hope that what $C^{\infty}_{0}$
means in each particular case will be clear
from the context.
 
 Talking about functions on domains
in a Euclidean space we always mean
 Borel measurable ones.

The summation   convention over repeated 
indices is enforced throughout the article.  

If 
\index{$N$@Norms!$"|a,b"|$}% 
$(a,b)=(a^{i...},b^{j...})$
 set   
$$
  |a,b|^{2}:=\sum_{i...}|a^{i...}|^{2}
+\sum_{j...}|b^{j...}|^{2},
$$
so that, for instance, if $u$ is a real-valued and $a^{i...}$ a tensor-valued
functions on $\bR^{d}$, then
$$
\|u,a\|_{L_{2}}^{2}=\int_{\bR^{d}}|u,a|^{2}\,dx=\int_{\bR^{d}}\Big(
|u|^{2}+\sum_{i...}|a^{i...}|^{2}\Big)\,dx.
$$
 
In the proofs of various results  we use
the symbol $N$  to denote finite 
nonnegative constants,
which may change from one occurrence to another and,
 if in the statement of a result there are constants
called $N$ which are claimed to depend only on certain
parameters, then in the proof of the result
the constants $N$ also depend only on the same
parameters unless specifically stated otherwise.
Of course, if we write 
$
N=N(...),
$
 this means that $N$ depends only
on what is inside the parentheses.

During the work on the articles the results of which constitute major part of this article I
was using advice, inspiration, and some criticism from
a few people. Hongjie Dong, A. Lerner, 
D. Kinzebulatov,  A.I. Nazarov,
 and I.~Gy\"ongy   deserve my special gratitude.

\mysection
{Some tools from  Real Analysis}

                                             \label{chapter 06.5.1.1}

\setcounter{equation}{0}
\setcounter{theorem}{0}

\subsection
{Partitions and stopping times}
                                           \label{section 06.5.1.1}

Fix some integers $k_{1},...,k_{d}\geq1$
and call any 
$$
\cB_{l,x}=x+[0,l ^{ k_{1}} )\times
...\times [0,l ^{ k_{d} }),\quad x\in\bR^{d},l>0
$$
 a ``box'' or half-closed ``box''
\index{$B$@Sets!$\cB_{l,x}$}% 
of ``size'' $l$.   
Observe that if 
$$
k_{1}=...=k_{d}=1,
$$
 then we deal just 
with usual boxes whose edges are parallel to the
coordinate axes and which are common
in the theory of elliptic equations. If 
$$
k_{1}=2,
 k_{2}=...=k_{d}=1,
$$
we are dealing with parabolic ``boxes'' often
used in the theory of second-order parabolic equations.

For $n\in\bZ=\{0,\pm1,...\}$
introduce the following families of ``dyadic boxes''
by  
$$
\frB_{n}=\{\cB_{n}(i ):
i \in\bZ^{d} \},
$$
 where, for 
$$i_{n}:=(i^{1}2^{-k_{1}n},
...,i^{d}2^{-k_{d}n}),\quad \bm1
=(1,...,1)
$$
we set 
$$
\cB_{n}(i)=[i_{n},(i+\bm1)_{n})=[i^{1}2^{-k_{1}n},(i^{1}+1)2^{-k_{1}n})\times...\times
[i^{d}2^{-k_{d}n},(i^{d}+1)2^{-k_{d}n}).
$$
Also set $\frB=\bigcup_{n\in \bZ}\frB_{n}$.

Define $\cF_{n}$ as the collection of subsets of $\bR^{d}$
consisting
of an empty set and of all unions
of some elements of $\frB_{n}$.
Obviously $\cF_{n}\subset\cF_{m}$
for $n\leq m$. We call $\cF_{n},n\in\bZ$, {\em a filtration
of partitions\/}. 
\index{filtration of partitions}%
One obvious but useful property of filtrations is that
if $n\leq m$ and $\cB_{n}\in\frB_{n}$, $\cB_{m}\in\frB_{m}$,
then $\cB_{n}\cap \cB_{m}$ is either $\cB_{m}$ 
($\cB_{n}$ is a ``parent'' of $\cB_{m}$) or empty.

If $\tau=\tau(x)$ is a function on $\bR^{d}$ with values in
$\{\infty,0,\pm1,\pm2,...\}$, we call $\tau$ {\em a stopping time\/}
\index{stopping time}%
(relative to $\{\cF_{n}\}$)
if, for each $n=0,\pm1,\pm2,...$,  
$$
\{x:\tau(x)=n\}\in\cF_{n}.
$$ 
The simplest example of a stopping time is given by
$\tau(x)\equiv0$. 

If $\tau$ is a stopping time we denote by $\cF_{\tau}$
the collection of 
  Borel $A$ such that, for any $n\in\bZ$
we have
$$
A\cap\{\tau=n\}\in \cF_{n}.
$$
Observe that
if we are given two stopping times $\tau$ and $\sigma$
and $\sigma\leq \tau$, then $\cF_{\sigma}
\subset\cF_{\tau}$ since for $A\in\cF_{\sigma}$
$$
 A\cap\{\tau=n\}=\bigcup_{k=-\infty}^{n }
 A\cap\{\sigma=k\}\cap\{\tau=n\}
$$
and $A\cap\{\sigma=k\}\in\cF_{k}\subset \cF_{n}$.
Obviously the intersection of two sets in $\cF_{\tau}$
belongs to $\cF_{\tau}$. An easy and useful
fact is that 
$$
\bR^{d},\{\tau<\infty\}\in \cF_{\tau}.
$$ 
Also a useful fact to remember is that if $A\in\cF_{\tau}$,
then $A\cap\{\tau<\infty\}$ is the disjoint union
of $A\cap\{\tau=n\}\in \cF_{n}$ and
each $A\cap\{\tau=n\}$ is either empty or is the disjoint union of
some $\cB\in\frB_{n}$
such that $\cB=\cB\cap\{\tau=n\}$, the latter showing that
$\cB\in \cF_{\tau}$.
 
We assume that we are given a measure $\mu$
on Borel subsets of $\bR^{d}$ finite on compact sets and such that, for any $x\in\bR^{d}$,
 \begin{equation}
                                                  \label{11.28.1}
\lim_{l\to\infty} \mu(\cB_{l,x} )=\infty.
\end{equation}

Whenever it makes sense we
  use the notation 
$$
  f_{A}=\dashint_{A}f\,\mu(dx):=\frac{1}{\mu(A)}\int_{A}f(x)\,\mu(dx)  
\quad \bigg(\frac{0}{0}:=0\bigg)
$$
 for  the 
\index{$S$@Miscelenea!$\dashint\phantom{m}$}%
\index{$S$@Miscelenea!$f_{A}$}%
average value of $f$ over $A$. 
 
Next, for each $x\in \bR^{d}$ and $n\in \bZ$ 
there exists (a unique)
$\cB\in\frB_{n}$ such that $x\in \cB$. We denote 
this $\cB$ by $\cB_{n}(x)$.

It turns out  that
for any $\lambda>0$ and $f\geq 0$, such that
$ f _{\cB_{n}(x)}\to0$ as $n\to-\infty$
for any $x$, 
\begin{equation}
                                            \label{11.20.1}
\tau_{\lambda}(x)=\inf\{n: f_{ |n}>\lambda\}
\quad   (\inf\emptyset:=\infty,\quad
f_{|n}(x)= f_{\cB_{n}(x)} ) 
\end{equation}
is a stopping time. Indeed,
  observe that 
$$
\cB_{n}(x)\subset \cB_{m}(x)
$$
for all $m\leq n$ since the partitions are nested.
It follows that, if $y\in \cB_{n}(x)$, then 
$$
\cB_{m}(y)=\cB_{m}(x),\quad  f_{|m}(y)
= f _{|m}(x),\quad
\forall m\leq n.
$$
 By adding that
$$
\tau(x)=n\Longleftrightarrow  f _{|n}(x)>\lambda,\quad
 f _{|m}(x)\leq\lambda\quad\forall m<n,
$$
we conclude that the set $\{\tau=n\}$ contains $\cB_{n}(x)$
 along with each $x$. Therefore, $\{\tau=n\}$ is indeed
the union of some elements of $\frB_{n}$.

 In our notation
$$
f_{|n}(x)= f_{\cB_{n}(x)}=\dashint_{\cB_{n}(x)}f(y)\,\mu(dy)
$$
we read $f_{|n}$ as ``$f$ given $\cF_{n}$'', continuing to
\index{$S$@Miscelenea!$f_{"|n}$}%
borrow  the terminology from probability theory. The most important property of $f_{|n}$ is that, obviously,
it is constant on each $\cB\in\frB_{n}$, so that, for any
$x\in \cB\in\frB_{n}$, we have
$$
\int_{\cB}f_{|n}(y)\,\mu(dy)=|\cB|f_{|n}(x)=\int_{\cB}f(y)\,\mu(dy).
$$
If we are also given a stopping 
\index{$S$@Miscelenea!$f_{"|\tau}$}%
time $\tau$, we let
$$
f_{|\tau}(x) =f_{|\tau(x)}(x)
$$
 for those $x$ for which
$\tau(x)<\infty$  and $f_{|\tau}(x)=f(x)$ otherwise.

We suppose that $\mu$ satisfies the ``doubling condition'':
  for any $n$, $\cB\in \frB_{n}$, and $\cB'\in\frB_{n-1}$ such that 
\index{doubling condition}% 
$\cB\subset \cB'$ we have
\begin{equation}
                                            \label{11.21.2}
\mu(\cB')\leq N_{0}\mu(\cB),
\end{equation}
 where $N_{0}$
is a constant independent of $n,\cB, \cB'$.
One of consequences of this condition is that
for $f\geq0$
on the set where $\tau_{\lambda}(x)=n$
we have $f_{|n-1}(x)\leq\lambda$ and
\begin{equation}
                                     \label{11.23.1}
f_{\tau}(x)=\frac{1}{\mu(\cB_{n}(x))}
\int_{\cB_{n}(x)}f\,\mu(dy)
\leq \frac{N_{0}}{\mu(\cB_{n-1}(x))}
\int_{\cB_{n-1}(x)}f\,\mu(dy)\leq N_{0}\lambda.
\end{equation}
 
Another consequence of \eqref{11.28.1} and \eqref{11.21.2}
is that $\mu(\cB_{l,x})>0$ for any $x\in \bR^{d}$ and $l>0$.
 
In the following lemma by $I_{A,\tau<\infty}$
we mean the indicator function of the set
$\{x\in A:\tau(x)<\infty\}$.
\begin{lemma}
                                 \label{lemma 06.5.30.1}

(i) Let   $f$ be Borel on $\bR^{d}$, $f\geq0$,
  let $\tau$ be a stopping time, and let $A\in \cF_{\tau}$.  
  Then
\begin{equation}
                                          \label{06.5.1.1}
\int_{\bR^{d}} f_{|\tau}(x) I_{A,\tau<\infty}\,\mu(dx)=
\int_{\bR^{d}} f(x)  I_{A,\tau<\infty}\,\mu(dx) ,
\end{equation} 
\begin{equation}
\label{07.9.27.2}
\int_{\bR^{d}} f_{|\tau}(x) I_{A }\,\mu(dx)=
\int_{\bR^{d}} f(x)  I_{A }\,\mu(dx) .
\end{equation}

(ii) Let   $f$ be Borel  on $\bR^{d}$, $f\geq0$,
and let $\lambda>0$ be a constant. 
Assume that $f_{|n}(x)\to0$ as $n\to-\infty$
at any $x$.
Then for $\tau=\tau_{\lambda}$ defined in \eqref{11.20.1} 
we have
\begin{equation}
                                          \label{11.18.6}
\lambda I_{\tau<\infty} <f_{|\tau}(x)
I_{\tau<\infty}\leq N_{0}\lambda,
\end{equation}
and for any $A\in\cF_{\tau}$
\begin{multline}
N_{0}^{-1}\lambda^{-1}\int_{\bR^{d}}f(x)
I_{A,\tau<\infty}\,\mu(dx)\leq 
\mu(\{x\in A:\tau(x)<\infty\})
\\                          \label{06.5.1.3}
\leq 
 \lambda^{-1}\int_{\bR^{d}}f(x)
I_{A,\tau<\infty}\,\mu(dx).
\end{multline}

\end{lemma}

 Proof. (i) Equation \eqref{07.9.27.2} follows immediately from
\eqref{06.5.1.1} since $f_{|\tau}=f$ on the set where $\tau=
\infty$.
Owing to the additivity of
the integral, it suffices to prove
\eqref{06.5.1.1} with $\tau=n$ in place
of $\tau<\infty$ and, since the set
$\{x\in A:\tau(x)=n\}$ is the disjoint
union of some $\cB\in \frB_{n}$, it only remains to observe that
for such $\cB$ we have $\cB\cap\{\tau=n\}=\cB$ and
$$
\int_{\bR^{d}} f_{|\tau}(x) I_{\cB,\tau=n}\,\mu(dx)
 =\int_{\bR^{d}} f_{\cB} I_{\cB}\,\mu(dx)
=\int_{\cB}f\,\mu(dx)=\int_{\bR^{d}}fI_{\cB,\tau=n}\,\mu(dx).
$$

(ii) Relations \eqref{11.18.6} follow
from the definition of $\tau$ and \eqref{11.23.1}.
The first inequality in \eqref{06.5.1.3}
follows from \eqref{06.5.1.1} and \eqref{11.18.6}.
The second one follows from Chebyshov's
inequality and \eqref{06.5.1.1} because
$$
\mu(\{x\in A:\tau(x)<\infty\})=
\mu(\{x\in \bR^{d}:f_{\tau}I_{A,\tau <\infty}>\lambda\}).
$$
The lemma is proved. \qed

It is worth clarifying the typical structure 
of $f_{|\tau}$ in \eqref{11.18.6}.
First we take any $m$ so largely negative that 
$$
A_{m}:=\{x:f_{|m}(x)\leq\lambda\}\ne\emptyset. 
$$

 Let $\cB_{j}$, $j=1,2,...$, be the set of all ``boxes'' in the family
$\frB_{m}$, so that 
$$
A_{m}=\bigcup_{j}\cB_{j}.
$$   
Then  we divide each ``box'' $\cB_{j}$
into smaller ``boxes'' $\cB_{jk}\in\frB_{m+1}$, so that 
$\cB_{j}=\bigcup_{k}\cB_{jk}$
and we look for those 
$$
\cB_{jk}\quad\text{on which}\quad f_{|m+1}>\lambda.
$$
 We set those
aside and set $\tau=m+1$ and accordingly $f_{|\tau}=f_{|m+1}$
on them.
With the {\em remaining\/} ``boxes'' on which $f_{|m+1}\leq\lambda$
we proceed in the same way  splitting each of them into
smaller ``boxes'' $\cB_{jkp}\in\frB_{m+2}$ and
 defining  $\tau=m+2$ and $f_{|\tau}=f_{|m+2}$
  on those $\cB_{jkp}$ on which $f_{|m+2}>\lambda$.
By continuing in this way, we define $\tau$ and $f_{|\tau}$
on a subset of $A_{m}$, which may not coincide with
$A_{m}$ (it is just empty if $f\leq\lambda$ everywhere).
The remaining set is the set of points $x\in A_{m}$ at which
 $f_{|n}(x)\leq\lambda$ for all $n$
and we let $\tau=\infty$ and $f_{|\tau}=f$ at those points. Going down to lower $m$ is done
in an obvious way.

What we obtain is called the Riesz-Cald\'eron-Zygmun decomposition of $\bR^{d}$.

\subsection  
{Maximal and sharp functions}  
                                             \label{section 06.6.9.1}
Here we keep using the notation
and the definitions from Subsection \ref{section 06.5.1.1} and 
the introduction to the section.

\subsection{Fefferman-Stein theorem}

The {\em ``dyadic'' maximal   function \/}
\index{``dyadic'' maximal  function}% 
of $f$ is defined by
\begin{equation}
                           \label{12.28.1}
\cM f(x)=\sup_{n<\infty}|f|_{|n}(x),
\end{equation}
 so that
$\cM f=\cM |f|$.
 Observe that, if $f\geq0$, then
$$
\{x:\cM f(x)>\lambda\}=\{x:\tau_{\lambda}<\infty\},
$$ 
where $\tau_{\lambda}$ is taken from
\eqref{11.20.1}. Therefore, Lemma \ref{lemma 06.5.30.1}  (ii) implies the following.

\begin{corollary}
                    \label{corollary 12.25.1}
Under the conditions of Lemma \ref{lemma 06.5.30.1}  (ii)
$$
(N_{0}\lambda)^{-1}\int_{\bR^{d}}f(x)
I_{\cM f>\lambda}\,\mu(dx)\leq 
\mu(\{x\in A:\cM f(x)>\lambda\}) 
$$
\begin{equation}
                       \label{12.25.2}
\leq 
 \lambda^{-1}\int_{\bR^{d}}f(x)
I_{\cM f>\lambda}\,\mu(dx).
\end{equation}
\end{corollary}

Here the right inequality is called
{\em the  maximal 
\index{maximal inequality}%
inequality\/}.

It is worth noting that both inequalities in
\eqref{06.5.1.3} are crucial
in the proof of the reverse H\"older's inequality
for Muckenhoupt's weights, and the right
inequality is used in a crucial way in proving
the Fefferman-Stein Theorem~\ref{theorem 11,23,1}.

The following remark will not be used in the future.
It was hard not to make it. 
\begin{remark}
                                 \label{remark 12.20.1}
The first inequality in \eqref{06.5.1.3} is instrumental
not only in proving
  the reverse H\"older
inequality for $A_{p}$-weights but its version
\eqref{12.25.2}
also
is crucial  in the proof of the first part
of a remarkable
Stein-Zygmund result that, for any 
Borel $f\geq0$ and $\lambda_{0}>0$,
\begin{equation}
                                      \label{12.20.3}
\int_{\bR^{d}}\cM fI_{\cM f>\lambda_{0}}\,\mu(dx)<\infty
\Longrightarrow 
\int_{\bR^{d}}  fI_{  f>\lambda_{0}}\log(f/\lambda_{0})\,\mu(dx)<\infty,  
\end{equation}
\begin{equation}
                                      \label{12.20.4}
\int_{\bR^{d}}  fI_{  f>\lambda_{0}}\log(f/\lambda_{0})\,\mu(dx)<\infty
\Longrightarrow 
\int_{\bR^{d}}\cM fI_{\cM f>2\lambda_{0}}\,\mu(dx)<\infty.
\end{equation}

Here \eqref{12.20.3} is obtained just by integrating
with respect $\lambda\in(\lambda_{0},\infty)$   the
first inequality in \eqref{12.25.2}, where on the left
$\cM f$ is replace with a smaller quantity $f$.
To prove \eqref{12.20.4}   use
the second inequality in \eqref{12.25.2}, which
after observing that 
$$
\cM (fI_{ f \leq \lambda})\leq\lambda
$$
implies that
$$
  \mu(
 x:\cM f(x)>2\lambda )\leq
 \mu(
 x:\cM (fI_{ f >\lambda})(x)> \lambda )
$$
$$
\leq \lambda^{-1}\int_{\bR^{d}}
f(x)I_{ f >\lambda} \,\mu(dx).
$$
Then again integrate
with respect $\lambda\in(\lambda_{0},\infty)$ the inequality between the extreme terms.
After that one will only need to estimate $2
\lambda_{0}\mu(
\{x:\cM f(x)>2\lambda_{0}\})$  and observe that
this quantity is less than
$$
2\lambda_{0}
\mu(
 x:\cM fI_{f>3\lambda_{0} /2}(x)>3\lambda_{0} /2 )
\leq(4/3)\int_{\bR^{d}}
fI_{f>3\lambda_{0} /2} \,\mu(dx),
$$
where the last term is finite due to  the condition in \eqref{12.20.4}.

At the same time, it is easy to see that,
 if $\mu$ is Lebesgue measure and $\cM f\in  L_{1}$, then $f=0$. Here and below
\index{$A$@Sets of functions!$L_{p}$}% 
\index{$N$@Norms!$"|"|f"|"|_{L_{p}}$}%
by $L_{p}$, $p\geq 1$, we mean the space of  functions $f$ on $\bR^{d}$ such 
that
$$
\|f\|^{p}_{L_{p}}:=\int_{\bR^{d}}|f|^{p}\,
\mu(dx)<\infty.
$$

\end{remark}

\begin{theorem}[Lebesgue]
                                       \label{theorem 9.3.1}
Let   $f\in L_{1,loc} $.  
\index{Lebesgue differentiating theorem}%
Then $f_{|n}\to f$  $\mu$-almost everywhere
as $n\to\infty$.
\end{theorem}
Proof. We may assume that $f\in L_{1} $.
If $f$ is the indicator of a box element of $\frB$, the assertion is obvious.
Then for any  $\phi$ which is a step-function
which is constant on each element of $\frB$ and 
$\varepsilon>0$
$$
\mu(\nlimsup_{n\to\infty}|f_{|n}- f|\geq 3\varepsilon)
\leq \mu(\cM(f-\phi)\geq\varepsilon)+
\mu(|f-\phi|\geq\varepsilon)
$$
$$
\leq \varepsilon^{-1}2\int_{\bR^{d}}|f-\phi|\,\mu(dx).
$$
 This yields the result since the last integral
 can be made arbitrarily small on account of
 choosing appropriate $\phi$. \qed

\begin{remark}
Our next interest in estimating $\mu(\cM g>\lambda)$ 
as in the right estimate in \eqref{12.25.2} is based on the
following formula valid for any $f\geq0$ in light of Fubini's
theorem:
$$
\int_{\bR^{d}} f(x)    \,\mu(dx)=\int_{\bR^{d}}\big(
\int_{0}^{ f(x) }dt\big)  \,\mu(dx)=
\int_{\bR^{d}}\big(
\int_{0}^{\infty}I_{ f(x) >t}\,dt\big)  \,\mu(dx)
$$
 \begin{equation}
                                                 \label{06.06.29.1}
=\int_{0}^{\infty}\big(\int_{\bR^{d}}I_{ f(x) >t}
 \,\mu(dx)\big)\,dt=\int_{0}^{\infty}\mu(  f  >t)\,dt.
\end{equation} 
\end{remark}

\begin{corollary}
                                          \label{corollary 06.5.31.9}
Let $p\in(1,\infty)$, $g\in L_{1} $,
$g\geq0$. Then
$$
 \|\cM g\|_{L_{p}}\leq q\|g\|_{L_{p}},
$$
where $q=p/(p-1)$.
\end{corollary}

Indeed, from \eqref{06.06.29.1}, \eqref{12.25.2}, and
Fubini's theorem we conclude
that, for any finite constant $\nu>0$,
$$
\|\nu\wedge \cM g\|_{L_{p}}^{p}=\int_{0}^{\infty}
\mu(\nu\wedge \cM g >\lambda^{1/p})\,d\lambda
$$
$$
=\int_{0}^{\nu^{p}}\mu(\cM g >\lambda^{1/p})\,d\lambda
\leq\int_{\bR^{d}}
g(x) \big(\int_{0}^{\nu^{p}}\lambda^{-1/p}
I_{\cM g(x)>\lambda^{1/p}}\,d\lambda\big)\,\mu(dx)
$$
$$
=\int_{\bR^{d}}
g(x)\big(\int_{0}^{(\nu\wedge \cM g(x))^{p}}
\lambda^{-1/p}\,d\lambda\big)\,\mu(dx)
=q\int_{\bR^{d}}
(\nu\wedge \cM g)^{p-1}g \,\mu(dx).
$$
This and $g\in L_{1}(\bR^{d})$ imply that $\|\nu\wedge \cM g\|_{L_{p}
  }
<\infty$. Then upon using  H\"older's inequality, we get
$$
\|\nu\wedge \cM g\|_{L_{p}}^{p}\leq
q\|g\|_{L_{p} }\|\nu\wedge \cM g\|_{L_{p}}^{p-1},
\quad
\|\nu\wedge \cM g\|_{L_{p}}\leq q\|g\|_{L_{p}}
$$ 
and it only remains to let $\nu\to\infty$ and
use Fatou's theorem.

Next we extend Corollary \ref{corollary 06.5.31.9} to
$g\in L_{p}$ .

\begin{theorem}
                                         \label{theorem 06.06.29.1}
For any
$p\in(1,\infty)$ and $g\in L_{p}$ 
$$
 \|\cM g\|_{L_{p} }\leq q\|g\|_{L_{p}}.
$$
\end{theorem}

Proof. Since
$$
\cM g=\cM |g| \quad\text{and}\quad\|g\|_{L_{p}}=
\|\,|g|\,\|_{L_{p} },
$$
we may concentrate on   $g\geq0$. 

In that case  take an increasing
sequence of $B_{m}\in \cF_{0}$
such that $\bR^{d}=\cup B_{m}$ and
$\mu(B_{m})<\infty$ for any $m$ and
 introduce
  $g^{m}(x)=g(x)I_{B_{m}}$. Then  
$g^{m}\in L_{1}$ and 
$$
\|\cM g^{m}\|_{L_{p}}\leq q\|g^{m}\|_{L_{p}}
\leq q\|g\|_{L_{p}}
$$
 by Corollary \ref{corollary 06.5.31.9}.
 It only remains to use Fatou's theorem
along with the observation that  
for any $x$ and $n$, since $ \mu(\cB_{n}(x))<\infty$, we have
$$
(g^{m})_{|n}(x)\to g_{|n}(x)\quad\text{as}\quad
m\to\infty,
 $$
 which implies
$$
 g_{|n}(x)\leq\nliminf_{m\to\infty}\sup_{r}(g^{m})_{|r}(x),
\quad
\cM g\leq\nliminf_{m\to\infty}\cM g^{m}.
$$
 \qed

Let   $f\in L_{1,loc} $.
 Define {\em the sharp function of\/}
\index{sharp function}% 
 $f$
\index{$S$@Miscelenea!$f^{\#}$}%
by
 \begin{equation}
                                                 \label{07.11.29.3}
f^{\#}(x)=\sup_{n<\infty}\dashint_{\cB_{n}(x)}|f(y)-
f_{|n}(y)|\,\mu(dy). 
\end{equation} 

We say that $f$ is a  {\em dyadic BMO function\/}
(bounded mean oscillation) 
\index{dyadic BMO function}%
if $f^{\#}$ is bounded.

\begin{remark}
                         \label{remark 9.3.2}
Observe that
$$
\dashint_{\cB_{n}(x)}|f(y)-
f_{|n}(y)|\,\mu(dy)=\dashint_{\cB_{n}(x)}\Big|f(y)-\dashint_{\cB_{n}(x)}f(z)\,
\mu(dz)\Big|\,\mu(dy)
$$
$$
\leq \dashint_{\cB_{n}(x)}\dashint_{\cB_{n}(x)}|f(y)-f(z)|\,\mu(dy)\mu(dz)
\leq 2\dashint_{\cB_{n}(x)}|f(y)-
f_{|n}(y)|\,\mu(dy),
$$
and, since 
$$
|\,|f(y)|-|f(z)|\,|\leq|f(y)-f(z)|,
$$
 we have
$|f|^{\#}\leq 2f^{\#}$.
\end{remark}

Obviously $f^{\#}(x)\leq 2\cM f(x)$. It turns out that $\cM f$ and
hence
$f$ are also controlled by $f^{\#}$.

The following result is false for
  locally integrable
functions.

\begin{lemma}
                                              \label{lemma 06.11,23,1}
For $\alpha=(2N_{0})^{-1}$,
 any constant $c>0$, and  $f\in L_{1}$, we have
 $$
\mu(|f|\geq c)\leq \frac{4}{c}\int_{\bR^{d}}I_{\cM f(x)>\alpha c}
f^{\#}(x)\,\mu(dx)
$$
and if $f\geq0$, then one can replace $4/c$ with $2/c$.
\end{lemma}

Proof. Remark \ref{remark 9.3.2} shows that
it suffices to prove the second assertion of the lemma.
Introduce
$$
\tau(x)=\inf\{n: f_{|n}(x)>c\alpha\}.
$$
Use Lemma \ref{lemma 06.5.30.1} (ii) to get that
$f_{|\tau}\leq c/2$ if $\tau<\infty$ and also use
  the fact that 
 $f_{|n }\to f$ (a.e.).
Then we find that (a.e.)
\begin{align*}
&\{x: f(x) \geq c\}=
\{x: f(x) \geq c,\tau(x)<\infty\}
\\ 
&\qquad=
\{x: f(x) \geq c,  f_{|\tau}(x)\leq c/2\}
\subset\{x:  
|f(x)-f_{|\tau}(x)| \geq c/2\} .
\end{align*} 
By
   Chebyshov's inequality  and Lemma \ref{lemma 06.5.30.1}
$$
\mu(|f| \geq c) \leq(2/c)\int_{\bR^{d}} 
|f(x)-f_{|\tau}(x)| \,\mu(dx)
$$
$$
=(2/c)\int_{\bR^{d}} 
|f-f_{|\tau}|_{|\tau}(x) \,\mu(dx)=(2/c)\int_{\bR^{d}} 
|f-f_{|\tau}|_{|\tau}(x)I_{\tau<\infty} \,\mu(dx),
$$ 
where in the  equalities we also used  the fact
that 
$$
|f-f_{|\tau}|_{|\tau}(x)=
|f(x)-f_{\tau}(x)|=0
$$
 when $\tau(x)=\infty$.

Finally, if at an $x$ we have $\tau(x)=n$, then
(recall that if $\tau(x)=n$ and $y\in \cB_{n}(x)$,
then $\tau(y)=n$)
\begin{align*}
|f-f_{|\tau} |_{|\tau}(x)&=\dashint
_{\cB_{n}(x)}|f(y)-f_{|\tau}(y)|\,\mu(dy)
\\
&= \dashint
_{\cB_{n}(x)}|f(y)-f_{|n}(y)|\,\mu(dy)\leq f^{\#}(x).
\end{align*}

Now it only remains to notice that
$$
\{\tau(x)<\infty\}
=\{\cM f(x)> c\alpha\}.
$$   \qed

\begin{theorem}[Fefferman-Stein]
                                              \label{theorem 11,23,1}
\index{Fefferman-Stein theorem}%
Let $p\in(1,\infty)$. Then  
for any $f\in L_{p}$ we have 
$$
\|f\|_{L_{p}}
\leq N\|f^{\#}\|_{L_{p}},
$$ 
where $N=(2q)^{p}N_{0}^{p-1}$, $q=p/(p-1)$.
\end{theorem}

Proof. As in the beginning of the
 proof of Corollary \ref{corollary 06.5.31.9} (with $\nu=\infty$)
we get from Lemma \ref{lemma 06.11,23,1} that
if $f\in L_{1}$, then
$$
\|f\|_{L_{p}}^{p}\leq
N\int_{\bR^{d}}f^{\#}(\cM f)^{p-1}\,\mu(dx).
$$
By using H\"older's inequality, we obtain
$$
\|f\|_{L_{p}}^{p}
\leq N\|f^{\#}\|_{L_{p}}\|\cM f\|^{p-1}_{L_{p}}.
$$
If in addition $f\in L_{p}$, then
  it only remains to use
Theorem \ref{theorem 06.06.29.1} and check that
the resulting constant is right.

If we only have $f\in L_{p}$, then it suffices
to take $f_{n}\in L_{1}\cap L_{p}$
converging to $f$ in $L_{p}$ and observe that
$f^{\#}_{n}\leq (f-f_{n})^{\#}+f^{\#}$ and
$$
\|(f-f_{n})^{\#}\| _{L_{p}}\leq 2
\|\cM (f-f_{n}) \| _{L_{p}}\leq 2q
\| f-f_{n}  \| _{L_{p}}\to0.
$$   \qed

\begin{remark}
This theorem is quite puzzling from the following point of view.
In the definition  \eqref{07.11.29.3} 
 only the integrals of
the {\em first\/} power of $f$ are involved
 and the sup there is certainly not attained
or approximated when $n\to\infty$. Then how can  $f^{\#}$ control the
integral of the $p$th power of $f$ may look mysterious since $p>1$
and generally the integrals of $|f|^{p}$ over sets of finite measure
 are {\em not\/}
controlled by the integrals of $|f|$.
\end{remark}

\subsection%[Comparing maximal and sharp functions]
{Comparing  the   maximal and sharp functions}
                                                \label{section 07.6.6.1}

The maximal and sharp functions introduced  
 in the previous subsection 
are related to the underlying filtration of partitions.
In particular applications the following
 classical maximal and sharp functions   
\index{classical maximal function}%
\index{classical sharp function}%
are used. Introduce  $\bE$ as the union of $\bE_{\rho}$,
$\rho>0$, where $\bE_{\rho}$ is
\index{$B$@Sets!$\bE$}% 
\index{$B$@Sets!$\bE_{\rho}$}% 
 the collection
 of {\em open\/} ellipsoids $\cE_{\rho}(x)$ which are centered at $x$, have main axes parallel
to the coordinate axes, and have the radiuses
parallel to the $x^{i}$-axis equal to
$\rho^{k_{i}} \nu_{0}^{k_{i}}$, where $\nu_{0}$
satisfies 
$$
\nu_{0}^{-2k_{1}}+...+
\nu_{0}^{-2k_{d}}=4.
$$ 
If $\cE=\cE_{\rho}(x)$ we call $\rho$
the ``size'' of $\cE$.
An easy exercise
is to show that
\begin{equation}
                       \label{10.29.1}
\cB_{n}(i)\subset \cE_{2^{-n}}(x_{n}),
\end{equation}
where $x_{n}=(1/2)(i_{n}+(i+\bm{1})_{n})$ is the center of $\cB_{n}(i)$.
Also notice that, if $k_{1}=...=k_{d}
=1$, then 
$$
\cE_{\rho}(x)=\{y:
|y-x|<\rho\nu_{0}\},\quad\nu_{0}=\sqrt d/2.
$$ 

Here we suppose that $\mu$ satisfies
a stronger doubling  condition
than in Subsection \ref{section 06.6.9.1},
namely, there exists a constant $N_{0}$
such that, for any $x\in\bR^{d}$,
$\rho\in (0,\infty)$, we have
$$
\mu(\cE_{2\rho}(x))\leq N_{0}
\mu(\cE_{\rho}(x)).
$$ 

We say that in such case the collection $(k_{1},...,k_{d},\mu,N_{0})$ defines
a (real analytic) {\em structure\/}
\index{real analytic structure}% 
 of $\bR^{d}$.

An example of real analytic structure
in $\bR^{d}$
is given by $k_{i}=1$
for all $i$, Lebesgue measure,
and $N_{0}=2^{d}$. In this case 
$\bE=\bB$, where $\bB=\cup_{\rho}\bB_{\rho}$
and $\bB_{\rho}$ is 
\index{$B$@Sets!$B_{\rho}(x)$}%
\index{$B$@Sets!$B_{\rho}$}%
\index{$B$@Sets!$\bB_{\rho}$}%
\index{$B$@Sets!$\bB$}%
the collection of balls $B_{\rho}(x)$ with
$$
B_{\rho}(x)=\{y\in\bR^{d}:|y-x|<\rho\},\quad B_{\rho}=B_{\rho}(0). 
$$

Another example of such structure is given
by $\mu(dx)=|x|^{\alpha}\,dx$ with $\alpha>-d$, $k_{i}=1$
for all $i$. To check that it is a valid 
example we need to show that
 for any ball
$B_{r}(x)$ we have 
$$
\mu(B_{2r}(x))
\leq N\mu(B_{r}(x)).
$$
 Because
of homogeneity it suffices to concentrate on $r=1$. For $|x|\leq4$,
clearly $\mu(B_{2}(x))$ is bounded
from above and $\mu(B_{1}(x))$ bounded
away from zero. But if
$|x|\geq 4$, then 
$$
\mu(B_{2}(x) )=\int_{B_{2} }|y-x|^{\alpha}\,dy=|x|^{\alpha}
\int_{B_{2} }|y/|x|-x/|x||^{\alpha}\,dy
\leq N |x|^{\alpha}
$$
and it is seen that $\mu(B_{2}(x))
\leq N\mu(B_{1}(x))$ again.

Next, for $\beta\geq0$ define
\index{classical maximal function}% 
the {\em classical maximal function\/}
$$
\bM_{\beta} g(x)=\sup_{\rho>0}\rho^{\beta}\sup_{\cE\in\bE_{\rho}
 }I_{\cE}(x)\dashint_{\cE}|g(y)|\,\mu(dy). 
$$
Since the $\cE$'s are open, $I_{\cE}(x)$ are lower 
semicontinuous, so is
$\bM_{\beta} g$ and, hence, it is Borel measurable.
In these notation we will drop $\beta$ if
$\beta=0$. Clearly, $\cM g\leq N\bM g$ with $N$
independent of $g$.

Also introduce the Euclidean
\index{sharp function}% 
\index{$S$@Miscelenea!$g^{\sharp}$}%
 sharp function 
  $$
g^{\sharp}(x):=\sup_{\cE\in\bE:x\in \cE}
\osc_{\cE}g,\quad \osc_{\cE}g:=\dashint_{\cE}|g(y)-g_{\cE}|\,\mu(dy).
$$
  The reader should
\index{$S$@Miscelenea!$\osc_{\cE}g$}%
notice the difference in shapes of the ``sharp" 
symbols in $g^{\#}$ and~$g^{\sharp}$
(much later for functions depending also on $t$ we introduce $g^{\shharp}$).

\begin{remark}
                   \label{remark 2.21.1}
If $g^{\sharp}$ is bounded, they say that
$g$ belongs to BMO (boun\-ded mean oscillation). Such function can have
rather wild discontinuities.
One example is $g(x)=\sin \ln|x|$.

\end{remark}

\begin{theorem}[Fefferman-Stein]
                                         \label{theorem 06.6.9.1}
Let $g\in L_{1,loc} $; then
\index{Fefferman-Stein theorem}% 
$$
g^{\#}\leq Ng^{\sharp},
$$ 
where the constant $N$ is independent of $g$. In particular,
owing to Theorem~\ref{theorem 11,23,1}, if $p\in(1,\infty)$,
$g\in L_{p} $,   then
$$
\|g\|_{L_{p}}\leq N\|g^{\sharp}\|_{L_{p}}
$$
where $N$   is independent of   $g$.  
\end{theorem}

Proof. For $x\in \cB_{n}(i)$
in light of \eqref{10.29.1} we have
$x\in \cE_{2^{-n}}(x_{n})=:\cE$ and
$$
 \dashint_{\cB_{n}(i)}\dashint_{\cB_{n}(i)}
|g( y)-g( z)|\,\mu(dy)\mu(dz)
$$
$$
 \leq N_{1}
 \dashint_{\cE}\dashint_{\cE}
|g( y)-g( z)|\,\mu(dy)\mu(dz)
\leq 2N_{1}g^{\sharp}(x),
$$
where
$$
N_{1}=\frac{\mu^{2}(\cE) }{\mu^{2}(\cB_{(n)}(i))}
$$
is   dominated by a constant
independent of $n$ and $x$.
 This proves the theorem. \qed

The following theorem is one of the
Hardy-Littlewood theorems.
\index{Hardy-Littlewood theorem}% 
In this theorem  assumption
\eqref{11.28.1} is not needed.

\begin{theorem}
                                            \label{theorem 7.7.1}
Let $p\in(1,\infty)$ and $g\in L_{p}$.
Then $\bM g\in L_{p}$  and
$$
\|\bM g\|_{L_{p}}\leq N\|g\|_{L_{p}},
$$ 
where the constants $N$ is independent of $g$.
\end{theorem}

Proof. 
We may certainly assume that $g\geq0$ and,
as in the case of Theorem \ref{theorem 06.06.29.1},
it suffices to prove that for any $\lambda>0$
and $g\in L_{1}$ we have
\begin{equation}
                                            \label{06.10.18.1}
\mu(A(\lambda))
\leq\frac{N}{\lambda}\int_{\bR^{d }}
I_{A(\lambda)}( x)g( x)\,\mu(dx),
\end{equation} 
where $N$ is independent of $\lambda$, $g$ and
$$
A(\lambda)=\{x:\bM g( x)>\lambda\}.
$$
 From the definition of $\bM g$ it follows easily 
that $A(\lambda)$ is an open set. Take a compact set $K\subset
A(\lambda)$. Then by the definition of $A(\lambda)$ for
any $x\in K$ there exists an $\cE\in\bE$ such that
$x\in  \cE$ and
\begin{equation}
                                            \label{06.10.18.2}
\int_{\cE}g\,\mu(dx)>\lambda\mu(E).
\end{equation} 
By the way, observe that $\cE\subset A(\lambda)$.

By the compactness of $K$, there is a finite collection
$\cE_{1},...,\cE_{n}\in\bE$ covering $K$
and such that for each $\cE=\cE_{i}$ equation \eqref{06.10.18.2} holds.

Now we use   Vitali  covering argument. If $\cE=\cE_{\rho}(x)\in\bE$, then
define $\cE^{*}=\cE_{3\rho}(x)$. Then denote by $\tilde{\cE}_{1}$
any of $\cE_{i}$ which has the largest size  and set it aside.
  Next,
introduce $\tilde{\cE}_{2}$ as one of the remaining $E_{i}$
which has the  largest size between those $\cE_{i}$ that have
{\em no\/} intersection with $\tilde{\cE}_{1}$. It may happen that
no such $\cE_{i}$ exists. Then it is almost obvious that
$\cE_{i}\subset\tilde{\cE}^{*}_{1}$ for any $i$.
If $\tilde{\cE}_{2}$ exists, we proceed further.

If we have already defined $\tilde{\cE}_{1},...,\tilde{\cE}_{k}$,
then we define $\tilde{\cE}_{k+1}$ as one of the ellipsoids in
the family
\begin{equation}
                                            \label{06.10.18.3}
\{\cE_{1},...,\cE_{n}\}
\setminus\{\tilde{\cE}_{1},...,\tilde{\cE}_{k}\}
\end{equation} 
 which is disjoint from $\tilde{\cE}_{1},...,\tilde{\cE}_{k}$ and has the
largest measure between those that
are  disjoint from $\tilde{\cE}_{1},...,\tilde{\cE}_{k}$.
In finitely many steps we will come to a $k$
 for which
any ellipsoid in family \eqref{06.10.18.3} intersects
one of $\tilde{\cE}_{1},...,\tilde{\cE}_{k}$ or else
the family is empty.
In the second case, obviously, for any $i$
 \begin{equation}
                                            \label{06.10.18.4}
\cE_{i}\subset\bigcup_{j=1}^{k}\tilde{\cE}_{j}^{*}.
\end{equation}  

It turns out that \eqref{06.10.18.4} also holds for any $i$
in the first case. Indeed, if, for a fixed $i$,
 $\cE_{i}$ has a nonempty intersection
with a $\tilde{\cE}_{j}$, then define $r=r(i)$ as the smallest such $j$
and observe that, if $r=1$, then as has been pointed out above,
$\cE_{i}\subset \tilde{\cE}_{1}^{*}$ and \eqref{06.10.18.4} holds.
If $r>1$, then the size of $ E_{i}$
is less than or equal to the size
of $\tilde{\cE}_{r}$
by the choice of $\tilde{\cE}_{r}$ and because $\cE_{i}$
has no intersection with $\tilde{\cE}_{1},...,\tilde{\cE}_{r-1}$
by the definition of $r$.
Now as above $\cE_{i}\subset \tilde{\cE}_{r}^{*}$, 
implying \eqref{06.10.18.4}.

  It follows that
 $$
K\subset\bigcup_{j=1}^{k}\tilde{\cE}_{j}^{*}.
 $$

By using \eqref{06.10.18.2} for $\tilde{\cE}_{j}$,
recalling that they are disjoint, and using that
$\cE_{i}\subset A(\lambda)$ for any $i$, we get
$$
\mu(K)\leq\sum_{j=1}^{k}\mu(\tilde{\cE}_{j}^{*})
=N\sum_{j=1}^{k}\mu(\tilde{\cE}_{j} )
$$
$$
\leq
N\lambda^{-1}\sum_{j=1}^{k}\int_{\tilde{\cE}_{j}}g\,\mu(dx)
\leq  \lambda^{-1}\int_{\bR^{d }}gI_{A(\lambda)}\,\mu(dx).
$$ 

We thus obtain \eqref{06.10.18.1} with $K$ in place of 
$A(\lambda)$.
By taking a sequence of compact sets $K_{m}\uparrow
A(\lambda)$ and passing to the limit, we get \eqref{06.10.18.1}.  \qed

Here is a result complementary to Theorem \ref{theorem 06.6.9.1}.
\begin{corollary}
For $p\in(1,\infty)$ there is a constant $N=N(d,p)$
such that, for any $g\in L_{p} $, we have
$$
\|g^{\sharp}\|_{L_{p} }
\leq N\|g^{\#}\|_{L_{p} }.
$$
\end{corollary}

Indeed, it suffices to observe that $g^{\sharp}\leq
2\bM g$  and to
use Theorems \ref{theorem 7.7.1} and~\ref{theorem 06.6.9.1}.

 \subsection
 {Application to the Laplacian}
 
 Let  $\mu$ be Lebesgue measure  and $k_{1}=...=k_{d}=1$. 
Introduce
\index{$C$@Operators!$D_{i}$}%
\index{$C$@Operators!$D_{ij}$}%
$$
D_{i}=\frac{\partial}{\partial x^{i}},
\quad D_{ij}=D_{i}D_{j} , \quad Du=(D_{i}u),
\quad D^{2}u=(D_{ij}u).
$$
 Generally, $D^{n}u$
is the collection of all deriavtives of $u$
of order $n$ with respect to $x$. 
\index{$C$@Operators!$Du$}%
\index{$C$@Operators!$D^{2}u$}%
\index{$C$@Operators!$D^{n}u$}%
By $ D^{n}u $ we always
mean Sobolev derivatives of $u$ (when  they exist). 
Define $W^{2}_{p}$ 
\index{$A$@Sets of functions!$W^{2}_{p}$}%
\index{$N$@Norms!$"|"|u"|"|_{W^{2}_{p}}$}%
as the set of functions $u$ on $\bR^{d}$ such that
$u,Du,D^{2}u\in L_{p}$. Define
$$
\|u\|_{W^{2}_{p} }:=\|u,Du,D^{2}u\|_{L_{p}}.
$$
Similarly we define $W^{1}_{p}$
as the set of 
\index{$A$@Sets of functions!$W^{1}_{p}$}%
\index{$N$@Norms!$"|"|u"|"|_{W^{1}_{p}}$}%
functions $u$ on $\bR^{d}$ such that
$u,Du \in L_{p}$ and set
$$
\|u\|_{W^{1}_{p} }:=\|u,Du \|_{L_{p}}.
$$

 Let $u\in C^{\infty}_{0}$.
 The consecutive integration by parts, first with respect to $i$
 and then with respect to $j$, proves that
 $$
 \int_{\bR^{d}}(D_{ii}u)D_{jj}u\,dx=\int_{\bR^{d}}(D_{ij}u)D_{ij}u\,dx.
 $$
 This yields
 \begin{equation}
                        \label{9.4.1}
 \int_{\bR^{d}}\sum_{i,j}|D_{ij}u|^{2}\,dx=\int_{\bR^{d}}f^{2}\,dx,\quad f=-\Delta u.
 \end{equation}

\begin{remark}
                         \label{remark 4.3,5}
Relation \eqref{9.4.1} is also true
for infinitely differentiable functions $u$
such that 
$$
|x|^{d-1}|Du|+|x|^{d-2}|u|
$$
is bounded.

Indeed, take a $\zeta\in C^{\infty}_{0}$
such that $\zeta=1$ in $B_{1}$, $\zeta=0$
outside $B_{2}$, $0\leq\zeta\leq1$, and set $\zeta_{n}(x)=\zeta(x/n)$. Then for each $n$
\begin{equation}
                           \label{2.24.5}
\|D^{2}(u\zeta_{n})\|_{L_{2}}=\|\Delta(u\zeta_{n})\|_{L_{2}}
\end{equation}
and
$$
\big|\|\Delta(u\zeta_{n})\|_{L_{2}}-
\|\zeta_{n}\Delta u\|_{L_{2}}\big|\leq Nn^{-1}
\||D\zeta|(\cdot/n)|\cdot|^{-(d-1)}\|_{L_{2}}
$$
$$
+Nn^{-2}
\||D^{2}\zeta|(\cdot/n)|\cdot|^{-(d-2)}\|_{L_{2}}= Nn^{1-d}
\||D\zeta| |\cdot|^{-(d-1)}\|_{L_{2}} 
$$
$$
+Nn^{-(d-2)}
\||D^{2}\zeta| |\cdot|^{-(d-2)}\|_{L_{2}}
\to 0 
$$
as $n\to\infty$. Similarly one treats the left-hand side of \eqref{2.24.5} and comes to
\eqref{9.4.1} (with, perhaps, both sides infinite).
 \end{remark}

To continue further 
  let $d\geq3$.
It is a classical fact that for
   an appropriate constant $c_{d}$ we have
 $$
 u(x)=Rf(x):=c_{d}\int_{\bR^{d}}\frac{1}{|x-y|^{d-2}}f(y)\,dy.
 $$
On the other hand,
\index{$C$@Operators!$Rf(x)$}%
 it is well known that,
if $g\in C^{\infty}_{0}$, then $G:=Rg$
satisfies 
$$
\Delta G=-g.
$$
We are going to estimate the integral oscillation of $D^{2}G$ in $B_{1}$ by
splitting $G$ into two parts, one of which
is harmonic in $B_{2}$ and, hence,
very smooth, and the other admits an
estimate of its second-order derivatives
on the basis of Remark \ref{remark 4.3,5}.

 Note that for any $g\in C^{\infty}_{0}$,
$Rg$ is infinitely differentiable and for $|x|$ large
 $$
 |D Rg(x)|\leq N\int_{\bR^{d}}\frac{1}{|x-y|^{d-1}}g(y)\,dy\leq N/|x|^{d-1},\quad|Rg(x)|
\leq N/|x|^{d-2}.
 $$

 Next, take $\zeta \in C^{\infty}_{0}$
such that $\zeta=1$ in $B_{2}$, $\zeta=0$
outside $B_{3}$, $1\leq\zeta\leq1$ and  
observe that, for   $g=f\zeta,h=f(1-\zeta)$
 and $(G,H)=R(g,h)$ in $B_{1}$, we have
 $$
 |D^{3}H(x)|\leq N\int_{|y|>2 }\frac{1}{|x-y|^{d+1}}|f(y)|\,dy
 \leq N\int_{|y|>1}\frac{1}{|y|^{d+1}}|f(y)|\,dy
 $$
 $$
 =N\int_{1}^{\infty}\frac{1}{r^{2}}\dashint_{B_{r}}|f(y)|\,dy dr-
 N \int_{B_{1}}|f(y)|\,dy,
 $$
where the equality is obtained by integrating by parts.
 It follows 
 that $ |D^{3}H(x)|\leq N\bM f(0)$. This yields
 \begin{equation}
                                          \label{9.4.2}
 \dashint_{B_{1}}\dashint_{B_{1}}
 |D^{2}H(y)-D^{2}H(z)|\,dydz\leq N \bM f(0)\leq N\big(\bM(f^{2})(0)\big)^{1/2}.
 \end{equation}
 
 As long as $D^{2}G$ is concerned, notice that in light of \eqref{9.4.1}
 $$
 \dashint_{B_{1}}\dashint_{B_{1}}
 |D^{2}G(y)-D^{2}G(z)|\,dydz\leq 2\dashint_{B_{1}}
 |D^{2}G |\,dx
 $$
 $$
 \leq N 
 \Big(\int_{\bR^{d}}|D^{2}G|^{2}\,dx\Big)^{1/2}
 \leq N\Big(\dashint_{B_{2}}|f|^{2}\,dx\Big).
 $$
 
 Upon combining this with \eqref{9.4.2} we find that for $\rho=1$
$$
 \dashint_{B_{\rho}}\dashint_{B_{\rho}}
 |D^{2}u(y)-D^{2}u(z)|\,dydz\leq  N\big(\bM(f^{2})(0)\big)^{1/2}.
 $$
For any other $\rho>0$ this is derived
by using scaling and yields
\begin{equation}
                      \label{10.31.1}
(D^{2}u)^{\sharp}\leq N\big(\bM(f^{2})\big)^{1/2}
\end{equation}
  at the origin. Of course, this
 inequality is valid also at any other point and then by applying Theorems
 \ref{theorem 06.6.9.1} and \ref{theorem 7.7.1}  we conclude that for any $p>2$
  \begin{equation}
                         \label{9.4.3}
 \|D^{2}u\|_{L_{p} }\leq N(p,d)\|\Delta u\|_{L_{p} }.
 \end{equation}
 
 \begin{theorem}
                    \label{theorem 9.6.1}
 For $d\geq1$, any $p\in(1,\infty)$, and $u\in W^{2}_{p}$ we have
 \eqref{9.4.3}.
 \end{theorem}

 Proof. First assume $d\geq3$. Then by using the fact
 that $C_{0}^{\infty}$ is dense in $W^{2}_{p}$
 we may assume that $u\in C_{0}^{\infty}$.
 Then 
$$
f:=-\Delta u\in C_{0}^{\infty},\quad u=Rf,
$$
 and
 \eqref{9.4.3} holds for $p>2$ by the above.
 
 This estimate also holds for $p=2$ in light 
 of \eqref{9.4.1}. To extend
 it to $p\in(1,2)$ observe that estimate \eqref{9.4.3} means that
 for any smooth $f,g$ with compact support and any $i,j$
\begin{equation}
                                          \label{9.4.4}
 \int_{\bR^{d}}g D_{ij}Rf\,dx\leq N\|g\|_{L_{q}}\|f\|_{L_{p}},
\end{equation}
 where $q=p/(p-1)$. Here on the left we can move first $D_{ij}$
 to $g$ and then $R$ to $D_{ij}g$. By observing that
 $$
RD_{ij}g=D_{ij}Rg
$$
 we come to the version of
 \eqref{9.4.4} with $f,g$ interchanged. Then the arbitrariness
 of $f$ leads to \eqref{9.4.3} with $q$ in place of $p$.
 This proves the theorem if $d\geq 3$.
 
If $d=1$, \eqref{9.4.3} is trivial.
For $d=2$ 
 take $u\in C^{\infty}_{0}(\bR^{2})$, $\zeta\in C^{\infty}_{0}
(\bR^{1})$, such that $\zeta(0)=1$, and for $n=1,2,...$
and $x\in \bR^{3}$ define $u_{n}(x)=u(x^{1},x^{2})\zeta(
x^{3}/n)$. Then \eqref{9.4.3}
yields
$$
\int_{\bR^{2}}  |D^{2}u|^{p}\,dx
\, n\|\zeta\|_{L_{p}(\bR )}^{p}
$$
$$
\leq N \Big(\int_{_{\bR^{2}}}
 |\Delta u|^{p}\,dx\,
n\|\zeta\|_{L_{p}(\bR )}^{p}
+\int_{\bR^{2}} |u|^{p}\,dx\,
n^{1-2p}\|\Delta \zeta\|_{L_{p}(\bR)}^{p}\Big).
$$
By dividing both sides by $n $ and sending
$n\to\infty$ we get \eqref{9.4.3} for $d=2$.  \qed

 We thus obtained \eqref{9.4.3} which  traditionally
 is attributed to Calder\'on-Zygmund and 
 is the basic
 estimate  in the Sobolev space theory of the {\em second-order\/}
 elliptic equations in  nondivergence form. 

\begin{remark} 
                    \label{remark 2.22.2}
 The above proof shows that the operators 
 $D_{ij}R+\delta_{ij}/d$  (which  are
 examples  of
  singular integral operators) extends from $C^{\infty}_{0}$
to $L_{p}$ as   continuous operators in $L_{p}$.
\end{remark}

\begin{remark}
                   \label{remark 1.25.7}
Estimate \eqref{9.4.3} also holds for
  functions not necessarily in $W^{2}_{p}$.
For instance, let $u$ be   smooth
and such that $u|x| ^{\beta}, |x| ^{\beta+1}Du$ are bounded for some $\beta>d/p-2$,   then it turns out that
\eqref{9.4.3}  holds for this $u$ (with perhaps both sides infinite).  This is proved
in the same way as the similar statement about 
\eqref{9.4.1}.

\end{remark}

To move forward we need some interpolation
inequalities.
\begin{lemma}
                       \label{lemma 2.1.1}
Let $u\in C^{\infty}_{0}$. Then 
for any $r>0,\varepsilon\in(0,r]$, $p\geq1$,
we have
\begin{equation}
                          \label{2.1.1}
\int_{B_{r}}|Du|^{p}\,dx\leq N \varepsilon ^{p}
\int_{B_{r}}|D^{2}u|^{p}\,dx+
N \varepsilon ^{-p}
\int_{B_{r}}|u|^{p}\,dx,
\end{equation}
\begin{equation}
                          \label{2.1.2}
(Du)^{\sharp}\leq N(\bM (D^{2}u))^{1/2}
(\bM u)^{1/2},
\end{equation}  
\begin{equation}
                          \label{2.1.3}
\|Du\|_{L_{p}}^{2}\leq N\|D^{2}u\|_{L_{p}}\| u\|_{L_{p}},
\end{equation}
where the constants $N$ are independent of
$r,\varepsilon,u$.

\end{lemma}

Proof. Using scalings allows us to take $r=2$ while proving \eqref{2.1.1}. Then for $\rho\in[1/2,1]$ introduce
$$
\phi(\rho)=\int_{B_{1}}|Du(\rho x)|^{p}\,dx,\quad
U=\int_{B_{1}}|u|^{p}\,dx,\quad W=\int_{B_{1}}|D^{2}u|^{p}\,dx
$$
and note that
$$
(\phi^{1/p}(\rho))'\leq\Big(\int_{B_{1}}|x|^{p}\,|D^{2}u(\rho x)|^{p}\,dx\Big)^{1/p}
$$
$$
=\rho^{-d/p-1}\Big(\int_{B_{ \rho}}|x|^{p}\,|D^{2}u( x)|\,dx\Big)^{1/p}\leq  \rho^{-d/p}W^{1/p}.
$$
Also
$$
\phi(\rho)=\rho^{-d}\int_{B_{ \rho}}|Du|^{p}\,dx
$$
and for $x\in B_{ \rho}$, $h= 1-\rho $,
  $i=1,...,d$, and $e_{i}$ defined as the $i$th
basis vector
$$
|D_{i}u(x)|\leq I_{1}(x)+I_{2}(x),
$$
where
$$
I_{1}(x)=|D_{i}u(x)-h^{-1}[u(x+e_{i}h)-u(x)]|
$$ 
$$
=h\big|\int_{0}^{1}\int_{0}^{1}sD^{2}_{i} 
u(x+tshe_{i})\,dtds\Big|,
$$
$$
I_{2}(x)=h^{-1}|u(x+e_{i}h)-u(x)|\leq h^{-1}[
|u(x+e_{i}h)|+|u(x)|].
$$
It follows that
$$
\phi(\rho)\leq N\rho^{-d}[hW+h^{-1}U],
$$
$$
\phi(1)\leq W\Big(\int_{\rho}^{1}\rho^{-d/p}\,d\rho\Big)^{p}
+N\rho^{-d}[hW+h^{-1}U].
$$
This yields \eqref{2.1.1} for $\varepsilon\leq1/2$ after setting
$\rho=1-\varepsilon=1-h$. The above restriction
on $\varepsilon$ is easily extended to 
$\varepsilon\leq1$.

While proving \eqref{2.1.2}, observe that
by Poincar\'e's inequality for $r\leq 1$
$$
\osc_{B_{r}}Du\leq Nr\bM(D^{2}u)\leq
N\bM(D^{2}u),
$$
whereas for $r>1$, in light of \eqref{2.1.1},
$$
\osc_{B_{r}}Du\leq 2\dashint_{B_{r}}
|Du|\,dx\leq N\bM(D^{2}u)+N\bM  u.
$$
Hence, for $\varepsilon=1$
$$
(Du)^{\sharp}\leq N\varepsilon\bM(D^{2}u)+
N\varepsilon^{-1}\bM  u.
$$
For any $\varepsilon>0$ this inequality is obtained by using scaling and then by minimizing with respect to $\varepsilon$
we come to \eqref{2.1.2}.

Estimate \eqref{2.1.3} follows
from \eqref{2.1.1} as $r\to\infty$.   \qed

The a priori estimate \eqref{9.4.3} leads
to an existence theorem.

\begin{theorem}
                      \label{theorem 12.22.1}
For any $p>1$, $\lambda>0$, and $f\in L_{p}$
there exists a unique $u\in W^{2}_{p}$
satisfying 
$$
\Delta u-\lambda u+f=0.
$$
Furthermore, there exists a constant
$N=N(d,p)$ such that for any $u\in W^{2}_{p}$
and $\lambda\geq0$
\begin{equation}
                          \label{12.22.1}
\|D^{2}u,\sqrt\lambda Du,\lambda u\|_{L_{p}}\leq
N\|\Delta u-\lambda u\|_{L_{p}}.
\end{equation}

\end{theorem}

Proof. We use Agmon's idea from \cite{Ag_62}.
 Consider the space
 $$
\bR^{d+1}=\{z=(x,y):y\in\bR,x\in\bR^{d}\}
 $$ 
and the function
 $$
u'(z)=u(x)\zeta(y)\cos(\mu y),
 $$
 where
 $\mu=\sqrt{\lambda}$ and $\zeta$ is a $C^{\infty}_{0}(\bR)$ function,
$\zeta\not\equiv0$.
Also introduce the operator 
 $$
\Delta'v(z)=\Delta_{x} v(z)+D_{y}^{2}v(z).
 $$

 In light of \eqref{9.4.3} applied to
$u'$ and $\Delta'$ we get
\begin{equation}
                                                     \label{06.6.27.9}
\|D^{2}_{x}u' ,D^{2}_{y}u' \|_{L_{p}(\bR^{d+1})}
\leq N\|\Delta'u' \|_{L_{p}(\bR^{d+1})}.
\end{equation}

It is an easy exercise to show that
$$
\int_{\bR}|\zeta(y)\cos(\mu y)|^{p}\,dy
$$
is bounded   away from zero for $\mu\in\bR$. 
It follows that

$$
\|u\|^{p}_{L_{p}}
=\mu^{-2p}\Big(\int_{\bR}|\zeta(y)\cos(\mu y)|^{p}\,dy\Big)^{-1}
\int_{\bR^{d+1}}\big|D^{2}_{y}u'(z) 
$$ 
$$
-u(x)[\zeta''(y)\cos(\mu y)-2\mu\zeta'(y)\sin(\mu y)
 ]\big|^{p}\,dz 
$$
$$
\leq N\mu^{-2p}\big(\|D^{2}_{y}u'\|^{p}_{L_{p}(\bR^{d+2})}
+(\mu^{p}+1)\|u \|^{p}_{L_{p} }\big). 
$$
This and \eqref{06.6.27.9} yield 
$$
\mu^{2} \|u \|_{L_{p}}\leq 
N\|\Delta'u' \|_{L_{p}(\bR^{d+1})}
+N(\mu+1)\|u \|_{L_{p}}.
$$
Since
$$
\Delta'u'=\zeta \cos(\mu y)
[\Delta u-\lambda u]+u[\zeta''\cos(\mu y)-2\mu\zeta'\sin(\mu y)],
$$
we have
$$
\|\Delta'u' \|_{L_{p}(\bR^{d+1})}\leq
N\|\Delta u-\lambda u\|_{L_{p}}
+  N(\mu+1)\|u\|_{L_{p}  },
$$
so that
$$
\lambda \|u \|_{L_{p} }\leq 
N\|\Delta u-\lambda u\|_{L_{p}}
+  N(\sqrt{\lambda}+1)\|u\|_{L_{p} }.
$$

Furthermore, \eqref{9.4.3} implies that
$$
\|D^{2}u\|_{L_{p}}\leq N\|\Delta u-\lambda u\|_{L_{p}}
+N\lambda \|u\|_{L_{p}},
$$
and by the multiplicative inequality
$$
\sqrt\lambda\|Du\|_{L_{p}}\leq N\sqrt\lambda\|D^{2}u\|^{1/2}_{L_{p}}\|u\|^{1/2}_{L_{p}}
\leq \|D^{2}u\|_{L_{p}}+N\lambda
\|u\|_{L_{p}}.
$$
Thus,
$$
\|D^{2}u\|_{L_{p}}+\sqrt\lambda\|Du\|_{L_{p}}+\lambda\|u\|_{L_{p}}
$$
$$
\leq N\|\Delta u-\lambda u\|_{L_{p}}+
N_{1}(\sqrt\lambda+1)\|u\|_{L_{p}}.
$$
For {\em fixed\/} $\lambda$ such that
$2N_{1}(\sqrt\lambda+1)\leq \lambda$ we obtain
\eqref{12.22.1}. After that for general
$\lambda$ estimate \eqref{12.22.1}
is obtained by using scaling.

This estimate, in particular, implies that  to
prove the first statement of the theorem
about the solvability of $\Delta u-\lambda u=f$,
it suffices to show that for $\lambda>0$ the set
$$
\{\Delta u-\lambda u:u\in W^{2}_{p}\}
$$
is dense in $L_{p}$. If it is not, then there exists a nonzero $g\in L_{q}$, $q=p/(p-1)$,
such that
$$
\int_{\bR^{d}}g(\Delta u-\lambda u)\,dx=0
$$
for any $u\in W^{2}_{p}$. Take a symmetric
$\zeta\in C^{\infty}_{0}$ and for $\varepsilon>0$ use the notation 
$$
\zeta_{\varepsilon}(x)=\varepsilon^{-d}
\zeta(x/\varepsilon),\quad
u^{(\varepsilon)}
=\zeta_{\varepsilon}*u.
$$ 
Then, for any $u\in C^{\infty}_{0}$, we have
$$
\int_{\bR^{d}}u(\Delta g^{(\varepsilon)}
-\lambda g^{(\varepsilon)})\,dx
=\int_{\bR^{d}}g^{(\varepsilon)}(\Delta u
-\lambda u)\,dx
$$
$$
=\int_{\bR^{d}}g(\Delta u^{(\varepsilon)}
-\lambda u^{(\varepsilon)})\,dx=0.
$$
It follows that 
$$
\Delta g^{(\varepsilon)}
-\lambda g^{(\varepsilon)}=0
$$
 and, since
$g^{(\varepsilon)}$ is obviously in $ W^{2}_{q} $,
$g^{(\varepsilon)}=0$ by \eqref{12.22.1}
and $g=0$ owing to $g^{(\varepsilon)}\to g$
in $L_{q}$ as $\varepsilon\downarrow 0$.
This proves the denseness and the theorem.
\qed

An important generalization of \eqref{10.31.1} is the following.

\begin{theorem}
                      \label{theorem 1.19.1}
For any $p>1$ there exists a constant $N=N(d,p)$
such that for any $u\in W^{2}_{p}$
\begin{equation}
                         \label{1.19.1}
(D^{2}u)^{\sharp}\leq N\big(\bM(|f|^{p})\big)^{1/p},
\end{equation} 
where $f=-\Delta u $.
\end{theorem}

Proof. We may assume that $u\in C^{\infty}_{0}$, so that $f\in C^{\infty}_{0}$ as well. The argument in the end of the proof of
Theorem \ref{theorem 9.6.1} along with the
obvious fact that
$$
\osc_{\{x\in \bR^{2}:|x|\leq 1\}}D^{2}u
=\lim_{n\to\infty}
\osc_{\{x\in \bR^{3}:|x|\leq 1\}}D^{2}(\zeta(\cdot/n)u) 
$$
allows us to concentrate
on $d\geq 3$. Then come back to the beginning of the subsection where we claimed that $u':=Rf$
satisfies $\Delta u'=-f$, so that 
$\Delta (u-u')=0$. Furthermore,
$u'$ is  infinitely differentiable and $u'$
goes to zero at infinity. By the maximum
principle   $u-u'$ is zero.
Thus, $u=Rf$ (which, by the way, is a well-known
fact which we have already used above).

 Next, let $\zeta\in C^{\infty}_{0}$
be such that $\zeta=1$ in $B_{2}$ and set
$v=R(f\zeta),w=R(f(1-\zeta))$, so that $u=v+w$.
The functions
$v,w$ are infinitely differentiable and they
and $v$ and $Dv$ decrease as $|x|^{2-d}$
and $|x|^{1-d}$, respectively, at infinity. By Remark \ref{remark 1.25.7}
$$
\|D^{2}v\|_{L_{1}(B_{1})}\leq N\|D^{2}v\|_{L_{p}(B_{1})}\leq N\|f\zeta\|_{L_{p}}\leq N\big(\bM(|f|^{p})(0)\big)^{1/p}.
$$
Furthermore, \eqref{9.4.2} implies that
$$
\dashint_{B_{1}}\dashint_{B_{1}}
 |D^{2}w(y)-D^{2}w(z)|\,dydz\leq N \bM (f(1-\zeta))(0)\leq N(\bM(f^{p})(0))^{1/p}.
$$
It follows that for $\rho=1$
$$
\dashint_{B_{\rho}}\dashint_{B_{\rho}}
 |D^{2}u(y)-D^{2}u(z)|\,dydz \leq N(\bM(f^{p})(0))^{1/p}.
$$
For any $\rho>0$ this is proved by using scaling and leads to \eqref{1.19.1}
at the origin, which is enough.   \qed

  \begin{remark}
                 \label{remark 1.16.1}
Next natural step is to add to
$\Delta$ lower order terms. Observe
that, if
\begin{equation}
            \label{10.31.5}
\|b^{i}D_{i}u\|_{L_{p}}+\|cu\|_{L_{p}}\leq (1/2)N^{-1}(p,d)\|D^{2}u\|_{L_{p}},
\end{equation}
where $N(p,d)$ is from \eqref{12.22.1},
then
  \begin{equation}
            \label{10.31.4}
 \|D^{2}u,\sqrt\lambda Du,\lambda u\|_{L_{p}}\leq 2N(p,d)\|\cL u-\lambda u\|_{L_{p}},
 \end{equation}
where 
$$
\cL u:=
\Delta u+b^{i}D_{i}u+cu ,
$$
and we can obtain the unique solvability
of $\cL u-\lambda u=f$
by using the method of continuity.
\end{remark}

Thus, we have the following result.
\begin{theorem}
              \label{theorem 1.10.1}
Let $p>1$. Assume that \eqref{10.31.5} holds
for any $u\in W^{2}_{p}$. Then
for any   $\lambda>0$  and $f\in L_{p}$
there exists a unique $u\in W^{2}_{p}$
satisfying $\cL u-\lambda u=f$.
Furthermore,   for any $u\in W^{2}_{p}$ 
and $\lambda\geq0$ estimate \eqref{10.31.4} holds.

\end{theorem}

One of the goals of the first subsection of the next section
is to give sufficient conditions for
\eqref{10.31.5} to hold.
We will see (Remark \ref{remark 9.5.3}) that it holds, in particular,
if $p<d$, $|b|=\kappa/|x|$,
$\kappa$ is sufficiently small, and
$c\equiv 0$, or $|c|=\kappa/|x|^{2}$
and $2p<d$.

\subsection{Elliptic equations with variable
main coefficients}
                     \label{section 4.5,1}

Another possible extension of the above
results is to consider operators with variable coefficients rather than just $\Delta$. 
\index{$B$@Sets!$\bS_{\delta}$}% 

Introduce $\bS_{\delta}$ as the set of symmetric   $d\times d$-matrices 
whose eigenvalues are in $[\delta,
\delta^{-1}]$, where $\delta\in(0,1]$.
It comes without saying that
linear nondegenerate change of coordinates
extends the above results to the operators
of the forms $a^{ij}D_{ij}$ as long as
{\em constant\/} 
$a\in\bS_{\delta}$. Of course, the above constants $N$ will change and start depending
on $\delta$. The possibility to make $a^{ij}$ variable and quite irregular
is provided by a
generalization of Theorem 
\ref{theorem 1.19.1}.

It is convenient to introduce the following notation. Set $|G|$ 
\index{$S$@Miscelenea!$"|G"|$}% 
to be the Lebesgue 
\index{$S$@Miscelenea!$\dashint\phantom{m}$}% 
\index{$N$@Norms!$\dashnorm f"|"|_{L_{p}}$}%
measure of $G\subset \bR^{d}$,
$$
\dashint_{G}f\,dx=\frac{1}{|G|}\int_{G}f\,dx,
\quad
\dashnorm f\| _{L_{p}(G)}=\Big(\dashint_{G}|f(x)|^{p}\,dx
\Big)^{1/p}.
$$

\begin{lemma} 
                  \label{lemma 2.19.2}
For any $p>1$, 
$\kappa>2$, $u\in C^{\infty}_{0}$, and $\rho>0$, at the origin we have
\begin{equation}
                         \label{2,19.1}
\osc_{B_{\rho}} D^{2}u
\leq N(d,p )\kappa^{d/p}\dashnorm \Delta u\|_{L_{p}(B_{\kappa\rho})}
+N(d,p)\kappa^{-1}  \bM |\Delta u|  .
\end{equation}   
\end{lemma}

Proof. As in the proof of Theorem 
\ref{theorem 1.19.1} we may concentrate
on $d\geq 3$. Scalings allow us to assume that $\rho=1$. Next, let $\zeta\in C^{\infty}_{0}$
be such that $\zeta=1$ in $B_{3\kappa/4}$
and $\zeta=0$ outside $B_{ \kappa}$  and set
$v=R(f\zeta),w=R(f(1-\zeta))$, so that $u=v+w$.
The functions
$v,w$ are infinitely differentiable and they
and $v$ and $Dv$ decrease as $|x|^{-(d-2)}$
and $|x|^{-(d-1)}$, respectively, at infinity. By Remark \ref{remark 1.25.7}
$$
\|D^{2}v\|_{L_{1}(B_{1})}\leq N\|D^{2}v\|_{L_{p}(B_{1})}\leq N\|f\zeta\|_{L_{p}}
$$
$$
\leq N\|f \|_{L_{p}(B_{\kappa})}=N\kappa^{d/p}\dashnorm f \|_{L_{p}(B_{\kappa})}.
$$
Furthermore, the computations before
Theorem \ref{theorem 9.6.1} show that,
for $|x|\leq 1$ we have 
$ |D^{3}w(x)|\leq N \kappa  ^{-1}\bM f(0)$ and therefore 
$$
\dashint_{B_{1}}\dashint_{B_{1}}
 |D^{2}w(y)-D^{2}w(z)|\,dydz \leq N \kappa  ^{-1}\bM f(0).
$$ 
This leads to \eqref{2,19.1} and proves the lemma.
\qed

Next, let $a(x)$ be an $\bS_{\delta}$-valued function on $\bR^{d}$. 
Set 
$$
\cL u=a^{ij}D_{ij}u
$$
 and
for $\rho>0$
define
$$
a^{\sharp}_{\rho}=\sup_{ r\leq \rho,x}\sup_{\,B\in\bB_{r}(x)}\osc_{B}a .  
$$
We follow the argument introduced in
\cite{Kr_07}. The same kind of results we are after (Theorem \ref{theorem 2.21.1})
was first obtained in \cite{CFL_91}
by using the Calder\'on-Zygmund singular
integral theory and a theorem
from \cite{CRW_76}
about commutators of singular integrals
and multiplications by a function.
 
The next result is derived by observing that Lemma \ref{lemma 2.19.2}
is valid also for $\bar\cL_{\rho} u:=\bar a^{ij}_{\rho}D_{ij}u$, where
$\bar a^{ij}_{\rho}$ are the averages of $a^{ij}$ over
$B_{\rho}$, and the fact that for
$f=\cL u$ we have 
$$
f=\bar\cL_{\rho} u
+[(a^{ij}-\bar a^{ij}_{\rho})D_{ij}u],
$$
 whereas
for $q>p$ by H\"older's inequality and the fact
that $a$ is bounded
$$
\dashnorm (a^{ij}-\bar a^{ij}_{\rho})D_{ij}u
\|_{L_{p}(B_{\rho})}\leq
N(d,\delta,p,q)a^{\sharp}_{\rho_{a}}
\big(\bM(|D^{2} u|^{q})\big)^{1/q} 
$$
for any $\rho>0$ provided that $u$
is supported in $B_{\rho_{a}}$.

\begin{lemma}  
                  \label{lemma 2.19.3}

For any $q> 1$,
$\kappa>2$,   $\rho_{a}>0$, and $u\in C^{\infty}_{0}(B_{\rho_{a}})$   we have
$$
( D^{2}u)^{\sharp}
\leq N(d,\delta,q,\kappa)\big(\bM(|\cL u|^{q})\big)^{1/q} 
$$
\begin{equation}
                         \label{2.20.1}
+\big(N(d,\delta,q,\kappa)a^{\sharp}_{\rho_{a}}+N(d,\delta, q)\kappa^{-1}\big)\big(\bM(|D^{2} u|^{q})\big)^{1/q}. 
\end{equation}
\end{lemma}

In light of
Theorems \ref{theorem 06.6.9.1} and \ref{theorem 7.7.1} this leads to the following
for any $p>1$ ($p>q>1$).
\begin{lemma}
                    \label{lemma 2.20.2}
For any $p> 1$,
$\kappa>2$,   $\rho_{a}>0$, and $u\in C^{\infty}_{0}(B_{\rho_{a}})$  we have
$$
\|D^{2}u\|_{L_{p}}\leq
  N_{0}(d,\delta,p,\kappa) \|\cL u\|_{L_{p}} 
$$
\begin{equation}
                         \label{2.20.3}
+\big(N_{1}(d,\delta,p,\kappa)a^{\sharp}_{\rho_{a}}+N_{2}(d,\delta, p)\kappa^{-1}\big)\|D^{2}u\|_{L_{p}}. 
\end{equation}
\end{lemma}

\begin{corollary}
                 \label{corollary 2.20.3}
Let $p>1$ and assume that
\begin{equation}
                         \label{2.20.4}
 a^{\sharp}_{\rho_{a}}\leq (1/4)
N_{1}^{-1}\big(d,\delta,p,4
N_{2} (d,\delta, p)\big).
\end{equation}
Then for any $u\in C^{\infty}_{0}(B_{\rho_{a}})$  we have
$$
\|D^{2}u\|_{L_{p}}\leq
  2N_{0}(d,\delta,p,4
N_{2} (d,\delta, p)) \|\cL u\|_{L_{p}}.
$$
\end{corollary}
From this ``local'' a priori estimate
we have a global one.
\begin{lemma}
                    \label{lemma 2.20.5}
Let $p> 1$ and assume \eqref{2.20.4}.  
Then there exist $N_{1}=N_{1}(d,\delta,p)$ and $N_{2}=N_{2}(d,\delta,p,\rho_{a})$ such that for any $u\in C^{\infty}_{0} $ we have
\begin{equation}
                         \label{2.20.6}
\|D^{2}u\|_{L_{p}}\leq
  N_{1} \|\cL u\|_{L_{p}}+N_{2}\|u\|_{L_{p}}.
\end{equation}
\end{lemma}

Proof. Take a nonnegative $\zeta\in C^{\infty}_{0}(B_{\rho_{a}})$ with $\zeta^{p}$ having unit integral and for
$y,x\in\bR^{d}$ define $\zeta_{y}(x)
=\zeta(x-y), u_{y}=u\zeta_{y}$.
Then 
$$
\|D^{2}u_{y}\|_{L_{p}}^{p}\leq N\|\cL u_{y}\|_{L_{p}}^{p}\leq N\|\zeta_{y}\cL u\|_{L_{p}}
$$
$$
+N\big(\|(D\zeta_{y})Du\|_{L_{p}}^{p}+
\|(D^{2}\zeta_{y}) u\|_{L_{p}}^{p}\big).
$$
Furthermore,
$$
2^{ p}\|D^{2}u_{y}\|_{L_{p}}^{p}\geq 
\|\zeta_{y}D^{2}u \|_{L_{p}}^{p}-
N\big(\|(D\zeta_{y})Du\|_{L_{p}}^{p}+
\|(D^{2}\zeta_{y}) u\|_{L_{p}}^{p}\big).
$$
Hence,
$$
\|\zeta_{y}D^{2}u \|_{L_{p}}^{p} \leq N\|\zeta_{y}\cL u\|_{L_{p}}+N\big(\|(D\zeta_{y})Du\|_{L_{p}}^{p}+
\|(D^{2}\zeta_{y}) u\|_{L_{p}}^{p}\big), 
$$
and it only remains to integrate through this
inequality over $y\in\bR^{d}$.  \qed

Next, we repeat part of the proof of
Theorem \ref{theorem 12.22.1} and see that
for $\lambda\geq0$ and $u\in C^{\infty}_{0}$
$$
\|D^{2}u,\sqrt\lambda Du,\lambda u\|_{L_{p}}\leq
N\|\cL u-\lambda u\|_{L_{p}}+
N\|u\|_{L_{p}}. 
$$
For $\lambda$ large enough the last term can be dropped and for $u\in C^{\infty}_{0}$ we arrive at the first assertion in the following theorem.

\begin{theorem}
                 \label{theorem 2.21.1}
Let $p> 1$ and assume \eqref{2.20.4}.  
Then there exist $N =N (d,\delta,p)$  and
$\lambda_{0}=\lambda_{0}(d,\delta,p,\rho_{a})>0$ such that for any $u\in W^{2}_{p}$
and $\lambda\geq\lambda_{0}$ 
\begin{equation}
                            \label{2.21.2}
\|D^{2}u,\sqrt\lambda Du,\lambda u\|_{L_{p}}\leq
N\|\cL u-\lambda u\|_{L_{p}}. 
\end{equation}
Furthermore, for any $f\in L_{p}$
and $\lambda\geq \lambda_{0}$
there exists a unique $u\in W^{2}_{p}$
satisfying $\cL u-\lambda u+f=0$.
\end{theorem}

 A priori estimate \eqref{2.21.2} extends
from $C^{\infty}_{0}$ to $W^{2}_{p}$
due to the denseness of the former in the latter and the solvability result is obtained by the method of continuity.

\begin{remark}
                         \label{remark 4.3,3}
In the literature quite often Theorem
\ref{theorem 2.21.1} is proved for $a$
that is uniformly continuous. 
Then it is instructive to emphasize that for
general $a$ satisfying \eqref{2.20.4}
the quantity $\sup_{x}|a-a'|$, where $a'$
is uniformly continuous, is greater than
$(1/2)\sup_{x,y}|a(x)-a(y)|$, see, for instance, Remark 
\ref{remark 2.21.1}.
\end{remark}

  \mysection{Morrey-Sobolev spaces}
                \label{chapter 1.7.1}
  
 The first subsection of this section is about
the  solvability in the {\em Sobolev\/} spaces of the Laplace equation
with singular drift. The remaining ones are
about the solvability in Morrey-Sobolev spaces.

\subsection{Adams's theorem}

 This is a generalization to any $p>1$ of an earlier result
 of Ch. Fefferman, the latter being also later extended by D. Adams and
 F. Chiarenza-M. Frasca with  different proofs.

    Let $d\geq 1$, $0< \alpha<d$. Introduce the Riesz potentials by 
$$
R_{\alpha}f(x)=\int_{\bR^{d}}\frac{f(y)}{|x-y|^{d-\alpha}}
\,dy. 
$$
We will often use representations of such integrals in polar coordinates.
For instance, for bounded $f$ with compact support
  integrating by parts we find
$$
\int_{0}^{\infty}\frac{1}{r}\Big(
r^{\alpha}\dashint_{B_{r} }f(x)\,dx
\Big)\,dr=N\int_{0}^{\infty} \Big(
 \int_{B_{r} }f(x)\,dx
\Big)\,dr^{\alpha-d}
$$
$$
=N\int_{0}^{\infty} r^{\alpha-d}
 \int_{\partial B_{r} }f(x)\,d\sigma_{r}\,dr
 =NR_{\alpha}f(0),
$$
where $\sigma_{r}$ is the element of the surface measure on $\partial B_{r}$.

\begin{remark}
                                 \label{remark 9.5.1}
 If $u$ is a smooth function with compact support,
 then for any $x$ with $|x|=1$ we have
 $$
 u(0)=-\int_{0}^{\infty}x^{i}(D_{i}u)(rx)\,dr,\quad
 |u(0)|\leq \int_{0}^{\infty}|D_{i}u|(rx)\,dr.
 $$
 By integrating this with respect to $x$
 over the unit sphere we find that at the origin
 \begin{equation}
                    \label{9.5.5}
 |u|\leq N(d)R_{1}|Du|.
 \end{equation}
 Obviously, \eqref{9.5.5} holds 
 at any other point as well.

Estimate \eqref{9.5.5} also follows from the well-known formula (proved, for instance,
by using the Fourier transform)
$$
u=N(d)
\sum_{i=1}^{d}R^{i}D_{i}u,\quad
R^{i}f(x)=\int_{\bR^{d}}\frac{y^{i}-x^{i}}{|y-x|^{d}}f(y)\,dy,
$$
which allows us to treat $|u(x)-u(0)|$.
Indeed, this representation implies
that
$$
|u(x)-u(0)|\leq N\int_{\bR^{d}}
\Big|\frac{y -x }{|y-x|^{d}}
-\frac{y }{|y |^{d}}\Big||Du(y)|\,dy.
$$
We split  the domain of integration
according to $|y|\leq 2|x|$ and $|y|>
2|x|$ and observe that for $|y|>2|x|$,
 $\varepsilon\in(0,1),|x|\leq 1$, we have
$$
\Big|\frac{y -x }{|y-x|^{d}}
-\frac{y }{|y |^{d}}\Big|\leq N|x|\frac{1}{|y|^{d }}\leq N|x|^{\varepsilon}
\frac{1}{|y|^{d-1+\varepsilon }},
$$
which follows after substitution
$|y|=|x|/s$ from
$$
\Big|\frac{l_{1} -sl_{2}}
{|l_{1} -sl_{2}|^{d}}- l_{1} \Big|\le Ns,\quad
|s|\leq 1/2,|l_{i}|=1.
$$
Then for $f=|Du|$, $\varepsilon\in(0,1),|x|\leq 1$ we find
$$
|u(x)-u(0)|\leq N\Big(R_{1}(I_{B_{2|x|}}f)(x)+R_{1}(I_{B_{2|x|}}f)(0)
$$
\begin{equation}
                      \label{1.5.6}
+|x|^{\varepsilon}R_{1-\varepsilon}
(I_{B^{c}_{2|x|}}f)(0)\Big).
\end{equation}
This result will be used in the proof 
of Theorem \ref{theorem 1.5.1}.
\end{remark}

The next result is   Proposition
3.3 of \cite{Ad_75} with a similar proof. 

\begin{lemma}
                       \label{lemma 9.5.1}
Let  $f\geq 0$ and assume that $R_{\alpha}f$ is locally
summable. Then
\begin{equation}
                                \label{5,30.1}
R^{\sharp}_{\alpha}f\leq N(d,\alpha)\bM_{\alpha}f.
\end{equation}

\end{lemma}

Proof. By using scalings and shifts we see that
it suffices to prove that
$$
I (f) :=\int_{B_{1}}\int_{B_{1}}|R_{\alpha}f(x)-R_{\alpha}f(y)
|\,dxdy\leq N\sup_{r\geq1}r^{\alpha}\dashint_{B_{r}}f\,dx=:J.
$$
Set
$f=g+h$, $g=fI_{B_{2}}$, $h=fI_{B^{c}_{2}}$
and observe that
$$
I (g) \leq 2\int_{B_{1}}R_{\alpha}g(x)\,dx=\int_{\bR^{d}}
g \phi\,dy\leq N\int_{B_{2}}f\,dy\leq NJ,
$$
where
$$
\phi(y)=\int_{B_{1}}\frac{2}{|x-y|^{d-\alpha}}\,dx.
$$
Furthermore,
$$
\int_{B_{1}}\int_{B_{1}}|R_{\alpha}h(x)-R_{\alpha}h(y)
|\,dxdy\leq N\sup_{B_{1}}|DR_{\alpha}h|,
$$
where for $|x|<1$
$$
|DR_{\alpha}h( x)|\leq NR_{\alpha-1}h(x)
=N\int_{0}^{\infty}\frac{1}{r^{2}}\Big(
r^{\alpha}\dashint_{B_{r} }h(x+y)\,dy
\Big)\,dr
$$
$$
\leq N\int_{1}^{\infty}\frac{1}{r^{2}}\Big(
r^{\alpha}\dashint_{B_{r} }f(y)\,dy
\Big)\,dr\leq NJ.
$$
This yields the desired estimate of $J$ and proves the lemma.
\qed

The following is extended to the case
of $L_{p}$-spaces with Muckenhoupt
weights in Theorem \ref{theorem 5.28,1}.

\begin{theorem} [D. Adams]
                            \label{theorem 10.7.1} 
Assume   $1<r<p $, $b\geq0$, $f\in L_{r}$, and
\index{Adams's theorem}
(the notation will be explained
\index{$S$@Miscelenea!$"|"b"|"_{\dot E_{p,\alpha}}$}% 
later)
\begin{equation}
                                \label{9.6.1}
\|b\|_{\dot E_{p,\alpha}}:=\sup_{\rho>0}\rho ^{\alpha}\sup_{B\in\bB_{\rho}}  
\Big(\dashint_{B}b^{p}\,dx\Big)^{1/p}\leq A . 
\end{equation}
Then
\begin{equation}
                                \label{9.25.1}
I:=\int_{\bR^{d}}b^{r}|R_{\alpha}f|^{r}\,dx\le N(\alpha,d,r,p)A^{r}
 \int_{\bR^{d}} |f|^{r}\,dx.
\end{equation}
\end{theorem}

\begin{corollary}
                       \label{corollary 10.7.1}
If $u\in C_{0}^{\infty}(\bR^{d})$ and $\alpha=1$
(and $d\geq2$), then
 owing to \eqref{9.5.5} (the Chiarenza-Frasca result)
$$
\int_{\bR^{d}}b^{r}|u|^{r}\,dx\le  NA^{r}
 \int_{\bR^{d}}  |Du|^{r}\,dx,
$$
\begin{equation}
                                \label{9.6.2} 
\int_{\bR^{d}}b^{r}|Du|^{r}\,dx\le  NA^{r}
 \int_{\bR^{d}}  |D^{2}u|^{r}\,dx,
\end{equation}
and if $\alpha=2$ (and $d\geq3$), then
\begin{equation}
                                \label{9.6.3}
\int_{\bR^{d}}b^{r}|u|^{r}\,dx\le  NA^{r}
 \int_{\bR^{d}} \big|\Delta u\big|^{r}\,dx.
\end{equation}
Generally (for $0<\alpha<d$),
$$
\int_{\bR^{d}}b^{r}|u|^{r}\,dx\le  NA^{r}
 \int_{\bR^{d}} \big|(-\Delta)^{\alpha/2}u\big|^{r}\,dx.
$$
 \end{corollary}
Indeed, $f:=(-\Delta)^{\alpha/2}u$ satisfies
$f\in L_{r}$ and $NR_{\alpha}f=u$.

 \begin{remark} 
                                    \label{remark 9.6.6}
 
 One may ask if \eqref{9.6.2} and \eqref{9.6.3}
  hold true for $d=1$ and $d=2$, respectively.
  The answer is ``yes'' by a very trivial reason that
 if $\alpha=d=1$ and $p>1$  
  and $b\not\equiv0$, then  $A=\infty$.
  In particular, to have a meaningful result,
   one can take $\alpha=1$ only if
  $1<p\leq d$ and $\alpha=2$ only if $1<p\leq d/2$.
 \end{remark}

\begin{remark}
                                 \label{remark 9.5.3}
The reader will easily check that
$b=1/|x|$ satisfies  the assumption
of Theorem \ref{theorem 10.7.1} with $A<\infty$ for $\alpha=1$
and any $p\in[1,d)$,
and $b=1/|x|^{2}$ satisfies it for $\alpha=2$ ($d\geq 3$)
and any $p\in[1,d/2)$. As
an application we have the following Hardy's inequalities:
$$
\int_{\bR^{d}}\frac{1}{|x|^{r}}|u|^{r}\,dx\leq
N\int_{\bR^{d}}  |Du|^{r}\,dx,\quad 1<r<d,
$$
$$
\int_{\bR^{d}}\frac{1}{|x|^{2 r}}|u|^{r}\,dx\leq
N\int_{\bR^{d}}  |\Delta u|^{r}\,dx,\quad 1<r<d/2.
$$
\end{remark}

{\bf Proof of Theorem \ref{theorem 10.7.1}}.
In light of Remark \ref{remark 9.6.6} we may assume that
$\alpha p\leq d$. In that case
we may also assume that $b$ and $f$ are bounded and have
compact support. This guarantees that $I<\infty$.
Then   set
$u=R_{\alpha}f$ and write   
$$
I=\int_{\bR^{d}}\big(b^{r}u^{r-1}\big) R_{\alpha}f
\,dx=\int_{\bR^{d}}R_{\alpha} \big(b^{r}u^{r-1}\big)  f
\,dx\leq \|f\|_{L_{r}}\big\|
R_{\alpha} \big(b^{r}u^{r-1}\big)\big\|_{L_{r'}},
$$
where $r'=r/(r-1)$. Observe that the last norm is finite because
$R_{\alpha} \big(b^{r}u^{r-1}\big)$ is locally bounded and at infinity it decreases as $1/|x|^{d-\alpha}$, whereas
$$
\int_{|x|\geq 1}\frac{1}{|x|^{(d-\alpha)r'}}\,dx
<\infty
$$
owing to $\alpha r<\alpha p\leq d$.

Now by the Fefferman-Stein theorem   
$$
\big\|
R_{\alpha} \big(b^{r}u^{r-1}\big)\big\|_{L_{r'}}
\leq N\big\|
R_{\alpha}^{\sharp} \big(b^{r}u^{r-1}\big)\big\|_{L_{r'}}.
$$

By Lemma \ref{lemma 9.5.1}
$$
\|
R_{\alpha}^{\sharp}\big(b^{r}u^{r-1}\big)\big\|_{L_{r'}}
\leq N\|
\bM_{\alpha}\big(b^{r}u^{r-1}\big)\big\|_{L_{r'}},
$$
where by H\"older's inequality
$$
\bM_{\alpha}\big(b^{r}u^{r-1}\big)
=\bM_{\alpha} \big(b(b^{r-1}u^{r-1})\big)\leq
\|b\|_{\dot E_{p,\alpha}} [\bM\big((bu)^{s}\big)]^{1/p'}
$$
with $s=(r-1)p'$ and $p'=p/(p-1)$.
Since $r<p$, we have $p'<r'$ and by the Hardy-Littlewood theorem
$$
\|
\bM_{\alpha} \big(b^{r}u^{r-1}\big)\big\|_{L_{r'}}
\leq N\|b\|_{\dot E_{p,\alpha}}\Big(\int_{\bR^{d}}
\big[(bu)^{s}\big]^{r'/p'}\,dx\Big)^{(r-1)/r}
$$
$$
=N\|b\|_{\dot E_{p,\alpha}}I^{(r-1)/r}.
$$
Thus, 
$$
I\leq N\|f\|_{L_{r}}\|b\|_{\dot E_{p,\alpha}}I^{(r-1)/r}.
$$
This, obviously, proves the theorem. \qed

As a direct consequence of
Corollary \ref{corollary 10.7.1}
and Theorem \ref{theorem 1.10.1} we have the following.

\begin{corollary}
               \label{corollary 1.16.1}
Let $1<p<q <\infty$. Then
there is a constant $N=N(d,p,q)$
such that, if 
$$
N\|b\|_{\dot E_{q,1}}
\leq 1,
$$ then
for any   $\lambda>0$  and $f\in L_{p}$
there exists a unique $u\in W^{2}_{p}$
satisfying $\cL u-\lambda u=f$, where 
$$
\cL u=\Delta u+b^{i}D_{i}u.
$$
Furthermore,   for any $u\in W^{2}_{p}$ estimate \eqref{10.31.4} holds.
\end{corollary}

 \begin{remark}
               \label{remark 1.17.1}
If $\|b\|_{\dot E_{q,1}}$ is not small enough, the assertion of Corollary
\ref{corollary 1.16.1} may fail for any
$\lambda>0$. For instance, if 
$$
b=-(d-1)x/|x|^{2},\quad 1< p<q<d,\quad  u(x)=\exp(-\sqrt\lambda|x|),
$$
 then as is easy
to check $\|b\|_{\dot E_{q,1}}<\infty$
(see, however, Remark \ref{remark 4.3,7}), 
 $u\in W^{2}_{p}$ and $\cL u-\lambda u=0$, so that estimate \eqref{10.31.4} fails.

\end{remark}

\begin{remark}
               \label{remark 1.17.2}
From the point of view of the theory
of elliptic equations the most desirable version of \eqref{9.6.2}
would be, of course,
\begin{equation}
                     \label{1.17.1}
\|b^{i}D_{i}u\|_{L_{p}}
\leq \varepsilon \|D^{2}u\|_{L_{p}}
+N(\varepsilon)\|u\|_{L_{p}}
\end{equation}
for a $p\in(1,d)$ any $u\in W^{2}_{p}$ and any $\varepsilon\in(0,1]$ with
$N(\varepsilon)$ independent of $u$
(may depend on $b$).

Indeed, if true for $\varepsilon$
small enough, this estimate would lead
to the unique solvability
in $W^{2}_{p}$ of 
$$
\cL u-\lambda u=f\in L_{p} 
$$
 for $\lambda$, large enough
to absorb $N(\varepsilon)\|u\|_{L_{p}}$
into the left-hand side of \eqref{10.31.4}. However, generally, 
there is no unique solvability for
$b$ such that $\|b\|_{\dot E_{q,1}}<\infty$ 
and, therefore, estimate \eqref{1.17.1} is 
impossible if $\varepsilon$ is too small.
For instance, \eqref{1.17.1} is 
impossible if $\varepsilon$ is too small if $|b|=|x|^{-1}$. 
 The example in Remark
\ref{remark 1.17.1} shows this.
\end{remark}

A useful observation is that \eqref{9.6.1}
is satisfied
with finite $A$ in case $b\in L_{p\alpha}$. 
Also, H\"older's inequality shows that 
the assumption that $A$ in \eqref{9.6.1} is finite
(and $p\alpha\leq d$)
becomes stronger if  we replace $p $ with any
$q\geq p $ and, if $p\alpha=d$,  it is satisfied
 in case $b\in L_{d}$.

It is instructive to know
that   it might happen
that $\|b\|_{\dot E_{p,\alpha}}<\infty$ but $b\not\in L_{p+\varepsilon,\loc}$,
no matter how small $\varepsilon>0$ is.

\begin{example}
                        \label{example 9.25.1}

 Take $p \in[d-1,d)$ and take $r_{n}>0$, $n=1,2,...$, such that
the  sum of $\rho_{n}:=r_{n}^{d-p}$ is $1/2$. Let $e_{1}$ be the first
basis vector  and set $b(x)=|x|^{-1}
I_{|x|<1}$, $x_{0}=1$,
$$
x_{n}=1-  2\sum_{1}^{n}r_{i}^{d-p},\quad 
c_{n}=(1/2)(x_{n}+x_{n-1})
$$
$$
  b_{n}(x)=r_{n}^{-1}b\big(r_{n}^{-1}
(x-c_{n}e_{1})\big),\quad b=\sum _{1}^{\infty}b_{n}.
$$
Since $r_{n}\leq 1$ and $d-p\leq 1$,  the supports of $b_{n}$'s are disjoint and
for $q>0$
$$
\int_{B_{1}}b^{q}\,dx=\sum _{1}^{\infty}\int_{\bR^{d}}b_{n}^{q}\,dx=N(d,p)\sum_{1}^{\infty}r_{n}^{d-q}.
$$
According to this we take the $r_{n}$'s so that
the last sum diverges for any $q>p$.
Then observe that for any $n\geq 1$ and any ball $B$
of radius $\rho$
$$
 \int_{B }  b_{n} ^{p}dx \leq N(d) \rho^{d-p} .
$$
Also, if the intersection of $B$ with $\bigcup B_{r_{n}}(c_{n})$
is nonempty, the intersection
 consists of some $B_{r_{i}}(c_{i})$, $i=i_{0},...,i_{1}$, and $B\cap B_{r_{i_{0}-1}}(c_{i_{0}-1})$ if $i_{0}\geq 2  $ and 
$B\cap B_{r_{i_{1}+1}}(c_{i_{1}+1})$.
In this situation $ c_{i_{0}}-c_{i_{1}}
\leq 2\rho$, $c_{i_{0}}\leq x_{i_{0}},
c_{i_{1}}\geq x_{i_{i-1}}$ and
$$
 \sum_{i=i_{0} +1}^{i_{1}-1 }
r^{d-p}_{i} \leq    \rho.
$$
Therefore,
$$
 \int_{B }  b  ^{p}\,dx\leq N(d)\sum_{i=i_{0}+1}^{i_{1}-1}r_{i}^{d-p}
$$
$$
+\int_{B }  [I_{i_{0}\geq2} b_{i_{0}-1}  ^{p} 
+b^{p}_{i_{0}}+ b_{i_{1} }  ^{p}+ b_{i_{1}+1}  ^{p}]dx \leq N(d)(\rho
+\rho^{d-p}),
$$
where the last term is less than $N(d)\rho^{d-p}$
for $\rho\leq 1$. This domination
also holds for $\rho>1$,
since $b\in L_{p}$ and $d\geq p$.
Thus, $b\in \dot E_{p,1}$ but $b\not\in
L_{q}$ for $q>p$.
\end{example}

 \subsection{Definition and simplest properties
of Morrey-Sobolev spa\-ces}
                  \label{section 2.16.1}

Take $p\in(1,\infty)$,
   $\beta\geq 0$, $r\in(0,\infty)$ and domain $\cO\subset
\bR^{d}$ and introduce the
Morrey space $E_{p,\beta;r}(\cO) $
as the set of $g\in  L_{p,\loc}(\cO)$ 
\index{$A$@Sets of functions!$E_{p,\beta;r}(\cO)$}%
\index{$N$@Norms!$"|"|f"|"|_{E_{p,\beta;r}(\cO)}$}%
such that  
\begin{equation}
                             \label{8.11.02}
\|g\|_{E_{p,\beta;r}(\cO) }:=
\sup_{\rho\leq r,x\in \cO }\rho^{\beta}
\dashnorm g I_{\cO} \|_{ L_{p}(B_{\rho}(x) )} <\infty .
\end{equation}  
Define the Morrey-Sobolev space
$E^{2}_{p,\beta;r}(\cO)$ by
$$
E^{2}_{p,\beta;r}(\cO) =\{u:u,Du,D^{2}u 
\in E_{p,\beta;r}(\cO) \},
$$
where $Du,D^{2}u$ are Sobolev derivatives,
and 
provide $E^{2}_{p,\beta;r}(\cO) $ with 
\index{$A$@Sets of functions!$E^{2}_{p,\beta}$}%
\index{$A$@Sets of functions!$E^{2}_{p,\beta;r}(\cO)$}%
\index{$A$@Sets of functions!$E^{1}_{p,\beta;r}(\cO)$}%
\index{$A$@Sets of functions!$E^{1}_{p,\beta}$}%
an obvious norm. 
Similarly $E^{1}_{p,\beta;r}(\cO)$
is introduced. If $\cO=\bR^{d}$ or $r=1$
we drop ``$(\cO)$'' or ``$;r$'' in the above notation.
For instance, $E^{2}_{p,\beta}=E^{2}_{p,\beta;1}(\bR^{d}).$

Let $r,s>0$.
Since $u\in E_{p,\beta;rs}$ if and only if $u(r\cdot)\in E_{p,\beta;s}$ and
$$
\|u\|_{E_{p,\beta;rs}}=r^{\beta}
\|u(r\cdot)\|_{E_{p,\beta;s}},
$$
we mainly concentrate on $E_{p,\beta}$.

\begin{remark}
                        \label{remark 4.3,7}
It turns out that $|x|^{-1}$ belongs
to $E_{p,1}$ for any $p<d$. Indeed, if
$|x|\leq 2\rho$, then
$$
\dashint_{B_{\rho}(x)}|y|^{-p}\,dy\leq N(d)
\dashint_{B_{2\rho} }|y|^{-p}\,dy=N^{p}_{0}(d,p)(2\rho)^{-p},
$$
and if
$|x|> 2\rho$, then $|y|^{-p}\leq N(p)\rho^{-p}$
on $B_{\rho}(x)$ and again
$$
\dashint_{B_{\rho}(x)}|y|^{-p}\,dy\leq N(p)\rho^{-p}.
$$
Also, obviously, if $x'\ne x''$ and $u_{x}(y)
=|x-y|^{-1}$, then 
$$
\|u_{x'}-u_{x''}\|_{E_{p,1}}\geq N_{0}(d,p),
$$
so that $E_{p,1}$ is a nonseparable Banach
space. The same is true for $E_{p,\beta}$ 
if $\beta<d/p$. This is similar to $C^{\alpha}$
H\"older spaces and, as in that case, the set
of smooth functions is not dense in $E_{p,1}$.
However, for $\beta\geq d/p$ the picture is quite
different, see Remark \ref{remark 7,2.1}.
\end{remark}
 
Observe that in \cite{Kr_22} the $E_{p,\beta}$-norm
is defined differently form \eqref{8.11.02}
without the restriction on $\rho$, just $\rho<\infty$.
In connection with this note that if $\beta\leq
d/p$ and the support of $g$ is in some $B\in\bB_{1}$,
then the $E_{p,\beta}$-norm of $g$ is the same:
taken from \eqref{8.11.02} or from \cite{Kr_22}.
One more useful observation is that
$$
\|g\|_{E_{p,\beta} }=\sup_{B\in \bB_{1}}
\|gI_{B}\|_{E_{p,\beta} }.
$$

 \begin{remark} 
                        \label{remark 7,2.1}
If $\beta\geq d/p$ then
$$
\|f\|_{E_{p,\beta}}=\sup_{B\in \bB_{1}}\|f\|_{L_{p}(B)}.
$$
Generally,  H\"older's inequality easily shows that
$$
\|f\|_{E_{p,\beta}}\leq N
\sup_{B\in \bB_{1}}\|f\|_{L_{r}(B)},
$$
where $N$ is independent of $f$ and 
$r=p\vee(d/\beta)$. 
\end{remark}

\begin{remark}
If $u$ vanishes outside some $B\in\bB_{1}$  and  if $d/p   \geq\beta $, one easily sees that
for any $r\geq 1$ and $B'\in \bB_{r}$
$$
r^{\beta }\dashnorm u\|_{L_{p}(B')}\leq N(d)
\|u\|_{L_{p}(B)}.
$$
\end{remark}

\begin{remark}
                \label{remark 10.3.1}
Quite often it is useful to know
that, if $f(x)\geq 0$ is Borel nonnegative, $r\in(0,\infty),\mu
\in(0,1)$, then
\begin{equation}
                     \label{6,30.60}
\dashint_{B_{r}}f\,dx
\leq  (1+\mu)^{d}\sup_{B\in\bB_{\mu r}}\dashint_{B}f\,dx.
\end{equation}
This inequality is obtained in the following way. Let 
$\zeta= |B_{\mu r}|^{-1} 
 I_{B_{\mu r}}$. Then
$$
\dashint_{B_{r}}f\,dx= |B_{r}|^{-1}
\int_{B_{r}}\int_{B_{r(1+\mu)}}
f(x)\zeta(x-y)\,dy\,dx
$$
$$
\leq |B_{r}|^{-1}\int_{B_{r(1+\mu)}}\int_{\bR^{d}}
f(x)\zeta(x-y)\,dx\,dy,
$$
which immediately leads to \eqref{6,30.60}.

In particular, spaces $E_{p,\beta;r}$
for different $r$ coincide and, if $r\leq 1$,
$$
\|u\|_{E_{p,\beta}}\leq 2^{d/p}r^{-\beta}\|u\|_{E_{p,\beta;r}}.
$$
\end{remark} 

\begin{remark}
              \label{remark 12.26.2}
If $R\leq 1 $ and $\beta\leq d/p$, then
$$
\|f\|_{E_{p,\beta}(B_{R})}
=  \sup_{r\leq  R,x\in B_{R}}r^{\beta}
\dashnorm fI_{B_{R}}\|_{L_{p}(B_{r}(x))}.
$$
This is true because, if $r\in[ R,1]$,
then (with $N(d,p)=|B_{1}|^{-1/p}$) 
$$
r^{\beta}
\dashnorm fI_{B_{R}}\|_{L_{p}(B_{r}(x))}\leq N(d,p)r^{\beta-d/p}
\|fI_{B_{R}}\|_{L_{p} }
$$
$$
\leq N(d,p) R ^{\beta-d/p}
\|fI_{B_{R}}\|_{L_{p}(B_{R})}
=  R^{\beta}
\dashnorm fI_{B_{R}}\|_{L_{p}(B_{R})}
$$
$$
\leq   \sup_{\rho\leq  R} \rho^{\beta}
\|fI_{B_{R}}\|_{L_{p}(B_{\rho}(0))}.
$$
\end{remark}

\begin{remark}
             \label{remark 12.26.4}
For any ball $B$ (note the range of $x$)
$$
\|f\|_{E_{p,\beta}(B)}=\sup_{x\in\bR^{d},\rho\leq 1 }\rho^{\beta}\dashnorm fI_{B}\|_{L_{p}(
B_{\rho}(x))}.
$$
Indeed, if $x\not\in B$ on can move
$B_{\rho}(x)$ to the  center of $B$ until
$x$ falls on $\partial B$ all the way
increasing (by inclusion) $B\cap B_{\rho}(x)$ and
$\dashnorm fI_{B}\|_{L_{p}(
B_{\rho}(x))}$.
\end{remark}

For functions $f(x)$ and $\varepsilon>0$ we define
\begin{equation}
                        \label{7,5.1}
f^{(\varepsilon)}=f*\zeta_{\varepsilon},
\end{equation}
 where
  $\zeta_{\varepsilon}=\varepsilon^{-d}\zeta
(x/\varepsilon)$ with a nonnegative
$\zeta\in C^{\infty}_{0} $ which has unit integral and support
in $B_{1}$. Observe that
owing to Minkowski's inequality  
$$
\|f^{(\varepsilon)}\|_{E_{p,\beta}}
\leq \|f \|_{E_{p,\beta}}
$$
for any $f\in E_{p,\beta}$.
Maximal function boundedness is
another useful property of Morrey spaces. One of the first results of that kind
 appeared in \cite{CF_88}.
\begin{lemma}
                 \label{lemma 7.5.1}
Set $M^{[c]}f=\sup_{\varepsilon\in(0,c]}|f|^{(\varepsilon)}$, then 
$$
\|M^{[1]}f\|_{E_{p,\beta}}\leq N
\|f\|_{E_{p,\beta}}
$$
 for any $f\geq0$, where $N=N(d,p,\beta)$.
\end{lemma}

Proof. We have to prove that, for any
$r\leq 1$ and $B\in \bB_{r}$,
\begin{equation}
                      \label{7,5.2}
r^{\beta}\dashnorm M^{[1]}f\|_{L_{p}(B)}
\leq N\|f\|_{L_{p,\beta}}.
\end{equation}
We may assume that $B=B_{r}$. In that
case first note that by the
Hardy-Littlewood maximal function theorem  
$$
r^{\beta}\dashnorm M^{[1]}(fI_{B_{3r}})\|_{L_{p}(B_{r})}\leq Nr^{\beta-d/p}
\|M^{[1]}(fI_{B_{3r}})\|_{L_{p}}
$$
\begin{equation}
                     \label{7,5.4}
\leq Nr^{\beta-d/p}\|fI_{B_{3r}}\|_{L_{p}}=N(3r)^{\beta}
\dashnorm f \|_{L_{p }(B_{3r})}
\leq N\|f\|_{E_{p,\beta}},
\end{equation}
where to derive the last inequality we used
Remark \ref{remark 10.3.1}.
Furthermore, for $x\in B_{r}$ and
$\varepsilon\leq r$ obviously
$(fI_{B^{c}_{3r}})^{(\varepsilon)}(x) 
=0$. Also for $\varepsilon>r$
$$
(fI_{B^{c}_{3r}})^{(\varepsilon)}(x)
\leq N\dashint_{B_{\varepsilon}} (fI_{B^{c}_{3r}})(x+y)\,dy,
$$
$$
r^{\beta}(fI_{B^{c}_{3r}})^{(\varepsilon)}(x)
\leq N\varepsilon^{\beta}\dashint_{B_{\varepsilon}} f(x+y)\,dy
$$
$$
\leq
N\varepsilon^{\beta} \dashnorm f\|
_{L_{p}(B_{\varepsilon}(x))}\leq N\|f\|_{E_{p,\beta}}.
$$
 Hence,
$
r^{\beta}M^{[1]}(fI_{B^{c}_{3r}})\leq
N\|f\|_{E_{p,\beta}}$ on $B_{r}$
and
$$
r^{\beta}\dashnorm M^{[1]}(fI_{B^{c}_{3r}})
\|_{L_{p}(B_{r})}\leq N\|f\|_{E_{p,\beta}}.
$$
By combining this with \eqref{7,5.4}
we come to \eqref{7,5.2} and the lemma is proved. \qed

\begin{remark}
                 \label{remark 7,7.1}
In Real Analysis traditionally
instead of $E_{p,\beta}$ one considers its ``homogeneous'' version
$\dot E_{p, \beta}$ (cf.~\eqref{9.6.1})
defined 
\index{$A$@Sets of functions!$\dot E_{p,\beta}$}%
\index{$N$@Norms!$"|"|f"|"|_{\dot E_{p,\beta}}$}%
as the set of $g$ such that
$$
\|g\|_{\dot E_{p,\beta} }:=
\sup_{\rho<\infty,B\in\bB_{\rho}}\rho^{\beta}
\dashnorm g  \|_{ L_{p}(B)} <\infty 
$$
(note $\rho<\infty$). 
\index{$A$@Sets of functions!$\dot E_{p,\beta}$}%

By applying scaling one derives form
Lemma \ref{lemma 7.5.1} that 
\begin{equation}
                     \label{7.7.2}
\|M^{[\infty]}f\|_{\dot E_{p, \beta}}\leq N
\|f\|_{\dot E_{p, \beta}}
\end{equation}
with the same $N$ as in the lemma. This estimate is found in 
\cite{CF_88} for $0<\beta<d/p$.
By the way, observe that  
estimate  \eqref{7.7.2} becomes
trivial   
if $\beta>d/p$. The reason for that
is that the right-hand sides in 
\eqref{7.7.2} become infinite 
unless $f=0$. For $\beta=d/p$ we have
$\dot E_{p, \beta}=L_{p}$ and \eqref{7.7.2}
is known from Theorem \ref{theorem 7.7.1} (observe that
$\bM f\leq NM^{[\infty]}f\leq N\bM f$ with
the constants independent of~$f$).
\end{remark}

Here is a useful approximation result.
\begin{lemma}
                               \label{lemma 3.14.3}
Let $0\leq\beta'<\beta$.
If   $\|f_{n}\|_{E_{p,\beta'}}$,
$n=0,1,...$, is a {\em bounded\/} sequence
and 
$$
\|f_{n}- f_{0}\|_{L_{p}(B)}\to0 
$$
for any ball $B$, then for any ball $B$
$$
\lim_{n\to\infty}\|f_{n}-f_{0}\|_{E_{p,\beta}(B)}=0.
$$
In particular, if
  $f\in E_{p,\beta'}$, then for any $B\in \bB$
\begin{equation} 
                            \label{3.14.10}
\lim_{\varepsilon\downarrow0}
\|f^{(\varepsilon)}-f \|_{E_{p,\beta}(B) }=0.
\end{equation}
\end{lemma}

Proof. Clearly, for any $r\in(0,1]$  
$$
\|f_{n}-f_{0}\|_{E_{p,\beta}(B)}
\leq 2r^{\beta -\beta'}
\sup_{n}\sup_{\rho\leq r}\rho^{\beta'}
\sup_{B ' \in\bB_{\rho}}\dashnorm f_{n}
\|_{L_{p}(B ' )}
$$
$$
+N(d,r,p,q, B) 
\|f_{n}-f_{0}\|_{L_{p}(B)}.
$$
Here the first term on the right can be made as small as
we like   on   account of $r$ and the second
term tends to zero as $n\to\infty$.
  \qed

Another result of approximation type
is the following.    
\begin{lemma}
                \label{lemma 7,4.3}
Let   $g(x)\geq 0$ be a Borel function such that for any smooth bounded $f(x)$ we have
\begin{equation}
                       \label{7,4.5}
\int_{\bR^{d}}g|f|\,dx\leq
\|f\|_{E_{p,\beta}}
\end{equation}
(bounded linear functional on $E_{p,\beta}$
of special type).
Then, for any $f\in E_{p,\beta}$,
\eqref{7,4.5} holds and, moreover,
$$
\lim_{\varepsilon\downarrow 0}
\int_{\bR^{d}}
g|f-f^{(\varepsilon)}|\,dx=0.
$$
\end{lemma}

Proof. First note that \eqref{7,4.5}
holds for any $f\in E_{p,\beta}$
in light of Fatou's lemma and the fact
that $\|f^{(\varepsilon)}\|_{E_{p,\beta}}
\leq \|f\|_{E_{p,\beta}}$.
After that it only remains to use
the dominated convergence theorem
along with the facts that, for $\varepsilon\leq 1$,
$$
|f-f^{(\varepsilon)}|\leq 2\sup_{ 
\gamma\leq1}|f|^{( \gamma)}= 2M^{[1]}f
$$
and  $\|M^{[1]}f\|_{E_{p,\beta}}
\leq N \|f\|_{E_{p,\beta}}$.   \qed

The following is similar to H\"older's
inequality.

\begin{lemma}
                \label{lemma 6,11.1}
If  
$
p\geq 1, \beta> 1, 
$ and
$ p_{0}=\beta p = s (\beta-1) $, then for any
$b,g$
$$
\|bg\|_{\dot E_{p,\beta}}
\leq \|b\|_{\dot E_{p_{0},1}}\|g\|_{\dot E_{s,\beta-1}},
$$ 
$$
\|bg\|_{E_{p,\beta}}
\leq \|b\|_{E_{p_{0},1}}\|g\|_{E_{s,\beta-1}}.
$$

\end{lemma}

Proof. Observe that $p/p_{0}+(p_{0}-p)/p_{0}=1,
pp_{0}/(p_{0}-p)=s$.
By H\"older's inequality, for any $\rho>0$ and $B\in\bB_{\rho}$,
$$
\rho^{\beta}\dashnorm I_{B }bg\|_{L_{p}}
\leq \rho\dashnorm b I_{B}\|_{L_{p_{0}}}\cdot \rho^{\beta-1}
\dashnorm I_{B} g \|_{L_{s}}
$$
and our assertion follows. \qed

\subsection{Interpolations and embeddings}

The interpolation inequality 
relating the $E_{p,\beta}$-norm
of $Du$ to the $E_{p,\beta}$-norms
of $u$ and $D^{2}u$ is quite natural. It looks like
it first appeared in \cite{FHS_17}
(2017), 
which makes the author wonder how in the past people
claiming that flat boundary and interior
Morrey-Sobolev estimates lead to global estimates
in smooth domains using flattening the boundary 
and partitions of unity. These procedures
unavoidably lead to  the  appearance of the first
order terms, the way to deal with which was not exhibited before \cite{FHS_17}. Unfortunately,
the proof in \cite{FHS_17} contains  an error 
(see Lemma 4.2 there). We give  a different proof mentioning that yet
another proof is given in Remark
\ref{remark 1.27.3}.
Here is Lemma 4.4 of \cite{Kr_22}.  

\begin{lemma}
                \label{lemma 12.26.1}
Let $p\in[1,\infty)$, $ \beta> 0$,
$R\in(0,\infty]$, $u\in W^{2}_{p}(B_{ R})$,
and $x_{0}\in B_{R}$. Then there is a constant
$N=N( d,p,\beta)$ such that, if
$\rho<\infty$,  $r\leq \rho/2\leq R/2 $, $\varepsilon
\in(0,1\wedge\rho] $, then
$$
 f(r):=r^{\beta}
\dashnorm I_{B_{R}}Du\|_{L_{p}(B_{r}(x_{0}))}
\leq 
 N\varepsilon   \sup_{r\leq s\leq \rho}s^{\beta}
\dashnorm I_{B_{ R}}D^{2}u\|_{L_{p}(B_{s}(x_{0}))}
$$
\begin{equation}
                      \label{7.30.1}
+N\varepsilon^{-1}  
\sup_{r\leq s\leq \rho}s^{\beta}\dashnorm I_{B_{R}}( u-c)\|_{L_{p}(B_{s}(x_{0}))},
\end{equation}
where $c$ is any constant. 
\end{lemma}

Proof.   Obviously, we may  assume that $c=0$. Then denote $v= Du$,
$w= D^{2}u$, $G_{s}=B_{s}(x_{0}) \cap B_{R} $,
$$
U=\sup_{r\leq s\leq \rho }s^{\beta}
\dashnorm I_{B_{R}}u\|_{L_{p}(B_{s}(x_{0}))},\quad
W=\sup_{r\leq s\leq \rho }s^{\beta}
\dashnorm I_{B_{R}}D^{2}u\|_{L_{p}(B_{s}(x_{0}))} . 
$$
 
By Poincar\'e's inequality, for $ r\leq s\leq  \rho$,
$$
\dashnorm v-v_{G_{s}}\|_{L_{p}(G_{s})}\leq
N(d)s\dashnorm w\|_{L_{p}(G_{s})}  
\leq Ns^{1-\beta}W.
$$
Also by the interpolation inequalities
for Sobolev spaces, 
   there
exists a constant $N=N( d,p)$
such that, for  
 $\varepsilon \in(0,1]$ and $\varepsilon\leq s\leq  \rho$ ,
$$
\dashnorm v-v_{G_{s}}\|_{L_{p}(G_{s})}
\leq 2\dashnorm v\|_{L_{p}(G_{s})}
$$
$$
\leq N \dashnorm w\|^{1/2}_{L_{p}(G_{s})}
\dashnorm  u\|^{1/2}_{L_{p}(G_{s})}
+Ns^{-1}\dashnorm  u\| _{L_{p}(G_{s})}
$$
\begin{equation}
                         \label{7.31.1}
\leq N \dashnorm w\|^{1/2}_{L_{p}(G_{s})}
\dashnorm  u\|^{1/2}_{L_{p}(G_{s})}
+N\varepsilon^{-1}\dashnorm  u\| _{L_{p}(G_{s})},
\end{equation}
which for $\rho\geq s\geq \varepsilon\vee r$    yields
$$
s^{\beta}\dashnorm v-v_{G_{s}}\|_{L_{p}(G_{s})}
\leq N W^{1/2}U^{1/2}   
+N \varepsilon^{-1} U.
$$ 
  Hence, for any $\varepsilon\in(0,1]$ 
(at this point the restriction $\varepsilon\leq\rho$ is not  needed)
and $r\leq s\leq \rho$
$$
g(s):=
s^{\beta}\dashnorm v-v_{G_{s}}\|_{L_{p}(G_{s})}
\leq N_{1}\varepsilon W 
+N_{2}\varepsilon^{-1}U=:h,
$$
where $N_{1}=N_{1}(d,p)$, 
$N_{2}=N_{2}( d,p)$.

Then we use Campanato's argument
 (see, for instance, Proposition 5.4 in \cite{MM_12}). Obviously,
\begin{equation}
                    \label{1.1.1}
f(s)\leq g(s)+s^{\beta}|v_{G_{s}}|
\leq h+s^{\beta}|v_{G_{s}}|
\end{equation}
and if $2s\leq \rho$ we have
$$
|v_{G_{s}}|\leq |v_{G_{s}}-v_{G_{2s}}|
+|v_{G_{2s}}|,
$$
where
$$
|v_{G_{s}}-v_{G_{2s}}|\leq
\dashint_{G_{s}}|v-v_{G_{2s}}|\,dx
$$
$$
\leq N(d)\dashint_{G_{2s}}|v-v_{G_{2s}}|\,dx\leq N(2s)^{-\beta}h.
$$
It follows for $r_{n}=r2^{n}$ in place of $s$  that, if $r_{n+1}\leq \rho$,
then
$$
 |v_{G_{r_{n}}}|\leq  N hr^{-\beta}2^{-(n+1)\beta} + |v_{G_{r_{n+1}}}|.
$$
By applying this inequality
for $n=0,1,...,n_{0}$, where $n_{0}=
\lfloor\log_{2}(\rho/r)\rfloor-1$, we find
$$
r^{\beta}|v_{G_{r }}|\leq Nh\sum_{n=0}^{n_{0}}
2^{-(n+1)\beta}+Nr^{\beta}|v_{G_{\rho}}|\leq Nh+N\rho^{\beta}|v_{G_{\rho}}|.
$$
This and \eqref{1.1.1} yields \eqref{7.30.1}
with the following term added to its
right-hand side:
$
N\rho^{\beta}|v_{G_{\rho}}|
$,
which we estimate as in \eqref{7.31.1}
and come to \eqref{7.30.1} as is.
The lemma is proved. \qed

\begin{theorem}
                \label{theorem 12.26.1}
Let $p\in[1,\infty)$, $0<\beta\leq d/p$,
$R\in(0,\infty]$, $u\in W^{2}_{p}(B_{ R})$.
 Then there is a constant
$N=N( d,p,\beta)$ such that, for any $\varepsilon
\in(0,1\wedge R] $,  
\begin{equation}
                      \label{8.4.3}
\|Du\|_{E_{p,\beta}(B_{R})}\leq
N \big(\varepsilon   \|D^{2}u\|_{E_{p,\beta}(B_{R})}
+ \varepsilon^{-1} 
\| u\|_{E_{p,\beta}(B_{R})}\big). 
\end{equation}
\end{theorem}

Proof. By 
Remarks \ref{remark 12.26.2}, 
\ref{remark 10.3.1}, and \ref{remark 12.26.4}
$$
\|Du\|_{E_{p,\beta}(B_{R})}
=\sup_{r\leq  R\wedge1,x\in B_{R}}r^{\beta}
\dashnorm I_{B_{R}}Du\|_{L_{p}(B_{r}(x))}
$$
$$
\leq N(d,p)
\sup_{r\leq  (R\wedge1)/2,x\in \bR^{d}}r^{\beta}
\dashnorm I_{B_{R}}Du\|_{L_{p}(B_{r}(x))}
$$
$$
=N(d,p)\sup_{r\leq (R\wedge1)/2,x\in B_{R}}
r^{\beta}\dashnorm I_{B_{R}}Du\|_{L_{p}(
B_{r}(x))}.  
$$
Now Lemma \ref{lemma 12.26.1} with $\rho=R\wedge1$
implies \eqref{8.4.3}.  \qed

To continue, we need the following
estimates going back to Adams~\cite{Ad_75}.

\begin{lemma}
                           \label{lemma 9.28.1}
Let $0\leq \gamma<\alpha<\beta$, $ \sigma<\beta\wedge d $.
Then there exist  a constant $N$
   such that for any $f\geq0$
and $\rho\in(0,\infty)$ we have
\begin{equation}
                             \label{1.17.2}
R_{\alpha}(I_{B _{ \rho}}f)(0)
\leq N\rho^{\alpha -\gamma}\bM_{\gamma} 
f(0) ,
\end{equation}
\begin{equation}
                        \label{1.3.1}
R_{\sigma}(I_{B^{c}_{ \rho}}f)(0) 
\leq N\rho^{\sigma-\beta}\bM_{\beta}
f (0),
\end{equation}
\begin{equation}
                       \label{9.28.1}
 R_{\alpha}f \leq    
N(\bM_{\beta}f )^{(\alpha-\gamma)/(\beta-\gamma)}
(  \bM_{\gamma} f)^{(\beta-\alpha)/(\beta-\gamma)}.
\end{equation}
In particular (by H\"older's inequality and $\gamma=0$), for   any $p\in[1,\infty]$,
$q\in(1,\infty]$, and measurable $\Gamma$
\begin{equation}
                                \label{9.28.2}
 \|R_{\alpha}f\|_{L_{r}(\Gamma)} \leq    
N\|\bM_{\beta}f \|_{L_{p}(\Gamma)}^{\alpha/\beta}
\|f\|_{L_{q}}^{1-\alpha/\beta},
\end{equation}
provided that
$$
\frac{1}{r}=\frac{\alpha}{\beta}\cdot\frac{1}{p}+
\Big(1-\frac{\alpha}{\beta}\Big)\frac{1}{q}.
$$
\end{lemma}

Proof. We basically repeat the proof
of Proposition 3.1 of \cite{Ad_75}.
Observe that \eqref{9.28.1} at the origin is easily obtained from summing up the inequalities in 
\eqref{1.17.2}  
  and \eqref{1.3.1} with $\sigma=\alpha$ and minimizing with respect 
to $\rho$. At any other point it is obtained
by changing the origin.
 Furthermore
clearly, we may assume that $f$ is bounded
with compact support. 

Then integrating by parts we find
$$
R_{\alpha}(I_{B _{ \rho}}f)(0)=N(d)
\int_{0}^{\rho}r^{\alpha-d}\frac{d}{dr}\Big(\int_{B_{r}}
f \,dx\Big)\,dr=N\rho^{\alpha}
\dashint_{B_{\rho}}f \,dx
$$
$$
+N\int_{0}^{\rho}r^{\alpha-\gamma-1}r^{\gamma}\dashint_{B_{r}}
f \,dxdr\leq N\rho^{\alpha-\gamma}\bM_{\gamma} f(0),
$$
which yields \eqref{1.17.2}.

Next,
$$
R_{\sigma}(I_{B^{c} _{ \rho}}f)(0)=N(d)
\int_{\rho}^{\infty}r^{\sigma-d}\frac{d}{dr}\Big(\int_{B_{r}}
f \,dx\Big)\,dr
\leq N\int_{\rho}^{\infty}r^{\sigma-1}
\dashint_{B_{r}}
f \,dxdr
$$
$$
\leq N\bM_{\beta}f(0)\int_{\rho}^{\infty}r^{\sigma-\beta-1}\,dr=
N\rho^{\sigma-\beta}\bM_{\beta}f(0),
$$
and the lemma is proved. \qed

 Here is Theorem 3.1 of \cite{Ad_75}.

\begin{theorem}
                       \label{theorem 9.29.1}
For any $0<\alpha<\beta$, 
$p\in(1,\infty)$, and $r$ such that
$$
r(\beta-\alpha)=p\beta ,
$$
there is a constant $N$ such that for any
   $f\geq0$ we have
\begin{equation}
                          \label{9.29.3}
 \|R_{\alpha}f\|_{\dot  E_{r,\beta-\alpha}}
\leq N \|f\|_{\dot E_{p,\beta}}.
\end{equation}
\end{theorem}

Proof. It suffices to prove that for any $\rho>0$
$$
\rho^{\beta-\alpha}\Big(\dashint_{B_{\rho}}|R_{\alpha}f|^{r}\,dx\Big)^{1/r}
\leq N\|f\|_{\dot E_{p,\beta}},
$$
that is
\begin{equation}
                          \label{9.29.4}
\rho^{\beta-\alpha-d/r}\Big(\int_{B_{\rho}}|R_{\alpha}f|^{r}\,dx\Big)^{1/r}
\leq N\|f\|_{\dot E_{p,\beta}},
\end{equation}

Observe that by H\"older's inequality
$$
\bM_{\beta}(I_{B_{2\rho}}f)\leq N\|f\|_{\dot E_{q,\beta}}
$$
 and by definition
$$
\Big(\int_{\bR^{d }}I_{B_{2\rho}}f^{p}\,dx
\Big)^{1/p}
\leq N\rho^{d/p-\beta  }\|f\|_{\dot E_{p,\beta}}.
$$
It follows from  \eqref{9.28.2}
with $p=\infty$ that
$$
\Big(\int_{B_{\rho}}|R_{\alpha}(I_{B_{2\rho}}f)|^{r}
dx\Big)^{1/r}\leq N\rho^{(d/p-\beta)(1-\alpha/\beta) }\|f\|_{\dot E_{p,\beta}}
$$
$$
=N\rho^{d/r-\beta+\alpha}\|f\|_{\dot E_{p,\beta}}.
$$
Furthermore, by \eqref{1.3.1}
$$
\Big(\int_{B_{\rho}}|R_{\alpha}(I_{B^{c}_{2\rho}}f)|^{r}
dx\Big)^{1/r}\leq N\rho^{d/r} \sup_{B_{\rho}}R_{\alpha}(I_{B^{c}_{2\rho}}f)
\leq N\rho^{d/r+\alpha-\beta}\|f\|_{\dot E_{p,\beta}}.
$$
By combining these estimates we come to 
\eqref{9.29.4} and the theorem is proved. \qed

\begin{remark}
                   \label{remark 1.4.4}
If $u\in \dot E^{1}_{p,\beta}$ has compact support, $p\in[1,\infty)$, then (a.e.)
\begin{equation}
                      \label{1.4.3}
|u|\leq NR_{1}|Du|.
\end{equation}
This fact is pointed out in Remark
\ref{remark 9.5.1} for smooth $u$
with compact support. In the general case it suffices to note that (a.e.) 
$$
|u|=\lim_{\varepsilon\downarrow 0}
|u^{(\varepsilon)}|\leq N\lim_{\varepsilon\downarrow 0}R_{1}|(Du)^{(\varepsilon)}|
\leq N\lim_{\varepsilon\downarrow 0}(R_{1}|Du|)^{(\varepsilon)}=NR_{1}|Du|, 
$$
where the last equality follows from the fact that $R_{1}|Du|$ is locally integrable.

Similar argument shows that if $u\in E^{2}_{p,\beta}$ has compact support,
then $|u|\leq NR_{2}|\Delta u|$.
\end{remark} 
\begin{corollary}[Morrey-Sobolev space embedding theorem]
                \label{corollary 1.4.1}
If $d>\beta>1$ and 
\index{Morrey-Sobolev space embedding theorem}%
$p\in(1,d/\beta]$,   then for any $u\in \dot E^{1}_{p,\beta}$, $v \in  E^{1}_{p,\beta}$
$$
\|u\|_{\dot E_{r,\beta-1}}
\leq N \|Du\|_{\dot E_{p,\beta}},
\quad \|v\|_{E_{r,\beta-1}}
\leq N (\|Dv\|_{E_{p,\beta}}+
\|v\|_{E_{p,\beta}}),
$$
where $r(\beta-1)=p\beta$ and $N$
is independent of $u,v$.
\end{corollary}

Indeed, if $u$ has compact support,
the first estimate follows from Theorem
\ref{theorem 9.29.1} and Remark \ref{remark 1.4.4}. In the general case
it suffices to take $\zeta\in C^{\infty}_{0}$ such that $\zeta=1$
on $B_{1}$, $\zeta=0$
outside $B_{2}$, $0\leq\zeta\leq1$, denote $\zeta_{\varepsilon}(x)=\zeta(\varepsilon x)$ and observe
that 
$$
\|u\|_{\dot E_{r,\beta-1}}\leq\nlimsup_{\varepsilon\downarrow 0}
\|u\zeta_{\varepsilon}\|_{\dot E_{r,\beta-1}}
\leq N  \nlimsup_{\varepsilon\downarrow 0}\|\zeta_{\varepsilon}Du\|_{\dot E_{p,\beta}}+N\nlimsup_{\varepsilon\downarrow 0}\varepsilon\|u(D\zeta)(\varepsilon\cdot)\|_{\dot E_{p,\beta}}.
$$

The second estimate follows from the first one after observing that
$$
\|v\|_{E_{p,\beta}(B_{1})}
\leq \|v\zeta\|_{\dot E_{p,\beta}
}\leq N(\|Dv\|_{E_{p,\beta}(B_{2})}
+\|v\|_{E_{p,\beta}(B_{2})}).
$$

\begin{corollary}
              \label{corollary 1.8.1}
If $1<\beta<q\leq d$, $b\in \dot E_{q,1}$,
$c\in  E_{q,1}$,  
  $p\in(1,q/\beta]$,  then
for any $u\in \dot E^{2}_{p,\beta}$,
 $v\in   E^{2}_{p,\beta}$
\begin{equation}
                       \label{1.8.1}
\|bDu\|_{\dot E_{p,\beta}}\leq N
\|b\|_{\dot E_{q,1}}\|D^{2}u
\|_{\dot  E_{p,\beta}}, 
\end{equation}
\begin{equation}
                       \label{1.9.1}
 \|cDv\|_{  E_{p,\beta}}\leq N
\|c\|_{  E_{q,1}}(\|D^{2}v
\|_{  E_{p,\beta}}+\| v
\|_{  E_{p,\beta}}) 
\end{equation}
where $N$ is independent of $b,c,u$.
\end{corollary}

Indeed, by Lemma \ref{lemma 6,11.1} 
$$
\|bDu\|_{\dot E_{p,\beta}} \leq 
\|b\|_{\dot E_{\beta p,1}}\|Du\|_{\dot E_{r,\beta-1}}\leq 
\|b\|_{\dot E_{q,1}}\|Du\|_{\dot E_{r,\beta-1}},
$$
where $r(\beta-1)=p\beta$ and by
Corollary \ref{corollary 1.4.1}
the last norm is dominated by $N\|D^{2}
u\|_{\dot E_{p,\beta}}$. This proves
\eqref{1.8.1}. To prove \eqref{1.9.1}
it suffices to use the same line of arguments combined at the last step
with the interpolation inequality
from Theorem \ref{theorem 12.26.1}.

Corollary \ref{corollary 1.8.1} implies the  following 
elliptic version of Lemma
5.1 of \cite{Kr_25} if $\rho_{b}=1$. 
The general case is obtained by
using scaling. 

\begin{theorem} 
              \label{theorem 8,23.3}
Let $b$ be an $\bR^{d}$-valued function
on $\bR^{d}$.
Let $p_{b}\in(1,d]$, $\rho_{b}\in(0,\infty)$ be such 
\index{$S$@Miscelenea!$\bar b_{p_{b},\rho_{b}}$}%
that 
$$
\bar b_{p_{b},\rho_{b}}:=
\sup_{\rho\leq \rho_{b}}\rho
\sup_{B\in \bB_{\rho}}\dashnorm
b\|_{L_{p_{b}}(B)}<\infty.
$$
Let  $ \beta>1$, $1<p\leq p_{b}/\beta$.  
Then the operator
$$
\cL u:=b^{i}D_{i}u
$$
 is a bounded operator
from $E^{ 2}_{p, \beta;\rho_{b}}$ to $E_{p, \beta;\rho_{b}}$ and for any $u\in E^{ 2}_{p, \beta;\rho_{b}}$ 
\begin{equation}
                         \label{1.17.3}
\|b^{i}D_{i}u \|_{E_{p,\beta;\rho_{b}}}
\leq N(d,p,p_{b},\beta)\bar b_{p_{b},\rho_{b}}\big(\|D^{2}u\|_{E_{p,\beta};\rho_{b}}
+ \rho^{-2}_{b}\|u\|_{E_{p,\beta;\rho_{b}}}).
\end{equation}
\end{theorem}

\begin{remark}
                \label{remark 1.17.4}

The most desirable version 
of \eqref{1.17.3} 
would be, of course, 
\begin{equation}
                     \label{1.17.4}
\|b^{i}D_{i}u\|_{E_{p,\beta}(B_ 1)}
\leq \varepsilon \|D^{2}u\|_{E_{p,\beta}(B_1)}
+N(\varepsilon)\|u\|_{E_{p,\beta}(B_1)} 
\end{equation}
for a $p\in(1, p_{b}/\beta )$ any $u\in E_{p,\beta}$ and any $\varepsilon\in(0,1]$ with
$N(\varepsilon)$ independent of $u$.
This fact is, actually, claimed
in Theorem 5.4 of \cite{FHS_17}.
However, it is wrong.

To show that
let 
$$
\beta\in(1,2),\quad p\in(1,d/\beta),
$$
 and $f(t)$, $t\geq0$, be a smooth function,
equal to zero for $t\leq 1$, equal
to $t^{2-\beta}$ for $t\in [2,3]$
and equal to zero for $t\geq4$. Clearly,
$$
|f'(t)|\leq Nt^{1-\beta},\quad |f''(t)|\leq Nt^{-\beta}
$$
 for some
constants $N$ and 
$$
f'(t)=(2-\beta)t^{1-\beta}
$$
 for $t\in [2,3]$. Now, for $\kappa\in(0,1/4]$ define 
$$
u_{\kappa}(x)=f(|x|/\kappa).
$$ 
Then
$$
D_{i}u_{\kappa}=f'\frac{x^{i}}{
\kappa|x|},\quad D_{ij}u_{\kappa}=
f''\frac{x^{i}x^{j}}{\kappa^{2}|x|^{2}}+f'\frac{1}{\kappa|x|}
\Big(\delta^{ij}-\frac{x^{i}x^{j}}{
|x|^{2}}\Big).
$$
It follows that 
$$
|D^{2}u_{\kappa}
|\leq N\kappa^{\beta-2}|x|^{-\beta}\quad\text{in}\quad
B_{1},\quad |Du|=(2-\beta)
\kappa^{\beta-2}|x|^{1-\beta}\quad\text{in}\quad
B_{3\kappa}\setminus B_{2\kappa}.
$$
With $b=|x|^{-1}$ this yields 
$$
\dashnorm bDu_{\kappa}\|_{L_{p}(B_{3\kappa})}
\geq N_{1}(3\kappa)^{-2},\quad
\dashnorm D^{2}u_{\kappa}\|_{L_{p}(B_{r})}
\leq N_{2}\kappa^{\beta-2}r^{-\beta},
\quad r<1. 
$$
By adding to this that $\dashnorm u_{\kappa}\|_{L_{p}(B_{r})}
\leq N\sup |f|$, we see that \eqref{1.17.4}
can only hold if $N_{1}3^{\beta-2}\leq\varepsilon N_{2}$.

\end{remark}

Let us 
\index{$A$@Sets of functions!$L_{p}^{\loc}$}% 
\index{$N$@Norms!$"|"|f"|"|_{L_{p}^{\loc}}$}%
write $f\in L_{p}^{\loc}$ if
$$
\|f\|_{L_{p}^{\loc}}:=\sup_{B\in\bB_{1}}\|f\|_{L_{p}(B)}<\infty.
$$
Here is the classical Sobolev embedding
\index{Sobolev embedding  theorem}%
theorem obtained from Corollary
\ref{corollary 1.4.1}
as $p\beta=d$.

\begin{corollary}
                \label{corollary 1.4.01}
If  $1<p <d$, then
for any $u\in W^{1}_{p}$   
\begin{equation}
                        \label{1.4.2}
 \|u \|_{L_{r} }\leq N\|Du\|_{L_{p} },\quad  \|u \|_{L_{r}^{\loc}}\leq N(\|Du\|_{L_{p}^{\loc} }+\|u\|_{L_{p}^{\loc} }),
\end{equation}
where $r(d/p-1)=d$, $N$
is independent of $u$.
\end{corollary}
\begin{remark}
              \label{remark 1.4.1}
Usually the first estimate in \eqref{1.4.2} is derived by means of the
Gagliardo-Nirenberg inequality,
which says that it holds for $p=1$.

\end{remark}

For $\alpha=2$ in Theorem \ref{theorem 9.29.1} similarly to
Corollary \ref{corollary 1.4.1}
we get the first estimate in \eqref{1.9.3}.
  
\begin{corollary}
                \label{corollary 1.4.100}
If $d\geq \beta>2$, 
$p\in(1,d/\beta]$, then for any $u\in \dot E^{2}_{p,\beta}$, 
$v\in   E^{2}_{p,\beta}$  
\begin{equation}
                      \label{1.9.3}
\|u\|_{\dot E_{r,\beta-2}}
\leq N \|\Delta u\|_{\dot E_{p,\beta}},
\quad \|v\|_{ E_{r,\beta-2}}
\leq N( \|D^{2} v\|_{  E_{p,\beta}}
+ \| v\|_{  E_{p,\beta}}) 
\end{equation}
where $r(\beta-2)=p\beta$ and $N$
is independent of $u$.
\end{corollary}

The second estimate is obtained from the first one as in Corollary~\ref{corollary 1.4.1}.

\begin{corollary}
                \label{corollary 1.4.10}
If $1<p<d/2$, then for any $u\in W^{2}_{p}$
\begin{equation}
                        \label{1.4.20}
\|u \|_{L_{r} }\leq N\|\Delta u\|_{L_{p} },\quad \|u \|_{L_{r}^{\loc}}\leq N(\|D^{2}u\|_{L_{p}^{\loc} }+
\| u\|_{L_{p}^{\loc} }), 
\end{equation}
where $r(d/p-2)=d$, $N$
is independent of $u$.
\end{corollary}

Now we are going to present a few results about embeddings of $E^{i}_{p,\beta}$ into $C^{\alpha}$ spaces.

As usual, when we say that a function
$u$,
uniquely defined only up to almost everywhere, is, say continuous, we mean
that there is a continuous
\index{Morrey-Sobolev space embedding theorem}% 
 function
equal to $u$ almost everywhere.

\begin{theorem}[Morrey-Sobolev space embedding theorem]
                             \label{theorem 1.5.1} 
(i) Let $0<\beta<2$, $\beta\leq d/p$,
 and $p\in[1,\infty)$.
\index{Morrey-Sobolev space embedding theorem}%
\index{Morrey theorem}%
Then any $u\in E^{2}_{p,\beta}$
is bounded and continuous  
 and for any $\varepsilon
\in(0,1]$  
\begin{equation}
                                    \label{3.20.06}
|u|
\leq \varepsilon^{2-\beta} \| D^{2}u\|_{E_{p, \beta}}+N(d, \beta)\varepsilon^{-\beta} \| u\|_{E _{p ,\beta}}.
\end{equation}  

(ii) 
In case $0<\beta<1 $, $\beta\leq d/p$, 
 $p\in[1,\infty)$, $\varepsilon\in(0,1]$, any
$u\in E^{1}_{p, \beta}$ is continuous, 
\begin{equation}
                 \label{1.5.5}
|u | \leq
N(d, \beta)(\varepsilon^{1-\beta} \|Du\|_{E _{p, \beta}}+ \varepsilon^{ -\beta} \| u\|_{E _{p, \beta}}),
\end{equation}
and for any
$ x,y \in B_{1}$ we have (Morrey theorem)
\begin{equation}
               \label{4.6.2}
| u(x)- u(y)| \leq
N(d, \beta)|x-y|^{1-\beta}  \|u\|_{E^{1}_{p, \beta}}.
\end{equation}

(iii) Let $1<\beta  <2$, $\beta\leq d/p$, $p\in[1,\infty)$. Then for any
$u\in E^{2}_{p, \beta}$,  
$ x,y \in B_{1}$, we have
\begin{equation}
                           \label{4.6.1}
|u( x )-u( y)|\leq
N(d,p ,\beta)|x-y|^{2-\beta}\|u\|_{E^{ 2}_{p, \beta}}.
\end{equation}

\end{theorem}

Proof. (i) Take nonnegative $\zeta\in C^{\infty}_{0}$ such that $\zeta=1$ on
$B_{1/2}$ and $\zeta=0$ outside $B_{1}$.
For $\varepsilon\in (0,1]$ introduce
$\zeta_{\varepsilon}(x)=\zeta(x/\varepsilon)$. By Remark \ref{remark 1.4.4} and \eqref{1.17.2} ($B_{1/2}$-a.e.)
$$
|u|\leq NR_{2}|\Delta(\zeta_{\varepsilon}u)|
\leq N\varepsilon^{2-\beta}
\|\Delta(\zeta_{\varepsilon}u)\|_{E_{p,\beta}}\leq N\varepsilon^{2-\beta}
\|\Delta u \|_{E_{p,\beta}}
$$
$$
+N\varepsilon^{1-\beta}\|Du\|_{E_{p,\beta}}
+N\varepsilon^{-\beta}\|u\|_{E_{p,\beta}}.
$$
After that the interpolation theorem
yields \eqref{3.20.06} almost everywhere in $B_{1/2}$ and, hence,
almost everywhere. To prove the continuity of $u$ it suffices to observe that
$u^{(\varepsilon)}\to u$ uniformly in light
of \eqref{3.20.06} and Lemma \ref{lemma 3.14.3},
in which we can take any $\beta'>\beta,\beta'<2$.
This proves
(i).

Estimate \eqref{1.5.5} is proved in the same way
as  \eqref{3.20.06} (without using interpolation).

To prove \eqref{4.6.2} observe that
for $\varepsilon\in(0,1)$ such that $1-\varepsilon<\beta$, $|x|\leq 1/2$, smooth $u$, 
$v=u\zeta$, and $f=|Dv|$,   we have
\eqref{1.5.6}, where by \eqref{1.17.2}
$$
R_{1}(I_{B_{2|x|}}f)(x)+R_{1}(I_{B_{2|x|}}f)(0)\leq N|x|^{1-\beta}\|f\|_{E_{p,\beta}}
$$
and by \eqref{1.3.1} ($f$ has support
in $B_{1}$ and $\beta\leq d/p$) 
$$
R_{1-\varepsilon}
(I_{B^{c}_{2|x|}}f)(0)\leq N|x|^{1-\varepsilon-\beta}\|f\|_{E_{p,\beta}}.
$$
The combination of these estimates
shows that
$$
|u(x)-u(0)|\leq N|x|^{1-\beta}
\|D(u\zeta)\|_{E_{p,\beta}}\leq
N|x|^{1-\beta}
\|u\|_{E^{1}_{p,\beta}},
$$
which yields \eqref{4.6.2} for smooth $u$. For general $u\in E^{1}_{p,\beta}$
we argue as in the case of \eqref{3.20.06}.

To prove \eqref{4.6.1} it suffices to
observe that $u\in E^{1}_{p,\beta-1}$
by Corollary \ref{corollary 1.4.1}
and after that use \eqref{4.6.2}. \qed

Notice that, if we want to have functions from
$W^{ 2}_{p }$ to be H\"older continuous
(as in \eqref{4.6.1}), we need $d/p <2$.
For H\"older continuity of the gradients of
$W^{ 2}_{p }$-functions we need $d/p <1$.
There are no restrictions on $d/p $ from above
in Theorem \ref{theorem 1.5.1} 
(instead, $\beta$ controls the local behavior of $u,Du$).
\begin{remark}
               \label{remark 1.5.4}
The restriction $\beta\leq d/p$
in Theorem \ref{theorem 1.5.1} 
is almost irrelevant. Without
it one need only replace $\beta$
with $\beta\wedge(d/p)$. This is because $E_{p,\beta}=L_{p}^{\loc}$
for $\beta\geq d/p$.

\end{remark}

\begin{corollary}
               \label{corollary 1.6.1}
Let $0<\beta  <1$, $\beta\leq d/p$, $p\in[1,\infty)$. Then for any
$u\in E^{2}_{p, \beta}$,  
$Du$
is bounded and continuous  and for any $\varepsilon
\in(0,1]$ 
\begin{equation}
                                    \label{3.20.60}
|Du|
\leq \varepsilon^{1-\beta} \| D^{2}u\|_{E_{p, \beta}}+N(d, \beta)\varepsilon^{-\beta-1} \| u\|_{E _{p, \beta}}.
\end{equation}

\end{corollary}
Indeed, $Du\in E^{1}_{p,\beta}$ and by 
\eqref{1.5.5}
$$
|Du | \leq
N(d, \beta)(\varepsilon^{1-\beta} \|D^{2}u\|_{E _{p, \beta}}+ \varepsilon^{ -\beta} \| Du\|_{E _{p, \beta}}),
$$
where the last term admits estimate
\eqref{8.4.3}.

\subsection{Solvability of the Laplace
equation in Morrey-Sobolev spaces} 
\begin{lemma}
                 \label{lemma 2.22.3}
Let $\beta>0$ and $p>1$. Then there is
a constant $N$ such that for any
$u\in C^{\infty}_{0}$
\begin{equation}
                     \label{2.22.4}
\|D^{2}u\|_{\dot E_{p,\beta}}\leq N \|f\|_{\dot E_{p,\beta}},
\end{equation}
where $f=\Delta u$.
\end{lemma}

Proof. We may assume that $d\geq3$,
in which case we observe that to prove
\eqref{2.22.4} it suffices to show that
\begin{equation}
                     \label{2.22.3}
\|D^{2}u\|_{L_{p}(B_{1})}\leq N \|f\|_{\dot E_{p,\beta}}.
\end{equation}
Indeed, then changing scales and the origin will yield 
$$
\rho^{\beta}\|D^{2}u\|_{L_{p}(B)}\leq N \|f\|_{\dot E_{p,\beta}}
$$
for any $\rho>0$ and $B\in\bB_{\rho}$, which is equivalent to \eqref{2.22.4}.

In case of \eqref{2.22.3}
let 
$$
g= fI_{B_{2}},h=fI_{B_{2}^{c}},\quad G=Rg,\quad  H=R h
$$
 We know that $u=-F-H$. By Theorem \ref{theorem 9.6.1}
 (see also Remark \ref{remark 2.22.2})
$$
\|D^{2}G\|_{L_{p}(B_{1})}\leq
\|D^{2}G\|_{L_{p} }\leq N \|f\|_{L_{p}
(B_{2}) }\leq N \|f\|_{\dot E_{p,\beta}}.
$$
In addition, by \eqref{1.3.1} with $\delta=0$ on $B_{1}$ we have
$$
|D^{2}H|\leq NR_{0}|h|\leq 
N\bM_{\beta}f\leq N\|f\|_{\dot E_{p,\beta}}.
$$
This and the above estimate lead to
\eqref{2.22.3}.  \qed

Next, define
$$
\sfG_{\lambda}(x)=\int_{0}^{\infty}\frac{1}{(4\pi t)^{d/2}}
e^{-|x|^{2}/(4t)-\lambda t}\,dt
$$
and   set
$u =\sfG_{\lambda}*f$. It is easy
\index{$S$@Miscelenea!$\sfG_{\lambda}$}%
 to check
that, if $f\in C^{\infty}_{0}$, then
$u$ is infinitely differentiable, each
of its derivatives decays as $|x|\to\infty$ exponentially fast and 
$u\in W^{2}_{p} $. Furthermore,
  it is
well known (see, for instance,
Exercise 1.6.5 of \cite{Kr_96}) that
$$
\lambda u-\Delta u=f
$$
 and, if $v\in C^{\infty}_{0}$, then 
$$v=\sfG_{\lambda}*(\lambda v-\Delta v).
$$

\begin{theorem}
                 \label{theorem 1.20.3}
Let $p\in(1,\infty),0<\beta\leq d/p$. Then
there exists $N$,
depending only on $d,p,\beta$, such that for any $\lambda>0$
and $f\in E_{p,\beta}$ there exists
a unique solution $u\in E^{2}_{p,\beta}$ of the equation 
$$
\Delta u-\lambda u=f.
$$
 Furthermore, for any $u\in E^{2}_{p,\beta}$ we have
\begin{equation}
                    \label{1.20.8}
 \|D^{2}u,\sqrt\lambda Du,\lambda u\|_{E_{p,\beta}}\leq N(1+\lambda^{-1})\|\Delta u-\lambda u\|_{E_{p,\beta}}.
\end{equation}
\end{theorem}

Proof. First, we prove 
\eqref{1.20.8} for $u\in C^{\infty}_{0}$. Let  $\zeta\in C^{\infty}_{0}$
be such that $\zeta=1$ in $B_{1}$,
$\zeta=0$ outside $B_{2}$,
and $1\leq\zeta\leq1$, and let $u\in C^{\infty}_{0}$. Then by Lemma \ref{lemma 2.22.3}
$$
\|D^{2}u\|_{E_{p,\beta}(B_{1})}\leq
\|D^{2}(u\zeta)\|_{\dot E_{p,\beta}}
$$
$$ 
\leq N\|\Delta(u\zeta)\|_{\dot E_{p,\beta}}\leq N\|\Delta u, Du,u\|_{E_{p,\beta}}.
$$
Using shifts of $B_{1}$ and the interpolation Theorem \ref{theorem 12.26.1}
($\beta\leq d/p$), we come to  
\begin{equation}
                      \label{1.21.1}
\|D^{2}u\|_{E_{p,\beta} } 
\leq N\|\Delta u \|_{E_{p,\beta}}
+N\| u\|_{E_{p,\beta}}.
\end{equation}

Next,   set $\lambda u-\Delta u=:f$.
Then
$u=\sfG_{\lambda}*f$, 
which implies that (Minkowski: the norm of a ``sum'' is less than the ``sum'' of norms)
\begin{equation}
                      \label{1.21.2}
\|u\|_{E_{p,\beta}}\leq\|\sfG_{\lambda}\|_{L_{1}}
\|f\|_{E_{p,\beta}}=\lambda^{-1}
\|f\|_{E_{p,\beta}}.
\end{equation}
This and \eqref{1.21.1} yield
$$
\|D^{2}u\|_{E_{p,\beta} }
\leq N\|\Delta u -\lambda u\|_{E_{p,\beta}}
+N\| u\|_{E_{p,\beta}} 
\leq N(1+\lambda^{-1})\|\Delta u -\lambda u\|_{E_{p,\beta}}, 
$$
which along with \eqref{1.21.2}
and interpolation lead to \eqref{1.20.8} for $u\in C^{\infty}_{0}$.

To prove it for general $u\in E^{2}_{p,\beta}$ define $\zeta_{\varepsilon}(x)
=\varepsilon^{-d}(x/\varepsilon)$,
$\zeta^{\varepsilon}(x)=\zeta(\varepsilon x)$, $u^{\varepsilon}=
(\zeta_{\varepsilon}*u),v^{\varepsilon}=u^{\varepsilon}\zeta^{\varepsilon}$. Since for each ball $B$,
for instance,
$$
\|D^{2}u\|_{L_{p}(B)}=\lim_{\varepsilon\downarrow0}
\|D^{2}v^{\varepsilon}\|_{L_{p}(B)},
$$
we have
$$
\|D^{2}u,\sqrt\lambda Du,\lambda u\|_{E_{p,\beta}}
\leq\nliminf_{\varepsilon\downarrow 0}
\|D^{2}v^{\varepsilon},\sqrt\lambda Dv^{\varepsilon},\lambda v^{\varepsilon}\|_{E_{p,\beta}}.
$$
Also observe that 
$$
\Delta v^{\varepsilon}-\lambda v^{\varepsilon}
=(\Delta u-\lambda u )^{\varepsilon}
\zeta^{\varepsilon}+2\varepsilon
(D_{i}\zeta)(\varepsilon x)(D_{i}u)^{\varepsilon}+\varepsilon^{2}(\Delta\zeta)(\varepsilon x)u ^{\varepsilon}
$$
which implies that
$$
\nlimsup_{\varepsilon\downarrow 0}
\|\Delta v^{\varepsilon} -\lambda v^{\varepsilon}\|_{E_{p,\beta}}
\leq \|\Delta u -\lambda u\|_{E_{p,\beta}}
$$
and finishes the proof of \eqref{1.20.8}.

To prove the solvability
in $E^{2}_{p,\beta}$ of 
$$
\lambda
u-\Delta u=f\in E_{p,\beta},
$$
 observe
that we have it with $f^{\varepsilon}\zeta^{\varepsilon}$ in place of  
$f$. Call $u_{\varepsilon}$ the corresponding
solution. Then the $E_{p,\beta}$-norms
of $D^{2} u_{\varepsilon},D u_{\varepsilon}, u_{\varepsilon}$
are uniformly bounded for $\varepsilon
\in(0,1]$. Their $L_{p}(B)$-norms
are also bounded for any ball $B$ and,
hence, there is a sequence $\varepsilon_{n}\downarrow0$ such that
these functions converge weakly
in $L_{p}(B)$ for any ball $B$. One knows that the weak limits of $D^{2}u_{\varepsilon_{n}},D  u_{\varepsilon_{n}} $ are $D^{2}u,Du$, where $u$
is the weak limit of $u_{\varepsilon_{n}}$. Since the norm of a weak limit
is less than the inflimit of the norms, $u\in E^{2}_{p,\beta}$, and by 
passing to the weak limit in
$$
\lambda
u_{\varepsilon_{n}}-\Delta u_{\varepsilon_{n}}= f^{\varepsilon_{n}}\zeta^{\varepsilon_{n}}
$$
 we get
$
\lambda u-\Delta u=f,
$
 thus proving the theorem. \qed

 \begin{theorem}
                 \label{theorem 1.21.2}
Let $p_{b}\in(1,d]$, $\rho_{b}\in(0,\infty)$, $p_{b}> \beta>1$, $1<p\leq p_{b}/\beta$. Under these conditions
there exist $N_{0}$ and $\lambda_{0}>0$,
depending only on $d,p,\beta$ and
there is  $N_{1}=N(d,p,p_{b},\beta)$ such that, if
\begin{equation}
                          \label{4.3,8}
N_{1}\bar b_{p_{b},\rho_{b}}\leq 1,
\end{equation}
then  for any $\lambda\geq \lambda_{0}
\rho_{b}^{-2}$
and $f\in E_{p,\beta;\rho_{b}}$ there exists
a unique   $u\in E^{2}_{p,\beta;\rho_{b}}$
satisfying 
$$
\cL u-\lambda u=f,
$$
 where 
$\cL u=\Delta u+b^{i}D_{i}u$.  Furthermore, for any $u\in E^{2}_{p,\beta;\rho_{b}}$ we have
\begin{equation}
                    \label{1.20.80}
 \|D^{2}u,\sqrt\lambda Du,\lambda u\|_{E_{p,\beta;\rho_{b}}}\leq N_{0}\|\cL  u-\lambda   u\|_{E_{p,\beta;\rho_{b}}}.
\end{equation}
In particular, in light of Remark
\ref{remark 10.3.1}, for $\rho_{b}\leq 1$
$$
 \|D^{2}u,\sqrt\lambda Du,\lambda u\|_{E_{p,\beta }}\leq N_{0}2^{d/p}\rho_{b}^{-\beta}\|\cL  u-\lambda   u\|_{E_{p,\beta }}.
$$

\end{theorem}

Proof. Changing scales shows that
it suffices to concentrate on $\rho_{b}=1$, in which case \eqref{1.20.80}
follows immediately from
\eqref{1.20.8} and \eqref{1.17.3}, provided $N(d,p,p_{b},\beta)\bar b_{p_{b},\rho_{b}}$
is small enough. Once the a priory 
estimate \eqref{1.20.80} is obtained
the existence of solutions is derived
by the method of continuity on the basis of
 Theorem \ref{theorem 1.20.3}.   \qed

\begin{remark}
                      \label{remark 4.8,1}
The reader have probably noticed that
adding the drift term in the Morrey space setting
only required Theorem \ref{theorem 8,23.3},
which is a simple corollary of H\"older's
inequality, whereas in the Sobolev space
setting the derivation of
Corollary \ref{corollary 1.16.1} was based
on rather deep Adams theorem.

\end{remark}
\begin{remark}
                \label{remark 1.22.5}
Note that if $b$ is bounded, one can satisfy \eqref{4.3,8} just by taking $\rho_{b}$
small enough. The same is true if $b
\in L_{d}$ since 
$$
\bar b_{p_{b},\rho_{b}}\leq \bar b_{d,\rho_{b}}
=\sup_{B\in \bB_{\rho_{b}}}\|b\|_{L_{d}(B)}.
$$

In Remark \ref{remark 1.17.1} the
quantity $\bar b_{p_{b},\rho_{b}}$
is independent of $\rho_{b}$ and is
finite for any $p_{b}<d$. This remark shows that Theorem \ref{theorem 1.21.2} is false
if $\bar b_{p_{b},\rho_{b}}$ is not sufficiently small.

\end{remark}

In \eqref{4.3,8} 
we have
$b\in L_{p_{b},\loc}$. It might happen
the $b$ is not summable to a much higher
power. We saw this in Example \ref{example 9.25.1}.
Here is another one.
\begin{example}
                 \label{example 7,29.1}
Let $d=d'+d''\geq3$, $d'\geq2$, and $|b(x)|\leq
1/|x'|$, where $x'=(x^{1},...,x^{d'})$.
Then $b\in E_{p,1}$ for any $p\in[1,d')$,
but $b\not\in L_{d',\loc}$.

\end{example}

\begin{remark}
                \label{remark 1.21.3}
If we want the solutions of 
$$
\Delta u
+b^{i}D_{i}u-\lambda u=f
$$
 with $b$ from a 
Morrey class to be bounded
and want to use Corollary \ref{corollary 1.16.1} from the Sobolev space theory, then we have to assume that $d/2<p_{b}\leq d$ ($d/2<p$ is
needed to have the functions from $W^{2}_{p}$ bounded) and, as always,
we need $\|b\|_{ E_{p_{b},1}}$ to be small enough.
In Theorem \ref{theorem 1.21.2}
about the solvability in the Morrey-Sobolev spaces
all solutions are bounded as long as
$1<\beta<2$ and $p_{b}$ can be any 
number~$>\beta $.
\end{remark}

\mysection
{Muckenhoupt's weights}

              \label{chapter 1.17.1}

\subsection{Muckenhoupt's weights}
                                 \label{section 11.20.1}

Here we, basically, follow \cite{GR_85}
and \cite{C-UMP_11}.
We use the setting of Subsection
\ref{section 06.5.1.1}, but the doubling condition
here is stronger than \eqref{11.21.2}.
Introduce $\frQ$
as the collection of (open)
$$
\cQ_{l}(x)= x+(-l^{k_{1}}/2,l^{k_{1}}/2)\times...\times
(-l^{k_{d}}/2,l^{k_{d}}/2)  ,\quad x\in\bR^{d},l>0,
$$
and suppose that
\index{doubling condition}%
\index{$B$@Sets!$\frQ$}%
\index{$B$@Sets!$\cQ_{l}(x)$}% 
for any $\cQ_{l}(x)\in \frQ$ 
we have 
\begin{equation}
                          \label{11.21.3}
0<\mu( \cQ_{2l}(x) )\leq N_{0}\mu( \cQ_{l}(x) ).
\end{equation}
It is not hard to see that \eqref{11.21.3} implies \eqref{11.21.2}
(with a different $N_{0}$).
 Thus, the real analytic structure of $\bR^{d}$
is fixed and in this subsection
 we will not mention it while saying
what our constants called $N$ depend on.

\begin{definition}[Muckenhoupt's weights]
                              \label{definition 11.17.1}
Let $w(x)$ be a  function
on $\bR^{d}$ such that $0< w<\infty$ ($\mu$-a.e.)
and $w_{\cQ}<\infty$ for any $\cQ\in\frQ$.
We call it an $A_{1}$-weight 
(relative to $\frQ$ or relative to $(\frQ,\mu)$) if there is a constant
$N$ such that
\begin{equation}
                        \label{11.19.01}
w_{\cQ}\leq Nw(x)\quad\forall x\in \cQ, \forall \cQ\in\frQ . 
\end{equation}
  The least constant $N$ satisfying 
\eqref{11.19.01} is called the $A_{1}$-constant of $w$
  denoted by $[w]_{A_{1}}$. For $p\in(1,\infty)$
we call $w$ an $A_{p}$-weight
\index{$A_{p}$-weight}%
 if there is a constant
$N$ such that
\begin{equation}
                           \label{11.17.02}
w_{\cQ}\Big(\big(w^{-1/(p-1)}\big)_{\cQ}\Big)^{p-1}\leq N\quad \forall \cQ\in\frQ. 
\end{equation}
  The least constant $N$ satisfying 
\eqref{11.17.02} is called the $A_{p}$-constant of $w$ and 
 denoted by $[w]_{A_{p}}$.

\end{definition}

In more explicit form \eqref{11.17.02} means
that
$$
\dashint_{\cQ}w\,\mu(dx)\Big(\dashint_{\cQ}w^{-1/(p-1)}\,\mu(dx)\Big)^{p-1}\leq N.
$$

 \begin{remark}
                   \label{remark 1.20.1}
It is known that, if $\mu$
is  Lebesgue measure and $k_{1}=...=k_{d}=1$, then $|x|^{\alpha}$ is an $A_{p}$-weight 
($p\in[1,\infty)$)
 if and only if $-d<\alpha<(p-1)d$. 
As is easy to see, the minimum of two
$A_{1}$-weights is also an $A_{1}$-weight.
In particular, $|x|^{\alpha}$ is an $A_{1}$-weight if and only if $-d<\alpha<0$
and then for any constant $c>0$
the function $w(x)=|x|^{\alpha}\wedge c$
is also an $A_{1}$-weight.
\end{remark}

Set 
$$
 w(A)=\int_{A}w\,\mu(dx), 
$$
and 
\index{$A$@Sets of functions!$L_{p}(w)$}% 
denote $L_{p}(w)=L_{p}(\bR^{d},w(dx))$.

\begin{remark}
                                       \label{remark 11.25.1}
H\"older's inequality implies that, if $p\in(1,\infty)$,
   $S $ is Borel, and $\infty>w>0$ ($\mu$-a.e.) on $S$, then  
\begin{equation}
                                    \label{11.25.2}
1=\dashint_{S}w^{1/p}w^{-1/p}\,\mu(dx)
\leq (w_{S})^{1/p}\Big(\big(w^{-1/(p-1)}\big)_{S}\Big)^{(p-1)/p}. 
\end{equation}
Therefore, $[w]_{A_{p}}\geq1$ if $w\in A_{p}$.
 This  obviously holds
for $p=1$ as well.
 
Also note that one can replace what is raised to the
power $p-1$ in \eqref{11.17.02} with
$\big(w^{-1/(p-1)}\big)_{S}\mu(S)/\mu(\cQ)$
for any Borel $S\subset \cQ$. Then \eqref{11.25.2}
implies that 
\begin{equation}
                           \label{11.17.4}
[w]_{A_{p}}\frac{w(S)}{w(\cQ)}
\geq\Big(\frac{\mu(S)}{\mu(\cQ)}\Big)^{p}.
\end{equation}

We obtained this result for $p>1$. It is also true
for $p=1$ since, for $w\in A_{1}$ we have $[w]_{A_{1}}
w\geq w_{\cQ}$
on $\cQ$ and hence
$$
\mu(S)=\int_{S}\frac{1}{w}w\,\mu(dx)
\leq [w]_{A_{1}} \frac{1}{w_{\cQ}}w(S),
$$
which is \eqref{11.17.4} for $p=1$.

\begin{corollary}
                 \label{corollary 11.17.1}
The measure $w$ satisfies the doubling
condition: 
$$
w( \cQ_{2l}(x) )\leq  w( \cQ_{l}(x) )[w]_{A_{p}}
\Big(\frac{\mu(\cQ_{l}(x))}{\mu(\cQ_{2l}(x) )}\Big)^{-p}\leq w( \cQ_{l}(x) )[w]_{A_{p}}N_{0}^{ p}.
$$

\end{corollary}

Another important consequence of \eqref{11.17.4}
is that for $p\in[1,\infty)$, $w\in A_{p}$, and
 any $\alpha\in(0,1)$ there exists $\beta\in(0,1)$
depending only on $\alpha$, $p$, and $[w]_{A_{p}}$, such that, 
for any   Borel $S\subset \cQ\in\frQ$ 
\begin{equation}
                                    \label{11.24.6}
\mu(S)\leq\alpha\mu(\cQ)\Longrightarrow
 w(S)\leq\beta w(\cQ).
\end{equation}
On gets \eqref{11.24.6} by inspecting the following
version of \eqref{11.17.4}
$$
[w]_{A_{p}}\Big(1-\frac{w(S)}{w(\cQ)}\Big)
\geq\Big(1-\frac{\mu(S)}{\mu(\cQ)}\Big)^{p}\geq (1-\alpha)^{p}.
$$
\end{remark}

\begin{remark}
                                      \label{remark 12.22.1}
In the future we several times say that
a constant entering an estimate
 depends only on ..., and $[w]_{A_{p}}$.
It is important to emphasize that, given a 
$K_{0}$, in these situations
the constant can be chosen  to
depend  only on ..., and $K_{0}$, provided
that $[w]_{A_{p}}\leq K_{0}$.

\end{remark}

As a simple exercise based on   H\"older's inequality one proves that
$A_{p}\subset A_{q}$ if $1< p\leq q$.
Here is an extension of this for $A_{1}$-weights
(when $v\equiv1$).

\begin{lemma}
                                 \label{lemma 11.19.1}
If $w ,v\in A_{1}$ and $p\in(1,\infty)$, then
$w v^{1-p}\in A_{p}$ and
$[w v^{1-p}]_{A_{p}}\leq[w] _{A_{1}}[v]^{p-1}_{A_{1}}$.
\end{lemma}

Proof. In light of \eqref{11.19.01} for any $\cQ\in\frQ$ 
$$
(w v ^{1-p})_{\cQ}\leq [v]_{A_{1}}^{p-1}(w (v_{\cQ})^{1-p})_{\cQ}
=[v]_{A_{1}}^{p-1}w_{\cQ}(v_{\cQ})^{1-p},
$$
$$
(w^{-1/(p-1)}v)_{\cQ}\leq [w]_{A_{1}}^{ 1/(p-1)}
((w_{\cQ})^{ -1/(p-1)}v)_{\cQ}
=[w]_{A_{1}}^{ 1/(p-1)}(w_{\cQ})^{-1/(p-1)}v_{\cQ}, 
$$
and our assertions follow. 
\qed

Here is a deep result which we prove at the end of the section.

\begin{theorem}[P. Jones]
                   \label{theorem 2.1.1}
For $p\in(1,\infty)$ any $A_{p}$-weight 
\index{Jones theorem}%
admits a representation
as $w v^{1-p}$ with $w,v\in A_{1}$
such that $[v,w]_{A_{1}}\leq N(d,p)[w]_{A_{p}}$.

\end{theorem}
 
\begin{theorem}[reverse H\"older's inequality]
                                 \label{theorem 11.17.3}
If $p\in[1,\infty)$ 
\index{reverse H\"older's inequality}%
and $w \in A_{p}$, then
there exist constants
$N=N(p,N_{0},[w]_{A_{p}})<\infty$ and
 $\varepsilon=\varepsilon(p,N_{0},[w]_{A_{p}})>0$
such that
for any  
 $\cQ \in  \frQ $  we have
\begin{equation}
                                          \label{11.25.4}
  (w^{1+\varepsilon})_{\cQ}\leq N
(w_{\cQ})^{1+\varepsilon}.
\end{equation}
\end{theorem}

Proof. Since $A_{1}\subset A_{p}$
for any $p>1$, we may assume that $p>1$.
We get an equivalent statement if
we take
$$
\cB^{o}_{l,x}:=x+(0,l ^{ k_{1}} )\times
...\times (0,l ^{ k_{d} })
$$
in place of $\cQ$ in \eqref{11.25.4}. 

  Then we observe that  $I_{\cB^{o}_{l,x}}$ is the pointwise
limit as $\varepsilon\downarrow 0$ of the indicators
of $\cB_{l,x(1+\varepsilon)}$, so that we only need
to prove \eqref{11.25.4} for $\cQ=\cB_{l,x(1+\varepsilon)}$. 
In that case 
 we can take $x(1+\varepsilon)$ as the new origin
and concentrate on proving \eqref{11.25.4}
for $\cQ=\cB_{l,0}$. 
Finally, by perhaps differently changing scales on different coordinate axes, we reduce
the general case to the one with $l=1$.
Then we are going to use
the results of Subsection \ref{section 06.5.1.1}.

Denote $\lambda=2N_{0}$   and for $k=0,1,...$
  define $\bar w=wI_{\cQ}$ and
$$
\tau_{k}(x)=\inf\{n\in\bZ:\bar w_{|n}(x)>\lambda^{k}w_{\cQ}
   \}.
$$
 Note that $\tau_{k+1}\geq\tau_{k}$ 
 and, obviously, $\tau_{0}=\infty$ outside $\cQ$.

Also observe that if $\cB\in \cF_{\tau_{k}}\cap \frB$,  
 then by the first inequality in
\eqref{06.5.1.3}
$$
I:=\int_{\bR^{d}}\bar w I_{\cB,\tau_{k+1}<\infty}\,\mu(dx)
\leq \int_{\bR^{d}}\bar w I_{\cB,\tau_{k }<\infty}\,\mu(dx)
\leq N_{0}\lambda^{k}w_{\cQ}\mu(\cB).
$$
On the other hand, by the second part of \eqref{06.5.1.3} 
$$
I\geq  \lambda^{k+1}w_{\cQ}\mu(\cB\cap\{\tau_{k+1 }<\infty\}).
$$
Thus,
$$
\frac{\mu(\cB\cap\{\tau_{k+1 }<\infty\})}{\mu(\cB)}\leq\frac{1}{2}.
$$
By \eqref{11.24.6}
$$
w(\cB\cap\{\tau_{k+1 }<\infty\})\leq \beta w(\cB) 
$$
and, since $\{\tau_{k}<\infty\}\in \cF_{\tau_{k}}$,
the set $\{\tau_{k}<\infty\}$ is represented as the
disjoint union of some $\cB\in \cF_{\tau_{k}}\cap \frB$, so that
$$
w( \{\tau_{k+1 }<\infty\}) \leq \beta w
(\{\tau_{k }<\infty\}).
$$

It follows that  
$$
w( \{\tau_{k  }<\infty\})\leq \beta^{k}
w( \{\tau_{0}<\infty\})=\beta^{k} w(\cQ)
=\mu(\cQ)\beta^{k}w_{\cQ}. 
$$

Now, observe that ($\mu$-a.e.) on the set
$\{\tau_{k}<\infty,\tau_{k+1}=\infty\}$
we have $w\leq \lambda^{k+1}w_{\cQ}$. Also
($\mu$-a.e.)
$$
\bigcup_{k=0}^{\infty}\{\tau_{k}<\infty,\tau_{k+1}=\infty\}
=\{\tau_{0}<\infty \}
$$
since $w<\infty$ ($\mu$-a.e.). Obviously, 
the terms in the above union
are disjoint. Therefore,
$$
\int_{\cQ}w^{1+\varepsilon}\,\mu(dx)
=\int_{\cQ}w^{1+\varepsilon}I_{\tau_{0}=\infty}\,\mu(dx)
+\sum_{k=0}^{\infty}\int_{\cQ}
w^{1+\varepsilon}I_{\tau_{k}<\infty,\tau_{k+1}=\infty}\,\mu(dx)
$$
$$
\leq \mu(\cQ) (w_{\cQ})^{1+\varepsilon}
+\sum_{k=0}^{\infty}\lambda^{(k+1) \varepsilon }
 (w_{\cQ})^{\varepsilon}
\int_{\cQ}
w I_{\tau_{k}<\infty,
\tau_{k+1}=\infty}\,\mu(dx)
$$
$$
\leq \mu(\cQ) (w_{\cQ})^{1+\varepsilon}
\Big[1+\sum_{k=0}^{\infty}\lambda^{(k+1) \varepsilon }
\beta^{k}\Big].
$$
We see that it only remains to choose $\varepsilon>0$
so small that the last series converges. \qed

The following somewhat sharper statement than  
\eqref{11.24.6} and converse to \eqref{11.17.4} 
will be needed to prove
the Fefferman-Stein theorem in $L_{p}(w)$
spaces.

\begin{corollary}
                    \label{corollary 11.21.1}
For any $p\in[1,\infty)$
and $w \in A_{p}$,  there exists
$\beta\in(0,1)$ and $N$,
both depending only on $ p,N_{0},[w]_{p} $, such that
for any Borel $S\subset \cQ\in\frQ$ we have
\begin{equation}
                            \label{11.24.2}
\frac{w(S)}{w(\cQ)}\leq N\Big(\frac{\mu(S)}{\mu(\cQ)}\Big)^{\beta}.
\end{equation}

\end{corollary}

Indeed, by H\"older's inequality
$$
w(S)=\int_{S}w \,\mu(dx)
\leq \Big(\int_{\cQ}w^{1+\varepsilon}
\,\mu(dx)\Big)^{1/(1+\varepsilon)}\,
\mu^{\varepsilon/(1+\varepsilon)}(S) 
$$
$$
=  \mu^{1/(1+\varepsilon)}(\cQ)\big[\big(w^{1+\varepsilon}\big)_{\cQ}\big]
^{1/(1+\varepsilon)}\mu^{\varepsilon/(1+\varepsilon)}(S) 
$$
$$
\leq N\mu^{1/(1+\varepsilon)}(\cQ)
 w_{\cQ} \mu^{\varepsilon/(1+\varepsilon}(S)
=Nw(\cQ)\mu^{-\varepsilon/(1+\varepsilon)}(\cQ)
\mu^{\varepsilon/(1+\varepsilon)}(S), 
$$
which is what is claimed with $\beta=\varepsilon/(1+\varepsilon)$.

\begin{corollary}[cf.~\eqref{11.28.1}]
                  \label{corollary 1.19.5}
Let $\mu$ be Lebesgue measure,   $p\in[1,\infty)
$, $w\in A_{p}$. Then, for any $x\in\bR^{d}$
$$
\lim_{l\to\infty} w(\cB_{l,x} )=\infty.
$$

\end{corollary}

Indeed, \eqref{11.24.2} implies that for 
any integers $r,n$
$$ w(\cB_{r^{n+1},x} )\geq w(\cB_{r^{n},x} )
N^{-1}\Big(\frac{\mu(\cB_{r^{n+1}})}{\mu(\cB_{r^{n}})}
\Big)^{\beta}=w(\cB_{r^{n},x} )
N^{-1}r^{ \beta \gamma},
$$
where $\gamma=k_{1}+...+k_{d}$ and for large $r$ the factor of $w(\cB_{r^{n},x} )$ becomes greater than 2.

For the reader's orientation we point out that
the weights satisfying \eqref{11.24.2}
\index{$A_{\infty}$-weight}%
are called $A_{\infty}$-weights.
 We are not going to use the remarkable
fact that 
$$
A_{\infty}=\bigcup_{p\in[1,\infty)}A_{p}.
$$ 
This fact follows from the observation that, actually,
in the proof of Theorem \ref{theorem 11.17.3}
and of \eqref{11.24.2}
only the implication
\eqref{11.24.6} was used, which turns out to be
(almost trivially) reversible.

The following result is used not only to prove
the Hardy-Littlewood weighted theorem but also
in many places in applications to partial
differential equations with almost VMO
leading coefficients.

\begin{theorem}[Self improving property]
                                 \label{theorem 11.18.3}
If $p\in(1,\infty)$
and $w \in A_{p}$, then there exists 
$$
q=q(p,N_{0},
[w]_{A_{p}})\in (1,p) 
$$
 such that $w\in A_{q}$.
Furthermore, $[w]_{A_{q}}$ is dominated by a constant
$N=N(p,N_{0},
[w]_{A_{p}})$.
\end{theorem}

Proof. Note that for $p\in(1,\infty)$
the condition $w\in A_{p}$ is equivalent to
$w^{-1/(p-1)}\in A_{p/(p-1)}$ and
$$
[w^{-1/(p-1)}]_{A_{p/(p-1)}}=[w]_{A_{p}}^{1/(p-1)}.
$$
By Theorem \ref{theorem 11.17.3} 
there exist $\varepsilon>0$
and $N_{1}$ depending only on $p,N_{0},
[w]_{A_{p}}$ such that
$$
\big(w^{-(1+\varepsilon)/(p-1)}\big)_{\cQ}
\leq N_{1}\big(\big(w^{-1/(p-1)})_{\cQ}\big)^{1+\varepsilon}.
$$
Obviously, $(1+\varepsilon)/(p-1)=1/(q-1)$
for some $q\in (1,p)$ and the above inequality
means that
$$
\big(\big(w^{-1/(q-1)})_{\cQ}\big)^{q-1}
\leq N^{q-1}_{1}\big(\big(w^{-1/(p-1)})_{\cQ}\big)^{p-1}.
$$
By multiplying both parts by $w_{\cQ}$, we get that
$w\in A_{q}$ and $[w]_{A_{q}}\leq N^{q-1}_{1}[w]_{A_{p}}$. \qed

Define the weighted Hardy-Littlewood maximal
\index{$S$@Miscelenea!weighted Hardy-Littlewood maximal operator}% 
\index{$C$@Operators!$\bM_{w}$}%
 operator. For an $A_{p}$-weight $w$ set
$$
\bM_{w} f(x)=\sup_{\cQ\in\frQ }I_{\cQ}(x)\frac{1}{w(\cQ)}
\int_{\cQ}|f|\,w(dy). 
$$

\begin{lemma}
                                 \label{lemma 11.17.2}
Let $p\in[1,\infty)$ and $w \in A_{p}$. Then
for any $\cQ\in\frQ$ and measurable $f\geq 0$ we have
\begin{equation}
                           \label{11.17.03}
 (f_{\cQ})^{p}w_{\cQ}\leq [w]_{A_{p}}(f^{p}w)_{\cQ}. 
\end{equation}
\end{lemma}

Proof. If  $p=1$, \eqref{11.17.03} follows from 
\eqref{11.17.4} and the linearity of the integral.
If $p>1$, by H\"older's inequality
$$ 
(f_{\cQ})^{p}w_{\cQ}=\big((fw^{1/p})w^{-1/p}\big)_{\cQ}^{p}w_{\cQ}
\leq (f^{p}w)_{\cQ}\Big(\big(w^{-1/(p-1)}\big)_{\cQ}\Big)^{p-1}
w_{\cQ}
$$
and \eqref{11.17.03} follows in light of our definitions. \qed

Observe that \eqref{11.17.03} implies the following.
\begin{corollary}
                                     \label{corollary 11.21.4}
For $p\in[1,\infty)$, $w\in A_{p}$, and Borel $f\geq0$
we have
$$
(\bM f)^{p}\leq  N(p) [w]_{A_{p}}\bM_{w}( f ^{p}).
$$  

\end{corollary}

Here is one of the fundamental results of the theory
of weights. Recall that
$q=q(p,N_{0},
[w]_{A_{p}})\in (1,p)$ is introduced
\index{$S$@Miscelenea!$\tau_{0}(p,N_{0},
[w]_{A_{p}})$}%
in Theorem \ref{theorem 11.18.3} and set
$$
\tau_{0}=\tau_{0}(p,N_{0},
[w]_{A_{p}})=q/p\in(0,1)
$$

\begin{theorem}[Muckenhoupt]
                      \label{theorem 11.21.3}
  If $p\in(1,\infty)$, $w\in A_{p}$
 and $\tau\in[\tau_{0},1]$, then  
$ \bM f $ is  a bounded operator in $L_{ \tau p}(w)$,
 that is there exists
a constant $N$, depending only on $
 \tau ,p,N_{0}$,  and 
$[w]_{A_{p}}$, such that for any   
$f\in L_{ \tau p}(w)$ we have
$$
\| (\bM f)^{\tau} \|_{L_{p}(w)}\leq N\|   |f|^{\tau} \|_{L_{p}(w)}.
$$

\end{theorem}

Proof.    According to Corollary \ref{corollary 11.17.1}
the measure $w$ satisfies the doubling condition and  
by Theorem \ref{theorem 11.18.3} we have 
$w\in A_{\tau_{ 0} p}$   and
then $w\in A_{\tau_{1}p}$ with $\tau_{1}<\tau_{0},
\tau_{1}p>1$.
 Theorem \ref{theorem 7.7.1}
applied to $w$ in place of $\mu$ implies that
for any   $f\geq0$,
$$
\|\bM_{w} f  \|_{L_{ \tau/\tau_{1}}(w)}
\leq N\|f \|_{L_{ \tau/\tau_{1}}(w)}. 
$$
Also $ (\bM f)^{\tau p}\leq N
[\bM_{w}(|f|^{\tau_{1} p})]^{\tau/\tau_{1}}  $
by Corollary \ref{corollary 11.21.4}
and   the desired result follows.
\qed
 
It is worth noting that in \cite{Le_08}
one can find an elegant and short proof of Theorem
\ref{theorem 11.21.3}  
not using Theorem \ref{theorem 11.18.3}.

Here is one more version of the
Fefferman-Stein theorem.

\begin{theorem}[Dong-Kim \cite{DK_18}]
                        \label{theorem 1.19.2} 

If $p\in(1,\infty)$ and $w\in A_{p}$, then    there exists
a constant $N$, depending only on $p,N_{0}$, and 
$[w]_{A_{p}}$, such that for any   
$f\in L_{p}(w)$ we have
\begin{equation}
                               \label{1.19.7}
\| f\|_{L_{p}(w)}\leq N\|f^{\#}\|_{L_{p}(w)}.
\end{equation}
\end{theorem}

To prove the theorem we need  a lemma from \cite{DK_18} similar to Lemma 
\ref{lemma 06.11,23,1} in the proof of which
we follow \cite{DK_18}. The constant $\beta
\in(0,1)$
is taken from \eqref{11.24.2}.
\begin{lemma}
                  \label{lemma 1.19.3}
For $\alpha=(2N_{0})^{-1}$,
 any constant $c>0$, and  $f\in L_{1} $, we have
 \begin{equation}
                             \label{1.19.6}
w(|f|\geq c)\leq \frac{4^{\beta}}{c^{\beta}}\int_{\bR^{d}}I_{\cM f(x)>\alpha c}
(f^{\#}(x))^{\beta}\,w(dx)
\end{equation}
and, if $f\geq0$, then one can replace $4/c$ with $2/c$.
\end{lemma}

Proof. Remark \ref{remark 9.3.2} shows that
it suffices to prove the second assertion of the lemma. For a while we follow the proof
of Lemma 
\ref{lemma 06.11,23,1}.
Introduce
$$
\tau(x)=\inf\{n: f_{|n}(x)>c\alpha\}.
$$
Use Lemma \ref{lemma 06.5.30.1} (ii) to get that
$f_{|\tau}\leq c/2$ if $\tau<\infty$ and also use
  the fact that 
 $f_{|n }\to f$ (a.e.).
Then we find that (a.e.)
$$
\{x: f(x) \geq c\}=
\{x: f(x) \geq c,\tau(x)<\infty\}
$$
$$=
\{x: f(x) \geq c,  f_{|\tau}(x)\leq c/2\}
\subset\{x:  
|f(x)-f_{|\tau}(x)| \geq c/2\} .
$$

Observe that for each $n$ the set
$$
A:=\{x:  \tau(x)=n,
|f(x)-f_{|\tau}(x)| \geq c/2\}
$$
 is the union
of some $\cB\in \frB_{n}$ and by Chebyshev's inequality for such $\cB$
$$
\mu(A\cap \cB)\leq (2/c)\int_{\cB}|f -f_{|n} |
\mu(dx)= (2/c)\mu(\cB)\dashint_{\cB}|f -f_{|n} |
\mu(dx),
$$
which by \eqref{11.24.2} implies that
$$
w(A\cap \cB)\leq  N(2/c)^{\beta}
\Big(\dashint_{\cB}|f -f_{|n} |\,\mu(dx)\Big)^{\beta}
w(\cB)
\leq  N(2/c)^{\beta}\int_{\cB}(f^{\#})^{\beta}\,w(dx).
$$
Hence,
$$
w(\{x:  \tau=n,
|f -f_{|\tau} | \geq c/2\})
\leq  N(2/c)^{\beta}\int_{\bR^{d}}
I_{\tau=n   }(f^{\#})^{\beta}\,w(dx),
$$
$$
w(\{x:   
|f -f_{|\tau} | \geq c/2\})
\leq  N(2/c)^{\beta}\int_{\bR^{d}}
I_{ \tau<\infty}(f^{\#})^{\beta}\,w(dx),
$$

Now it only remains to notice that
$$
\{\tau(x)<\infty\}
=\{\cM f(x)> c\alpha\}.
$$    \qed  

{\bf Proof of Theorem \ref{theorem 1.19.2}}.
We multiply \eqref{1.19.6} by $c^{p-1}$ and
integrate with respect to $c\in(0,\infty)$.
Then we get
$$
\int_{\bR^{d}}|f|^{p}\,w(dx)
\leq N\int_{\bR^{d}}(\cM f)^{p-\beta}
(f^{\#}(x))^{\beta}\,w(dx).
$$
By applying H\"older's inequality,
the estimate $\cM f\leq N \bM f$, and
Theorem \ref{theorem 11.21.3} we come
to \eqref{1.19.7} provided $f\in L_{1}
\cap L_{p}(w)$. In the general case take
$f_{n}\in L_{1}
\cap L_{p}(w)$ such that $f_{n}\to f$ in
$L_{p}(w)$. Then \eqref{1.19.7} holds
for $f_{n}$ in place of $f$. By passing to the limit while observing that 
$$
|f^{\#}-f_{n}^{\#}|\leq 2\cM(f-f_{n})
\leq N\bM(f-f_{n})
$$
in light of $||a-b|-|c-d||\leq |a-c|+|b-d|$
and applying Theorem \ref{theorem 11.21.3},
we come to \eqref{1.19.7} for $f\in L_{p}(w)$.\qed

The most general version of Theorem
\ref{theorem 1.19.2} the reader can find in \cite{DK_19}.
Here is a somewhat unexpected result.

\begin{corollary}
                     \label{corolary 1.19.8}
For any $p\in(1,\infty)$ and $w\in A_{p}$
we have $w(\bR^{d})=\infty$.
\end{corollary}

Indeed, otherwise $1\in L_{p}(w)$ and $1^{\#}=0$.

For applications in PDEs the following version
of the Fefferman-Stein theorem
\index{Fefferman-Stein theorem}%
is important, in which the case that
$\beta=d/p$ is taken care of
by Theorem \ref{theorem 11,23,1} in light of the fact that
then
$\dot E_{p, \beta}=L_{p}$. 

\begin{theorem}
                    \label{theorem 1.20.2}
Let $0<\beta\leq d/p$, $1<p<\infty$. Then,
for any $f\in \dot E_{p,\beta}$,
\begin{equation}
                          \label{1.20.3}
\|f\|_{\dot E_{p,\beta}}
\leq N(d,p,\beta)\|f^{\#}\|_{\dot E_{p,\beta}}.
\end{equation} 
\end{theorem}

Proof. We use an observation from \cite{CF_88}.
We only need to concentrate on $\beta<d/p$. Then
note that the right-hand side of \eqref{1.20.3}
is finite, since $f^{\#}\leq 2\cM f\leq N\bM f$
and we can refer to Remark \ref{remark 7,7.1}.
To finish with preparations,
we may assume that $f\geq 0$.

Take an $\alpha\in(0,d)$, such that $ 2\alpha>d$ and $\alpha>d-p\beta$, set $v(r)=r^{-\alpha}\wedge 1$,  
$w(x)=v(|x|)$, and by using Remark
\ref{remark 1.20.1} and Theorem \ref{theorem 1.19.2}  write  
$$
\dashint_{B_{1}}f^{p}\,dx\leq
N\int_{\bR^{d}}f^{p}w\,dx
\leq N\int_{\bR^{d}}(f^{\#})^{p}w\,dx.
$$
The last integral equals
$$
N\int_{0}^{\infty}v(r)
\frac{d}{dr}\Big(\int_{B_{r}}(f^{\#})^{p} \,dx\Big)
\,dr=NI.
$$
We want to perform  the  integration by parts.
To do that observe that
$$
\int_{B_{r}}(f^{\#})^{p} \,dx\leq Nr^{d-p\beta}
\Big(r^{p\beta}\dashint_{B_{r}}(f^{\#})^{p} \,dx
\Big)\leq Nr^{d-p\beta}\|f^{\#}\|^{p}_{\dot E_{p,\beta}}
$$
and since $d-p\beta>0$ the integrated out term
for $r=0$ vanishes. At infinity it also
vanishes since $d-p\beta-\alpha<0$.
Thus
$$ 
I=\alpha\int_{1}^{\infty}r^{d-\alpha-1-p\beta}
\Big(r^{p\beta}\dashint_{B_{r}}(f^{\#})^{p} \,dx
\Big)\,dr  
$$
$$
\leq 
N\|f^{\#}\|^{p}_{\dot E_{p,\beta}}
\int_{1}^{\infty}r^{d-\alpha-1-p\beta}\,dr.
$$
Since $d-\alpha-p\beta<0$, the integral converges and for $r=1$ we get
$$
r^{\beta}\dashint_{B_{r}}f^{p}\,dx\leq
N\|f^{\#}\|^{p}_{\dot E_{p,\beta}}.
$$
For any $r>0$ one obtains this inequality
by using scaling and then using shifts
of the origin one arrives at \eqref{1.20.3}.  \qed

{\bf Proof of Theorem \ref{theorem 2.1.1}}.
We follows the proofs of Lemma 5.1
and Theorem 5.2 of \cite{GR_85} and, first, observe
that
we need to prove the theorem only for
$ p\in(1,2]$ because for any $w\in A_{p}$,
$w^{-1/(p-1)}\in A_{p/(p-1)}$. So let
$p\in(1,2]$ and $w\in A_{p}$. We need to find
$v>0$ (a.e.) such that, for a constant $N$,
\begin{equation}
                      \label{2.6.1}
\bM (vw)\leq Nvw,\quad \bM (v^{1/(p-1)})\leq Nv^{1/(p-1)}.
\end{equation}
Indeed, then for $w_{2}=vw$ and $w_{1}=
v^{1/(p-1)}$ we have $w_{1},w_{2}\in A_{1}$
and $w=w_{1}^{1-p}w_{2}$. We will achieve
\eqref{2.6.1} if we find $v$ such that
$$
Sv:=w^{-1}\bM (vw)+[\bM (v^{1/(p-1)})]^{p-1}\leq Nv.
$$
Here the first and the second term behave
like an $L_{1}$- and an $L_{1/(p-1)}$-norm,
respectively. Therefore, they are
sublinear with respect to $v$. Furthermore, $S$ is a bounded operator
in $L_{p/(p-1)}(w)$ because $w\in A_{p}$,
$w^{-1/(p-1)}\in A_{p/ (p-1)} $ and 
$$
\int_{\bR^{d}}|w^{-1}\bM (vw)|^{p/(p-1)}
w\,dx=\int_{\bR^{d}}|\bM (vw)|^{p/(p-1)}
w^{-1/(p-1)}\,dx
$$
$$
\leq N\int_{\bR^{d}}(vw) ^{p/(p-1)}w^{-1/(p-1)}\,dx=N\int_{\bR^{d}}v^{p/(p-1)}w\,dx,
$$
$$
\int_{\bR^{d}}[\bM (v^{1/(p-1)})]^{p}
w\,dx\leq \int_{\bR^{d}}v^{p/(p-1)}w\,dx.
$$

Now the most natural candidate for
satisfying $Sv\leq Nv$ is obtained
by taking any $f\geq 0$ in
$L_{p/(p-1)}(w)$  and setting
$$
v=\sum_{n=0}^{\infty}\frac{S^{n}f}{2^{n}\|S\|^{n}}.
$$
Indeed, then $Sv\leq 2\|S\|v$.   \qed

This theorem will be used in the proof of
the Rubio De Francia extrapolation Theorem
 \ref{theorem 11.21.6}. For the same purpose we also need the following.

\begin{lemma}
                     \label{lemma 2.7.1}
Let $p\in(1,\infty)$ and $w_{p}\in A_{p}$.
Then there exists a sublinear bounded
operator $S_{p',w_{p}}:L_{p'}(w_{p})\to L_{p'}(w_{p})$,
where $p'=p/(p-1)$,
such that, for any $u\in L_{p'}(w_{p})$, we have 
$$
|u|\leq S_{p',w_{p}}u ,\quad
\|S_{p',w_{p}}u\|_{L_{p'}(w_{p})}\leq
2\| u\|_{L_{p'}(w_{p})},
$$
\begin{equation}
                           \label{4.4,2}
\bM(w_{p}S_{p',w_{p}}u)\leq 2Nw_{p}S_{p',w_{p}}u,
\end{equation}
so that $w_{p}S_{p',w_{p}}u\in A_{1}$, where
$N$ is the norm in $L_{p'}(w_{p})$ 
of the operator $T$ given by
$T u:=w_{p}^{-1}\bM(uw_{p})$.
\end{lemma}

The proof of the lemma is easily
extracted from the proof of Theorem 
\ref{theorem 2.1.1}. Indeed, we saw that
 $T$ is  a sublinear bounded
operator in $L_{p'}(w_{p})$ and
$$
S_{p',w_{p}}u:=\sum_{n=0}^{\infty}\frac{T^{n}u}
{2^{n}\|T\|^{n}}
$$
is  a sublinear bounded
operator in $L_{p'}(w_{p})$
satisfying $|u|\leq S_{p',w_{p}}u$.
Then, obviously $TS_{p',w_{p}}u\leq 2\|T\|S_{p',w_{p}}u$, that is \eqref{4.4,2}. \qed

\begin{remark}
                    \label{remark 2.9.1}
It is worth emphasizing that the
$A_{1}$-constant of $w_{p}S_{p',w_{p}}u$
is {\em independent\/} of $u$.
\end{remark}

\subsection{Application to the Laplacian}  

In this section the real analytic structure of $\bR^{d}$ is given by
  Lebesgue
measure $\mu$  and $k_{1}=...=k_{d}=1$.
First we apply Theorem \ref{theorem 1.19.2}.

\begin{lemma}
               \label{lemma 1.21.4}
For any $p\in(1,\infty)$, $w\in A_{p}$,
and $u\in C^{\infty}_{0}$ we have
\begin{equation}
                    \label{1.21.5}
\|D^{2}u\|_{L_{p}(w)}\leq N(d,p,[w]_{p})\|\Delta u\|_{L_{p}(w)}.
\end{equation}
\end{lemma}

Indeed, for any $r\in(1,p)$ close
to $1$ such that $w$ is an $A_{p/r}$-weight
$$
\|D^{2}u\|_{L_{p}(w)}\leq N
\|(D^{2}u)^{\sharp}\|_{L_{p}(w)}
\leq N\|\bM^{1/r}(|\Delta u|^{r}) \|_{L_{p}(w)}\leq N\|\Delta u\|_{L_{p}(w)}.
$$

\begin{lemma}
               \label{lemma 1.22.3}
For any $u\in C^{\infty}_{0}$
and $\varepsilon>0$ we have
\begin{equation}
                    \label{1.22.1}
(Du)^{\sharp}\leq N( d) \big(\varepsilon \bM(D^{2}u)+\varepsilon^{-1}\bM u\big)
\end{equation}
and, if $p\in(1,\infty)$ and $w\in A_{p}$, then
\begin{equation}
                    \label{1.22.2}
\|Du\|_{L_{p}(w)}\leq  N(d,p,[w]_{p}) \big(\varepsilon\|D^{2}u\|_{L_{p}(w)}
+\varepsilon^{-1}\|u\|_{L_{p}(w)}\big).
\end{equation}

\end{lemma}

Here \eqref{1.22.1} follows from Lemma \ref{lemma 2.1.1} 
 and \eqref{1.22.1}
  implies \eqref{1.22.2} as a few times above.

\begin{remark}
                \label{remark 1.27.3}
Estimate \eqref{1.22.1} allows us
to obtain also the interpolation inequality
for Morrey spaces as in
 Theorem \ref{theorem 12.26.1} differently, albeit
only for $p>1$ and in the whole space
but with any $\varepsilon>0$.

Indeed, take an $\alpha\in(0,d)$, such that $ 2\alpha>d$ and $\alpha>d-p\beta$, set $v(r)=r^{-\alpha}\wedge 1$, 
$w(x)=v(|x|)$, and proceed as in the proof
of Theorem \ref{theorem 1.20.2}. Then we find
$$
\dashint_{B_{1}}|Du|^{p}\,dx
\leq N\int_{\bR^{d}}|Du|^{p}w\,dx
\leq N\int_{\bR^{d}}((Du)^{\sharp})^{p}w\,dx 
$$
$$
\leq N\int_{\bR^{d}}[\varepsilon \bM(D^{2}u)
+\varepsilon^{-1}\bM u
]^{p}w\,dx\leq N\int_{\bR^{d}}[\varepsilon  (D^{2}u)
 ^{p}+\varepsilon^{-1} |u|^{p}]w\,dx .
$$
As in the proof of Theorem \ref{theorem 1.20.2} the last term is dominated by
$$
N[\varepsilon \|D^{2}u\|_{\dot E_{p,\beta}}
+\varepsilon^{-1}\| u\|_{\dot E_{p,\beta}}]
$$
and  again arguing as in the proof of Theorem \ref{theorem 1.20.2} we get that 
$$
\|Du\|_{\dot E_{p,\beta}}\leq
N[\varepsilon \|D^{2}u\|_{\dot E_{p,\beta}}
+\varepsilon^{-1}\| u\|_{\dot E_{p,\beta}}]
$$
as long as  $p\in(1,\infty)$, $0<\beta\leq d/p$,
and $u\in W^{2}_{p} $. 

However, in the whole space the interpolation
inequalities admit much simpler proofs
outlined, for instance, in Lemma \ref{lemma 3.31.4}. 
 
\end{remark}
 
 For $w\in A_{p}$ introduce the space $W^{2}_{p}(w)$ as
\index{$A$@Sets of functions!$W^{2}_{p}(w)$}% 
 the set of
$u\in L_{p}(w)$ such that  (Sobolev)  $Du,D^{2}u
\in L_{p}(w)$ and provide it with
a natural norm. It is not hard to see
that $C^{\infty}_{0}$ is dense
in $W^{2}_{p}(w)$ and, therefore,
\eqref{1.21.5} and \eqref{1.22.2} hold  for any
$u\in W^{2}_{p}(w)$. In particular,
for obvious mollifiers $u^{\varepsilon}$
of $u$ we have $u^{\varepsilon}\to u$
(a.e.) and $\sup_{\varepsilon}| u^{\varepsilon}|\leq N\bM u\in L_{p}(w)$,
so that $u^{\varepsilon}\to u$
in $L_{p}(w)$ by the dominated convergence theorem.

Then follows an analog of part of Theorem
\ref{theorem 12.22.1} with similar proof.

\begin{theorem}
                      \label{theorem 1.21.4} 
For any $p\in(1,\infty)$, $w\in A_{p}$,     there exists a constant
$N=N(d,p,[w]_{A_{p}})$ such that for any $u\in W^{2}_{p}(w)$
and $\lambda\geq0$
\begin{equation}
                          \label{12.22.10}
\|D^{2}u,\sqrt\lambda Du,\lambda u\|_{L_{p}(w)}\leq
N\|\Delta u-\lambda u\|_{L_{p}(w)}.
\end{equation}

\end{theorem}

Proof.  
 In the space
 $$
\bR^{d+1}=\{z=(x,y):y\in\bR,x\in\bR^{d}\}
 $$
introduce   $w'(x,y)=w(x)$ and 
boxes
$$
C'_{l,z}=z+[0,l   )^{d+1} ,\quad z\in\bR^{d+1},l>0.
$$
Obviously, $w'$ is an $A_{p}$-weight
with the same $[w]_{A_{p}}$ and \eqref{1.21.5} is applicable to $\Delta'=
\Delta_{x}+\partial^{2}/(\partial y)^{2}$ 
 and the function
 $$
u'(z)=u(x)\zeta(y)\cos(\mu y),
 $$
 where
 $\mu=\sqrt{\lambda}$ and $\zeta$ is a $C^{\infty}_{0}(\bR)$ function,
$\zeta\not\equiv0$.

 In light of \eqref{9.4.3} applied to
$u'$ and $ \Delta'$ we get
\begin{equation}
                                                     \label{06.6.27.90}
\|D^{2}_{x}u' ,D^{2}_{y}u' \|_{L_{p}(w')}
\leq N\| \Delta'u' \|_{L_{p}(w')}.
\end{equation}

As we have mentioned before,
$$
\int_{\bR}|\zeta(y)\cos(\mu y)|^{p}\,dy
$$
is bounded   away from zero for $\mu\in\bR$. 
It follows that

$$
\|u'\|^{p}_{L_{p}(w')}
=\mu^{-2p} \Big( \int_{\bR}
|\zeta(y)\cos(\mu y)|^{p}\,dy \Big) ^{-1}
\int_{\bR^{d+1}}\Big|D^{2}_{y}u'(z) 
$$
$$
-u(x)[\zeta''(y)\cos(\mu y)-2\mu\zeta'(y)\sin(\mu y)
 ] \Big| ^{p}w(x)\,dxdy 
$$
$$
\leq N\mu^{-2p}\big(\|D^{2}_{y}u'\|^{p}_{L_{p}(w')}
+(\mu^{p}+1)\|u \|^{p}_{L_{p}(w) }\big).
$$
This and \eqref{06.6.27.90} yield 
$$
\mu^{2} \|u \|_{L_{p}(w) }\leq 
N\| \Delta'u' \|_{L_{p}(w') }
+N(\mu+1)\|u \|_{L_{p}(w) }.
$$
Since
$$
\Delta'u'=\zeta \cos(\mu y)
[\Delta u-\lambda u]+u[\zeta''\cos(\mu y)-2\mu\zeta'\sin(\mu y)],
$$
we have
$$
\|\Delta'u' \|_{L_{p}(w') }\leq
N\|\Delta u-\lambda u\|_{L_{p}(w) }
+  N(\mu+1)\|u\|_{L_{p}(w) },
$$
so that
$$
\lambda \|u \|_{L_{p}(w) }\leq 
N\|\Delta u-\lambda u\|_{L_{p}(w) }
+  N(\sqrt{\lambda}+1)\|u\|_{L_{p}(w) }.
$$

Furthermore, \eqref{1.21.5} implies that
$$
\|D^{2}u\|_{L_{p}(w)}\leq N\|\Delta u-\lambda u\|_{L_{p}(w)}
+N\lambda \|u\|_{L_{p}(w)},
$$
and by \eqref{1.22.2}
$$
\sqrt\lambda\|Du\|_{L_{p}(w)} 
\leq \|D^{2}u\|_{L_{p}(w)}+N\lambda
\|u\|_{L_{p}(w)}.
$$
Thus,
$$
\|D^{2}u\|_{L_{p}(w)}+\sqrt\lambda\|Du\|_{L_{p}(w)}+\lambda\|u\|_{L_{p}(w)} 
$$
$$
\leq N\|\Delta u-\lambda u\|_{L_{p}(w)}+
N_{1}(\sqrt\lambda+1)\|u\|_{L_{p}(w)}.
$$
For {\em fixed\/} $\lambda$ such that
$2N_{1}(\sqrt\lambda+1)\leq \lambda$ we obtain
\eqref{12.22.10}. After that for general
$\lambda$ estimate \eqref{12.22.10}
is obtained by using scaling
based on the fact that for any $c>0$
the weight $w(cx)$ is also an $A_{p}$-weight with the same $[w]_{A_{p}}$.
\qed

\begin{theorem}
                \label{theorem 1.22.2}
For any $p\in(1,\infty)$, $w\in A_{p}$,
$\lambda>0$, and 
$f\in L_{p}(w)$, there is a unique
$u\in W^{2}_{p}(w)$ such that $\Delta u
-\lambda u=f$.

\end{theorem}  

Proof. Recall that
$$
\sfG_{\lambda}(x)=\int_{0}^{\infty}\frac{1}{(4\pi t)^{d/2}}e^{-|x|^{2}/(4t)-\lambda t}\,dt
$$
and for $f\in C^{\infty}_{0}$ we set
$u= \sfG_{\lambda}*f$. It is easy to check
that
$u$ is infinitely differentiable, each
of its derivatives decays as $|x|\to\infty$ exponentially fast and,
since according to the doubling condition $w(B_{r})$ grows at most
polynomially as $r\to \infty$,
$u\in W^{2}_{p}(w)$. Furthermore,
$$
\lambda u-\Delta u=f.
$$
 After that the solvability of $\lambda u-\Delta u=f$
for general $f\in L_{p}(w)$ follows
from the denseness of $C^{\infty}_{0}$
in $L_{p}(w)$ and the estimates from Theorem \ref{theorem 1.21.4}.
  \qed
  
As in the case of $w=1$ to include the first-order terms it is convenient to use
Adams's theorem adjusted to the case
of general weights.

\begin{theorem}
                      \label{theorem 5.28,1}

Assume  $0< \alpha<d$,  $1<p<q $, $b\geq0$, $w\in A_{p}, 
w^{s}\in A_{r} $,  where $s=q/(q-p),r=p(q-1)(q-p)$, $f\in L_{p}(w)$,  and
\index{Adams's theorem}
\begin{equation}
                                \label{5.28,4}
\|b\|_{\dot E_{q,\alpha}}=\sup_{\rho>0}\rho ^{\alpha}\sup_{B\in\bB_{\rho}}  
\Big(\dashint_{B}b^{q}\,dx\Big)^{1/q}\leq A , 
\end{equation}
\begin{equation}
                                \label{5.30,1}
J:=\int_{1}^{\infty}\frac{1}{r^{\gamma}}
\dashint_{B_{r}}w^{1/(1-p)}(x)\,dx\,dr<\infty , 
\end{equation}
where
$$
\gamma=1+\frac{d-\alpha p}{p-1}.
$$
Then
\begin{equation}
                                \label{5.28,3}
I:=\int_{\bR^{d}}b^{p}|R_{\alpha}f|^{p}w\,dx\le N(\alpha,d,p,q,
[w]_{A_{p}}, [w^{s}]_{A_{r}})A^{p}
 \int_{\bR^{d}} |f|^{p}w\,dx.
\end{equation}

\end{theorem}

Proof. We may assume that $b\not\equiv0$ and
$A<\infty$. In that case necessarily $\alpha\leq d/q$, $\alpha p<d$ and $\gamma>1$.
In particular it follows that we may assume
that $b,f$ are bounded, have compact supports,
and $f\geq0$.   Let $H$ be the subset of the unit ball
 in $L_{p'}(w)$,  $(p'=p/(p-1))$, consisting of bounded 
  nonnegative function with compact support.
Then $I^{1/p}$
equals
$$
\sup_{h\in H}\int_{\bR^{d }}h|b|R_{\alpha}fw
\,dx=\sup_{h\in H}\int_{\bR^{d}}R_{\alpha }(h|b|w)f {}
\,dx
$$
\begin{equation}
                                \label{5.30,3}
\leq \sup_{h\in H}\|R_{\alpha} (h|b|w)\|_{L_{p'}(w')}\|f\|_{L_{p}(w)},
\end{equation} 
where   $w'=w^{1/( 1-p)}$ $(\in A_{p'})$. Observe that
$$
\|R_{\alpha} (h|b|w)\|_{L_{p'}(w')}
=\lim_{n\to\infty}\|R_{\alpha} 
(h|b|(w\wedge n))\|_{L_{p'}(w')}.
$$

Here $R_{\alpha} (h|b|(w\wedge n)\in L_{p'}(w')$ for each $n$, because this functions
is bounded in each ball and outside the ball $B_{R}$
containing the support of $b$ we have
$R_{\alpha} (h|b|(w\wedge n)(x)\leq N|x|^{\alpha-d}$, where $N$ is independent of $x$
with
$$
\int_{|x|\geq 1}|x|^{(\alpha-d)p'}
w^{1/(p-1)}(x)\,dx=N\int_{1}^{\infty}r^{(\alpha-d)p'-1}\int_{B_{r}\setminus B_{1}}
w^{1/(p-1)}(x)\,dxdr
$$
$$
\leq N\int_{1}^{\infty}r^{(\alpha-d)p'-1+d}\dashint_{B_{r}}
w^{1/(p-1)}(x)\,dxdr
=NJ<\infty.
$$

Hence,
by Theorem \ref{theorem 1.19.2} and
Lemma \ref{lemma 9.5.1}
$$
\|R_{\alpha} (h|b|(w\wedge n))\|_{L_{p'}(w')}
\leq N\|R^{\sharp}_{\alpha} (h|b|(w\wedge n))\|_{L_{p'}(w')}
$$
$$
\leq N\|\bM_{\alpha} (h|b|(w\wedge n))
\|_{L_{p'}(w')}\leq N\|\bM_{\alpha} (h|b|w)\|_{L_{p'}(w')}.
$$
We see that in \eqref{5.30,3}
$$
\|R_{\alpha} (h|b|w)\|_{L_{p'}(w')}
\leq N\|\bM_{\alpha} (h|b|w)\|_{L_{p'}(w')}.
$$
  After that
we observe that by H\"older's inequality
$$
\bM_{\alpha} (h|b|w)
\leq \|b\|_{\dot E_{q,\alpha}}
[\bM\big((hw)^{q'}\big)]^{1/q'},
$$
where  $q'=q/(q-1)$ and,  since 
the condition $w^{s}\in A_{r}$ means that
$w'\in A_{p'/q'}$,  
by Theorem \ref{theorem 11.21.3}
$$
\|[\bM\big((hw)^{q'}\big)]^{1/q'}\|
_{L_{p'}(w')}\leq
N\|hw\|_{L_{p'}(w')}=N\|h\|_{L_{r'}(w)}.
$$
This, obviously, proves the theorem. \qed

Similarly to Corollary \ref{corollary 10.7.1}
we have the following.

\begin{corollary}
                       \label{corollary 5.28,1}
If $u\in C_{0}^{\infty}(\bR^{d})$ and $\alpha=1$
(and $d\geq2$), then
 owing to \eqref{9.5.5}  
$$
\int_{\bR^{d}}b^{p}|u|^{p}w\,dx\le  NA^{p}
 \int_{\bR^{d}}  |Du|^{p}w\,dx,
$$
$$ 
\int_{\bR^{d}}b^{p}|Du|^{p}w\,dx\le  NA^{p}
 \int_{\bR^{d}}  |D^{2}u|^{p}w\,dx,
$$
and if $\alpha=2$ (and $d\geq3$), then
$$
\int_{\bR^{d}}b^{p}|u|^{p}w\,dx\le  NA^{p}
 \int_{\bR^{d}} \big|\Delta u\big|^{p}w\,dx.
$$
Generally,
$$
\int_{\bR^{d}}b^{p}|u|^{p}w\,dx\le  NA^{p}
 \int_{\bR^{d}} \big|(-\Delta)^{\alpha/2}u\big|^{p}w\,dx.
$$
 \end{corollary}

\begin{remark}
                     \label{remark 5.30,1}
In the model case when $\alpha=1$, 
$p< q<  d$, $|b|=1/|x|$ and
$w=1/|x|^{\beta}$ with $\beta\in(-(p-1)d, (q-p)d/q )$,
so that $w\in A_{p}$   and   $w^{s}\in A_{r}$,  it turns out that
to satisfy \eqref{5.30,1} it suffices to have
$p+\beta<d$.  Since $d-p>(q-p)d/q$, for $\beta\in
(-(p-1)d,(q-p)d/q)$  and $u\in C^{\infty}_{0}$
$$
\int_{\bR^{d}}\frac{|u(x)|^{p}}{|x|^{p}}
\frac{1}{|x|^{\beta}}\,dx
\leq N\int_{\bR^{d}}|Du(x)|^{p}
\frac{1}{|x|^{\beta}}\,dx,
$$
where $N$ is independent of $u$.
\end{remark}

Next, take an $\bR^{d}$-valued $b$ and introduce
$$
\cL u :=\Delta u+b^{i}D_{i}u.
$$
By using that
$$
\|\Delta u-\lambda u\|_{L_{p}(w)}
\leq \|\cL u-\lambda u\|_{L_{p}(w)}
+\|b^{i}D_{i}u\|_{L_{p}(w)},
$$
and using \eqref{12.22.10} and \eqref{1.22.2} we come
to \eqref{5.28,2} below.

\begin{theorem} 
               \label{theorem 5.28,2}
Let $1<p<q <\infty$, $w\in A_{p}$. 
Also assume \eqref{5.30,1}  
with $\alpha=1$. Then
there  are constants  $N=N(d,p,q , [w]_{A_{p} } )$  
  and $\lambda_{0}=\lambda_{0} 
(d,p,q,[w]_{A_{p}} )>0$ such that, if  
$$
N\|b\|_{\dot E_{q,1}(\bR^{d })}
\leq 1,
$$
 then
for any   $\lambda \geq \lambda_{0}$  and $f\in L_{p}(w) $
there exists a unique $u\in W^{ 2}_{p}(w) $
satisfying 
\begin{equation}
                           \label{5.28,1}
\cL u-\lambda u:=\Delta u+b^{i}D_{i}u-\lambda u=f.
\end{equation}

Furthermore,   for any $u\in W^{ 2}_{p}(w) $  and $\lambda 
 \geq \lambda_{0} $ 
\begin{equation}
                           \label{5.28,2}
\|\partial_{t}u,D^{2}u,\sqrt{\lambda}Du,
\lambda u\|_{L_{p}(w) }\leq N(d,p,q , [w]_{A_{p}})
\|\cL u-\lambda u\|_{L_{p}(w) }.
\end{equation}  
\end{theorem}

The existence part of this theorem,
as usual, is obtained by the method of continuity.

We can also reprove Lemma \ref{lemma 2.22.3}.
\begin{lemma}
                 \label{lemma 1.20.2}
For any $\beta>0,p>1$, $u\in \dot E^{2}_{p,\beta}$ we have
$$
\|D^{2}u\|_{\dot E^{2}_{p,\beta}}
\leq N(d,p,\beta)\|\Delta u\|_{\dot E _{p,\beta}}.
$$ 
\end{lemma}

Indeed, for $1<r<p$, $(D^{2}u)^{\sharp}\leq N\big(\bM(|\Delta u|^{r})\big)^{1/r}$ by Theorem \ref{theorem 1.19.1} and it only remains to apply Theorem \ref{theorem 1.20.2} and Remark \ref{remark 7,7.1}.

\mysection
{Parabolic equations}

In this section
 we deal with  $\bR^{d+1}=\{(t,x):
t\in \bR,x\in\bR^{d}\}$
and the real analytic structure 
defined by Lebesgue measure
and  $k_{0}=2,k_{1}=...=k_{d}=1$,
so that we are dealing with the ``parabolic boxes''
$$
\cB_{l,(t,x)}=(t,x)+[0,l^{2})
\times [0,l)^{d}.
$$
This affects such object
\index{$B$@Sets!$\cB_{l,(t,x)}$}% 
 as $\bM_{\alpha} u$, $u^{\sharp}$, 
and
Muckenhoupt weights
among others.

Introduce
$$
C_{T,R}=[0,T)\times B_{R}, 
\quad C_{R}=
C_{R^{2},R},
$$
$$
C_{T,R}(t,x)=(t,x)+C_{T,R},
 \quad C_{R}(t,x)=C_{R^{2},R}(t,x),
$$
\index{$B$@Sets!$C_{T,R}$}%
\index{$B$@Sets!$C_{R}$}%
\index{$B$@Sets!$C_{T,R}(t,x)$}%
\index{$B$@Sets!$C_{R}(t,x)$}%
\index{$B$@Sets!$\bC_{R}$}%
\index{$B$@Sets!$\bC$}%
and let $\bC_{R}$ be the collection of $C_{R}(t,x)$ and $\bC=\cup_{R}\bC_{R}$.
Also 
\index{$C$@Operators!$\partial_{t}$}%
define
$$
\quad \partial_{t}=\frac{\partial}{\partial t}.$$

\subsection{Application of Theorem
\protect\ref{theorem 11,23,1} to the heat equation}

                 \label{section 1.26.2} 
For $p\in[1,\infty)$ and domain $\cO\subset
\bR^{d+1}$ define
$W^{1,2}_{p}(\cO) $ 
\index{$A$@Sets of functions!$W^{1,2}_{p}(\cO)$}%
\index{$A$@Sets of functions!$W^{1,2}_{p}$}%
as the set of functions $u$ on $\cO$ such that
$u$ and its Sobolev derivatives $Du,D^{2}u,\partial_{t}u$ are in 
$L_{p}(\cO) $. Define
$$
\|u\|_{W^{1,2}_{ p}(\cO)}:=\|u,Du,D^{2}u,\partial_{t}u\|_{_{L_{p}(\cO)}} 
$$
\index{$N$@Norms!$"|"|u"|"| _{W^{1,2}_{p}(\cO)}$}%
and drop $\cO$ in the above notation if
$\cO=\bR^{d+1}$.

Let $u\in C^{\infty}_{0}$. 
Set $-f=\partial_{t}u+\Delta u$. Then
by recalling \eqref{9.4.1} and observing that
$$
\int_{\bR^{d+1}}\partial_{t}u\Delta u
\,dxdt=-\int_{\bR^{d+1}}\partial_{t}D_{i}uD_{i} u
\,dxdt
=-(1/2)\int_{\bR^{d+1}}\partial_{t}|Du|^{2}
\,dxdt=0,
$$
we obtain
\begin{equation}
                          \label{1.25.1}
\int_{\bR^{d+1}}\big(|\partial_{t}u|^{2}+|D^{2}u|^{2}\big)\,dxdt=\int_{\bR^{d+1}}f^{2}\,dxdt.
\end{equation}
 \begin{remark}
                         \label{remark 4.4,3}

Equality  \eqref{1.25.1} is also true for
infinitely differentiable functions $u$ on
$\bR^{d+1}$ such that, for some $\alpha>d/2$ 
\index{$S$@Miscelenea!$\rho(t,x)$}%
and $\rho(t,x)=\sqrt{|t|}+|x|,\xi 
=\rho^{ \alpha-1},\eta=\rho^{\alpha}$ the function 
$$
\xi|u|+\eta|Du| 
$$
 is bounded.

Indeed, take a $\zeta\in C^{\infty}_{0}$ such that $\zeta=1$ in $(-1,1)\times B_{1}$ and $\zeta=0$ outside $(-4,4)\times B_{2}$,
$0\leq \zeta\leq 1$, and set $\zeta_{n}(t,x)
=\zeta(t/n^{2},x/n)$. Then for each $n$
\begin{equation}
                          \label{2.24.1}
 \|\partial_{t}(\zeta_{n}u),D^{2}(\zeta_{n}u)\|_{L_{2}(\bR^{d+1})}
 =\|\partial_{t}(\zeta_{n}u)+\Delta(\zeta_{n}u)\|_{L_{2}(\bR^{d+1})}.
\end{equation}
Here
$$
\big|\|\partial_{t}(\zeta_{n}u),D^{2}(\zeta_{n}u)\|_{L_{2}(\bR^{d+1})}
-\|\zeta_{n}\partial_{t}u,\zeta_{n}D^{2}u\|_{L_{2}(\bR^{d+1})}\big|\leq I_{n}+J_{n},
$$
where 
$$
I_{n}=\big|\|\partial_{t}(\zeta_{n}u) \|_{L_{2}(\bR^{d+1})}
-\|\zeta_{n}\partial_{t}u \|_{L_{2}(\bR^{d+1})}\big|  
$$
$$
\leq Nn^{-2}\|(\partial_{t}\zeta )(\cdot/n^{2},\cdot/n)\xi^{-1} \|_{L_{2}(\bR^{d+1})}
=Nn^{d/2-\alpha}\|(\partial_{t}\zeta )  \xi^{-1} \|_{L_{2}(\bR^{d+1})}, 
$$
$$
J_{n}
=\big|\| D^{2}(\zeta_{n}u)\|_{L_{2}(\bR^{d+1})}
-\| \zeta_{n}D^{2}u\|_{L_{2}(\bR^{d+1})}\big| 
$$
$$
\leq Nn^{-1}\|(D\zeta)(\cdot/n^{2},\cdot/n)
\eta^{-1}\|_{L_{2}(\bR^{d+1})}+Nn^{-2}
\|(D^{2}\zeta)(\cdot/n^{2},\cdot/n)
\xi^{-1}\|_{L_{2}(\bR^{d+1})} 
$$
$$
\leq Nn^{d/2-\alpha}\big(\|(D\zeta)\eta^{-1}\|_{L_{2}(\bR^{d+1})}+\|(D^{2}\zeta)
\xi^{-1}\|_{L_{2}(\bR^{d+1})}\big). 
$$
Since $n^{d/2-\alpha}\to0$ as $n\to\infty$
and, owing to the fact that $\zeta=1$ near the origin 
$$
\|(\partial_{t}\zeta )  \xi^{-1} \|_{L_{2}(\bR^{d+1})}+\|(D\zeta)\eta^{-1}\|_{L_{2}(\bR^{d+1})}+\|(D^{2}\zeta)
\xi^{-1}\|_{L_{2}(\bR^{d+1})}<\infty, 
$$
we easily get \eqref{1.25.1} from \eqref{2.24.1}.
\end{remark}

Next,
set $c_{d}=(4\pi)^{-d/2}$ and define
$$
\sfp(t,x)=c_{d}t^{-d/2}e^{-|x|^{2}/(4t)}
I_{t>0} 
$$
which is the fundamental
\index{$S$@Miscelenea!$\sfp(t,x)$}%
\index{$C$@Operators!$R_{\lambda}f(t,x)$}%
 solution of  the heat equation.
It is a classical fact that for  
$u\in C^{\infty}_{0}$, $\lambda\geq0$, and 
$$
-f:=\partial_{t}u+\Delta u-\lambda u,
$$
 we
\index{$C$@Operators!$R_{\lambda}$}%
 have
$$
u(t,x)=R_{\lambda}f(t,x):= \int_{0}^{\infty}\int_{\bR^{d}}e^{-\lambda s}
\sfp(s,y)f(t+s,x+y)\,dsdy
$$
\begin{equation}
                            \label{4.1,3}
=\int_{t}^{\infty}\int_{\bR^{d}}e^{-\lambda (s-t)}
\sfp(s-t,x-y)f( s, y)\,dsdy.
\end{equation}

On the other hand, it is well known that, if
$g\in C^{\infty}_{0}$, then $G:=R_{0}g$
is infinitely differentiable and satisfies  the heat equation
\begin{equation}
                      \label{1.24.3}
-g =\partial_{t}G+\Delta G.
\end{equation}

Next, note that, for any integer  $n \geq0$
we have 
$$|D^{n}e^{-|x|^{2}}|\leq N(n)e^{-|x|^{2}/2},
$$
 which along with $\partial_{t} \sfp=\Delta \sfp$ implies that for any integer
$m\geq0$ and $t>0$ 
\begin{equation}
                          \label{1.25.2} 
| \partial_{t} ^{m}D^{n}\sfp(t,x)|\leq 
N(m,n)\frac{1}{t^{m+(d+n)/2} }e^{-|x|^{2}/(8t)}.
\end{equation}

Then we make a slight digression with
the goal of not repeating the same argument
twice in the future. 
For $k,s,r>0,\alpha\in \bR $, and appropriate $f(t,x)$'s
on $\bR^{d+1}$
\index{$S$@Miscelenea!$p_{\alpha,k}(s,r)$}%
\index{$C$@Operators!$P_{\alpha,k}f(t,x)$}%
 define
$$
p_{\alpha,k}(s,r)=\frac{1}{s^{(d+2-\alpha)/2}}e^{-r^{2}/(ks)}I_{s>0}, 
$$
$$
P_{\alpha,k}f(t,x)=\int_{\bR^{d+1} }p_{\alpha,k}(s,|y|)f(t+s,x+y)\,dyds.
$$
$$
=\int_{t}^{\infty}\int_{\bR^{d} }p_{\alpha,k}(s-t,|y-x|)
f(s,y)\,dsdy.
$$  

\begin{remark}
                  \label{remark 5.31,1}
For $\alpha\leq d+2$ there is a constant $N=N(d,\alpha)$
such that
\begin{equation}
                        \label{5.31,1}
p_{\alpha,k}(s,|x|)\leq N\frac{1}{\rho
^{d+2-\alpha}(s,x)}.
\end{equation}
Indeed, if $|x|\leq \sqrt{|s|}$, then
$$
p_{\alpha,k}(s,|x|)\leq 
\frac{1}{s
^{(d+2-\alpha)/2}}=\frac{1}{\rho
^{d+2-\alpha}(s,x)}\frac{\rho
^{d+2-\alpha}(s,x)}{s
^{(d+2-\alpha)/2}},
$$
where the last fraction is majorated
by $2^{d+2-\alpha}$. In case
$|x|\geq \sqrt{|s|}$,
$$
p_{\alpha,k}(s,|x|)\leq 
\frac{1}{|x|
^{d+2-\alpha}}\phi(|x|/\sqrt{|s|}),
$$
where $\phi(t)= t^{ d+2-\alpha}e^{-t^{2}/k} $ 
 is a bounded function on
$(0,\infty)$. In that case we get 
\eqref{5.31,1} again.

As a consequence of this estimate we
obtain that, if $f$ is a bounded function
with compact support, then 
\begin{equation}
                        \label{5.31,2}
|P_{\alpha,k}f(t,x)|
\leq N(1+\rho(t,x))^{-(d+2-\alpha)},
\end{equation}
where $N$ is independent of $(t,x)$.
In particular, if
$g\in C^{\infty}_{0}$, then $G=R_{0}g$,
which 
is infinitely differentiable and satisfies  the heat equation \eqref{1.24.3} also admits the estimates
$$
|G|\leq N(1+\rho )^{-d},\quad
|DG|\leq N(1+\rho )^{-(d+1)}.
$$
It follows that the conditions of
Remark \ref{remark 4.4,3} are satisfied,
for instance, with $\alpha=d$ and 
\eqref{1.25.1} holds for $u=G$.
\end{remark}

\begin{lemma}
                   \label{lemma 19.30.1}
For any $ \alpha<\beta $, $\beta\geq0$, and $k>0$, there exists a   constants $N$
  such that for any $f\geq0$, $\rho>0$,
 we have   on $C_{\rho}$ that
\begin{equation}
                     \label{1.17.20}
P_{\alpha,k}(I_{C^{c}_{2\rho}}f)  
\leq N \rho^{\alpha-\beta}\bM_{\beta}
f(0)  .
\end{equation}
\end{lemma}

Proof. First observe that if $\alpha\geq d+2$,
then $\beta>d+2$ and the right-hand side
of \eqref{1.17.20} is infinite unless $f=0$. Hence, we suppose that $\alpha<d+2$
Self-similar transformations 
allow us to concentrate on $\rho=1$.
In that case, clearly, we may assume that $f$ is bounded
with compact support.  Set 
$$
Q_{1}=\{(s,y):|y|\geq\sqrt s,s>0\},\quad
 Q_{2}=\{(s,y):|y|\leq \sqrt s,s>0\}.
 $$
Dealing with $P_{\alpha,k}(fI_{Q_{1}\cap C^{c}_{2}})$
we observe that $p_{\alpha,k}(s,r)\leq Nr^{-(d+2-\alpha)}$ since $\alpha\leq d+2$. Therefore,
$$
P_{\alpha,k}(fI_{Q_{1}\cap C^{c}_{2}})(0)
\leq N\int_{2}^{\infty}\frac{1}{r^{d+2-\alpha}}
\int_{0}^{r^{2}}\Big(\int_{|y|=r}f(s,y)\,d\sigma_{r}\Big)\,dsdr,
$$
where $d\sigma_{r}$ is the element of the surface area on $|y|=r$. By using that $\alpha<d+2$ and integrating by parts we get
\begin{align*}
P_{\alpha,k}(fI_{Q_{1}\cap C^{c}_{2}})&(0)
\\
&\leq N
\int_{2}^{\infty}\frac{1}{r^{d+3-\alpha}}
\int_{2}^{r}\Big(\int_{0}^{\rho^{2}}
\Big(\int_{|y|=\rho}f(s,y)\,d\sigma_{\rho}\Big)\,ds\Big)\,d\rho dr \\
&
\leq N
\int_{2}^{\infty}\frac{1}{r^{d+3-\alpha}}
\int_{0}^{r}\Big(\int_{0}^{r^{2}}
\Big(\int_{|y|=\rho}f(s,y)\,d\sigma_{\rho}\Big)\,ds\Big)\,d\rho dr \\
&
=N\int_{2}^{\infty}\frac{1}{r^{d+3-\alpha}}I(r)\,dr ,
\end{align*}
where  
$$
I(r)=\int_{C_{r}}f(s,y)\,dyds.
$$
Observe that $r\geq 2$ and for any $(t,x)\in C_{1}(-1,0)$
$$
I(r)\leq 
\int_{C_{2r}(t,x)}f(s,y)\,dyds
$$
$$
\leq 
Nr^{ d+2-\beta }\Big((2r)^{\beta}\dashint_{C_{2r}(t,x)}f(s,y)\,dyds\Big)\leq Nr^{ d+2-\beta }\bM_{\beta}f(t,x). 
$$
We also use that    $\alpha<\beta$. Then we see that  
\begin{equation}
                          \label{1.17.030}
P_{\alpha,k}(fI_{Q_{1}\cap C^{c}_{2}})(0)\leq N \bM_{\beta}f(t,x). 
\end{equation}

Next, by integrating by parts
we  obtain  that
$$
P_{\alpha,k}(fI_{Q_{2}\cap C^{c}_{2}})(0)
\leq
\int_{4}^{\infty}\frac{1}{s^{(d+2-\alpha)/2}}
\int_{|y|\leq\sqrt s}f(s,y)\,dyds 
$$
$$
\leq N\int_{4}^{\infty}\frac{1}{s^{(d+4-\alpha)/2}}I(\sqrt s)\,ds=
N\int_{2}^{\infty}\frac{1}{r^{d+3-\alpha}}I(r)\,dr. 
$$
This along with \eqref{1.17.030} prove  that
$P_{\alpha,k}(I_{C^{c}_{2}}f)(0) 
\leq N \bM_{\beta}
f (t,x) $, which is equivalent to our statement.  \qed
 
Then coming back to the setting of this
section when $u\in C^{\infty}_{0}$,
$\lambda=0$,
and $-f=\partial_{t}u+\Delta u$, take $\zeta \in C^{\infty}_{0}$
such that $\zeta=1$ in $(-4,4)\times B_{2}$, $\zeta=0$
outside $(-9,9)\times B_{3}$, $0\leq\zeta\leq1$ and  
note that for   $g=f\zeta,h=f(1-\zeta)$
 and $(G,H)=R_{0}(g,h)$ in $C_{1}$ we have
$u=G+H$ and
 $$
 \partial_{t} ^{m}D^{n}
H(t,x)= \int_{C_{2}^{c}}
 \partial_{t} ^{m}D^{n}
\sfp(s-t,x-y)h( s, y)\,dsdy.
$$

By using \eqref{1.25.2} and taking  $\beta=0$ and $\alpha=-1,-2$
in Lemma \ref{lemma 19.30.1}
we get that  in $C_{1}$
\begin{equation}
                         \label{1.26.1}
|\partial_{t}D^{2}H(t,x)|+|D^{3}H(t,x)|
\leq N\bM h(0),
\end{equation}
implying that
$$
\int_{C_{1}}\int_{C_{1}}
|D^{2}H(t,x)-D^{2}H(s,y)|\,dxdydtds\leq N\bM f(0).
$$

Regarding $G$ we have by \eqref{1.25.1}
that
$$
\int_{C_{1}}\int_{C_{1}}
|D^{2}G(t,x)-D^{2}G(s,y)|\,dxdydtds\leq N\int_{C_{1}} 
|D^{2}G(t,x) |\,dxdt
$$
$$
\leq N\|g\|_{L_{2}(\bR^{d+1})}\leq N(\bM (f^{2}))^{1/2}(0). 
$$

It follows that for $\rho=1$
$$
\dashint_{C_{\rho}}\dashint_{C_{\rho}}
|D^{2}u(t,x)-D^{2}u(s,y)|\,dxdydtds
\leq N(\bM (f^{2}))^{1/2}(0).
$$
For any other $\rho>0$ this is derived
by using scaling and yields 
\begin{equation}
                            \label{1.25.4}
(D^{2}u)^{\sharp}\leq N(\bM (f^{2}))^{1/2}
\end{equation}
at the origin. Of course, this estimate holds at any other point and by Theorems
\ref{theorem 06.6.9.1} and \ref{theorem 7.7.1} and using that $\partial_{t}u=-f-\Delta u$ we conclude that for any $p>2$ and $u\in C^{\infty}_{0}$ 
  \begin{equation}
                                          \label{1.25.5}
  \|\partial_{t}u\|_{L_{p}(\bR^{d+1})}+\|D^{2}u\|_{L_{p}(\bR^{d+1})}\leq N(p,d)\|\partial_{t}u+\Delta u\|_{L_{p}(\bR^{d+1})}.
 \end{equation}

 \begin{theorem}
                  \label{theorem 1.25.3}
 For $d\geq1$, any $p\in(1,\infty)$, and $u\in W^{1,2}_{p}$ we have
 \eqref{1.25.5}.
 \end{theorem}
 
 Proof.   By using the fact
 that $C_{0}^{\infty}(\bR^{d+1})$ is dense in $W^{1,2}_{p}$
 we may assume that $u\in C_{0}^{\infty}(\bR^{d+1})$.
 Then $f:=-\partial u-\Delta u\in C_{0}^{\infty}(\bR^{d+1})$ and $u=R_{0}f$.
In that case
 \eqref{1.25.5} holds for $p>2$ by the above.
 
 This estimate also holds for $p=2$ in light 
 of \eqref{1.25.1}. To extend
 it to $p\in(1,2)$ observe that estimate \eqref{1.25.5} implies that
 for any smooth $f,g$ with compact support and any $i,j$
\begin{equation}
                                          \label{9.4.40}
 \int_{\bR^{d+1}}g D_{ij}R_{0}f\,dx\leq N\|g\|_{L_{q}(\bR^{d+1})}\|f\|_{L_{p}(\bR^{d+1})},
\end{equation}
 where $q=p/(p-1)$. Here on the left we can move first $D_{ij}$
 to $g$ and then $R_{0}$ to $D_{ij}g$
which will hit it as $R^{*}_{0}$, the conjugate
to $R_{0}$. The operator $R^{*}_{0}$ is related
to $-\partial_{t} +\Delta  $ in the same way as $R_{0}$ is related to $ \partial_{t} +\Delta  $ and the one is obtained from the other just by reflecting the $t$-axis.
In short 
$$
R^{*}_{0}f(t,x)=R_{0}(f(-\cdot,\cdot))(-t,x).
$$
By observing that
 $R^{*}_{0}D_{ij}g=D_{ij}R^{*}_{0}g$ we come to the version of
 \eqref{9.4.40} with $f,g$ interchanged and $R_{0}$   replaced with $R^{*}_{0}$. Then the arbitrariness
 of $f$ and similar argument regarding
$\partial_{t}$ leads to \eqref{1.25.5} with $q$ in place of $p$ and $-\partial_{t}$
in place of $\partial_{t}$, with $(-)$ which is easily eliminated.
\qed   

\begin{remark}
               \label{remark 1.26.1}
As in Remark \ref{remark 1.25.7} we point   
out that 
estimate \eqref{1.25.5} also holds for
  functions not necessarily in $W^{1,2}_{p}$.
For instance, let $u$ be infinitely differentiable
and such that $\xi|u|+\eta|Du|$ is bounded,
where $\xi,\eta$ is introduced after \eqref{1.25.1} with $\alpha>(d+2)/p-1$. Then, by repeating the computations after \eqref{1.25.1}
we obtain \eqref{1.25.5} for such an $u$.

\end{remark}

Now we prove the interpolation inequality
similar to its elliptic counterpart
from Lemma \ref{lemma 2.1.1}.

\begin{lemma}
                          \label{lemma 2.1.3}
For any $\varepsilon>0$ and $u\in C^{\infty}_{0}$
\begin{equation}
                  \label{2.1.4}
(Du)^{\sharp}\leq N\varepsilon(\bM(\partial_{t}u)+\bM(D^{2}u))+N\varepsilon^{-1}\bM u,
\end{equation}
in particular, for any $p>1$,
\begin{equation}
                  \label{2.2.8}
\|Du\|_{L_{p}(\bR^{d+1})}\leq
N(\varepsilon \|\partial_{t}u,D^{2}u\|_{L_{p}(\bR^{d+1})}+\varepsilon^{-1}
\|u\|_{L_{p}(\bR^{d+1})}),
\end{equation}
where the constants $N$ are independent of
$\varepsilon,u$.
\end{lemma}

Proof. By Poincar\'e's inequality
(see, for instance, Lemma \ref{lemma 10.10.1}) for $r\leq 1$
$$
\osc_{C_{r}}Du\leq Nr (\bM(\partial_{t}u)+\bM(D^{2}u))\leq N(\bM(\partial_{t}u)+\bM(D^{2}u)),
$$
whereas for $r>1$ by Lemma \ref{lemma 2.1.1}
$$
\osc_{C_{r}}Du\leq 2\dashint_{C_{r}}|Du|\,dxdt
\leq N( \bM(D^{2}u)+\bM u).
$$
This implies \eqref{2.1.4} for $\varepsilon=1$.
For other $\varepsilon>0$ one obtains it
by using self-similarity. 
\qed

After that almost literally repeating the proof
of Theorem \ref{theorem 12.22.1} 
and using \eqref{2.2.8} we come to the following.

\begin{theorem}
           \label{theorem 1.25.4}
For any $p>1$, $\lambda>0$, and $f\in L_{p}(\bR^{d+1})$
there exists a unique $u\in W^{1,2}_{p}$
satisfying 
\begin{equation}
                           \label{4.4,6}
\partial_{t}u+\Delta u-\lambda u+f=0.
\end{equation}
Furthermore, there exists a constant
$N=N(d,p)$ such that for any $u\in W^{1,2}_{p}$
and $\lambda\geq0$
\begin{equation}
                          \label{1.25.60}
\|\partial_{t}u,D^{2}u,\sqrt\lambda Du,\lambda u\|_{L_{p}(\bR^{d+1})}\leq
N\|\partial_{t}u+\Delta u-\lambda u\|_{L_{p}(\bR^{d+1})}.     
\end{equation}

\end{theorem}

\begin{remark}
                       \label{remark 4.4.6}
Solutions of the parabolic equation
\eqref{4.4,6} for $\lambda>0$ possess a feature absent
for solutions of elliptic equations.
What we mean is that if, for a fixed $T\in\bR$, $f(t,x)=0$ for $t\geq T$ and all $x$, then for the solution $u$ the same is true: $u(t,x)=0$ for $t\geq T$ and all $x$. 
We call this {\em causality\/}.
\index{$S$@Miscelenea!causality}% 
This is a direct consequence of formula \eqref{4.1,3},
which is valid, as we know, for $u\in C^{\infty}_{0}$ and then by the continuity extends to $u\in W^{1,2}_{p}$.
Note that for $\lambda>0$ the $L_{p}(\bR^{d+1})$-norm of the operator $R_{\lambda}$ is dominated by $\lambda^{-1}$ thanks to Minkowski's inequality.

\end{remark}

This remark implies the following in
which 
\index{$B$@Sets!$\bR^{d+1}_{t}$}% 
$$
\bR^{d+1}_{t}=(t,\infty)\times \bR^{d}. 
$$

\begin{theorem}
           \label{theorem 4.26,1}
For any $p>1$, $\lambda>0$, $s\in[-\infty,\infty)$, and $f\in L_{p}(\bR^{d+1}_{s})$
there exists a unique $u\in W^{1,2}_{p}(\bR^{d+1}_{s})$
satisfying in $\bR^{d+1}_{s}$
\begin{equation}
                           \label{4.26,1}
\partial_{t}u+\Delta u-\lambda u+f=0.
\end{equation}
Furthermore, there exists a constant
$N=N(d,p)$ such that for any $u\in W^{1,2}_{p}(\bR^{d+1}_{s})$
and $\lambda\geq0$
\begin{equation}
                          \label{1.25.6}
\|\partial_{t}u,D^{2}u,\sqrt\lambda Du,\lambda u\|_{L_{p}(\bR^{d+1}_{s})}\leq
N\|\partial_{t}u+\Delta u-\lambda u\|_{L_{p}(\bR^{d+1}_{s})}.     
\end{equation}

\end{theorem}

Proof. If $s=-\infty$, the assertions
are proved in Theorem \ref{theorem 1.25.4}. General $s\in\bR$ is as good as $s=0$. The
solvability of \eqref{4.26,1} with
$fI_{(0,\infty)}$ in place of $f$
even in $W^{1,2}_{p}(\bR^{d+1} )$
follows from Theorem \ref{theorem 1.25.4}.

To prove \eqref{1.25.6} for $\lambda>0$,
take $u\in W^{1,2}_{p}(\bR^{d+1}_{0})$,
define $f$ by \eqref{4.26,1}, find
a $w\in W^{1,2}_{p}(\bR^{d+1} )$
satisfying \eqref{4.26,1} with $fI_{(0,\infty)}$ in place of $f$ and  
set $v(t,x)=u(|t|,x)$. Then $w-v$
satisfies \eqref{4.26,1} with the free term
which is zero for $t>0$. Hence, by  causality
$w=u$ for $t>0$ and \eqref{1.25.6}
follows from Theorem \ref{theorem 1.25.4}.
This implies uniqueness for \eqref{4.26,1}
and by sending $\lambda\downarrow0$
also shows that \eqref{1.25.6}
is true for $\lambda=0$ as well. \qed

Here is the parabolic analog
of Theorem \ref{theorem 1.19.1}
\begin{theorem}
                      \label{theorem 1.26.1}
For any $p>1$ there exists a constant $N=N(d,p)$
such that for any $u\in W^{1,2}_{p}$
\begin{equation}
                         \label{1.19.10}
(D^{2}u)^{\sharp}\leq N\big(\bM(|f|^{p})\big)^{1/p},
\end{equation} 
where $f=\partial_{t}u+\Delta u $.
\end{theorem}

The proof of this theorem is obtained
from the proof of Theorem \ref{theorem 1.19.1} by using Remark \ref{remark 1.26.1}
in place of Remark \ref{remark 1.25.7}
and using estimate \eqref{1.26.1}
in place of \eqref{9.4.2}.

\subsection{Parabolic Adams theorem
and the heat equation}

                 \label{section 1.26.1}

Observe that, if $f$ is independent of $t$ 
and $0<\alpha<d$,
then
$$
P_{\alpha,k}f(t,x)=
P_{\alpha,k}f(x)=N(\alpha,k)\int_{\bR^{d}}
\frac{1}{|y|^{d-\alpha}}f(x+y)\,dy= NR_{\alpha}f(x),
$$
where $R_{\alpha}$ is the Riesz potential.
In our presentation the most important values   
of $\alpha$ are 1 and 2.

Here is an analog of the Adams Theorem \ref{theorem 10.7.1}, which is extended to the mixed-norm case
in Theorem \ref{theorem 5.25,1}.

\begin{theorem}
                    \label{theorem 1.18.1}
Let $  q> p>1$, $\alpha>0$, let $b$ be defined on $\bR^{d+1}$, and
\index{$N$@Norms!$"|"|b"|"|_{\dot E_{q,\alpha}(\bR^{d+1})}$}%
set
\begin{equation}
                       \label{1.28.3}
\|b\|_{\dot E_{q,\alpha}(\bR^{d+1})}:=\sup_{\rho>0}\rho ^{\alpha}\sup_{C\in\bC_{\rho}}  
\Big(\dashint_{C}|b|^{q}\,dxdt\Big)^{1/q} . 
\end{equation}
    Then for any $s\in [-\infty,\infty)$ and $f\geq0$ on $\bR^{d+1}_{s}$
we have
\begin{equation}
                         \label{5.30.1}
I:= \int_{\bR^{d+1}_{s}}|b|^{p}(P_{\alpha,k}f)^{p}\,dz\leq
N\|b\|_{\dot E_{q,\alpha}(\bR^{d+1})}^{p}\|f\|_{L_{p}(\bR^{d+1}_{s})},
\end{equation}
where $N$ depends only on $d,p,q,\alpha,k$.
In particular, 
 for any $u\in C^{\infty}_{0}$
\begin{equation}
                       \label{1.18.3}
 \int_{\bR^{d+1}_{s} } |b|  ^{p}|Du |^{p}\,dz
\leq N\|b\|_{\dot E_{q,1}}^{p}K,
\quad
\int_{\bR^{d+1}_{s} }| b|  ^{p}  |u |^{p}\,dz
\leq N\|b\|_{\dot E_{q,2}}^{p}K,
\end{equation}
where
$
K=\|D^{2}u,\partial_{t}u\|^{p}_{L_{p}(\bR^{d+1}_{s})}
$
and $N$ depends only on $d,p,q$.

\end{theorem}
\begin{remark}
                     \label{remark 11.2.1}
The first estimate in \eqref{1.18.3} follows from \eqref{5.30.1} with $\alpha=1$   and the fact that
for $f=\partial_{t}u+\Delta u$ we have
$$
Du(t,x)=c\int_{\bR^{d +1 }_{0}}
\frac{y}{s^{(d+2)/2}}e^{-|y|^{2}/(4s)}
f(t+s,x+y)\,dyds,
$$
where $c$ is a constant  and
$(|y|/s^{1/2})e^{-|y|^{2}/(4s)}\leq
Ne^{-|y|^{2}/(8s)}$. The second estimate
follows when $\alpha=2$ since
$$
u(t,x)=c\int_{\bR^{d +1 }_{0}}
\frac{1}{s^{d/2}}e^{-|y|^{2}/(4s)}
f(t+s,x+y)\,dyds.
$$

\end{remark}

\begin{remark}
                      \label{remark 4.27,1}
It suffices to prove Theorem \ref{theorem 1.18.1} for $s=-\infty$. Indeed, if
it is true for such $s$, then plugging in
$f(t,x)I_{(s,\infty)}(t)$ in place of $f$,
we obtain the result as stated.
Below we argue in case $s=-\infty$ when
$\bR^{d+1}_{s}=\bR^{d+1}$.
\end{remark}

 In Theorem \ref{theorem 1.18.1}
we did not restrict $\alpha$ from above.
However, observe that if $\alpha>(d+2)/q$, then
automatically $b=0$ if $\|b\|_{\dot E_{q,\alpha}(\bR^{d+1})}<\infty$.

Before proving  Theorem \ref{theorem 1.18.1}
we need some preparations.
 
\begin{corollary}[of Lemma \ref{lemma 19.30.1}]
                        \label{corollary 10.30.1}
There is a constant $N$ such that for any $f$,
for which $P_{\alpha,k}|f|$ is finite on $C_{1}$,
on $C_{1}$ we have
\begin{equation}
                             \label{10.30.1}
|DP_{\alpha,k}(I_{C^{c}_{2}}f) |+|\partial_{t}
P_{\alpha,k}(I_{C^{c}_{2}}f)  |
\leq N \bM_{\alpha}
f (0) .
\end{equation}
\end{corollary}

Indeed, one can concentrate on $f$ with compact
support 
when
$$
|D_{i}P_{\alpha,k}(I_{C^{c}_{2}}f)(t,x)|\leq\int_{\bR^{d+1} }
|D_{i}p_{\alpha,k}(s,|y|)|\,|(I_{C^{c}_{2}}f)(t+s,x+y)|      \,dyds,
$$
where 
$$
|D_{i}p_{\alpha,k}(s,|y|)|\leq N
\frac{1}{s^{(d+2-(\alpha-1))/2}}e^{-r^{2}/(2ks)}I_{s>0}.
$$
Similarly, for $(t,x)\in C_{1}$
$$
|\partial_{t}P_{\alpha,k}(I_{C^{c}_{2}}f)(t,x)|\leq\int_{\bR^{d+1} }
|\partial_{s}p_{\alpha,k}(s,|y|)|\,|(I_{C^{c}_{2}}f)(t+s,x+y)|      \,dyds,
$$
where 
$$
|\partial_{s}p_{\alpha,k}(s,|y|)|\leq N
\frac{1}{s^{(d+2-(\alpha-2))/2}}e^{-r^{2}/(2ks)}I_{s>0}. 
$$

Corollary \ref{corollary 10.30.1} and the mean value theorem  yield
\begin{corollary}
                        \label{corollary 10.30.2}
There is a constant $N$ such that for any $f$,
for which $P_{\alpha,k}|f|$ is finite on  $C_{1}$,
on $C_{1}$ we have
\begin{equation}
                             \label{10.30.10}
\int_{C_{1}} |
P_{\alpha,k}(I_{C^{c}_{2}}f)(z )-P_{\alpha,k}(I_{C^{c}_{2}}f)(0 )|\,dz 
\leq  N  \bM_{\alpha}
f (0) .
\end{equation}
\end{corollary}

\begin{lemma}
                   \label{lemma 10.30.2}

For any $\alpha,k>0$ there exists a constant $N$ such that for any $f$ and $\rho>0$
\begin{equation}
                             \label{10.30.3}
\int_{C_{1}}P_{\alpha,k}(I_{C _{2}}|f|) \,dz
\leq N  \bM_{\alpha }f(0).
\end{equation}
\end{lemma}

Proof. We have
$$
\int_{C_{1}}P_{\alpha,k}(I_{C _{2}}|f|)\,dz
=P_{\alpha,k}g(0),
$$
where for any $(s,y)$
$$
g(s,y)=\int_{C_{1}}(I_{C_{2}}|f|)(s+t,x+y)\,dxdt
\leq  I_{C_{3}}(s,y)\bM_{\alpha}f(0)
$$
This proves \eqref{10.30.3} and the lemma.
\qed

The next result is similar to Lemma \ref{lemma 9.5.1} with   similar proof.
By the way, if $f$ is independent of $t$,
it yields Lemma \ref{lemma 9.5.1} .  

\begin{theorem}
                           \label{theorem 10.30.1}
Let $\alpha\in(0,d+2]$ and $k>0$. Then there is a constant $N$ such that for any $f$,
for which $P_{\alpha,k}|f|$ is
  locally summable,   we have
\begin{equation}
                                 \label{10.30.5}
P^{\sharp}_{\alpha,k}f\leq N \bM_{\alpha}f.
\end{equation}

\end{theorem}

Proof.   It suffices to prove \eqref{10.30.5}
at the origin, that is, it suffices to prove that
for any $\rho>0$
\begin{equation}
                                 \label{10.30.6}
\dashint_{C_{\rho}}\dashint_{C_{\rho}}|
P_{\alpha,k}f(z_{1})-P_{\alpha,k}f(z_{2})|\,dz_{1}dz_{2}
\leq N \bM_{\alpha}f(0).
\end{equation}
 By self-similarity, it suffices to concentrate
on $\rho=1$, in which case the left-hand side of
\eqref{10.30.6} is dominated by
$$
N\int_{C_{1}} |
P_{\alpha,k}(I_{C^{c}_{2}}f)(z )-P_{\alpha,k}(I_{C^{c}_{2}}f)(0)|\,dz 
+N\int_{C_{1}}P_{\alpha,k}(I_{C _{2}}|f|) \,dz.
$$
After that the result follows immediately from
Corollary \ref{corollary 10.30.2} and Lemma \ref{lemma 10.30.2}.
\qed
  
{\bf Proof of Theorem \ref{theorem 1.18.1}}.
As we have noticed
before, to avoid trivialities, we may
assume that $q\alpha\le d+2$ when $
p\alpha< d+2$.
We may also assume that $b$ and $f$ are
nonnegative, bounded, and have
compact support. This guarantees that $I<\infty$.
Then   set
$u=P_{\alpha,k}f$, and write  
$$
I=\int_{\bR^{d+1}}\big(b^{p}u^{p-1}\big) P_{\alpha}f
\,dz=\int_{\bR^{d+1}}P_{\alpha,k}^{*}\big(b^{p}u^{p-1}\big)  f
\,dz
$$ 
\begin{equation}
                           \label{5.30.02}
\leq \|f\|_{L_{p}(\bR^{d+1})}\big\|
P_{\alpha,k}^{*}\big(b^{p}u^{p-1}\big)\big\|_{L_{p'}(\bR^{d+1})},
\end{equation}
where $p'=p/(p-1)$ and $P^{*}_{\alpha,k}$ is the conjugate operator to $P_{\alpha,k}$, namely, for any $g\geq0$,   
\begin{equation}
                           \label{5.30.10}
(P^{*}_{\alpha,k}g)
(s,x)=\big(P_{\alpha,k}(g(-\cdot,-\cdot)\big)(-s,-x).
\end{equation}
An easy computation based
on Remark \ref{remark 5.31,1}
and the fact that $p\alpha<d+2$ shows that $ P^{*}_{\alpha,k}g \in L_{p'}(\bR^{d+1})$.
 By   using Theorem \ref{theorem 10.30.1}
and the Fefferman-Stein theorem we conclude
$$
\|
P_{\alpha,k}^{*}\big(b^{p}u^{p-1}\big)\big\|_{L_{p'}(\bR^{d+1})}
\leq N\|
\bM_{\alpha} \big(b^{p}u^{p-1}\big)\big\|_{L_{p'}(\bR^{d+1})},
$$
where by H\"older's inequality
$$
\bM_{\alpha} \big(b^{p}u^{p-1}\big)
=\bM_{\alpha} \big(b(b^{p-1}u^{p-1})\big)\leq
\|b\|_{\dot E_{q,\alpha}(\bR^{d+1})} [\bM \big((bu)^{r}\big)]^{1/q'}
$$
with $r=q(p-1)/(q-1)$ and $q'=q/(q-1)$.
Since $p<q$, we have $p'>q'$ and by the Hardy-Littlewood theorem
$$
\|
\bM_{\alpha} \big(b^{p}u^{p-1}\big)\big\|_{L_{p'}(\bR^{d+1})}
\leq N\|b\|_{\dot E_{q,\alpha}(\bR^{d+1})}\Big(\int_{\bR^{d+1}}
\big[(bu)^{r}\big]^{p'/q'}\,dz\Big)^{(p-1)/p}
$$
$$
=N\|b\|_{\dot E_{q,\alpha}(\bR^{d+1})}I^{(p-1)/p}.
$$
This, obviously, proves the theorem.\qed

\begin{corollary}
                 \label{corollary 1.29.2}
Let $  q> p>1$,   let $b$ be
an $\bR^{d}$-valued function defined on $\bR^{d+1}$, such that $\|b\|_{\dot E_{q,1}(\bR^{d+1})}<\infty$. Then 
\begin{equation}
                             \label{4.4,7}
\cL=\partial_{t}
+\Delta +b^{i}D_{i}
\end{equation}
 is a bounded operator
from $W^{1,2}_{p}$ to $L_{p}(\bR^{d+1})$.

\end{corollary}

\begin{theorem} 
               \label{theorem 4.27,1}
Let $1<p<q<\infty$. Then
there is a constant $N=N(d,p,q)$
such that, if 
$$
N\|b\|_{\dot E_{q,1}(\bR^{d+1})}
\leq 1,
$$
 then
for any $s\in[-\infty,\infty)$,  $\lambda>0$  and $f\in L_{p}(\bR^{d+1}_{s})$
there exists a unique $u\in W^{1,2}_{p}(\bR^{d+1}_{s})$
satisfying 
\begin{equation}
                           \label{4.5,1}
\cL u-\lambda u=f,
\end{equation}
 where $\cL$ is defined in \eqref{4.4,7}.
Furthermore,   for any $u\in W^{1,2}_{p}(\bR^{d+1}_{s})$  and $\lambda \geq0$ 
\begin{equation}
                           \label{4.5,2}
\|\partial_{t}u,D^{2}u,\sqrt{\lambda}Du,
\lambda u\|_{L_{p}(\bR^{d+1}_{s}  )}\leq N(d,p,q)
\|\cL u-\lambda u\|_{L_{p}(\bR^{d+1}_{s})}.
\end{equation}  
\end{theorem}

Proof. A priori estimate \eqref{4.5,2}
and uniqueness 
follow directly from \eqref{1.25.6}
 and the first estimate in \eqref{1.18.3}.
To prove the existence we use the method
of continuity introduced
by S.M. Bernstein defining
$\cL_{\nu}=\partial_{t}
+\Delta +\nu b^{i}D_{i}$ for $\nu\in[0,1]$.
If the result is true for a $\nu_{0}\in[0,1]$,
then we rewrite the equation $\cL_{\nu}u-\lambda u=f$ as
$$
\cL_{\nu_{0}}u-\lambda u=f+(\nu_{0}-\nu) b^{i}D_{i}u
$$
and solve it by successive approximations
by defining $u^{0}=0$ and finding $u^{n+1}$
from
$$
\cL_{\nu_{0}}u^{n+1}-\lambda u^{n+1}=f+(\nu_{0}-\nu) b^{i}D_{i}u^{n}.
$$
It is easy to check that under our assumptions, if $|\nu_{0}-\nu|$
is small enough, all $u^{n}$ are well defined
and converge in $W^{1,2}_{p}(\bR^{d+1}_{s})$
to the solution of 
$
\cL_{\nu}u-\lambda u=f.
$   \qed

\begin{remark}
                      \label{remark 4.27,3}
The reason we explained here the method
of continuity, which we have already alluded to before, is that we wanted to show in detail
that causality is preserved after adding
first order terms.

\end{remark}

Here we did not restrict $q$ from above.
However, observe that if $q> d+2 $, then
automatically $b=0$ if $\|b\|_{\dot E_{q,1}(\bR^{d+1})}<\infty$ and the corollary becomes just
Theorem \ref{theorem 1.25.4}.
 
In the proof of Theorem \ref{theorem 1.18.1}
we implicitly used the following parabolic analog
of the Muckenhoupt-Wheeden
\index{Muckenhoupt-Wheeden theorem}%
 theorem
(see Theorem 3.6.1 in \cite{AH_96}).

\begin{corollary}[of Theorem
\ref{theorem 10.30.1}]
                       \label{corollary 10.30.4}
Let $\alpha\in(0,d+2]$, $k>0$, and $r\in(1,\infty)$. Then there is a constant $N$
such that, for any $f\geq 0$,
$\|P_{\alpha,k}f\|_{L_{r}(\bR^{d+1})}\leq N\|\bM_{\alpha}f\|_{L_{r}(\bR^{d+1})}$.

\end{corollary}

\subsection{Application of the weight theory to the heat equation}   

 For $\alpha\geq 0$
and $p>1$ introduce $\dot E_{q,\alpha}(\bR^{d+1})$ as the space
\index{$A$@Sets of functions!$\dot E_{q,\alpha}(\bR^{d+1})$}%
\index{$N$@Norms!$"|"|b"|"|_{\dot E_{q,\alpha}(\bR^{d+1})}$}%
 of functions $f$ on
$\bR^{d+1}$ with finite
$$
\|f\|_{\dot E_{q,\alpha}(\bR^{d+1})}=
\sup_{\rho>0}\rho^{\alpha}\sup_{C\in\bC_{\rho}}
\dashnorm f\|_{L_{q}(C)}.
$$

\begin{lemma}
                          \label{lemma 1.29.3}
Define $w_{\alpha}=(|x|+ \sqrt{|t|})^{-\alpha}$. Then for any $\alpha\in[0,d+2)$ the function $w_{\alpha}$
 is an $A_{1}$-weight.
\end{lemma}

Proof. We need to show that for any
$\cQ_{l}(t,x)\in \frQ$
$$
I=\dashint_{\cQ_{l}(t,x)}w_{\alpha}\,dyds\leq N
\inf_{\cQ_{l}(t,x)}w_{\alpha}.
$$
Observe that self-similar transformations
allow us to concentrate on $l=2$, in which case
$$
I=\dashint_{(-2,2)\times(-1,1)^{d}}
(|y-x| + \sqrt{|s-t|})^{-\alpha }\,dyds.
$$

If $|x|+\sqrt{ |t| }>2(d+2)$, then for $(s,y)\in \cQ_{2} 
( 0 ) $ we have 
$$
|y|+\sqrt{|s|}\leq d+2,\quad
|y-x| + \sqrt{|s-t|}\geq |x|+\sqrt{ |t| }-(|y|+
\sqrt{|s|})
$$
$$
\geq |x|+\sqrt{ |t| }-(d+2)\geq (  |x|+\sqrt{ |t| })/2
$$
  and ($\alpha\geq0$)
\begin{equation}
                             \label{1.30.1}
I\leq Nw_{\alpha}(t,x)\leq N\inf_{\cQ_{2} (t,x)}w_{\alpha}.
\end{equation}
In case that $|x|+\sqrt{ |t| }\leq 2(d+2)$ as is easy to see, for a constant $R=R(d)$,
$$
I\leq \int_{B_{R}}\int_{-R^{2}}^{R^{2}}
(|y|+\sqrt{|s|})^{-\alpha}\,dsdy
$$
$$
=N\int_{B_{R}}|y|^{2-\alpha}
\int_{-R^{2}/|y|^{2}}^{R^{2}/|y|^{2}}
(1+\sqrt{|t|})^{-\alpha}\,dtdy.
$$
For $\alpha\in[0,2]$
$$
|y|^{2-\alpha}
\int_{-R^{2}/|y|^{2}}^{R^{2}/|y|^{2}}
(1+\sqrt{|t|})^{-\alpha}\,dt\leq N(1+I_{\alpha=2}|\ln|y|\,|),
$$
so that $I\leq N$ and \eqref{1.30.1} holds.
If $\alpha>2$ the integral with respect to $t$ is bounded and, since $2-\alpha>-d$,
the integral with respect to $y$ is finite,
so that again \eqref{1.30.1} holds. 
   \qed

\begin{lemma}
                        \label{lemma 1.30.2}
Let $\alpha\geq 0,\alpha+\beta>d+2 $. Then
for any $f\geq0$  
\begin{equation}
                            \label{1.30.4}
\int_{\bR^{d+1}}f(w_{\alpha}\wedge 1)\,dyds\leq  
\bM_{\beta}f(0).
\end{equation}

\end{lemma}

Proof. Without losing generality we suppose that $f$ is bounded and has compact support. Set 
$$
Q_{1}=\{(s,y):|y|\geq\sqrt {|s| }\},\quad Q_{2}=\{(s,y):|y|\leq \sqrt{ |s|}\}.
$$
 Note that
on $Q_{1}$ we have $w_{\alpha}(s,y)\leq N |y|^{-\alpha}$
and
$$
\int_{Q_{1}}f(w_{\alpha}\wedge 1)\,dyds\leq N\int_{\bR^{d}}
(1\wedge|y|^{-\alpha})\int_{-|y|^{2}}^{|y|^{2}}
f(s,y)\,dsdy 
$$
$$
\leq N\int_{0}^{\infty}(1\wedge r^{-\alpha})
\Big(\frac{\partial}{\partial r}\int_{B_{r}}
\int_{- |y|^{2}} ^{ |y|^{2}} 
f(s,y)\,ds\,dy\Big)\,dr
$$
$$
=N\int_{1}^{\infty}r^{-\alpha-1}
 \int_{B_{r}}
\int_{- |y|^{2}} ^{ |y|^{2}} 
f(s,y)\,ds\,dy \,dr
$$
$$
\leq N\bM_{\beta}f(0)\int_{1}^{\infty}r^{-\alpha-\beta+(d+2)-1}\,dr\leq N\bM_{\beta}f(0).
$$

Next, on $Q_{2}$ we have $w_{\alpha}(s,y)\leq N |s|^{-\alpha/2}$
and 
$$
\int_{Q_{2}}f(w_{\alpha}\wedge 1)I_{s>0}\,dyds\leq N\int_{0}^{\infty}
(1\wedge s^{-\alpha/2})\int_{B_{\sqrt s}}
f(s,y)\,dyds
$$
$$
 \leq N\int_{0}^{\infty}
(1\wedge s^{-\alpha/2})\Big(\frac{\partial}{\partial s}\int_{0}^{s}\int_{B_{\sqrt s}}
f(t,y)\,dtdy\Big)ds
$$
$$
\leq N\int_{1}^{\infty}
  s^{-\alpha/2-1}  \int_{0}^{s}\int_{B_{\sqrt s }}
f(t,y)\,dtdy ds
$$
$$
\leq N\bM_{\beta}f(0)\int_{1}^{\infty}
  s^{-\alpha/2+(d+2)/2-\beta/2-1}\,ds
\leq N\bM_{\beta}f(0).
$$
Similarly we estimate the part of the integral
of $f(w_{\alpha}\wedge1)$ over $Q_{2}\cap\{s<0\}$.   \qed

The following result essentially belongs  
to \cite{CF_88}.

\begin{theorem}[Hardy-Littlewood]
                  \label{theorem 1.30.8}
Let $p>1$, $\beta>0$. Then there exists
$N=N(d,p,\beta)$ such that for any $f\geq0$
we have
\begin{equation}
                         \label{1.30.7}
\|\bM f\|_{\dot E_{p,\beta}(\bR^{d+1})}\leq N\|f\|_{\dot E_{p,\beta}(\bR^{d+1})}.
\end{equation}

\end{theorem}

Proof. For $\cQ=\cQ_{1} (0,0) $ and $\alpha\in[0,d+2)$ such that $\alpha+p\beta>d+2$,
by using the fact that $(w_{\alpha}\wedge1)$
is an $A_{1}$- and, hence, an $A_{p}$-weight
and using Lemma \ref{lemma 1.30.2}, we get
$$
\dashint_{\cQ}(\bM f)^{p}\,dxdt\leq N\int_{\bR^{d+1}}
(\bM f)^{p}(w_{\alpha}\wedge 1)\,dxdt
$$
$$
\leq N
\int_{\bR^{d+1}}
f^{p}(w_{\alpha}\wedge 1)\,dxdt\leq N\bM_{p\beta}(f^{p})(0 )\leq N\|f\|_{\dot E_{p,\beta}(\bR^{d+1})}^{p}.
$$
By using self-similar transformations and shifts of the origin we arrive at
$$
l^{p\beta}\dashint_{\cQ_{l}(t,x)}(\bM f)^{p}\,dxdt\leq \|f\|_{\dot E_{p,\beta}(\bR^{d+1})}^{p},
$$
which is equivalent to \eqref{1.30.7}.
  \qed

The following is a parabolic version
of Theorem \ref{theorem 1.20.2}
similar to the Fefferman-Stein theorem.
\index{Fefferman-Stein theorem}%
\begin{theorem}
                    \label{theorem 1.31.2}
Let $p>1,\beta>0$. Then,
for any $f\in \dot E_{p,\beta}(\bR^{d+1})$,
\begin{equation}
                          \label{1.20.03}
\|f\|_{\dot E_{p,\beta}(\bR^{d+1})}
\leq N(d,p,\beta)\|f^{\#}\|_{\dot E_{p,\beta}(\bR^{d+1})}.
\end{equation}
\end{theorem}

Proof. With the notation from the above proof
by using Theorem \ref{theorem 1.19.2} we get
$$
\dashint_{\cQ}|f|^{p}\,dxdt\leq N\int_{\bR^{d+1}}
|f|^{p}(w_{\alpha}\wedge 1)\,dxdt
$$
$$
\leq N\int_{\bR^{d+1}}
(f^{\#}) ^{p}(w_{\alpha}\wedge 1)\,dxdt\leq N\bM_{p\beta}((f^{\#}) ^{p})(0 )\leq N\|f^{\#}\|_{\dot E_{p,\beta}(\bR^{d+1})}^{p}.
$$
After that we finish the proof as that
of Theorem \ref{theorem 1.30.8}. \qed

The arguments similar to the ones used above
also derive from Lemma \ref{lemma 2.1.3} the following interpolation inequality analogous to its elliptic counterpart.
\begin{lemma}
                          \label{lemma 2.1.5}
Let $p>1,\beta>0$, $u\in C^{\infty}_{0}$. Then for any $\varepsilon>0$ 
$$
\|Du\|_{\dot E_{p,\beta}(\bR^{d+1})}
\leq N\varepsilon\|D^{2}u\|_{\dot E_{p,\beta}(\bR^{d+1})}+N\varepsilon^{-1}\| u\|_{\dot E_{p,\beta}(\bR^{d+1})},
$$
where the constants $N$ are independent of $\varepsilon,u$.

\end{lemma}

\begin{theorem}
                       \label{theorem 1.30.1}
For $p>1$, $\beta>0$, and any $u\in C^{\infty}_{0}$ we have  
\begin{equation}
                             \label{1.30.9}
\|D^{2}u\|_{\dot E_{p,\beta}(\bR^{d+1})}\leq N\|f\|_{\dot E_{p,\beta}(\bR^{d+1})},
\end{equation}
where $f=\partial_{t}u+\Delta u$
and $N$ is independent of $u$.
\end{theorem}

Proof. For $1<r<p$  
by Theorem \ref{theorem 1.26.1}
$(D^{2}u)^{\sharp}\leq N\big(\bM(|f|^{r})\big)^{1/r}$. By Theorem \ref{theorem 1.30.8}
$$
\|\big(\bM(|f|^{r})\big)^{1/r}\|_{\dot E_{p,\beta}(\bR^{d+1})}\leq N\|f\|_{\dot E_{p,\beta}(\bR^{d+1})}. 
$$
After that it only remains to apply
Theorem \ref{theorem 1.31.2}.   \qed

Introduce the space $\dot E^{1,2}_{p,\beta}$
as the space of
\index{$A$@Sets of functions!$\dot E^{1,2}_{p,\beta}$}%
 functions $u(t,x)$ on $\bR^{d+1}$ such that $$u,Du,D^{2}u,\partial_{t} u
\in  \dot E _{p,\beta}(\bR^{d+1}).
$$
 We provide
this space with an obvious norm.
Observe that $\dot E^{1,2}_{p,\beta}$
consists of only zero function if $\beta>(d+2)/p$ and otherwise $W^{1,2}_{(d+2)/\beta}\subset \dot E^{1,2}_{p,\beta}$ (a consequence of H\"older's inequality).

The following result about the heat equation 
in Morrey spaces is analogous to
its elliptic counterpart from the Morrey-space theory of the Laplacian and Sobolev-space theory of the heat equation.

\begin{theorem}
                  \label{theorem 1.31.3}
Let $p>1,\beta>0$. Then  for any $\lambda>0$  and
$f\in \dot E_{p,\beta}(\bR^{d+1})$ there exists
a unique $u\in \dot E^{1,2}_{p,\beta}$ satisfying 
$$
\partial_{t}u+\Delta u-\lambda u=f.
$$
Furthermore, there exists $N=N(d,p,\beta)$ such that for any $\lambda\geq0$ and $u\in
\dot E^{1,2}_{p,\beta}$ we have
\begin{equation}
                           \label{1.31.9}
\|\lambda u,\sqrt\lambda Du,D^{2}u,\partial_{t} u\|_{\dot E_{p,\beta}(\bR^{d+1})}\leq N\|\partial_{t}u+\Delta u-\lambda u\|_{\dot E_{p,\beta}(\bR^{d+1})}.
\end{equation}

\end{theorem}

Proof. We may assume that $p\beta\leq d+2$
because otherwise any $g$ having finite $\dot E_{p,\beta}
(\bR^{d+1})$-norm is zero. Then for $u\in C^{\infty}_{0}$ estimate \eqref{1.31.9} is derived from \eqref{1.30.9}
as in the proof of Theorem \ref{theorem 12.22.1}. To any $u\in \dot E^{1,2}_{p,\beta}$
it is extended as in the proof of Theorem
\ref{theorem 1.20.3} (see the argument about
$u^{\varepsilon} ,v^{\varepsilon} $).
This, in particular, implies uniqueness
of solutions.

To prove the existence observe that for 
$f\in L_{(d+2)/\beta}$ we have it in light
of Theorem \ref{theorem 1.25.4}. After that
we finish the proof as that of Theorem
\ref{theorem 1.20.3}. \qed

\begin{remark}
             \label{remark 2.2.1}
For $\beta=(d+2)/p$ Theorem \ref{theorem 1.31.3} coincides with Theorem
\ref{theorem 1.25.4}. This is because
$\dot E_{p,(d+2)/p}(\bR^{d+1}) = L_{p}(\bR^{d+1}) $.
\end{remark}

In order to be able to include the term
$b^{i}D_{i}u$ with $b\in \dot E_{q,1}(\bR^{d+1})$ into the right-hand side
of \eqref{1.31.9} we need some embedding
theorems, which we prove in the next section. In this connection 
observe that for such $b$, $u\in \dot E^{1,2}_{p,\beta}$, $p\geq 1,\beta>1,q=\beta p$,
and $r>0$
by H\"older's inequality
$$
\dashnorm |b|Du\|_{L_{p}(C_{r})}
\leq \dashnorm b \|_{L_{\beta p}(C_{r})}
\dashnorm  Du\|_{L_{\beta'p}(C_{r})},
$$
where $\beta'=  \beta /(\beta-1) $, which
implies the following.

\begin{lemma}
                   \label{lemma 2.2.1}
Let $p\geq 1,\beta>1$,$u\in \dot E^{1,2}_{p,\beta}$, and
$b=|b|\in \dot E_{\beta p,1}(\bR^{d+1})$.
Then
\begin{equation}
                          \label{2.2.2}
\|bDu\|_{\dot E_{p,\beta }(\bR^{d+1})}
\leq \|b\|_{\dot E_{\beta p,1}(\bR^{d+1})}
\|Du\|_{\dot E_{\beta' p,\beta-1}(\bR^{d+1})}.
\end{equation}

\end{lemma}

\subsection{Embedding theorems for
parabolic Morrey-Sobolev spaces
and  applications} 
 
\begin{lemma}
                           \label{lemma 1.29.1}
For any $\alpha>\gamma\geq0$, $k>0$ there exists a   constant  $N$
 such that for any $f\geq0$
and $\rho\in(0,\infty)$ we have
\begin{equation}
                             \label{1.17.020}
P_{\alpha,k}(I_{C _{ \rho}}f)(0)
\leq N \rho^{\alpha-\gamma }\bM_{\gamma} f (0).
\end{equation}
\end{lemma}

Proof. We are going to use the notation
from the proof of Lemma \ref{lemma 19.30.1}. We have
$$
P_{\alpha,k}(fI_{Q_{1}\cap C _{\rho}})(0)
\leq N\int_{0}^{\rho}\frac{1}{r^{d+2-\alpha}}
\int_{0}^{r^{2}}\Big(\int_{|y|=r}f(s,y)\,d\sigma_{r}\Big)\,dsdr
$$
$$
=N\int_{0}^{\rho}\frac{1}{r^{d+2-\alpha}}
\Big(\frac{\partial}{\partial r}\int_{0}^{r}
\Big(\int_{0}^{\tau^{2}}\int_{|y|=\tau}f(s,y)\,d\sigma_{\tau}\,ds\Big)\,d\tau\Big)dr
$$
$$
=J_{1}+ N
\int_{0}^{\rho}\frac{1}{r^{d+3-\alpha}}
\int_{0}^{r}\Big(\int_{0}^{\tau^{2}}
\Big(\int_{|y|=\tau}f(s,y)\,d\sigma_{\tau}\Big)\,ds\Big)\,d\tau dr 
$$
$$
\leq J_{1}+N\int_{0}^{\rho}\frac{1}{r^{d+3-\alpha}}I(r)\,dr ,
$$
where
$$
J_{1}=N\frac{1}{\rho^{d+2-\alpha}}\int_{0}^{\rho}
\Big(\int_{0}^{\tau^{2}}\int_{|y|=\tau}f(s,y)\,d\sigma_{\tau}\,ds\Big)\,d\tau 
\leq N\frac{1}{\rho^{d+2-\alpha}}I(\rho).
$$
Here $I(r)
\leq Nr^{ d+2 -\gamma}  \bM_{\gamma} f(0)$ and $\alpha>\gamma$,
so that 
\begin{equation}
                          \label{1.17.04}
P_{\alpha,k}(fI_{Q_{1}\cap C _{\rho}})(0)
\leq N\rho^{\alpha-\gamma }\bM_{\gamma} f(0).
\end{equation}
Furthermore,
$$
P_{\alpha,k}(fI_{Q_{2}\cap C _{\rho}})(0)
\leq N\int_{0}^{\rho^{2}}\frac{1}{s^{(d+2-\alpha)/2}}
\int_{|y|\leq\sqrt s}f(s,y)\,dyds
$$
$$
\leq J_{2}+ N\int_{0}^{\rho^{2}}\frac{1}{s^{(d+4-\alpha)/2}}I(\sqrt s)\,ds= J_{2}+
N\int_{0}^{\rho}\frac{1}{r^{d+3-\alpha}}I(r)\,dr,
$$
where
$$
J_{2}=N\frac{1}{\rho^{ d+2-\alpha }}\int_{0}^{\rho^{2}}\int_{|y|\leq \sqrt\tau}f(\tau,y)\,dyd\tau  \leq N\frac{1}{\rho^{d+2-\alpha}}I(\rho).
$$
This and \eqref{1.17.04} prove   \eqref{1.17.020}. \qed

Here is one of embedding
\index{embedding theorem}%
 results.
\begin{corollary}
                  \label{corollary 2.2.8}
For any $u\in C^{\infty}_{0}$,
$\gamma\in[0,2)$, $\nu\in[0,1)$, $p\geq1$
\begin{equation}
                      \label{2.2.01}
|u|\leq N(d,p,\gamma)\|u\|_{\dot E^{1,2}_{p,\gamma}},\quad |Du|\leq N(d,p,\nu)\|u\|_{\dot E^{1,2}_{p,\nu}}.
\end{equation}
\end{corollary}
 
Indeed, take $\zeta\in C^{\infty}_{0}((-1,1)\times B_{1})$ such that $\zeta(0)=1$. Then for $f:=(\partial_{t}+\Delta)(\zeta u)$   we have
$$
u(0)=-R_{0}f(0)=NP_{2,4}f(0),\quad
|u(0)|\leq N\bM_{\gamma} f(0) 
$$
$$
\leq N\bM_{\gamma}(|u|+|Du|+|D^{2}u|
+|\partial_{t}u|)
$$
and the first estimate in \eqref{2.2.01} at the origin follows
by H\"older's inequality. The second estimate is obtained similarly starting with
$$
|Du(0)|=|DR_{0}f(0)|\leq NP_{1,8}|f|(0).
$$

\begin{remark}
                     \label{remark 4.27,4}
For $p\geq 1, \beta\geq0$ introduce $E_{p,\beta}(\bR^{d+1})$ 
as the space of functions $u(t,x)$ on $\bR^{d+1}$
\index{$A$@Sets of functions!$E_{p,\beta}(\bR^{d+1})$}%
 such
\index{$N$@Norms!$"|"|f"|"|_{E_{p,\beta}(\bR^{d+1})}$}%
 that
$$
\|u\|_{E_{p,\beta}(\bR^{d+1})}:=
\sup_{\rho\leq 1}\rho^{\beta}\sup_{C\in\bC_{\rho}}\dashnorm u\|_{L_{p}(C)}<\infty
$$
Introduce the space $  E^{1,2}_{p,\beta}$
as the space of
\index{$A$@Sets of functions!$E^{1,2}_{p,\beta}$}%
 functions $u(t,x)$ on $\bR^{d+1}$ such that $$u,Du,D^{2}u,\partial_{t} u
\in  E _{p,\beta}(\bR^{d+1}).
$$
 We provide
this space with an obvious norm.
It turns out that 
for any $u\in C^{\infty}_{0}$,
$\beta\in[0,2)$, 
$1\leq p\leq (d+2)/\beta$ 
$$
|u|\leq N(d,p,\beta)\|u\|_{ E^{1,2}_{p,\beta}}.  
$$

Indeed, to obtain this estimate at the origin,
it suffices to substitute $u\zeta$ in place of $u$ in the first estimate in \eqref{2.2.01},
where $\zeta\in C^{\infty}_{0}$
with $\zeta(0)=1$ and small support
($\beta\leq (d+2)/p$ is needed to guarantee
that $\|u\zeta\|_{\dot E^{1,2}_{p,\beta}}
\leq \|u\zeta\|_{  E^{1,2}_{p,\beta}}$).
\end{remark}

\begin{remark}
                     \label{remark 2.2.6}
In the Sobolev-space theory
of parabolic equations to have
functions in $W^{1,2}_{p}$ bounded we need
$p>(d+2)/2$. In Corollary \ref{corollary 2.2.8} this restriction is removed
on account of requiring $\gamma<2$.
\end{remark}

\begin{lemma}
                    \label{lemma 2.2.4}
For  $0<\alpha <\beta$ and $\kappa>0$
    there exist a  constant  $N$
   such that for any $f\geq0$
  we have
\begin{equation}
                        \label{2.2.4}
 P_{\alpha,k}f \leq    
N(\bM_{\beta}f )^{\alpha/\beta}
(  \bM f)^{1-\alpha/\beta}.
\end{equation}
In particular (by H\"older's inequality), for any $p\in[1,\infty]$,
$q\in(1,\infty]$, and measurable $\Gamma
\subset \bR^{d+1}$
\begin{equation}
                        \label{2.2.9}
 \|P_{\alpha,k}f\|_{L_{r}(\Gamma)} \leq    
N\|\bM_{\beta}f \|_{L_{p}(\Gamma)}^{\alpha/\beta}
\|f\|_{L_{q}(\bR^{d+1})}^{1-\alpha/\beta},
\end{equation}
provided that
$$
\frac{1}{r}=\frac{\alpha}{\beta}\cdot\frac{1}{p}+
\Big(1-\frac{\alpha}{\beta}\Big)\frac{1}{q}.
$$
\end{lemma}

Proof. Observe that by Lemmas \ref{lemma 1.29.1}
and \ref{lemma 19.30.1}
\begin{equation}
                        \label{2.2.3}
P_{\alpha,k}(I_{C _{ \rho}}f)(0)
\leq N\rho^{\alpha }\bM 
f(0) ,\quad
P_{\alpha,k}(I_{C^{c}_{ \rho}}f)(0) 
\leq N\rho^{\alpha-\beta}\bM_{\beta}
f (0).
\end{equation}
Then
 \eqref{2.2.4} at the origin is easily obtained from summing up the inequalities in 
\eqref{2.2.3} and minimizing with respect 
to $\rho$. At any other point it is obtained
by changing the origin.   \qed

 \begin{corollary}
                    \label{corollary 9.29.1}
If $\alpha\in(0,(d+2)/q)$, $q\in(1,\infty)$, and $k>0$,
  then
there exists a constant $N$
  such that for any
 $f\geq0$ we have
$$
\|P_{\alpha,k}f\|_{L_{r}(\bR^{d+1})}\leq N\|f\|_{L_{q}(\bR^{d+1})}
$$
as long as   
$$
\frac{d+2}{q}-\alpha=\frac{d+2}{r}.
$$
In particular, (a classical embedding 
\index{embedding theorem}%
result) if
  $1<q<d+2$ and $u\in 
C^{\infty}_{0} $, then
$$
\|Du\|_{L_{r}(\bR^{d+1})}\leq N\|\partial_{t}u+\Delta u\|_{L_{q}(\bR^{d+1})} 
$$
as long as
$$
\quad  \frac{d+2}{q}-1= \frac{d+2}{r}.
$$
\end{corollary}

Indeed, the first assertion follows from
H\"older's inequality and
\eqref{2.2.9} with $p=\infty$ and $\beta=(d+2)/q$ ($>\alpha$). The second assertion follows from the first one with $\alpha=1$ ($<\beta$) and the fact that
for $f=\partial_{t}u+\Delta u$ we have
$$
Du(t,x)=c\int_{\bR^{d+1}_{0}}
\frac{y}{s^{(d+2)/2}}e^{-|y|^{2}/(4s)}
f(t+s,x+y)\,dyds,
$$
where $c$ is a constant  and
$(|y|/s^{1/2})e^{-|y|^{2}/(4s)}\leq
Ne^{-|y|^{2}/(8s)}$.

Here is   Morrey space analog
of Corollary \ref{corollary 9.29.1}.
\begin{theorem}
                  \label{theorem 2.2.1}
For any $k>0$, $0<\alpha<\beta$, 
$q\in(1,\infty)$, and $r$ defined by
$$
r(\beta-\alpha)=q\beta ,
$$
there is a constant $N$ such that for any
   $f\geq0$ we have
\begin{equation}
                          \label{2.2.5}
 \|P_{\alpha,k}f\|_{\dot E_{r,\beta-\alpha}(\bR^{d+1})}
\leq N \|f\|_{\dot E_{q,\beta}(\bR^{d+1})}.
\end{equation}
\end{theorem}

Proof. It suffices to prove that for any $\rho>0$
$$
\rho^{\beta-\alpha}\Big(\dashint_{C_{\rho}}|P_{\alpha,k}f|^{r}\,dz\Big)^{1/r}
\leq N\|f\|_{\dot E_{q,\beta}(\bR^{d+1})}
=:NF,
$$
that is
\begin{equation}
                          \label{2.2.6}
\rho^{\beta-\alpha-(d+2)/r}\Big(\int_{C_{\rho}}|P_{\alpha,k}f|^{r}\,dz\Big)^{1/r}
\leq NF.
\end{equation}

Observe that by H\"older's inequality
$\bM_{\beta}f\leq N F$ and by definition
$$
\Big(\int_{\bR^{d+1}}I_{C_{2\rho}}f^{q}\,dz
\Big)^{1/q}
\leq N\rho^{(d+2)/q-\beta  }F.
$$
It follows from  Lemma \ref{lemma 2.2.4} 
with $p=\infty$ that
$$
\Big(\int_{C_{\rho}}|P_{\alpha,k}(I_{C_{2\rho}}f)|^{r}
dz\Big)^{1/r}\leq N\rho^{((d+2)/q-\beta)(1-\alpha/\beta) }F
=N\rho^{(d+2)/r-\beta+\alpha}F.
$$
Furthermore, by Lemma \ref{lemma 19.30.1}
$$
\Big(\int_{C_{\rho}}|P_{\alpha,k}(I_{C^{c}_{2\rho}}f)|^{r}
dz\Big)^{1/r}
$$
$$
\leq N\rho^{(d+2)/r} \sup_{C_{\rho}}P_{\alpha,k}(I_{C^{c}_{2\rho}}f)
\leq N\rho^{(d+2)/r+\alpha-\beta}F.
$$
By combining these estimates we come to 
\eqref{2.2.6} and the theorem is proved.
\qed

\begin{remark}
                    \label{remark 9.29.2}
As usual, Theorem \ref{theorem 2.2.1}
is nontrivial only if $q\beta\leq  d+2 $
because otherwise the spaces $\dot E_{r,\beta-\alpha}(\bR^{d+1})$ and $\dot E_{q,\beta}(\bR^{d+1})$ contain only
one function: the zero function. 
\end{remark}
 
Here is one more 
\index{embedding theorem}%
embedding theorem.
\begin{theorem}
                    \label{theorem 2.2.3}
Let $p>1,\beta>1$, and $u\in C^{\infty}_{0}$. Then
$$
\|Du\|_{\dot E_{r,\beta-1}(\bR^{d+1})}\leq N(d,p,\beta)\|\partial_{t}u+\Delta u\|_{\dot E_{p,\beta}(\bR^{d+1})}
$$
as long as  $r(\beta-1)=p\beta$.
 
\end{theorem}

This is, actually, a consequence of
Theorem \ref{theorem 2.2.1} because,
according to Remark \ref{remark 11.2.1},
for $f=|\partial_{t}u+\Delta u|$ we have 
$|Du |\leq NP_{1,8}f$.

\begin{remark}
                   \label{remark 2.2.9}
For $\beta=(d+2)/p,p<d+2$ Theorem \ref{theorem 2.2.3}
yields a classical result
$$
\|Du\|_{L_{r}(\bR^{d+1})}\leq N
\|\partial_{t}u+\Delta u\|_{L_{p}(\bR^{d+1})},
$$
whenever $p>1$ and $r((d+2)/p-1)=d+2$.
\end{remark}

A useful addition to Corollary 
\ref{corollary 2.2.8} is the following
embedding
\index{embedding theorem}%
 results.
\begin{theorem}
                  \label{theorem 2.3.1}
Let $1<\beta <2$, $p>1$. Then for any
$u\in C^{\infty}_{0}$, $\rho>0$,
$(t_{i},x_{i})\in C_{\rho}$, $i=1,2$, we have
\begin{equation}
                           \label{2.3.2}
|u(t_{1},x_{1})-u(t_{2},x_{2})|\leq
N(d, \beta)\rho^{\alpha}U,
\end{equation}
where $\alpha=2-\beta,U:=\|\partial_{t}u,D^{2} u\|_{\dot E_{p,\beta}(\bR^{d+1})}$.
In case $0<\beta <1$, we have
\begin{equation}
                           \label{2.3.3}
|Du(t_{1},x_{1})-Du(t_{2},x_{2})|\leq
N(d, \beta)\rho^{1-\beta}U.
\end{equation}

\end{theorem}

Proof. Use the notation $(f)_{r}$ for the average of $f$ over $C_{r}$
and observe that by Poincar\'e's inequality
$$
\dashint_{C_{r}}|u-(u)_{r}-x^{i}(D_{i}u)_{r}|
\,dxdt
$$
$$
\leq N r^{2}\dashint_{C_{r}}|\partial_{t}u|\,dxdt+Nr\dashint_{C_{r}}
|Du-(Du)_{r}|\,dxdt,
$$
\begin{equation}
                        \label{4.6.3}
\dashint_{C_{r}}
|Du-(Du)_{r}|\,dxdt\leq Nr\dashint_{C_{r}}(
|\partial_{t}u|+|D^{2}u|)\,dxdt.
\end{equation}
Also by Theorem  \ref{theorem 2.2.3} and H\"older's
inequality  
$$
\dashint_{C_{r}}| x^{i}(D_{i}u)_{r}|
\,dxdt\leq Nr\dashint_{C_{r}}|Du|\,dxdt
\leq Nr^{\alpha}U.
$$
By combining these estimates we see that 
\begin{equation}
                     \label{2.3.4}
\dashint_{C_{r}}\dashint_{C_{r}}|u(z_{1})-u(z_{2})|\,dz_{1}dz_{2}\leq 2\dashint_{C_{r}}|u-(u)_{r} |
\,dxdt\leq Nr^{\alpha}U.
\end{equation}

This estimate also holds for translates of $C_{r}$. Then set
$$
K=\sup_{\rho>0,C\in \bC_{\rho}}\rho^{-\alpha}\sup_{z_{1},z_{2}\in C}|u(z_{1})-u(z_{2})|.
$$
Obviously, $K<\infty$ and for any $z_{1},z_{2}\in C\in \bC_{\rho}$
$$
|u(z_{1})-u(z_{2})|\leq
|u(z_{1})-u(\tilde z_{1})|+|u(\tilde z_{1})-u(\tilde z_{2})|+|u(\tilde z_{2})-u(z_{2})|.
$$
Define $\varepsilon>0$ from $4\varepsilon
^{\alpha}=1$ and take the averages with respect to
$\tilde z_{1},\tilde z_{2}$ over cylinders
$C',C''$ of class $\bC_{\varepsilon \rho}$
containing $z_{1}$ and $z_{2}$, respectively.
Then we get
$$
|u(z_{1})-u(z_{2})|\leq (1/2)K \rho 
^{\alpha}+\dashint_{C'}\dashint_{C''}
|u(\tilde z_{1})-u(\tilde z_{2})|\,d\tilde
z_{1}d\tilde z_{2}.
$$
After that replace $C',C''$ in the last integral with $\tilde C$ of class
$\bC_{4\rho}$ containing all $C,C'$ and $C''$. Then in light of \eqref{2.3.4}
and the arbitrariness of $z_{1},z_{2}$
we conclude 
$$
|u(z_{1})-u(z_{2})|\leq (1/2)K \rho 
^{\alpha}+N\rho^{\alpha}U,\quad
K\leq(1/2)K+NU,\quad K\leq NU.  
$$
This proves \eqref{2.3.2}. Similarly
\eqref{2.3.3} is derived from
\eqref{4.6.3}.  \qed

The proof of the following result  is obtained by just  combining
Lemma \ref{lemma 2.2.1} and Theorem~\ref{theorem 2.2.3}.
\begin{lemma}
                   \label{lemma 2.2.5}
For  $p,\beta>1,q\geq \beta p $, and $u\in C^{\infty}_{0}$ we have 
\begin{equation}
                         \label{2.2.7}
\| |b|Du \|_{\dot E_{p,\beta}(\bR^{d+1})}
\leq N(d,p,q,\beta)\|b\|_{\dot E_{q,1}(\bR^{d+1})} 
\|D^{2}u,\partial_{t}u\| _{\dot E_{p,\beta}(\bR^{d+1})}. 
\end{equation}

\end{lemma}
 
Having this result, 
a familiar pattern
allows the reader to deduce from
Theorem \ref{theorem 1.31.3}
the following.

\begin{theorem}
                  \label{theorem 2.2.4}
Let $p>1,\beta>1,q\geq \beta p$ and let $b$ be an $\bR^{d}$-valued function. Then 
there exists a constant $N=N(d,p,q,\beta)$
such that, if 
$$
N\|b\|_{\dot E_{q,1}(\bR^{d+1})}\leq 1,
$$
 then
 for any $\lambda>0$  and  
$f\in \dot E_{p,\beta}(\bR^{d+1})$ there exists
a unique $u\in \dot E^{1,2}_{p,\beta}$ satisfying  \eqref{4.5,1}, where $\cL$ is defined in \eqref{4.4,7}. 

Furthermore, there exists $N_{0}=N(d,p,q,\beta)$ such that for any $\lambda\geq0$ and $u\in
\dot E^{1,2}_{p,\beta}$ we have
\begin{equation}
                           \label{2.2.08}
\|\lambda u,\sqrt\lambda Du,D^{2}u,\partial_{t} u\|_{\dot E_{p,\beta}(\bR^{d+1})}\leq N_{0}\|\cL u-\lambda u\|_{\dot E_{p,\beta}(\bR^{d+1})}.
\end{equation}

\end{theorem}

\begin{remark}
                  \label{remark 2.3.1}
By using the Sobolev-space theory of parabolic 
equations with $b\in \dot E_{q,1}(\bR^{d+1})$, $q<d+2$, we can only have solutions
in $W^{1,2}_{p}$ with $p<q$ no matter
how regular the free term is, because,
generally $b\not\in L_{q+\varepsilon}(
\bR^{d+1})$ even for small $\varepsilon>0$. These solutions are not guaranteed to be even bounded if $q\leq (d+2)/2$.

At the same time with the same $b$, if
the free term is bounded with compact support, the same solution $u$ turns out
to be bounded, H\"older continuous with
$Du$ locally summable to any power
(just in case, to have $Du$
H\"older continuous we need $\beta<1$
and not $\beta>1$ as in Theorem~\ref{theorem 2.2.4}).

\end{remark}

 \mysection{Rubio De Francia extrapolation theorem
and mixed norm inequalities}

\subsection{Rubio de Francia extrapolation theorem}
                                       \label{section 11.21.10}
We work in the framework of 
 Subsection 
\ref{section 11.20.1}, where the Muckenhoupt weights are introduced.
 Here is a surprising result found
by J.L. Rubio de Francia, showing how amazingly
many functions $w$ satisfy
  $[w]_{A_{1}}<\infty$. In our exposition we follow
some ideas from the proof of Theorem 2.5
of \cite{DK_18} or the proof of Theorem 1.4
of \cite{C-UMP_11}, which are streamlined versions
of the proof of Theorem IV.5.19 of \cite{GR_85}. A  different proof one can find in 
\cite{Du_11}. See also \cite{Sa_15} for
generalizations in the case of abelian groups.

\begin{theorem}[Rubio de Francia]
                   \label{theorem 11.21.6}
Let $p,q,K_{q}\in(1,\infty)$, $w_{q}\in A_{q}$,
$[w_{q}]_{A_{q}}\leq K_{q}<\infty$.
Let  $f,g$ be Borel
nonnegative functions on $\bR^{d}$.
Then there exists $K_{p}$, depending only on $p,q,K_{q}$, and the real-analytic structure of $\bR^{d}$, such that  if
\begin{equation}
                                     \label{11.21.6}
\int_{\bR^{d}}f^{p} w_{p}\,\mu(dx)
\leq \int_{\bR^{d}}g^{p} w_{p}\,\mu(dx)
\end{equation}
for any $w_{p}\in A_{p}$ with $[w]_{A_{p}}\leq K_{p}$,
then
\begin{equation}
                                     \label{11.21.7}
F^{q} =\int_{\bR^{d}}f^{q} w_{q}\,\mu(dx)
\leq 2^{q+1}\int_{\bR^{d}}g^{q} w_{q}\,\mu(dx)=:2^{q+1}G^{q}. 
\end{equation}

\end{theorem}

Proof.  Below   by $N,N',N''$ we denote generic constants depending
only on $p,q,K_{q}$ and the real-analytic structure of $\bR^{d}$. We will also say that
some quantity is {\em under control\/} if
it is dominated by a constant $N$.

First we suppose that $g\in A_{1}$
with $[g]_{A_{1}}$ under control, 
where the way of control is to be 
specified later.
Then, in particular, $g>0$ and we will prove
that
\begin{equation}
                         \label{2.9.1}
 F^{q}\leq 2 G^{q}.
\end{equation}

While doing so, we may assume that, for a constant $c$, 
we have $f\leq cg$. Indeed, if not,
replace $f$ with $f\wedge(ng)$ and after
getting the result set $n\to\infty$.

If $p=q$, there is nothing to prove.
Then first, let $p< q$. 
By the Jones theorem we have
$w_{q}=v^{1-q}w$ with $[v,w]_{A_{1}}$  
under control and for any $r\in[0,q]$
for $w_{r}=v^{1-r}w$ the constant
$[w_{r}]_{A_{r}}$ is also under control.

We need
to show that  
$$
\int_{\bR^{d}}f^{q}v^{1-q}w\,\mu(dx)
\leq 2\int_{\bR^{d}}g^{q}v^{1-q}w\,\mu(dx)
$$
in case, of course, the right-hand side
is finite and so is the left-hand side owing to $f\leq cg$.
 
Observe that $f^{q}v^{1-q}w=f^{p}v^{1-p}u$, where $u=
 f^{q-p}v^{p-q}w $. Introduce $r=q/p$  and note that $u/w_{r}\in L_{r'}(w_{r})$, $r'=r/(r-1)$, since
$$
 (f^{q-p}v^{p-q}w)^{r'}w_{r}^{-r'}
w_{r} =f^{q}v^{1-q}w.
$$
By Lemma \ref{lemma 2.7.1} we have
$$
u=w_{r}(u/w_{r})\leq w_{r}S_{r',w_{r}}(u/w_{r})=:u_{r}\in A_{1}
$$
 with
$[u_{r}]_{A_{1}}$  under control
(despite the fact that $f$ is entering
the definition of $u$). It follows that  
$$
 \hat w _{p}:=v^{1-p} u_{r} \in A_{p},\quad [ \hat w _{p}]_{A_{p}}\leq N'.
$$
We take $K_{p}\geq N'$ and note that 
$$
f^{q}v^{1-q}w\leq f^{p}v^{1-p} u_{r} .
$$ Then by assumption
$$
\int_{\bR^{d}}f^{q} w_{q}\,\mu(dx)
\leq\int_{\bR^{d}}g^{p}
v^{1-p} u_{r} \,\mu(dx).
$$
To estimate the last integral we use
H\"older's inequality with  
$v^{1-q}w\,\mu(dx)$ in place of $\mu$.
Then we see that it is less than $G ^{p }$
times the $[(q-p)/q]$th power of
$$
\int_{\bR^{d}}\big(v^{q-p}w^{-1} u_{r} \big)^{q/(q-p)}v^{1-q}w\,\mu(dx)
$$
$$
=\int_{\bR^{d}}\big(S_{r',w_{r}}(u/w_{r})\big)^{r'}w_{r}\,\mu(dx).
$$ 
By Lemma \ref{lemma 2.7.1} the last integral is dominated by  
$$
2^{r'}\int_{\bR^{d}} (u/w_{r})^{r'}w_{r}\,\mu(dx)=2^{q/(q-p)}F ^{q}.
$$
By combining these estimates we see
that $F^{q} \leq 2 G ^{p }F^{q-p} $,
which is equivalent to \eqref{2.9.1}.

In case $q<p$   we may assume that $G<\infty$.   Then observe that
$$
\int_{\bR^{d}}f^{q}w_{q}\,\mu(dx)
\leq\Big(\int_{\bR^{d}}f^{p}g^{q-p}w_{q}\,\mu(dx)\Big)^{q/p}G^{(p-q)q/p},
$$
where 
$$
g^{q-p}w_{q}=g^{1-p}(g^{q-1}w_{q}),\quad g^{q-1} \in L_{q'}(w_{q}),
$$
$$
 g^{q-1}w_{q}
\leq w_{q}S_{q',w_{q}}(g^{q-1}):=v\in A_{1}.
$$
 It follows that
$$
g^{q-p}w_{q}\leq g^{1-p}v
\in A_{p}
$$
with the $A_{p}$-constant of 
$g^{1-p}v$ dominated by $N''$ (recall that
$[g]_{A_{1}}$ is under control). We suppose that
$K_{p}\geq N''$ and then
 due to our assumption
$$
\int_{\bR^{d}}f^{p}g^{q-p}w_{q}\,\mu(dx)
\leq\int_{\bR^{d}}gv
\,\mu(dx)
$$
$$
\leq G \Big(
\int_{\bR^{d}} \big[S_{q',w_{q}}(g^{q-1})\big]^{q'}w_{q}\,\mu(dx)\Big)^{1/q'}
\leq 2G^{q}.
$$
By combining the above estimates we come 
to \eqref{2.9.1} in both cases: $p<q,p>q$ provided
$K_{p}\geq N'\vee N''$ and $[g]_{A_{1}}$
is under control. 

In the   case of general $g$ we still may assume that $G<\infty$ and then 
define 
$$
\hat g=
\hat w_{q'}S_{q,\hat w_{q'}}(g/\hat w_{q'})\in A_{1},
$$
where 
$\hat w_{q'}=w_{q}^{-1/(q-1)}$ ($\in A_{q'}$).
Notice that $\hat g$ is well defined and
$g\leq \hat g$ because
$$
\int_{\bR^{d}}(g/\hat w_{q'})^{q}\hat w_{q'}\,\mu(dx)=G^{q}<\infty.
$$
Furthermore, since $[\hat w_{q'}]_{A_{q'}}
=[w_{q}]_{A_{q}}^{1/(q-1)}$, $[\hat g]_{A_{1}}$ is under control (see Lemma
\ref{lemma 2.7.1}).
 Of course, assuming that
\eqref{11.21.6}   holds ,
whenever $[w_{p}]_{A_{p}}\leq K_{p}$
with the above $K_{p}$ implies
that it also holds with $\hat g$ in place of $g$.
Then by the above 
result   
$$
F^{q}\leq 2\int_{\bR^{d}}\big[S_{q,w_{q'}}(g/\hat w_{q'})\hat w_{q'}\big]^{q}w_{q}\,\mu(dx)=2
\int_{\bR^{d}}\big[S_{q,w_{q'}}(g/\hat w_{q'})\big]^{q}\hat w_{q'} \,\mu(dx)
$$
$$
\leq 2^{q+1}\int_{\bR^{d}}(g/\hat w_{q'})^{q}
\hat w_{q'}\,\mu(dx)=2^{q+1}G^{q} .
$$
This proves the theorem.
\qed

\subsection
{Mixed-norms inequalities, Dong-Kim theorem}

Here we consider the product of
two Euclidean spaces: $d'$-dimen\-sional $\bR'=\{x'\}$
and $d''$-dimensional $\bR''=\{x''\}$. We assume that
each of them is provided with
  a  real analytic   structure 
related to measures $\mu'$ on $\bR^{'}$ and $\mu''$ on $\bR^{''}$.
Suppose that these measures satisfy
the doubling condition \eqref{11.21.3}
with constants $N'_{0}$ for $\mu'$
and $N''_{0}$ for $\mu''$ and with $\cQ'_{l}(x')$ and $\cQ''_{l}(x'')$ having
obvious meanings, respectively.
In the space $\bR^{d}=\bR^{'}\times\bR^{''}=\{x=(x',x'')\}$   we introduce the measure $\mu
=\mu'\times \mu''$ and $\cQ_{l}(x)=
\cQ'_{l}(x')\times\cQ''_{l}(x'')$. We are, actually, in the same setting as in Subsection \ref{section 11.20.1} with the only difference that $\mu$ is the product measure (of special type) with doubling constant $N_{0}=N_{0}'N_{0}''$.

The set of Muckenhoupt weights of class $A_{p}$ on $\bR'$ and $\bR''$
are denoted by $A_{p}(\bR')$ and
$A_{p}(\bR'')$, respectively.
As is easy to see, if $w'\in A_{p}(\bR')$ and $w''\in A_{p}(\bR'')$,
then $w'w''\in A_{p}$ and $[w]_{A_{p}}
=[w']_{A_{p}(\bR')}[w'']_{A_{p}(\bR'')}$.

Finally, for  
 $p',p''\in(1,\infty)$, and
weights $w',w''$ given on $\bR',\bR''$ introduce
\begin{equation}
                                       \label{4.9.60}
\| f
\|_{L_{p' , p''}(w',w'')} 
:=\Big(\int_{\bR''}\Big( 
 \int_{\bR' }
  |f|^{p' }w'  \,\mu '(d x ')
\Big)^{p'' /p ' }w''\,\mu''(d x'')\Big)^{1/p'' }.
\end{equation}

 \begin{theorem}[Dong-Kim \cite{DK_18}]
                \label{theorem 12.2.1}
Let $K',K'',p,p',p''\in(1,\infty)$, 
$$
w'\in
A_{p'}( \bR'),\quad w''\in
A_{p''}( \bR''),\quad [w']_{A_{p'}(\bR')}\leq K',\quad [w'']_{A_{p''}(\bR'')}\le K^{''},
$$  and
let
$f,g$ be nonnegative measurable functions on $\bR^{d}$. Then there exists
a constant $K_p$, depending only
on $K', K'',p,p',p''$ and the real-analytic
structures of $\bR',\bR''$  such that, if
\begin{equation}
                     \label{2.9.4}
\|f\|_{L_{p }(w )}
\le  \|g\|_{L_{p}(w )}
\end{equation}
  for every $w\in
A_{p } $ with
 $[w]_{p }\le K_p$, then  
\begin{equation}
                     \label{2.9.3}
\| f\|_{L_{p' , p''}(w',w'')}\le
2^{\alpha} \| g\|_{L_{p' , p''}(w',w'')},
\end{equation}
where $\alpha=2+1/p'+1/p''$.
\end{theorem}

Proof. We follow the proof of Corollary 2.7 in \cite{DK_18} or Theorem 8.1 in \cite{DK_19} that is about multiple repeated norms. 
First define $K_{p''}$
as in Theorem
\ref{theorem 11.21.6} on the basis of
$d'', p',p'' ,K''$ and then define   $K_{p}$ as in Theorem
\ref{theorem 11.21.6} on the basis of
$d,p,p ',K ' \cdot K_{p''}$ and assume that \eqref{2.9.4}  holds for every $w\in
A_{p } $ with $[w]_{p }\le K_p$. Then
take $  v_{p' }(x'')$ with $[v_{p '}]_{A_{p '}(\mu '' )}\leq K_{p''}$   and
introduce
$$
\hat f(x'')=\Big(\int_{\bR '}f^{p '}(x',x'')w ' ( x ' )
\mu '  (dx ' )\Big)^{1/p '},
$$
$$
\hat g(x'')=\Big(\int_{\bR '}g^{p '}(x',x'')w ' ( x ')
\mu '(dx ')\Big)^{1/p '}.
$$
In light of Theorem
\ref{theorem 11.21.6} we have
$$
\int_{\bR''}\hat f^{p' }v_{p '}
\,\mu''(dx'')\leq 2^{p' +1}
\int_{\bR''}\hat g^{p' }v_{p '}
\,\mu''(dx'')=  
\int_{\bR''} (2^{1+1/p '}\hat g)^{p '}v_{p '}
\,\mu''(dx'').
$$
This holds for any $v_{p' }$ with
$[v_{p '}]_{A_{p '}(\mu'')}
\leq K_{p''}$. Hence, again by 
Theorem~\ref{theorem 11.21.6}  
$$
\int_{\bR''}\hat f^{p'' }w_{p ''}
\,\mu''(dx'')\leq 2^{\alpha p''}
\int_{\bR''}\hat g^{p'' }w_{p ''}
\,\mu''(dx''),
$$
which is \eqref{2.9.3}. 
\qed

\begin{remark}
                \label{remark 4.2,1}
Here is a typical example of applications of Theorem \ref{theorem 12.2.1} when $w'=w''=1$.
Assume that we have two functions $f,g\geq0$
on $\bR^{d}$ and a finite increasing function $N(t)\geq0$ on $[0,\infty)$
such that for a fixed $p\in(1,\infty)$ and any
$w \in A_{p}$
\begin{equation}
                                 \label{4.2,6}
\int_{\bR^{d}}f^{p}w \,\mu(dx)\leq
N([w]_{A_{p}})\int_{\bR^{d}}g^{p}w \,\mu(dx).
\end{equation}
Let $p',p''\in(1,\infty)$. Then it turns out
that there exists a constant $N$ depending
only on $p,p',p''$,  the function $N(t)$,
and the real-analytic
structures of $\bR',\bR''$ such that
\begin{equation}
                                 \label{4.2,9}
\| f\|_{L_{p' , p''} (1,1) }\le
N\| g\|_{L_{p' , p''} (1,1)} .
\end{equation}

Indeed, here $K'=K''=1$ and   for every $w\in
A_{p } $ with
 $[w]_{p }\le K_p$ we have
 $$
\|f\|_{L_{p }(w )}
\le  \|\bar g\|_{L_{p}(w )},
$$
where $\bar g=N^{1/p}(K_{p})g$.
\end{remark}

\begin{remark}
                \label{remark 12.21.1}
Obviously, Theorem \ref{theorem 12.2.1}
and \eqref{4.2,9} hold true also if we interchange
$\bR'$ and $\bR''$.
 
\end{remark}

\mysection{Parabolic equations
in mixed norms Sobolev and Morrey-Sobolev spaces}
                \label{chapter 2.11.1}

In $\bR^{d+1}=\{t\in\bR,x\in\bR^{d}\}$
we introduce the real analytic structure
given by $(2,1,...,1)$ (parabolic boxes)
and Lebesgue measure $dz=dtdx$. Suppose
that we are given an $\bR^{d}$-valued
function $b$ and as usual set
$$
\cL=\partial_{t}+\Delta +b^{i}D_{i}
$$

For $p ,q\in[1,\infty]$
\index{$N$@Norms!$"|"|f"|"|_{L_{p,q}}$}%
\index{$N$@Norms!$\dashnorm f"|"|_{L_{p,q}}$}%
 and measurable $f$ and $\Gamma\subset \bR^{d+1}$ introduce
$$
\| f\|_{L_{p,q }(\bR^{d+1}) } 
= \Big( \int_{\bR}\Big(\int_{\bR^{d}}|f(t,x)|^{p }
\,dx\Big)^{q/p }\,dt \Big)^{1/q},
$$
$$
\| f \|_{L_{p,q}(\Gamma) }=\| fI_{\Gamma}\|_{L_{p,q}(\bR^{d+1}) },
$$
$$
\dashnorm f\|_{L_{p,q}(\Gamma)}
 =\|  I_{\Gamma}\|_{L_{p,q}(\bR^{d+1}) }^{-1}\| f \|_{L_{p,q}(\Gamma) }.
$$

 We also introduce the 
\index{$A$@Sets of functions!$L_{p,q}$}%
spaces $L_{p,q}(\Gamma)$ as the spaces of
functions whose $L_{p,q}(\Gamma)$-norms are finite.
We abbreviate $L_{p,q} =L_{p,q}
( \bR^{d+1})$.

Similarly we 
\index{$A$@Sets of functions!$\L_{p,q}$}%
introduce $\L_{q,p}(\Gamma)$,
$\L_{q,p} $,
$\dashnorm\cdot\|_{\L_{q,p}(\Gamma)}$ reversing
 the order of 
\index{$N$@Norms!$"|"|f"|"|_{\L_{p,q}}$}%
\index{$N$@Norms!$\dashnorm f"|"|_{\L_{p,q}}$}%
integration in $L_{p,q}$:
$$
\| f\|_{\L_{q,p }(\bR^{d+1}) } 
= \Big(\int_{\bR^{d}} \Big(\int_{\bR}|f(t,x)|^{q }
\,dt\Big)^{p/q }\,d x  \Big)^{1/p}.
$$
\subsection{Mixed norms parabolic inequalities in Sobolev spaces}

For $p ,q\in[1,\infty)$ define
$W^{2,1}_{p,q} $ 
\index{$A$@Sets of functions!$W^{2,1}_{p,q}$}%
as the set of functions $u$ on $\bR^{d+1}$ such that
$u$ and its Sobolev derivatives $Du,D^{2}u,\partial_{t}u$ are in 
$L_{p,q} $. Define
$$
\|u\|_{W^{2,1}_{p,q} }:=\|u,Du,D^{2}u,\partial_{t}u\|_{L_{p,q} } .
$$
\index{$N$@Norms!$"|"|u"|"|_{W^{2,1}_{p,q}}$}%
Similarly, we define $\W^{1,2}_{q,p} $ 
\index{$A$@Sets of functions!$\W^{1,2}_{q,p}$}%
on the basis of $\L_{q,p}$.

We start with the following version
of the Hardy-Littlewood theorem  
found in \cite{DK_18}.
\index{Hardy-Littlewood theorem}%

\begin{theorem}
                  \label{theorem 2.11.5}
Let $f$ be a nonnegative function on $\bR^{d+1}$, $p,q\in(1,\infty)$. Then 
\begin{equation}
                             \label{10.6.4}
\|\bM f \|_{L_{p,q}}
\leq N \| f \|_{L_{p,q}},\quad
\|\bM f \|_{\L_{q,p}}
\leq N \| f \|_{\L_{q,p}}
\end{equation}
where $N$  depends only on $d,p,q$.

\end{theorem}

For the proof it suffices to refer
to Theorems \ref{theorem 11.21.3}
\index{Fefferman-Stein theorem}%
and \ref{theorem 12.2.1} (and, perhaps,
Remark \ref{remark 4.2,1}).

  Here is the mixed norm 
version of Fefferman-Stein theorem.

\begin{theorem}[Dong-Kim \cite{DK_18}]
                  \label{theorem 2.14.2}   
Let $p,q\in(1,\infty)$. Then there
is a constant $N=N(d,p,q)$ such that
for any $f\in L_{p,q}$
\begin{equation}
                           \label{2.14.4}
\|f\|_{L_{p,q}}\leq N\|f^{\#}\|_{L_{p,q}}
\end{equation}

\end{theorem}

For the proof it suffices to refer
to Theorems \ref{theorem 1.19.2}
and \ref{theorem 12.2.1}. Similar statement
holds for $f\in\L_{q,p}$.

Applying Lemma \ref{lemma 2.1.1} for each
fixed $t$ and then integrating we get
the first part of the following
interpolation inequalities.

\begin{lemma}
                    \label{lemma 2.14.1}
Let $p,q\in[1,\infty)$ and $u\in C^{\infty}_{0}$. Then there is a constant $N$,
depending only on $d,p,q$, such that
for any $r,\varepsilon>0$
\begin{equation}
                         \label{2.14.2}
\| Du \|_{L_{p,q}(C_{r})}\leq N((\varepsilon\wedge r)
\|D^{2} u \|_{L_{p,q}(C_{r})}
+(\varepsilon\wedge r)^{-1}\|u\|_{L_{p,q}(C_{r})}).  
\end{equation}
\begin{equation}
                         \label{2.14.20}
\| Du \|_{L_{p,q}}\leq N(\varepsilon
\|D^{2} u \|_{L_{p,q}}
+\varepsilon^{-1}\|u\|_{L_{p,q}}).   
\end{equation}
\end{lemma}
Estimate \eqref{2.14.20} is  obtained
from \eqref{2.14.2} by sending
$r\to\infty$, of course, if the support of $u$ is in $\bigcup_{r}C_{r}$. The latter is irrelevant and can be taken care of
by shifts of the origin.

Lemma \ref{lemma 3.31.4} proved later shows that
\eqref{2.14.20} holds with $\L_{q,p}$
in place of $L_{p,q}$ if $\|D^{2} u \|_{L_{p,q}}$ is replaced with $\|\partial_{t}u,D^{2} u \|_{\L_{q,p}}$.

\begin{theorem}
                 \label{theorem 2.14.1}
Let $p,q\in(1,\infty)$, $s\in[-\infty,\infty)$. Then there is
a constant $N$, depending only on
$d,p,q$, such that for any
$\lambda \geq 0$ and $u\in W^{2,1}_{p,q}(\bR^{d+1}_{s})$
\begin{equation}
                         \label{2.14.1}
\|\lambda u,\sqrt\lambda Du,D^{2}u,
\partial_{t}u\|_{L_{p,q}(\bR^{d+1}_{s})}\leq N
\|\partial_{t}u+\Delta u-\lambda u\|_{L_{p,q}(\bR^{d+1}_{s})}.
\end{equation}
Furthermore, for any $\lambda>0$ and
$f\in L_{p,q}(\bR^{d+1}_{s})$, there exists a unique
$u\in W^{2,1}_{p,q}(\bR^{d+1}_{s})$ satisfying in $\bR^{d+1}_{s}$
\begin{equation}
                       \label{4.29,1}
\lambda u-\Delta u-\partial_{t}u=f.
\end{equation}
\end{theorem}

Proof. We start with the main case that $s=-\infty$. First, in light of Theorems
\ref{theorem 1.26.1} (with small $p$), \ref{theorem 2.11.5}, and \ref{theorem 2.14.2} we obtain \eqref{2.14.1}
for $\lambda=0$ and $u\in C^{\infty}_{0}$. Then, as a few times before,
we use Agmon's idea and the interpolation
inequality from Lemma \ref{lemma 2.14.1}
to extend \eqref{2.14.1} to all $\lambda
\geq0$ and $u\in C^{\infty}_{0}$. Since $C^{\infty}_{0}$ is dense in $W^{2,1}_{p,q}$, a priori
estimate \eqref{2.14.1} is proved.

The solvability for $\lambda>0$ is proved,
for instance, by almost literally repeating the corresponding part of
the proof of Theorem \ref{theorem 12.22.1}
(assume $(\partial_{t}+\Delta -\lambda)
W^{2,1}_{p,q}$ is not dense in $L_{p,q}$
...).  

In the case of general $s\in\bR$
we only need to prove uniqueness.
So let $f(t,x)=0$ for $t>s$
and $u\in W^{2,1}_{p,q}(\bR^{d+1}_{s})$ satisfy \eqref{4.29,1} in $\bR^{d+1}_{s}$. Continue $u$ across $\{t=s\}$
in an even way, keep the notation $u$ for the function thus obtained, and
call $f$ the left-hand side of \eqref{4.29,1}. We claim that (a.e.)
 \begin{equation}
                       \label{4.29,2}
u=R_{\lambda}f,
\end{equation}
where $R_{\lambda}$ is introduced
in \eqref{4.1,3}. 

Indeed,
if $v\in C^{\infty}_{0}$,
we have $v=R_{\lambda}(\lambda v-\Delta v-\partial_{t}v)$ by \eqref{4.1,3}. Since by Minkowski's
inequality the norm
of $R_{\lambda}$ in  $L_{p,q}$ is less than $1/\lambda$, \eqref{4.29,2} follows
by continuity.

Now from \eqref{4.29,2} we see
that $u(t,x)=0$ for $t\geq s$ and
the theorem is proved.
\qed 

\begin{remark}
                     \label{remark 2.25.1}
For $\lambda=0$ estimate \eqref{2.14.1} reads
\begin{equation}
                         \label{2.25.1}
\| D^{2}u,
\partial_{t}u\|_{L_{p,q}}\leq N
\|\partial_{t}u+\Delta u \|_{L_{p,q}}.
\end{equation}
It is important to know that \eqref{2.25.1} holds
for not only $u\in W^{2,1}_{p,q}$. For
instance, let $u$ be infinitely differentiable
and such that 
$$
\xi |u|+\eta|Du|
$$
 are bounded
with $\xi,\eta$ introduced in Remark
\ref{remark 4.4,3}   with $\alpha>d/p+2/q-1$, then, by repeating the computations 
in Remark \ref{remark 4.4,3} with insignificant changes,
we obtain \eqref{2.25.1} for such an $u$.

One of typical examples is $u=P_{2,4}f$
with $f\in C^{\infty}_{0}$. Indeed,
according to \eqref{1.25.2} the quantity
$\xi |u|+\eta|Du|$ is bounded for $\alpha\leq d+1$
and there are such $\alpha$ which also
satisfy $\alpha>d/p+2/q-1$ ($p,q>1$).

\end{remark}

Theorem \ref{theorem 2.11.5}, Lemma \ref{lemma 2.2.4}, and H\"older's inequality immediately lead to the following.

\begin{lemma}
                           \label{lemma 10.7.1}
For $k>0$, $0<\alpha<\beta$, 
 $p\in[1,\infty]$,
$q_{1},q_{2}\in(1,\infty]$, there exists a constant $N$ such that for any $f\geq0$    and measurable $\Gamma$ we have
\begin{equation}
                                \label{10.7.1}
 \|P_{\alpha,k}f\|_{L_{r_{1},r_{2}}(\Gamma)} \leq    
N\|\bM_{\beta}f \|_{L_{p}(\Gamma)}^{\alpha/\beta}
\|f\|_{L_{q_{1},q_{2}}}^{1-\alpha/\beta},
\end{equation}
provided that
$$
\frac{1}{r_{i}}=\frac{\alpha}{\beta}\cdot\frac{1}{p}+
\Big(1-\frac{\alpha}{\beta}\Big)\frac{1}{q_{i}},\quad i=1,2.
$$
\end{lemma}  

Similarly to Corollary \ref{corollary 9.29.1}
we have
 \begin{corollary}
                    \label{corollary 10.7.01}
Let    
$q_{1},q_{2}\in(1,\infty]$,
$$
\beta:=\frac{d}{q_{1}}+\frac{2}{q_{2}}>0,\quad \alpha\in(0,\beta ),
\quad k>0.
$$
Then for any $f\geq0$
we have
$$
\|P_{\alpha,k}f\|_{L_{r_{1},r_{2}}}\leq N\|f\|_{L_{q_{1},q_{2}}}
$$
as long as $q_{i}\beta=r_{i}(\beta-\alpha)$,
$i=1,2$.

In particular, $($almost follows from
Theorem 10.2 of \cite{BIN_75}$)$
if    $\beta>1$, and    $u\in C^{\infty}_{0} $, then
\begin{equation}
                            \label{10.8.4}
\|Du\|_{L_{r_{1},r_{2}}}\leq N\|\partial_{t}u+\Delta u\|_{L_{q_{1},q_{2}}}
\end{equation}
as long as
$q_{i}\beta=r_{i}(\beta-1)$,
$i=1,2$.
\end{corollary}

\begin{corollary}
                     \label{corollary 10.8.2}
Under the assumptions of Corollary
\ref{corollary 10.7.01}, if $\beta>1$, there is a constant $N$ such that, for any $b=(b^{i})
\in L_{\beta q_{1} ,\beta q_{2}}$ and $u\in C^{\infty}_{0} $,
\begin{equation}
                         \label{2.29.2}
\|b^{i}D_{i}u\|_{L_{q_{1},q_{2}}}\leq N
\|b  \|_{L_{\beta q_{1},\beta q_{2}}}
\|\partial_{t}u+\Delta u\|_{L_{q_{1},q_{2}}}.
\end{equation}

\end{corollary}

Indeed, by H\"older's inequality
$$
\|b^{i}D_{i}u\|_{L_{q_{1},q_{2}}}
\leq \|b  \|_{L_{\beta q_{1},\beta q_{2}}}
\|Du\|_{L_{r_{1},r_{2}}}. 
$$
\begin{remark}
                       \label{remark 5.11,1}
Theorem \ref{theorem 2.14.1}, Lemma \ref{lemma 10.7.1} and Corollaries \ref{corollary 10.7.01}
and \ref{corollary 10.8.2} remain true
if we replace $L$ with $\L$ and in the corollaries in the definition of $\beta$
exchange $q_{1}$ and $q_{2}$.

\end{remark}
 
\begin{remark}
                   \label{remark 2.17.1}

Corollary \ref{corollary 10.8.2} and Theorem \ref{theorem 2.14.1} immediately lead to the unique solvability in $W^{2,1}_{p,q}$ of 
$$
\partial_{t}u+\Delta u
+b^{i}D_{i}u-\lambda u=f\in L_{p,q}
$$
for $\lambda>0$ provided $\|b\|_{L_{\beta p ,\beta q}}$ is small enough.
Observe that in the corollary 
$$
\frac{d}{ \beta q_{1}}+\frac{2}{ \beta q_{2}}=1.
$$
 In particular, $\beta p> d$, which does not allow singular $b$ like
$|x|^{-1}$ or even $b\in L_{d}$. Still allowing
such singular $b$ is possible in the Sobolev space theory. 
\end{remark}

Here is mixed-norm analog of the parabolic Adams  theorem containing Theorem~\ref{theorem 1.18.1}.
\begin{theorem}
                        \label{theorem 5.25,1}
Let $ \alpha>0 ,q_{1},q_{2}\in(1,\infty),q>\max(q_{1},q_{2})$, $k>0$, $b(x)\geq0$. Then for any $f(t,x)\geq0$
\begin{equation}
                            \label{5.25,1}
\|bP_{\alpha,k}f\|_{L_{q_{1},q_{2}}}
\leq N\|b\|_{\dot E_{q,\alpha}(\bR^{d+1})}
\|f\|_{L_{q_{1},q_{2}}},
\end{equation}
where $N$ depends only on $d,q_{i},q,\alpha,k$.
In particular, for any $u\in C^{\infty}_{0}$
$$
\|b Du \|_{L_{q_{1},q_{2}}}
\leq N\|b\|_{\dot E_{q,1}(\bR^{d+1})}K,
\quad \| b u\|_{L_{q_{1},q_{2}}}
\leq N\|b\|_{\dot E_{q,2}(\bR^{d+1})}K,
$$
where $K=\|D^{2}u,\partial_{t}u\|_{L_{q_{1},q_{2}}}$
and $N$ depends only $d,q_{i},q$.
\end{theorem}

Proof. We may assume that  $ d+2\geq \alpha q $ because otherwise either
$b=0$ or $\|b\|_{\dot E_{q,\alpha}(\bR^{d+1})}=\infty$.
We may also assume that $b$ is bounded
and has compact support. Define
$q_{i}'=q_{i}/(q_{i}-1)$ and let $H$ be the subset
of the unit ball in $L_{q_{1}',q_{2}'}$
consisting 
of nonnegative bounded functions with compact support.
Then the left-hand side of \eqref{5.25,1}
equals
$$
\sup_{h\in H}\int_{\bR^{d+1}}h b P_{\alpha,k}f
\,dz=\sup_{h\in H}\int_{\bR^{d+1}}P_{\alpha,k}^{*}(h b )f
\,dz
$$
$$
\leq \sup_{h\in H}\|P_{\alpha,k}^{*}(h b )\|_{L_{q_{1}',q_{2}'}}\|f\|_{L_{q_{1},q_{2}}}.
$$ 
It is not hard to check that 
 $P_{\alpha,k}^{*}(h b )\in
L_{q_{1}',q_{2}'}$ by using Remark
\ref{remark 5.31,1} and the fact that
$$
\alpha\leq \frac{d+2}{q}<\frac{ d}{q_{1}}
+\frac{2}{q_{2}}.
$$
 
Then owing to Theorems \ref{theorem 2.14.2} and \ref{theorem 10.30.1}  
$$
\|P_{\alpha,k }^{*}(h b )\|_{L_{q_{1}',q_{2}'}}
\leq N 
\|(P_{\alpha,k }^{*}(h b ))^{\sharp}\|_{L_{q_{1}',q_{2}'}}\leq
N \|\bM_{\alpha}(h b )\|_{L_{q_{1}',q_{2}'}}.
$$
 
Here with $q'=q/(q-1)$ by H\"older's inequality
$$
\bM_{\alpha}(h|b|)\leq \big(\bM(|h|^{q'})\big)^{1/q'}\|b\|_{\dot E_{q,\alpha}(\bR^{d+1})}
$$
and it only remains to observe that
in light of Theorem \ref{theorem 2.11.5}
$$
\|\big(\bM(|h|^{q'})\big)^{1/q'}\|_{L_{q_{1}',q_{2}'}}=\| \bM(|h|^{q'}) \|_{L_{q_{1}'/q',q_{2}'/q'}}^{1/q'}
$$
$$
\leq N\|  |h|^{q'}  \|_{L_{q_{1}'/q',q_{2}'/q'}}^{1/q'}=N\|h\|_{L_{q_{1}',q_{2}'}}.
$$
This proves the theorem. \qed

The following is proved in the same way 
as Theorem \ref{theorem 4.27,1}.

\begin{theorem} 
               \label{theorem 5.25,3}
Let $q_{1},q_{2}\in(1,\infty),q>\max(q_{1},q_{2})$. Then
there is a constant $N=N(d,q_{1},q_{2}, q )$
such that, if 
\begin{equation}
                      \label{5.31,4}
N\|b\|_{\dot E_{q,1}(\bR^{d+1})}
\leq 1,
\end{equation}
 then
for any   $\lambda>0$  and $f\in L_{q_{1},q_{2}} $
there exists a unique $u\in W^{2,1}_{q_{1},q_{2}} $
satisfying 
$$
\cL u-\lambda u=f,
$$
 where $\cL$ is defined in Remark
\ref{remark 2.17.1}.
Furthermore,   for any $u\in W^{2,1}_{q_{1},q_{2}} $  and $\lambda \geq0$ 
$$
\|\partial_{t}u,D^{2}u,\sqrt{\lambda}Du,
\lambda u\|_{L_{q_{1},q_{2}} }\leq N(d, q_{1},q_{2} )
\|\cL u-\lambda u\|_{L_{q_{1},q_{2}} }.
$$  
\end{theorem}

Just in case, observe that \eqref{5.31,4}
implies that either $b=0$ or $q\leq d+2$.

The following result turns out to be also 
 useful in various circumstances.  
We will be using it only for $p=q=1$.

\begin{lemma}[Poincar\'e's inequality] 
                     \label{lemma 10.10.1}
Let $1\leq q,p<\infty$,
$u\in C^{\infty}_{0}$, $\rho\in(0,\infty)$.
\index{Poincar\'e's inequality}%  
Then   
\begin{equation}
                         \label{10.10.1}
\dashnorm Du-(Du)_{C_{\rho}}\|_{L_{p,q}(C_{\rho})}^{q}\leq
N(d,q,p)\rho^{q}\dashnorm \, \partial_{t}u, D^{2}u \|_{L_{p,q}(C_{\rho})}^{q}.
\end{equation}

\end{lemma}

Proof. We follow the usual way (see, for instance,
Lemma 4.2.2 of \cite{Kr_08}). First, due
to self-similar transformations, we may take
$\rho=1$. In that case,  for
a $\zeta\in C^{\infty}_{0}(B_{1})$ with unit integral, introduce
$$
v(t)=\int_{B_{1}}\zeta(y)Du(t,y)\,dy.
$$
Then by the usual Poincar\'e inequality
$$
\int_{B_{1}}|Du(t,x)-v(t)|^{p}\,dx=
\int_{B_{1}}\big|\int_{B_{1}}[Du(t,x)-Du(t,y)]\zeta(y)\,dy\big|^{p}\,dx
 $$
 \begin{equation}                
                              \label{7.8.3}
\leq N\int_{B_{1}}\int_{B_{1}}
|Du(t,x)-Du(t,y)|^{p}\,dxdy
\leq N\int_{B_{1}}|D^{2}u(t,x)|^{p}\,dx.
\end{equation}

Next, observe that for any constant vector $v$
the left-hand side of
\eqref{10.10.1} is less than a constant times (recall that $\rho=1$)
 $$
\int_{0}^{1}\Big(\int_{B_{1}}|Du(t,x)-v|^{p}\,dx\Big)^{q/p}dt
$$
$$
\leq N
\int_{0}^{1}\Big(\int_{B_{1}}|Du(t,x)-v(t)|^{p}\,dx
\Big)^{ q /p}dt+N
\int_{0}^{1}|v(t)-v|^{q}\,dt. 
 $$
By \eqref{7.8.3} the first term on the right is less
than the right-hand side of \eqref{10.10.1}. 
To estimate the second term,
take 
 $$
v=\int_{0}^{1}v(t)\,dt.
 $$
Then by Poincar\'e's inequality
 $$
\int_{0}^{1}|v(t)-v|^{q}\,dt\leq N
\int_{0}^{1}\big|\int_{B_{1}}\zeta 
\partial_{t} Du\,dx\big|^{q}\,dt
$$
$$
=N
\int_{0}^{1}\big|\int_{B_{1}}(D\zeta) 
\partial_{t}u\,dx\big|^{q}\,dt 
$$
and to finish the proof it only remains
to use  
 H\"older's inequality.  \qed

The following is a corollary of
Theorem 10.2 of \cite{BIN_75} which
we give with a different proof.

\begin{theorem}
                 \label{theorem 2.16.2}
Let  
$r\geq p$  and assume
that 
$$
\gamma:= \frac{d}{p}
+\frac{2}{q} -\frac{d}{r} <1
$$
 Then for any
$u\in W^{2,1}_{p,q}$ and $\varepsilon>0$ we 
have
\begin{equation}
                            \label{6,17.2}
\|D u(0,\cdot)\|_{L_{r}}\leq
N\varepsilon \|\partial_{t}u,D^{2}u
\|_{L_{p,q} }+N\varepsilon^{-(1+\gamma)/(1-\gamma)}
\|u
\|_{L_{p,q} }.
\end{equation}
\end{theorem}

Proof. The case of arbitrary $\varepsilon>0$
is reduced to that of $\varepsilon=1$
by using   self-similarity.
To treat $\varepsilon=1$ take $\zeta\in
C^{\infty}_{0}(\bR)$ such that $\zeta(t)=1$
for $|t|\in[0,1]$, $\zeta(t)=0$ for $|t|\geq 2$, and define 
$$
-f=
\partial_{t}(\zeta u)+\Delta(\zeta u).
$$ 
We know  that $\zeta u=R_{0}f$, where
$$
R_{0}f(t,x)=N(d)\int_{0}^{\infty}s^{-d/2}
\int_{\bR^{d}}e^{-|x-y|^{2}/(4s)}
f(t+s,y)\,dyds.
$$
It follows that  
$$
D u(0,x)= -N(d)\int_{0}^{\infty}t^{-d/2}
\int_{\bR^{d}} \frac{y}{2t} e^{-| y|^{2}/(4t)}f(t,x-y)\,dydt.
$$
By observing that $|y|/\sqrt t|
e^{-| y|^{2}/(4t)}\leq Ne^{-| y|^{2}/(8t)}$ we conclude that
$$
|D u(0,x)|\leq NF(x),
$$
where
$$
 F(x) =
\int_{0}^{2} 
\int_{\bR^{d}}p_{1,8}(t,y)|f(t,x-y)|\,dydt .
$$

By Minkowski's inequality
$$
\|F\|_{L_{r} }
\leq \int_{0}^{2}\Big(
\int_{\bR^{d}}\Big(\int_{\bR^{d}}
p_{1,8}(t,x-y)|f(t,y)|\,dy\Big)^{r}
\,dx\Big)^{1/r}
\,dt,
$$ 
where inside the integral with respect
to $t$ we have the norm of   convolution, so that by Young's
inequality this expression is dominated by
$$
\|f(t,\cdot)\|_{L_{p}}\|p_{1,8}(t,\cdot)
\|_{L_{s}},
$$
where $1/s=1+1/r-1/p$ ($ \leq 1$ since $r\geq p$). An easy computation
shows that $\|p_{1,8}(t,\cdot)
\|_{L_{s}}=N(d)t^{\alpha}$ with
$\alpha=-1/2+(d/2)(1/r-1/p)$, which  yields
$$
\|F\|_{L_{r} }
\leq N\int_{0}^{2}\|f(t,\cdot)\|_{L_{p} }t^{\alpha}\,dt.
$$
Now use H\"older's inequality along with the observation that
$\alpha q/(q-1)>-1$, due to 
the assumption that $\gamma<1$, that is,  $1 +d/r>
d/p+2/q$, to conclude that
$$
\|F\|_{L_{r} }
\leq N \|f \|_{L_{p,q}  }.
$$  \qed

\begin{corollary}
           \label{corollary 6,19.1}
For any $\rho>0,\varepsilon>0$ and
$u\in W^{2,1}_{p,q}(C_{2\rho})$ we have
$$  
\dashnorm D  u(0,\cdot)\|_{L_{r}(B_{\rho})}
\leq N\varepsilon  
\dashnorm \partial_{t}u,D^{2}u
\|_{L_{p,q}(C_{2\rho})}
$$
$$
+N
\big(\varepsilon \rho^{-2}+\varepsilon^{-(1+\gamma)/(1-\gamma)} \rho^{2\gamma/(1-\gamma)}\big)
\dashnorm u
\|_{L_{p,q}(C_{2\rho})}.
$$
\end{corollary}

Indeed, the case of arbitrary $\rho>0$ is reduced
to $\rho=1$ by means of parabolic dilation. In the latter case
take $\zeta\in C^{\infty}_{0}
(\bR^{d+1})$ such that $\zeta=1 $ on $C_{1}$ and $\zeta=0$ in $\bR^{d+1}_{0}
\cap C_{2}$ and $\zeta(-t,x)=\zeta(t,x)$. Then  apply
\eqref{6,17.2} to $u(|t|,x)$ in place of
$u(t,x)$ to see that 
$$
\|D u(0,\cdot)\|_{L_{r}(B_{1})}
\leq N\varepsilon  
\|\partial_{t}(\zeta u),D^{2}(\zeta u)
\|_{L_{p,q} }+N\varepsilon^{-(1+\gamma)/(1-\gamma)}
\|u
\|_{L_{p,q}(C_{2}) } 
$$
$$
\leq N\varepsilon  
\|\partial_{t}u,D^{2}u
\|_{L_{p,q}(C_{2}) }+N\varepsilon \|u,Du\|_{L_{p,q}(C_{2}) }
+N\varepsilon^{-(1+\gamma)/(1-\gamma)}
\|u
\|_{L_{p,q}(C_{2}) }.
$$
After that it only remains to use
the interpolation inequality
$$
\|Du(t,\cdot)\|_{L_{p}(B_{2})}\leq
N \|D^{2}u(t,\cdot)\|_{L_{p}(B_{2})}+
N\| u(t,\cdot)\|_{L_{p}(B_{2})}. 
$$  

\subsection{Mixed norms parabolic inequalities in Morrey-Sobolev spaces}

The above results concern Sobolev space
theory with mixed norms.
To extend them
to Morrey spaces with mixed norms
take $p,q\in(1,\infty)$ and
   $\beta> 0$ and introduce the
Morrey space $\dot E_{p,q,\beta} $
as the set of $g\in  L_{p,q,\loc}$ 
\index{$A$@Sets of functions!$\dot E_{p,q,\beta}$}%
\index{$N$@Norms!$"|"|f"|"|_{\dot E_{p,q,\beta}}$}%
such that  
\begin{equation}
                             \label{2.16.5}
\|g\|_{\dot E_{p,q,\beta} }:=
\sup_{\rho >0,C\in\bC_{\rho}}\rho^{\beta}
\dashnorm g  \|_{ L_{p,q}(C)} <\infty .
\end{equation}  
Define
$$
\dot E^{2,1}_{p,q,\beta} =\{u:u,Du,D^{2}u, 
\partial_{t}u\in \dot E_{p,q,\beta} \},
$$
where $Du,D^{2}u,
\partial_{t}u$ are Sobolev
\index{$A$@Sets of functions!$\dot E^{1,2}_{p,q,\beta}$}%
 derivatives,
and 
provide $\dot E^{2,1}_{p,q,\beta} $ with an obvious norm.
Recall that the subsets of these spaces
consisting of functions independent of
$t$ are denoted by $\dot E_{p,\beta}$
and $\dot E^{2}_{p,\beta}$, respectively.
\index{$A$@Sets of functions!$\dot E_{p,\beta}$}%
\index{$A$@Sets of functions!$\dot E^{2}_{p,\beta}$}%
It is a good idea for the
 reader to keep in mind that
if $\beta>d/p+2/q $ and $
f\in \dot E_{p,q,\beta} $, then
$f=0$ (a.e.).

By replacing above $L_{p,q}$ with $\L_{q,p}$
we introduce $\dot \L_{q,p,\beta}$,
$\dot \L^{1,2}_{q,p,\beta}$ and the respective
norms.

For functions $f(t,x)$ and $\varepsilon>0$ we define
\begin{equation}
                        \label{2.16.6}
f^{(\varepsilon)}=f*\zeta_{\varepsilon},
\end{equation}
 where
  $\zeta_{\varepsilon}=\varepsilon^{-d-2}\zeta
(t/\varepsilon^{2},x/\varepsilon)$, with a nonnegative
$\zeta\in C^{\infty}_{0}$ which has unit integral and support
in $C_{1}(-1,0)$. Observe that
owing to Minkowski's inequality  
$$
\|f^{(\varepsilon)}\|_{\dot E_{q,p,\beta}}
\leq \|f \|_{\dot E_{q,p,\beta}} 
$$
for any $f\in \dot E_{q,p,\beta}$.
Maximal function boundedness is
another notable property of Morrey spaces. One of the first results of that kind
in elliptic setting appeared in 
\cite{CF_88}. The following result also holds
if we replace $\dot E_{p,q,\beta}$ with
$\dot \E_{p,q,\beta}$.
 
\begin{theorem}
              \label{theorem 2.15.1}
Let   $p,q\in(1,\infty)$  $\beta>0$. Then there is a constant
$N $ such that, for any $f\geq0$

 \begin{equation}
                   \label{2.15.1}
\|\bM f\|_{\dot E_{p,q,\beta} }\leq
N   \| f\|_{\dot E_{p,q,\beta} }.
\end{equation}
In particular,
$$
\|\sup_{\varepsilon>0}f^{(\varepsilon)}
\|_{\dot E_{p,q,\beta} }\leq
N   \| f\|_{\dot E_{p,q,\beta} }.
$$
\end{theorem}   

Proof. It suffices to prove that
\begin{equation}
                    \label{2.15.2}
 I:=\|\bM f\|_{L_{p,q}(C_{1}) }
\leq
N   \| f\|_{\dot E_{p,q,\beta} }.
\end{equation}
Indeed, in that case \eqref{2.15.2}
holds with any $C\in\bC_{1}$ in place
of $C$, by homogeneity it extends
to any $\rho>0$ as
$$
 \rho^{\beta}\dashnorm\bM f\|_{L_{p,q}(C ) }
\leq
N   \| f\|_{\dot E_{p,q,\beta} }
$$
for any $C\in \bC_{\rho}$ and this is equivalent to \eqref{2.15.1}.

While proving \eqref{2.15.2}, introduce
$C$ as $2C_{1}$, that is the twice dilated $C_{1}$ 
concentric with $C_{1}$. Then
$I\leq I_{1}+I_{2}$, where
$$
I_{1}=\|\bM (fI_{C})\|_{L_{p,q}(C_{1}) }\leq \|\bM (fI_{C})\|_{L_{p,q} } 
$$
$$
\leq N\| fI_{C} \|_{L_{p,q} }
\leq N\| f\|_{\dot  E_{p,q,\beta} },
$$
and
$$
I_{2}=\|\bM (fI_{C^{c}})\|_{L_{p,q}(C_{1}) }\leq  \sup_{C_{1}}
\bM (fI_{C^{c}})
$$
$$
\leq\sup_{z\in C_{1}}
\sup_{\rho\geq 1/2}\sup_{C'\in\bC_{\rho},C'\ni z}\dashint_{C'} f
\,dy
\leq 2^{\beta}\| f\|_{\dot E_{p,q,\beta} }.
$$  \qed

Here are   useful approximation results,
which are proved in the same way 
as similar results in Subsection 
\ref{section 2.16.1}.    Set
$$
\|f\|_{E_{p,q,\beta}}:=\sup_{\rho \leq 1,C\in\bC_{\rho}}\rho^{\beta}
\dashnorm f  \|_{ L_{p,q}(C)},\quad \|f\|_{E_{p,q,\beta}(C)}:=
\|fI_{C}\|_{E_{p,q,\beta}}
$$   
and introduce the space $E_{p,q,\beta}$ as the set
of functions with finite norm introduced above.
\begin{lemma}
                               \label{lemma 2.16.2}
Let $q,p\in(1,\infty)$, $0\leq\beta'<\beta$.
If   $\|f_{n}\|_{E_{p,q,\beta'}}$,
$n=0,1,...$, is a {\em bounded\/} sequence
and $f_{n}\to f_{0}$ in $L_{p,q}(C)$
for any
$C\in \bC$, then for any
$C\in \bC$
$$
\lim_{n\to\infty}\|f_{n}-f_{0}\|_{E_{p,q,\beta}(C)}=0. 
$$
In particular, if
  $f\in E_{p,q,\beta'}$, then for any $C\in \bC$
\begin{equation} 
                            \label{2.16.3}
\lim_{\varepsilon\downarrow0}
\|f^{(\varepsilon)}-f \|_{E_{p,q,\beta}(C) }=0.
\end{equation}
\end{lemma}

\begin{lemma}
                \label{lemma 2.16.3}
Let $p,q\in(1,\infty)$,  $g(t,x)\geq 0$ be a Borel function such that for any smooth bounded $f(t,x)$ we have
(bounded linear functional  of special type)
\begin{equation}
                       \label{2.16.4}
\int_{\bR^{d+1}}g|f|\,dxdt\leq
\|f\|_{E_{p,q,\beta}}. 
\end{equation}
Then, for any $f\in E_{p,q,\beta}$,
\eqref{2.16.4} holds and, moreover,
$$
\lim_{\varepsilon\downarrow 0}
\int_{\bR^{d+1}}
g|f-f^{(\varepsilon)}|\,dxdt=0.
$$
\end{lemma}
 
The following is similar to Lemma
\ref{lemma 6,11.1}.

\begin{lemma}
                \label{lemma 2.16.5}
If $p,q\in(1,\infty)$, $\beta>1$, 
and
$(p_{0},q_{0})=(p,q)\beta=(r,s)(\beta-1) $, then for any
$f,g$
$$
\|fg\|_{\dot E_{p,q,\beta}}
\leq \|f\|_{\dot E_{p_{0},q_{0},1}}\|g\|_{\dot E_{r,s,\beta-1}}.
$$

\end{lemma}

Next comes a parabolic analog of Theorem
\ref{theorem 9.29.1}, which is nontrivial
only if $\beta\leq d/p+2/q $.
\begin{theorem}
                       \label{theorem 2.11.2}
Let    
$q_{1},q_{2}\in(1,\infty]$, $k>0$,
$0<\alpha<\beta $. Then 
there is a constant $N$ such that for any $f\geq0$
we have
\begin{equation}
                          \label{10.7.4}
 \|P_{\alpha,k}f\|_{\dot E_{r_{1},r_{2},\beta-\alpha}}
\leq N \|f\|_{\dot E_{q_{1},q_{2},\beta}},
\end{equation}
where $r_{i}(\beta-\alpha)=q_{i}\beta$,
$i=1,2$.
\end{theorem}

Proof. It suffices to prove that for any $\rho>0$
$$
\rho^{\beta-\alpha}\dashnorm P_{\alpha,k}f
\|_{L_{r_{1},r_{2}}(C_{\rho})}\leq N\|f\|_{\dot E_{q_{1},q_{2},\beta}},
$$
that is  
\begin{equation}
                          \label{2.11.3}
\rho^{\beta-\alpha-(d/r_{1}+2/r_{2})}\|P_{\alpha,k}f
\|_{L_{r_{1},r_{2}}(C_{\rho})}
\leq N\|f\|_{\dot E_{q_{1},q_{2},\beta}}.
\end{equation}

Observe that by H\"older's inequality
$\bM_{\beta}f\leq N\|f\|_{\dot E_{q_{1},q_{2},\beta}}$ and by definition
$$
\|I_{C_{2\rho}}f\|_{L_{q_{1},q_{2}}}
\leq N\rho^{d/q_{1}+2/q_{2}-\beta} \|f\|_{\dot E_{q_{1},q_{2},\beta}} 
$$
$$
=N\rho^{(d/r_{1}+2/r_{2})\beta/(\beta-\alpha)-\beta} \|f\|_{\dot E_{q_{1},q_{2},\beta}}.
$$
 
It follows from  Lemma \ref{lemma 10.7.1} 
with $p=\infty$ that \eqref{2.11.3} holds with
$I_{C_{2\rho}}f$ in place of $f$ on the left.
 
Furthermore,  by Lemma \ref{lemma 19.30.1} we have 
 $|P_{\alpha,k}(I_{C_{2\rho}^{c}}f)|
\leq N\rho^{\alpha-\beta}\bM_{\beta}f$ in $C_{\rho}$.
Therefore,
$$
\rho^{\beta-\alpha}\dashnorm P_{\alpha,k}
(I_{C_{2\rho}^{c}}f)\|_{L_{r_{1},r_{2}}(C_{\rho})}
\leq N\sup \bM_{\beta}f
\leq N\|f\|_{\dot E_{q_{1},q_{2},\beta}}.
$$
By combining these results we come to 
\eqref{2.11.3} and the theorem is proved.
\qed
 
As usual Theorem \ref{theorem 2.11.2}
leads to the following 
\index{embedding theorem}%
embedding result.

\begin{corollary}
                      \label{corollary 10.8.1}
Under the assumptions of Theorem \ref{theorem 2.11.2}, if $\beta>1$,  for any $u\in C^{\infty}_{0}$  we have
$$
\|Du\|_{\dot E_{r_{1},r_{2},\beta-1}}
\leq N\|\partial_{t}u+\Delta u\|_{\dot E_{q_{1},q_{2},\beta }},
$$
where $r_{i}(\beta-1)=q_{i}\beta$,
$i=1,2$. This coincides with \eqref{10.8.4} if
$\beta=d/q_{1}+2/q_{2}$.
\end{corollary}

\begin{remark}
                      \label{remark 10.8.2}
Corollary \ref{corollary 10.8.1} opens up the
possibility to treat the terms like
$b^{i}D_{i}u$ as perturbation term
in operators like 
$$
\partial_{t}u+\Delta u
+b^{i}D_{i}u
$$
 with   low 
summability properties of $b=(b^{i})$.
To show this observe that  for  
$q_{1},q_{2},\beta$ as in Theorem \ref{theorem 2.11.2} in the nontrivial case with $\beta>1$ and
$s_{i}=\beta q_{i}\in(1,\infty]$, $i=1,2$, we have
$$
\rho^{\beta}\dashnorm I_{C_{\rho}}b^{i}D_{i}u\|_{L_{q_{1},q_{2}}}
\leq \rho\dashnorm b I_{C_{\rho}}\|_{L_{s_{1},s_{2}}}\cdot \rho^{\beta-1}
\dashnorm I_{C_{\rho}} Du \|_{L_{r_{1},r_{2}}}
$$
implying that
$$
\| b^{i}D_{i}u\|_{\dot E_{q_{1},q_{2},\beta}}
\leq \|b  \|_{\dot E_{s_{1},s_{2},1}} \|   Du  \|_{\dot E_{r_{1},r_{2},\beta-1}}
$$
\begin{equation}
                         \label{2.29.1}
\leq N\|b  \|_{\dot E_{s_{1},s_{2},1}}
\|\partial_{t}u+\Delta u\|_{\dot E_{q_{1},q_{2},\beta }},
\end{equation}
where  
$d/s_{1}+2/s_{2}\geq 1$ if we want the norm of $b\not\equiv0$
to be finite.

Under the famous
\index{Ladyzhenskaya-Prodi-Serrin condition}
 Ladyzhenskaya-Prodi-Serrin condition
$$
b\in L_{s_{1},s_{2}},\quad \frac{d}{s_{1}}+\frac{2}{s_{2}}=1
$$
we have $\|b  \|_{\dot E_{s_{1},s_{2},1}}
=N\|b\|_{L_{s_{1},s_{2}}}$. It is worth noting that, if $d/s_{1}+2/s_{2}>1$ and
we look for a condition on $\tau$ and
$b\in L_{\tau s_{1},\tau s_{2}}$
to guarantee that $\|b  \|_{\dot E_{s_{1},s_{2},1}}<\infty$, then $\tau$ should satisfy 
$\tau=d/s_{1}+2/s_{2}$, so that
$d/(\tau s_{1})+2/(\tau s_{2})=1$,
thus $b$ should be in the Ladyzhenskaya-Prodi-Serrin class. Still
functions in $\dot E_{s_{1},s_{2},1}$ may have much higher singularities than those in
$L_{\tau s_{1},\tau s_{2}}$.
For instance, $|x|^{-1}$ is in any 
$\dot E_{s_{1},s_{2},1}$ for $s_{1}<d$
and $\tau s_{1}>d$.

  Another advantage of 
\eqref{2.29.1} in comparison with \eqref{2.29.2} is seen when $b$ depends only
on $t$ or $|b(t,x)|\leq \hat b(t)$. In that case \eqref{2.29.1}
becomes 
$$
\| b^{i}D_{i}u\|_{\dot E_{q_{1},q_{2},\beta}}
\leq N\|\hat b  \|_{\dot E_{ \beta q_{2},1/2}(\bR)}
\|\partial_{t}u+\Delta u\|_{\dot E_{q_{1},q_{2},\beta }},
$$
and if $\beta q_{2}=2$, then
$$
\|\hat b  \|_{\dot E_{ \beta q_{2},1/2}(\bR)}=\|\hat b\|_{L_{2}(\bR)}.
$$
Thus for any $q_{1}\in(1,\infty]$ and  $q_{2}\in(1,2)$
$$
\| b^{i}D_{i}u\|_{\dot E_{q_{1},q_{2},2/q_{2}}}
\leq N \|\hat b\|_{L_{2}(\bR)}
\|\partial_{t}u+\Delta u\|_{\dot E_{q_{1},q_{2},2/q_{2}}}.
$$
In case  $q_{1}\in(1,d)$, $q_{2}\in(1,\infty]$, $1<\beta\leq d/q_{1}$, and 
$$
\|b(\cdot,t)\|_{\dot E_{\beta q_{1},1}(\bR^{d})}\leq \hat b<\infty
$$
 for any $t$,
we also have
$$
\| b^{i}D_{i}u\|_{\dot E_{q_{1},q_{2},\beta}}
\leq N\hat b
\|\partial_{t}u+\Delta u\|_{\dot E_{q_{1},q_{2},\beta }}.
$$
An application of the last inequality
in case $u,b$ are independent of $t$,
  $\beta=d/q_{1}$, $q_{1}\in(1,d)$, and $q_{2}=\infty$, yields the well-known estimate
$$
 \| b^{i}D_{i}u \|_{L_{q_{1}}(\bR^{d})}\leq N \| b \|_{L_{d}(\bR^{d})}
 \|  \Delta u \|_{L_{q_{1}}(\bR^{d})}.
$$
 
\end{remark}

We also have an interpolation result.

\begin{theorem}
              \label{theorem 2.11.4}
Let   $q_{1},q_{2}\in[1,\infty)$, $\beta>0$. Then 
there is a constant
$N $ such that, for any
   $\varepsilon
>0 $   
and $u\in C^{\infty}_{0}$,
 \begin{equation}
                              \label{10.10.4}
\|Du\|_{\dot E_{q_{1},q_{2},\beta} }\leq
N\varepsilon   \| \partial_{t}u,D^{2}u\,\|_{\dot E_{q_{1},q_{2},\beta} }
+N\varepsilon^{-1} 
\| u\|_{\dot E_{q_{1},q_{2},\beta} }.
\end{equation}
\end{theorem}

This theorem is a particular case
of the following general result. 
Recall that $R_{\lambda}f$ is introduced
in \eqref{4.1,3} and according to
\eqref{1.25.2} 
$$
|DR_{\lambda}f(t,x)|\leq N
\int_{0}^{\infty}\int_{\bR^{d}}e^{-\lambda s}
\frac{|f(s+t,x+y)|}{s^{(d+1)/2}}
e^{-|y|^{2} /(8s)}\,dyds,
$$
which estimates  $|DR_{\lambda}f|$ in term of the ``sum'' with respect to $(s,y)$
with weights depending only on $(s,y)$
of functions $|f(s+t,x+y)|$ as functions of $(t,x)$. Minkowski's inequality that the norm
of a sum is less than the sum of norms
now leads to the first assertion in the following. The second assertion follows
from \eqref{4.1,3}.

\begin{lemma}
                      \label{lemma 3.31.4}
 Let $s\in[-\infty,\infty)$ and let $\frB$ be a Banach space of
functions on $\bR^{d+1}_{s}$ such that if $f\in \frB$
and $(s,y)\in\bR_{0}^{d+1}$, then $f_{s,y}\in\frB$ and
$\|f_{s,y}\|_{\frB}\leq \|f\|_{\frB}$, where
 $f_{s,y}(t,x)=f(s+t,x+y)$. Also assume that  
 if $g\in\frB$ and a measurable $f$ satisfies $|f|\leq g$, then $f\in\frB$ and  $\|f\|_{\frB}\leq \|g\|_{\frB}$. Then for any $\lambda>0$
and $f\in\frB$ we have 
$$
 \lambda^{1/2}\|DR_{\lambda}f\|_{\frB}\leq N(d)\|
f\|_{\frB}.
$$
Furthermore, if $C^{\infty}_{0}
\subset \frB$ then for any
$u\in C^{\infty}_{0}$ and $\lambda>0$ 
$$
\|Du\|_{\frB}\leq N\lambda^{-1/2}
\|\partial_{t}u+\Delta u-\lambda u\|_{\frB}
\leq N\big(\lambda^{-1/2}
\|\partial_{t}u,D^{2} u \|_{\frB}
+\lambda^{1/2}
\| u\|_{\frB}\big), 
$$
where $N=N(d)$.

\end{lemma}

\begin{remark}
\label{remark 2.16.1}
Another way to prove Theorem \ref{theorem 2.11.4} 
 is to follow the proof
of Theorem  \ref{theorem 12.26.1}.
By the way, 
if $u$ is independent of $t$, \eqref{10.10.4}
provides an  interpolation inequality in the ``elliptic"
case.
 \end{remark}

\begin{theorem}[embedding theorem]
                  \label{theorem 2.16.6}
(i) Let $ 0< \beta <2 $. 
\index{embedding theorem}%  
Then any $u\in \dot E^{2,1}_{p,q,\beta}$
is bounded and continuous  and for any $\varepsilon>0$ 
\begin{equation}
                    \label{2.16.7}
|u|
\leq \varepsilon^{2-\beta} \|\partial_{t}u,D^{2}u\|_{\dot E_{p,q,\beta}}+N(d, \beta)\varepsilon^{-\beta} \| u\|_{\dot E _{p,q,\beta}}. 
\end{equation} 
 (ii)
Let $1<\beta<2$. Then for any
$u\in \dot E^{2,1}_{p,q,\beta}$, $\rho>0$,
$(t_{i},x_{i})\in C_{\rho}$, $i=1,2$, we have
\begin{equation}
                           \label{2.16.8}
|u(t_{1},x_{1})-u(t_{2},x_{2})|\leq
N(d,p,q,\beta)\rho^{2-\beta}\|\partial_{t}u,D^{2}u\|_{\dot E _{p,q,\beta}}. 
\end{equation}
 
\end{theorem}

Proof.   Take $\zeta\in C^{\infty}_{0}$
such that $\zeta=1$ in $C_{1/2}$, $\zeta(t,x)=0$
in $\{t\geq0\}\cap C_{1}^{c}$, and $\zeta$ depends only
on $t$ and $|x|$.  For $\varepsilon
\in(0,1]$ define $\zeta_{\varepsilon}(t,x)=\zeta(t/\varepsilon^{2},x/\varepsilon)$. Introduce
  $$
f=-\partial_{t}u- \Delta u. 
$$
 Then
   $u(0)=u^{\varepsilon}_{1} 
+u_{2}^{\varepsilon} $, where 
with a constant $\gamma=
\gamma(d)$ and $\bR^{d+1}_{t}=(t,\infty)\times\bR^{d}$
$$
\gamma^{-1}u^{\varepsilon}_{1}=
\int_{\bR^{d+1}_{0} }\zeta_{\varepsilon}(t,x)t^{-d/2}
e^{- |x|^{2}/(4t)}f(t ,x )\,dxdt,
$$
$$
\gamma^{-1}u_{2}^{\varepsilon} =
\int_{\bR^{d+1}_{0} }(1-\zeta_{\varepsilon}(t,x))t^{-d/2}
e^{-|x|^{2}/(4t)}f(t ,x )\,dyds
$$
$$
=\int_{\bR^{d+1}_{0} }u(t ,x )\big(\partial_{t}- \Delta )\big[(1-\zeta_{\varepsilon}(t,x))t^{-d/2}
e^{- |x|^{2}/(4t)}\big]\,dxdt
$$
\begin{equation}
                      \label{3.31.2}
=  \int_{\bR^{d+1}_{0} }u(t ,x ) \eta _{\varepsilon}(t,x) t^{-d/2}
e^{-|x|^{2}/(4t)}\,dxdt,
\end{equation}
where
$$
\eta_{\varepsilon}(t,x)= -(\partial_{t}- \Delta ) 
\zeta_{\varepsilon}(t,x))-x^{i}D_{i}\zeta_{\varepsilon}(t,x) t^{-1}. 
$$
Note that $\eta_{\varepsilon}=0$ outside the part of $C_{\varepsilon}$ lying in $\bR^{d+1}_{0}$, $D\zeta(t,0)=0$, 
and 
$$
|x^{i}D_{i}\zeta_{\varepsilon}(t,x)| t^{-1}
\leq N\varepsilon^{-2}|x|^{2}/t,\quad e^{-|x|^{2}/(4t)}|x|^{2}/t
\leq Ne^{-|x|^{2}/(8t)},
$$
$$
|(\partial_{t}- \Delta ) 
\zeta_{\varepsilon}(t,x))|\leq \varepsilon^{-2}
\sup|(\partial_{t}- \Delta ) 
\zeta|.
$$

Therefore,
by   Lemma \ref{lemma 1.29.1}
 with  $\alpha=2,\gamma=\beta$, $\rho=\varepsilon$ we get
that
$$
\gamma^{-1}|u(0)|\leq P_{2,8}( |f\zeta_{\varepsilon}|I_{C_{\varepsilon}})
+N P_{2,8}( |u\eta _{\varepsilon}|I_{C_{\varepsilon}}) 
$$
$$
\leq N\varepsilon^{2-\beta}\bM_{\beta}f(0)
+N\varepsilon^{-\beta}\bM_{\beta}u(0). 
$$
Now \eqref{2.16.7} at the origin follows by H\"older's
inequality.

The proof of \eqref{2.16.8} is obtained
from the proof of Theorem \ref{theorem 2.3.1} by replacing there $U$ with 
$\|\partial_{t}u,D^{2}u\|_{\dot E _{p,q,\beta}}$.
  \qed

Observe that, if we want to have functions from
$W^{2,1}_{p,q}$ to be H\"older continuous
(as in \eqref{4.6.1}), we need $d/p+2/q<2$.
There are no restrictions on $d/p+2/q$ from above
in Theorem  \ref{theorem 2.16.6}.

\begin{remark}
                  \label{remark 2.16.2}
Assertion (ii) of Theorem \ref{theorem 2.16.6} holds also if $0<\beta\leq1$
and $1<d/p+2/q$
with only  few  corrections: $\rho\leq 1$
and $2-\beta$ in \eqref{2.16.8}
should be replaced with $2-\gamma$ for any 
   $\gamma\in (1,d/p+2/q]$,  $\dot E^{2,1}_{p,q,\beta}$
   shoud be replaced wih $E^{2,1}_{p,q,\beta}$ 
   and the norm in \eqref{2.16.8} should be replaced with
   the $E^{2,1}_{p,q,\beta}$-norm of $u$. This follows from the
fact that, for any $\zeta\in C^{\infty}_{0}$, $u\zeta\in \dot E_{p,q,\gamma}$ for any such $\gamma$
(we need $\gamma\leq d/p+2/q$ in order to have the new right-hand side of 
\eqref{2.16.8} finite).

\end{remark}

\begin{remark}
                      \label{remark 5.11,5}
Up to this point in this
subsection all results with identical proofs
hold for the $\L$-spaces in place of the $L$-spaces with only one changes that   we  should interchange
1 and 2, $p$ and $q$ in $\dot E^{1,2}_{p,q,\beta}$ and in similar situations and
in the expressions like $d/p+2/q$ use $p$
as the power of summability with respect to $x$ and $q$ as that to~$t$.

\end{remark}

Here is a theorem about 
\index{$t$-traces
of $\dot E^{1,2}_{q,p,\beta}$-functions}%
$t$-traces
of $\dot E^{2,1}_{p,q,\beta}$-functions.

\begin{theorem}
               \label{theorem 2.18.2}
Let   $p,q\in(1,\infty)$,
$r\geq p$  and let
$$
1< \beta  \leq  \frac{d}{p}+\frac{2}{q}
<1+\frac{d}{r}.
$$ 
Then for any
$u\in \dot E^{2,1}_{p,q,\beta}$ its trace $u(0,\cdot)$ is uniquely defined
and  
\begin{equation}
                     \label{2.18.3}
\|D u(0,\cdot)\|_{\dot E_{r,\beta -1} }\leq N \|\partial_{t}u,
D^{2}u\|
_{ \dot E _{p,q,\beta}}, 
\end{equation} 
where the constant  $N$ depends only
on $d,p,q,r,\beta$.

\end{theorem}

Proof. Take $\zeta\in C^{\infty}_{0} (\bR^{d+1})$, such that $\zeta(0)=1\geq \zeta\geq0$, define $\zeta_{n}(t,x )=
\zeta(t/n^{2},x/n )$ and observe that,
as $n\to\infty$,
$$
 \big|\|
\partial_{t}(\zeta_{n}u)\|_{\dot E_{p,q,\beta}} -\|\zeta_{n}
\partial_{t}u\|_{\dot E_{p,q,\beta} }\big|\leq n^{-2}\sup |\partial_{t}\zeta|
\|
u\|_{\dot E_{p,q,\beta}} \to 0.
$$
Also
$$
 \big|\|
D(\zeta_{n}u)\|_{\dot E_{p,q,\beta}} -\|\zeta_{n}
Du\|_{\dot E_{p,q,\beta} }\big|\leq n^{-1}\sup |D\zeta|
\|
u\|_{\dot E_{p,q,\beta}} \to 0, 
$$
$$
 \big|\|
D^{2}(\zeta_{n}u)\|_{\dot E_{p,q,\beta}} -\|\zeta_{n}
D^{2}u\|_{\dot E_{p,q\beta} }\big|\leq n^{-2}\sup |D^{2}\zeta|
\|
u\|_{\dot E_{p,q,\beta}} 
$$
$$
+2n^{-1}\sup |D\zeta|\|
Du\|_{\dot E_{p,q,\beta} } \to 0.
$$

It follows that it suffices to concentrate on $u$ that  
vanishes
for large $|t|+|x|^{2}$. In that case set 
$$
-f=\partial_{t}u+ \Delta u.
$$
To further reduce our
problem observe that using translations  
shows that it suffices
to prove that for any
  $\rho>0$,
$$
\rho^{\beta+\gamma-2}\dashnorm D^{\gamma}u(0,\cdot)\|_{L_{r}( B_{\rho})}\leq N\sup_{\rho_{1}\geq \rho}
\rho^{\beta }_{1}\dashnorm f\|
_{L_{p,q} (C_{\rho_{1}}    )}
$$
\begin{equation}
                     \label{2.18.4}
=  N\sup_{\rho_{1}\in[\rho,\rho+\rho_{2}]}
\rho^{\beta }_{1}\dashnorm f\|
_{L_{p,q} (C_{\rho_{1}})   },
\end{equation}
where $\rho_{2}$ is such that
$u(t,x)=0$ for $|t|+|x|^{2}\geq \rho_{2}^{2}$  
and the last equality is due to
$\beta \leq d/p+2/q$.

It is easy to pass to the limit
in \eqref{2.18.4} from smooth functions to arbitrary ones in
$W^{2,1}_{p,q}(C_{\rho_{2}})\supset
E^{2,1}_{p,q,\beta}(C_{\rho_{2}})$.
Therefore, we may assume that $u$ is smooth.
We thus reduced the general case to the task of proving the first estimate
in  \eqref{2.18.4} for smooth $u$ with compact support. One more reduction
is achieved by using the self-similarity
which shows that we only need to concentrate on $\rho=1$, that is
we only need to prove
\begin{equation}
                    \label{2.18.5}
 \dashnorm D u(0,\cdot)\|_{L_{r}( B_{1})}\leq N\sup_{\rho \geq 1}
\rho^{\beta} \dashnorm f\|
_{L_{p,q} (C_{\rho })   )}
\end{equation}
for smooth $u$ with compact  support.

  Now   define
  $
g=|f|I_{C_{2}},h=|f|I_{C_{2}^{c}}$.
As we know,
$$
|D u(0,x)|\leq NG (x)+NH (x),
\quad
 (G ,H )(x) =P_{1,8}(g,h)(0,x).
$$

By Lemma \ref{lemma 19.30.1} for $|x|\leq 1$
$$
H (x)\leq N\sup_{\rho>1}\rho^{\beta}
\dashint_{(0,x)+C_{\rho}}h\,dydt
\leq N\sup_{\rho>1}\rho^{\beta}
\dashint_{ C_{2\rho}}h\,dydt
$$
$$
\leq N\sup_{\rho \geq1}
\rho^{\beta} \dashnorm f\|
_{L_{p,q}  (C_{\rho} )},
$$
where the last inequality is
due to H\"older's inequality.
Hence,
\begin{equation}
                     \label{2.18.6}
 \| H \|_{L_{r}(  B_{1})}\leq N\sup_{\rho \geq1}
\rho^{\beta} \dashnorm f\|
_{L_{p,q}  (C_{\rho} )}. 
\end{equation}

By Minkowski's inequality
$$
\|G \|_{L_{r}(B_{1})}
\leq \int_{0}^{\infty}\Big(
\int_{B_{1}}\Big(\int_{\bR^{d}}
p_{1,8}(t,x-y)g(t,y)\,dy\Big)^{r}\,dx\Big)^{1/r}
\,dt,
$$ 
where inside the integral with respect
to $t$ we have the norm of a convolution, so that by Young's
inequality this expression is dominated by
$$
\|g(t,\cdot)\|_{L_{p}}\|p_{1,8}(t,\cdot)
\|_{L_{s}},
$$
where $1/s=1+1/r-1/p$ ($ \leq 1$ since $r\geq p$). We know that $\|p_{1,8}(t,\cdot)
\|_{L_{s}}=N(d)t^{\alpha}$ with
$\alpha=-1/2+(d/2)(1/r-1/p)$, which after taking into account that 
$\|g(t,\cdot)\|_{L_{p}}=0$ for $t\geq4$, yields
$$
\|G\|_{L_{r}(B_{1})}
\leq N\int_{0}^{4}\|f(t,\cdot)\|_{L_{p}(B_{2})}t^{\alpha}\,dt.
$$
Now use H\"older's inequality along with the observation that
$\alpha q/(q-1)>-1$, due to 
the assumption that   $1+d/r>
d/p+2/q$, to conclude that
$$
\|G\|_{L_{r}(B_{1})}
\leq N \|f \|_{L_{p,q}(C_{2})}.
$$
This and \eqref{2.18.6} prove
\eqref{2.18.5} and the theorem. \qed

\begin{corollary}
               \label{corollary 2.18.1}
If $q>2$, then
for any
$u\in \dot E^{2,1}_{p,q,\beta}$  and any $\varepsilon>0$  (there is no dot
on the left)
\begin{equation}
                     \label{2.18.7}
\|D u(0,\cdot)\|_{  E_{p,\beta } }\leq N\varepsilon \|\partial_{t}u,
D^{2}u\|
_{\dot E _{p,q,\beta}}
+N \varepsilon^{-(q+2)/(q-2)} 
\|u\|
_{ \dot E _{p,q,\beta}}. 
\end{equation} 
\end{corollary}
  
Indeed, if $\rho\leq \varepsilon$, then by \eqref{2.18.3} for $r=p$ (and $q>2$)
$$
\rho^{\beta}\dashnorm D(0,\cdot)\|_{L_{p}(B_{\rho})}\leq \varepsilon\rho^{\beta-1}\dashnorm D(0,\cdot)\|_{L_{p}(B_{\rho})}
\leq N \varepsilon\|\partial_{t}u,
D^{2}u\|
_{ \dot E _{q,p,\beta}}. 
$$
If $\varepsilon<\rho\leq1$, in light
of Corollary \ref{corollary 6,19.1} and of $(q+2)/(q-2)\geq1$ , 
$$ 
\rho^{\beta}\dashnorm D  u(0,\cdot)\|_{L_{r}(B_{\rho})}
\leq N\varepsilon \|\partial_{t}u,D^{2}u
\|_{\dot E_{q,p,\beta} }
+N \varepsilon^{-(q+2)/(q-2)}\| u
\|_{\dot E_{q,p,\beta} }. 
$$

\subsection
{Application to the heat
equation}

The following is the basic a priori estimate.

\begin{lemma}
              \label{lemma 2.15.3}
Let   $p,q\in(1,\infty)$, $\beta>0$. Then there is a constant
$N $ such that, for any  $u\in C^{\infty}_{0}$,
 \begin{equation}
                   \label{2.15.4}
\|D^{2}u\|_{\dot E_{p,q,\beta} }\leq
N  \| \partial_{t}u+\Delta u\,\|_{\dot E_{p,q,\beta} }.
\end{equation}
\end{lemma}

Proof. It suffices to prove that
\begin{equation}
                     \label{2.19.1}
\|D^{2}u\|_{L_{p,q}(C_{1})}\leq 
N  \| f\|_{\dot E_{p,q,\beta} },
\end{equation}
where 
$$
f=-\partial_{t}u-\Delta u.
$$
Take a $\zeta\in C^{\infty}_{0}$
such that $\zeta=1$ on $(-1,1)\times B_{1}$
and $\zeta=0$ outside $(-4,4)\times B_{2}$.
Let $g=\zeta f$, $h=(1-\zeta)f$,
$ (G,H)=P_{2,4}(g,h)$. We know that
$u=N(F+G)$. By Remark \ref{remark 2.25.1}
$$
\|D^{2}G\|_{L_{p,q}(C_{1})}\leq 
N  \| f\|_{\dot E_{p,q,\beta} }  .
$$
By Lemma \ref{lemma 19.30.1}
and H\"older's inequality
$$
\|D^{2}H\|_{L_{p,q}(C_{1})}
\leq N\bM_{\beta}f(0)
  \leq 
N  \| f\|_{\dot E_{p,q,\beta} } . 
$$
This proves \eqref{2.19.1} and the lemma. \qed

By combining this lemma with Lemma \ref{lemma 3.31.4}
 and using Agmon's idea,
we come to the a priori estimate
\eqref{2.25.3} if $u\in C^{\infty}_{0}
(\bR^{d+1})$. 

\begin{theorem}
               \label{theorem 2.25.1}
Let   $p,q\in(1,\infty)$, $\beta>0$. Then there is a constant
$N $ such that, for any $\lambda\geq0$ and  $u\in \dot E^{2,1}_{p,q}$,
 \begin{equation}
                   \label{2.25.3}
\|D^{2}u,\partial_{t}u,\sqrt\lambda Du,
\lambda u\|_{\dot E_{p,q,\beta} }\leq
N  \| \partial_{t}u+\Delta u-\lambda u
\|_{\dot E_{p,q,\beta} }.
\end{equation}
Furthermore, for any $\lambda>0$
and $f\in \dot E_{p,q,\beta}(\bR^{d+1})$
there exists a unique $u\in 
\dot E^{2,1}_{p,q,\beta}$ satisfying
$$
\partial_{t}u+\Delta u-\lambda u=f.
$$

\end{theorem}

The extension of \eqref{2.25.3} to
$u\in \dot E^{2,1}_{p,q}$ and the existence part of the theorem is done
as in the case of Theorem \ref{theorem 1.31.3}.

After that we use 
Remark \ref{remark 10.8.2} 
to add the term $b^{i}D_{i}u$ and in an already familiar way we come to the following, that generalizes Theorem
\ref{theorem 2.25.1} if $\beta>1$.

\begin{theorem}
               \label{theorem 2.25.2}
Let   $p,q\in(1,\infty)$, $\beta>1$ 
and let $b$ be a $\bR^{d}$-valued function. Then 
there exists a constant $N=N(d,p,q,\beta)$
such that, if 
$$
N\|b\|_{\dot E_{\beta p,\beta q,1}(\bR^{d+1})}\leq 1,
$$
 then
 for any $\lambda>0$  and
$f\in \dot E_{p,q,\beta}(\bR^{d+1})$ there exists
a unique $u\in \dot E^{2,1}_{p,q,\beta}$ satisfying $\cL u-\lambda u=f$,
where 
$$
\cL u=\partial_{t}u+\Delta u+b^{i}D_{i}u.
$$
Furthermore, there exists $N_{0}=N(d,p,q,\beta)$ such that for any $\lambda\geq0$ and $u\in
\dot E^{2,1}_{p,q,\beta}$ we have
 \begin{equation}
                   \label{2.25.4}
\|D^{2}u,\partial_{t}u,\sqrt\lambda Du,
\lambda u\|_{\dot E_{p,q,\beta} }\leq
N  \| \cL u-\lambda u
\|_{\dot E_{p,q,\beta} }.
\end{equation}
 
\end{theorem}
 
More information about parabolic equations
in Morrey spaces with mixed norms
can be found in \cite{Kr_25} and in
Section \ref{chapter 4.5,1}.

\mysection{Parabolic equations with variable
main part}

                       \label{chapter 4.5,1}

In the previous sections we have given
applications of the results from Real Analysis
to only model equations with the main part
being the Laplacian. The only exception
was short Subsection \ref{section 4.5,1}
where the main part was given by 
$\bS_{\delta}$-valued $a$ of
BMO class. In this section we concentrate
only on applications of Real Analysis results to parabolic operators with $a(t,x)$ which is measurable in $t$
and BMO in $x$ and drift term in
the mixed-norm local Morrey class. Here we consider mixed-norm
spaces when the inside integration is performed
with respect to $t$ and not to $x$ as is customary. This is mainly caused by an application, the author has in mind,
to the strong uniqueness problem for
stochastic It\^o equations and also, of course,
has an independent interest in its own,
especially that in the literature one can find 
claims that the theory in such spaces
is quite parallel to the one when the
interior integration is done with respect to $x$. In this respect
  the main  results presented here
appear for the first time.  However, it is worth
noting that equations in
the ``unusual'' mixed normed Sobolev (not Morrey-Sobolev)
spaces were earlier considered
in \cite{KN_14} (see also the references therein).
By the way, the reader will easily see that
all our results hold true also in the usual mixed-norm spaces, where they are known.
The methods and results presented here mimic those
from \cite{DK_18} and \cite{Kr_25}, where
the mixed norms are defined in the traditional way. 

In this section, as usual in the theory
of parabolic equations, in $\bR^{d+1}=\{t\in\bR,x\in\bR^{d}\}$
we introduce the real analytic structure
given by $(2,1,...,1)$ (parabolic boxes)
and Lebesgue measure $dz=dtdx$.

\subsection{Main coefficients depending  only on $t$}

Fix a $\delta\in(0,1]$ and let $a=a(t)$
be an $\bS_{\delta}$-valued function on $\bR$.

\begin{lemma}
                   \label{lemma 3.29.1}
Let $a \in \bS_{\delta}$. Then for any $d\times d $ 
symmetric matrix $u$
we have
\begin{equation}
                              \label{3.29.2}
\{a,u\}:=a^{ij}a^{kr}u_{ik}u_{jr}\geq \delta^{2}\sum_{i,j}u_{ij}^{2},
\end{equation}
\begin{equation}
                              \label{3.29.3}
 (1-\delta^{2})^{2} \{a,u\}\geq \{a-\delta(\delta^{ij}),u\} \geq 0. 
\end{equation}
\end{lemma}

Proof. If $\lambda_{p},\ell_{p}$ are the eigenvalues and eigenvectors of $a$, then
the left-hand side of \eqref{3.29.2} is written
as
$$
\sum_{p,q}\lambda_{p}\lambda_{q}\Big(\sum_{i,j}
u_{ij}\ell_{p}^{i}\ell_{q}^{j}\Big)^{2}
\geq\delta^{2}\sum_{p,q} \Big(\sum_{i,j}
u_{ij}\ell_{p}^{i}\ell_{q}^{j}\Big)^{2}
=\delta^{2}u_{ij}\ell_{p}^{i}\ell_{q}^{j}
u_{kr}\ell_{p}^{k}\ell_{q}^{r}
$$
and the latter equals the right-hand side of
\eqref{3.29.2} because $\ell^{i}_{p}\ell^{k}_{p}=\delta^{ik},\ell^{j}_{q}\ell^{r}_{q}=\delta
^{jr}$.  

Since $a\in\bS_{\delta}$, $(1-\delta^{2})\lambda_{p}\geq \lambda_{p}-\delta\geq0$.
This yields \eqref{3.29.3}. \qed

Let $u\in C^{\infty}_{0}$. 
Set 
$$
-f=\partial_{t}u+a^{ij}(t)D_{ij}u.
$$
 Then
$$
\int_{\bR^{d+1}}f^{2}\,dz=I_{1}+ I_{2}+I_{3},
$$
where integrating by parts 
and observing that
$$
\int_{\bR^{d+1}}\partial_{t}u \delta^{ij}D_{ij}u\,dz=-(1/2)\int_{\bR^{d+1}}\partial_{t}(|Du|^{2})\,dz=0
$$

we find
$$
I_{1}=\int_{\bR^{d+1}}|\partial_{t}u|^{2}\,dz,
$$
$$
I_{3}=\int_{\bR^{d+1}}a^{ij}D_{ij}ua^{kr}D_{kr}u\,dz=\int_{\bR^{d+1}}\{a,D^{2} u\}\,dz
$$
$$
I_{2}=2\int_{\bR^{d+1}}\partial_{t}u a^{ij}D_{ij}u\,dz=2\int_{\bR^{d+1}}\partial_{t}u [a^{ij}-\delta\delta^{ij}]D_{ij}u\,dz 
$$
$$
\geq-2I_{1}^{1/2}\Big(\int_{\bR^{d+1}}\{a
-\delta (\delta^{ij}),D^{2} u\}\,dz\Big)^{1/2}
\geq-2(1-\delta^{2})I_{1}^{1/2}I_{3}^{1/2}.
$$
It follows that
$$
I_{1}+I_{3}\leq\delta^{-2}\int_{\bR^{d+1}}f^{2}\,dz, 
$$
\begin{equation}
                          \label{3.29,3}
\int_{\bR^{d+1}}\big(|\partial_{t}u|^{2}+\delta^{2}|D^{2}u|^{2}\big)\,dz\leq
\delta^{-2}\int_{\bR^{d+1}}f^{2}\,dz. 
\end{equation}

Introduce $B^{0,\infty}$ as the space
\index{$A$@Sets of functions!$B^{0,\infty}$}%
of function $u(t,x)$ on $\bR^{d+1}$
which are infinitely differentiable in $x$
for every $t$ with each derivative
locally bounded on $\bR^{d+1}$. By $B^{1,\infty}$ we mean the 
\index{$A$@Sets of functions!$B^{1,\infty}$}%
\index{$A$@Sets of functions!$B^{1,\infty}_{0}$}%
subspace of $B^{0,\infty}$ consisting of functions
$u(t,x)$ such that $\partial_{t}u\in B^{0,\infty}$. By $B^{i,\infty}_{0}$
we mean the subspaces of $B^{i,\infty} $  consisting of functions
with compact support, $i=0,1$.

 Note that, like \eqref{1.25.1},
  estimate  \eqref{3.29,3}  is also 
  true by the same reasons as in Remark \ref{remark 4.4,3}
for
  $u\in B^{1,\infty}$  such that, for some $\alpha>d/2$ and $\rho(t,x)=\sqrt{|t|}+|x|,\xi 
=\rho^{ \alpha-1},\eta=\rho^{\alpha}$ the function 
\begin{equation}    
                             \label{3.30.2}
\xi|u|+\eta|Du| 
\end{equation}
 is bounded.
 
Next, define 
$$
A_{t,t+ s}=\int_{t}^{t+s }a(r)\,dr
=s\int_{0}^{1}a(t+rs)\,dr,
\quad
\sigma_{t,t+ s}= s ^{-1/2}A^{1/2}_{t,t+s} ,
$$
recall that for $c_{d}=(4\pi)^{-d/2}$ 
$$
\sfp(t,x)=c_{d}t^{-d/2}e^{-|x|^{2}/(4t)}
I_{t>0}
$$
and for $\lambda\geq0$ introduce
$$
R_{\lambda}f(t,x)=\int_{0}^{\infty}
\int_{\bR^{d}}e^{-\lambda s}f(t+s, x+
\sigma_{t,t+ s}y)\sfp(s,y)\,dyds
$$
$$
=e^{\lambda t}\int_{t}^{\infty}
\int_{\bR^{d}}e^{-\lambda s}f( s, x+
\sigma_{t,  s}y)\sfp(s-t,y)\,dyds 
$$
$$
=\int_{t}^{\infty}e^{-\lambda s}
\int_{\bR^{d}}f( s,y )\sfp(t, s,y-x)\,dyds,
$$
where 
\index{$S$@Miscelenea!$\sfp(t,s,x)$}%
$\sfp(t,s, x)=\sfp(s-t ,\sigma_{t, s}^{-1}x)\det \sigma_{t, s }^{-1}$.

Introduce
\index{$C$@Operators!$\cL^{0}$}%
$$
\cL^{ 0} u(t,x)=\partial_{t}u(t,x)+a^{ij}(t)D_{ij}u(t,x).
$$

\begin{theorem}
                     \label{theorem 3.30.1}
 Let $ f\in B^{0,\infty}_{0}$. Then

(i) we have $u:=R_{\lambda}f\in B^{1,\infty}$
and   
\begin{equation}
                             \label{3.30,3}
\cL^{ 0} u-\lambda u =-f   ;
\end{equation}
 
(ii) there is $\alpha>d/2$ such that the function \eqref{3.30.2} is bounded.
\end{theorem}

Proof. Obviously, $D^{n}u $ is bounded
on $\bR^{d+1}$ for any $n\geq0$ and
is continuous with respect to $\lambda\geq0$.
In addition, we can rewrite \eqref{3.30,3} by using
integration with respect to $t$ and, therefore while proving \eqref{3.30,3}, we may concentrate on $\lambda>0$. In that case
$$
\sup_{t}\int_{\bR^{d}}|D^{n}u(t,x)|\,dx\leq
\sup_{t}\int_{0}^{\infty}e^{-\lambda s}
\int_{\bR^{d}}|D^{n}f(t+s, x ) |\,dxds<\infty.
$$
This allows us to use the Fourier transform:
$$
\tilde u(t,\xi)=\int_{\bR^{d}}u(t,x)e^{i(\xi,x)}\,dx 
$$
and, in light of the Fubini theorem and the fact that
$\tilde \sfp(t,\xi)=e^{- |\xi|^{2}t}$,
find that  
$$
\tilde u(t,\xi)=\int_{0}^{\infty}
\int_{\bR^{d}}e^{-\lambda s}\tilde f(t+s,\xi)
e^{-i(\xi , \sigma_{t,s}y)} \sfp(s,y)\,dyds 
$$
$$
= \int_{0}^{\infty}e^{-\lambda s}\tilde f(t+s,\xi)\exp\big(- ( A _{t,t+s}\xi,\xi)\big) \, ds. 
$$
For $\tilde v(t,\xi)=e^{-\lambda t}\tilde u(t,\xi)$ this yields
\begin{equation}
                                \label{3.31.3}
\tilde v(t,\xi)=\int_{t}^{\infty}e^{-\lambda r}
\tilde f(r,\xi)\exp\big(- ( A _{t,r}\xi,\xi)\big) \, dr,
\end{equation}
$$
\partial_{t}\tilde v(t,\xi)=-e^{-\lambda t}\tilde f(t,\xi)
+(a(t)\xi,\xi)\int_{t}^{\infty}e^{-\lambda r}
\tilde f(r,\xi)\exp\big(- ( A _{t,r}\xi,\xi)\big) \, dr 
$$
$$
=-e^{-\lambda t}\tilde f(t,\xi)+(a(t)\xi,\xi)
\tilde v(t,\xi). 
$$
Coming back to the original $\tilde u(t,\xi)$
we get
$$
e^{-\lambda(t+s)}\tilde u(t+s,\xi)-e^{-\lambda t }\tilde u(t,\xi)
$$
\begin{equation}
                                \label{3.30.4}
=-\int_{0}^{s}e^{-\lambda r}\big[\tilde f(t+r,\xi)+
( a(t+r)\xi,\xi)\tilde u(t+r,\xi)\big]\,dr.
\end{equation}

We intend to apply the inverse Fourier
transform to \eqref{3.30.4}. Recall  that since $f\in B^{0,\infty}_{0}$ the function $u$ is infinitely differentiable in $x$ and each of its derivatives is in $L_{1}$. 
In particular, $|\tilde u(t,\xi)|$
decays as $|\xi|\to \infty$ faster that $|\xi|$
to any negative power. This and Fubini's
theorem allow  us to
insert $e^{-i(x,\xi)}$ and the integration
with respect to $\xi$ inside the integral in
\eqref{3.30.4}. After doing so we obtain
$$
e^{-\lambda(t+s)}  u(t+s,x)-e^{-\lambda t } u(t,x) 
$$
$$
=-\int_{0}^{s}e^{-\lambda r}\big[  f(t+r,x)-
 a^{ij}(t+r)D_{ij} u(t+r,x)\big]\,dr,
$$
which is equivalent to \eqref{3.30,3}.
Along with what is said in the beginning of the proof this implies that $D^{n}\partial_{t}u$
is bounded for any $n\geq0$ so that
$u\in B^{1,\infty}$. This proves (i).

Remark \ref{remark 5.31,1} shows that 
$|D^{n}u|\rho^{d}$ is bounded for any $n\geq0$
and that this fact implies the boundedness
of \eqref{3.30.2}. The theorem is proved.\qed

\begin{lemma}
                        \label{lemma 3.31.1}
Let $u\in C^{\infty}_{0}$, $\lambda
\geq0$. Then
$$
u=R_{\lambda}(\lambda u-\cL^{ 0} u).
$$
\end{lemma}

Proof. Introduce $f$ by \eqref{3.30,3} and define $w=R_{\lambda}f$. Then according to
\eqref{3.31.3}
$$
e^{-\lambda t}\tilde w(t,\xi)=I_{1}+I_{2},
$$
where by using integration by parts we have
$$
I_{1}=-\int_{t}^{\infty}e^{-\lambda r}
\partial\tilde w(r,\xi)\exp\big(- ( A _{t,r}\xi,\xi)\big) \, dr=e^{-\lambda t}\tilde u(t,\xi)-I_{2},
$$
$$
I_{2}=\int_{t}^{\infty}e^{-\lambda r}
(a(r)\xi,\xi)\tilde u(r,\xi)\exp\big(- ( A _{t,r}\xi,\xi)\big) \, dr. 
$$
Thus, $\tilde w(t,\xi)=\tilde u(t,\xi)$
and $w=u$. \qed

\begin{remark}
                   \label{remark 3.30.2}
Observe that $\sigma_{t, s}\geq\delta^{1/2}$, which similarly to \eqref{1.25.2} implies that
for $n\geq0$ and $s>t$
\begin{equation}
                            \label{3.31.1}
|D^{n}\sfp(t,s, x)|\leq N(d,n,\delta)
\frac{1}{(s-t)^{(d+n)/2}}e^{-|x|^{2}\delta/(8  (s-t))}.
\end{equation}
Furthermore, from
$\sfp(t,s, x)=\sfp(s-t ,\sigma_{t, s}^{-1}x)\det \sigma_{t, s }^{-1}$ we get that 
$$
\tilde p(t,s,\xi)=\exp(-(A_{t,s}\xi,\xi)),\quad
\partial_{t}\tilde p(t,s,\xi)
-(a(t)\xi,\xi)\tilde p(t,s,\xi)=0.
$$
It follows that
$$
\partial_{t} p(t,s,x)+a^{ij}(t)D_{ij}
p(t,s,x)=0 
$$
and \eqref{3.31.1} implies that for $n\geq0$ and $s>t$
\begin{equation}
                            \label{3.31,2}
|\partial_{t}D^{n}\sfp(t,s, x)|\leq N(d,n,\delta)
\frac{1}{(s-t)^{1+(d+n)/2}}e^{-|x|^{2}\delta/(8  (s-t))}.
\end{equation}

By the way, differentiating
$\tilde p(t,s,\xi)=\exp(-(A_{t,s}\xi,\xi))$
with respect to $s$ we see that
$$
\partial_{s} p(t,s,x)=a^{ij}(s)D_{ij}
p(t,s,x) . 
$$
\end{remark}

\begin{theorem}
                        \label{theorem 3.31.1}
Let $f\in B^{0,\infty}_{ 0}$ and $u=R_{0}f$. Then
for any $\kappa\geq 2$ and $\rho>0$ 
$$
\osc_{C_{\rho}}D^{2}u\leq N(d,\delta)\kappa^{(d+2)/2}
\dashnorm f\|_{L_{2}(C_{\kappa\rho})}+
N(d,\delta)\kappa^{-1} \bM  f  (0).
$$
In particular, if $u\in C^{\infty}_{0}$ then 
$$
\osc_{C_{\rho}}D^{2}u\leq N(d,\delta)\kappa^{(d+2)/2}
\dashnorm \cL^{ 0} u\|_{L_{2}(C_{\kappa\rho})}+
N(d,\delta)\kappa^{-1} \bM (\cL^{ 0} u)  (0). 
$$

\end{theorem}

Proof. We follow the proof of \eqref{1.25.4}
after observing that the case of arbitrary $\rho>0$ is reduced to $\rho=1$ by using
scale changes. In that case
take $\zeta \in C^{\infty}_{0}$
such that $\zeta=1$ in $C_{ \kappa}$, $\zeta=0$ 
outside $C_{2 \kappa}^{c}\cap \bR^{d+1}_{0}$, $0\leq\zeta\leq1$ and  
note that for   $g=f\zeta,h=f(1-\zeta)$
 and $(G,H)=R_{0}(g,h)$   we have
$u=G+H$, and for $m=0,1,n\geq0$    
 $$
 \partial_{t} ^{m}D^{n}
H(t,x)= \int_{C_{2}^{c}}
 \partial_{t} ^{m}D^{n}
\sfp(t, s,y-x)h( s, y)\,dsdy.   
$$

By using   \eqref{3.31.1} we get that in $C_{1}$  
$$
|\partial_{t}D^{2}H |+|D^{3}H |
\leq P_{-2,\delta/8 }|h|+P_{-1,\delta/8 }|h|, 
$$
which by Lemma \ref{lemma 19.30.1} 
(with $\beta=0$) implies that 
\begin{equation}
                         \label{1.26,1}
|\partial_{t}D^{2}H(t,x)|+|D^{3}H(t,x)|
\leq N\kappa^{-1}\bM h(0).
\end{equation}
It follows  that 
$$
\int_{C_{1}}\int_{C_{1}}
|D^{2}H(t,x)-D^{2}H(s,y)|\,dxdydtds\leq N
\kappa^{-1}\bM f(0). 
$$

Regarding $G$ we have in light of   \eqref{3.29,3} 
 and Theorem \ref{theorem 3.30.1}
that 
$$
\int_{C_{1}}\int_{C_{1}}
|D^{2}G(t,x)-D^{2}G(s,y)|\,dxdydtds\leq N\int_{C_{1}} 
|D^{2}G(t,x) |\,dxdt 
$$
$$
\leq N\|g\|_{L_{2}(\bR^{d+1})}\leq N
\kappa^{(d+2)/2}\dashnorm f\|_{L_{2}(C_{\kappa })}.
$$
This proves the first assertion of the theorem. 
The second one is a consequence of the first one and Lemma \ref{lemma 3.31.1}.  \qed

 Recall that  for  $p,q\in(1,\infty)$ we set
$\L_{q,p}$ to be the space of functions on $\bR^{d+1}$ with finite norm defined by
\begin{equation}
                        \label{4.29,5}
\|f\|^{p}_{\L_{q,p}}=\int_{\bR^{d}}\Big(
\int_{\bR}|f(t,x)|^{q}\,dt\Big)^{p/q}\,dx. 
\end{equation}
The space $\W^{1,2}_{q,p}$ is defined as 
\index{$A$@Sets of functions!$\W^{1,2}_{q,p}$}%
the collection of functions $u(t,x)$ on $\bR^{d+1}$
such that $u,Du,D^{2}u,\partial_{t} u
\in \L_{q,p}$. The space $\W^{1,2}_{q,p}$ is
 provided with a natural norm.  
 The folowing is part of Proposition 2 of \cite{KN_14},
 where the reference to the proof based 
 on singular integrals is given. 
\begin{theorem}
                         \label{theorem 3.31.2}
There is a constant $N=N(d,\delta,p,q)$
such that for any $u\in \W^{1,2}_{q,p}$
\begin{equation}
                               \label{4.1,5}
\|\partial_{t}u,D^{2}u\|_{\L_{q,p}}
\leq N\|\cL^{ 0} u\|_{\L_{q,p}}.
\end{equation}
\end{theorem}

Proof. By Lemma \ref{lemma 3.31.1} and Theorem \ref{theorem 3.31.1}
 with $\kappa=2$ for $u\in C^{\infty}_{0}
(\bR^{d+1})$ we have
\begin{equation}
                               \label{9.18.10}
(D^{2}u)^{\sharp}\leq N(d,\delta)(\bM((\cL^{ 0} u)^{2}))^{1/2}, 
\end{equation}
 which is similar to \eqref{1.25.4} and leads
to the following in the same way \eqref{1.25.4} leads to \eqref{1.25.5}.
 \begin{theorem}
                  \label{theorem 1.25.33}
 For $d\geq1$, any $p\in(1,\infty)$, and $u\in W^{1,2}_{p}$ we have
 \begin{equation}
                                          \label{9.18,0}
  \|\partial_{t}u\|_{L_{p}(\bR^{d+1})}+
  \|D^{2}u\|_{L_{p}(\bR^{d+1})}\leq N(p,d,\delta)
  \|\cL^{0} u\|_{L_{p}(\bR^{d+1})}.
 \end{equation}
 \end{theorem}

 In turn this theorem allows us to get that
 for $u\in C^{\infty}_{0}
(\bR^{d+1})$, $p>1$ we have
\begin{equation}
                                          \label{9.18,11}
(D^{2}u)^{\sharp}\leq N(p,d,\delta)(\bM((\cL^{ 0} u)^{p}))^{1/p} 
\end{equation}
repeating the derivation of \eqref{9.18.10} along the lines
in the derivation of \eqref{1.19.1}.
As in the proof of Lemma \ref{lemma 1.21.4} this yields that
for $p>1 $ and $w\in A_{p}$  
\begin{equation}
                                          \label{9.18,12}
\int_{\bR^{d+1}}|D^{2}u|^{p}w\,dxdt
\leq N(d,\delta,p,[w]_{p})\int_{\bR^{d+1}}|\cL^{0} u|^{p}w\,dxdt. 
\end{equation}

After that the estimate for $D^{2}u$ follows from
Theorem \ref{theorem 12.2.1}. 
The term $\partial_{t}u$ as usual
is estimated from $\partial_{t}u=\cL^{ 0 } u-a^{ij}D_{ij}u$. \qed

Estimate \eqref{4.1,5} for $p=q$ (and variable
$a$) was first obtained in \cite{Jo_64},
for $p\ne q$ it is proved in \cite{Kr_01}.
The technique used here is quite different
from both \cite{Jo_64} and \cite{Kr_01}
and is the same as in \cite{DK_18},
where mixed norm estimate are obtained
with the integration on $x$ and $t$ in any order.

Next, we introduce the term $\lambda u$
by using Agmon's method and adding Lemma
\ref{lemma 3.31.4} conclude that
the a priori estimate \eqref{4.1,1}
with $s=-\infty$
is true for $u\in C^{\infty}_{0}$.
For any $u\in \W^{1,2}_{q,p}$ it holds
because $C^{\infty}_{0}$ is dense
in $\W^{1,2}_{q,p}$.

\begin{theorem}
                         \label{theorem 4.2,1}
There is a constant $N=N(d,\delta,p,q)$
such that for any $s\in[-\infty,\infty)$,
 $u\in \W^{1,2}_{q,p}(\bR^{d+1}_{s})$ and
$\lambda\geq0$ 
\begin{equation}
                               \label{4.1,1}
\|\lambda u,\sqrt\lambda Du,\partial_{t}u,D^{2}u\|_{\L_{q,p}(\bR^{d+1}_{s})}
\leq N\|\cL^{ 0 } u-\lambda u\|_{\L_{q,p}(\bR^{d+1}_{s})}.  
\end{equation}
Furthermore, for any $\lambda>0$ and
$f\in \L_{q,p}(\bR^{d+1}_{s})$ there is a unique
$u\in \W^{1,2}_{q,p}(\bR^{d+1}_{s})$ such that
$\cL^{ 0 } u-\lambda u=f$ in $\bR^{d+1}_{s}$.  
\end{theorem}

The existence result for $s=-\infty$
is obtained by the method
of continuity starting from
Theorem \ref{theorem 2.14.1}, 
in which the order of integration
with respect to $t,x$ is different
but the proof for the needed version
is easily obtained from
its proof with only  using Lemma \ref{lemma 3.31.4} in place of 
Lemma \ref{lemma 2.14.1}.

The uniqueness of $\W^{1,2}_{q,p}(\bR^{d+1}_{s})$-solutions follows from the fact that it holds in Theorem \ref{theorem 2.14.1} and that the method
of continuity preserves causality.

Theorem \ref{theorem 1.25.33} was instrumental
in proving \eqref{9.18,11}. It also implies the following
generalization of Theorem \ref{theorem 3.31.1} obtained just by repeating the proof
of it having in mind Theorem \ref{theorem 1.25.3}.
\begin{theorem}
                        \label{theorem 9.18,00}
Let $p>1$, $f\in B^{0,\infty}_{ 0}$ and $u=R_{0}f$. Then
for any $\kappa\geq 2$ and $\rho>0$ 
$$
\osc_{C_{\rho}}D^{2}u\leq N(p,d,\delta)\kappa^{(d+2)/p}
\dashnorm f\|_{L_{p}(C_{\kappa\rho})}+
N(p,d,\delta)\kappa^{-1} \bM  f  (0).
$$
In particular, if $u\in C^{\infty}_{0}$ then 
$$
\osc_{C_{\rho}}D^{2}u\leq N(p,d,\delta)\kappa^{(d+2)/p}
\dashnorm \cL^{ 0} u\|_{L_{p}(C_{\kappa\rho})}+
N(p,d,\delta)\kappa^{-1} \bM (\cL^{ 0} u)  (0). 
$$

\end{theorem}

\subsection{The case of variable main coefficients and Sobolev spaces}
 \label{section 4.29.1}
 
Now we are going to extend Theorem
\ref{theorem 3.31.2} to operators with 
$\bS_{\delta}$-valued $a=a(t,x)$.
Define 
$$
a^{\shharp} =
 \sup
_{ C\in\bC }\dashint_{C}|a(t,x)-\tilde a_{C}(t)|\,dxdt,\quad \tilde a_{C}(t)=\dashint_{C}a(t,x)\,dxds
$$
(note $t$ and $ds$). 
\index{$S$@Miscelenea!$a^{\shharp}$}%
\index{$S$@Miscelenea!$a^{\shharp}$@$\tilde a_{C}$}%
Observe that if $a$
is independent of $x$, then $a^{\shharp}  =0$. 

Set
\index{$C$@Operators!$\cL_{0}$}%
$$
\cL _{ 0 } u=\partial_{t}u+a^{ij}D_{ij}u.
$$

For  $ C\in\bC$  and
$$
 \cL_{ C} u:= \partial_{t}u+
\tilde a^{ij}_{ C}D_{ij}u,\quad f=\cL_{ 0 } u
$$
  we have 
$$
f= \cL_{ C} u
+[(a^{ij}- \tilde a^{ij}_{C})D_{ij}u],
$$
 whereas
for $q> p>1 $
$$
\dashnorm (a^{ij}- \tilde a^{ij}_{ C})D_{ij}u
\|_{L_{p} C )}\leq
N(d,\delta, q,p)a^{\shharp} 
\big(\bM(|D^{2} u|^{q})\big)^{1/q}. 
$$
This and Theorem  \ref{theorem 9.18,00} 
 lead  to the following  (in which we keep $p$
 to facilitate its derivation from Theorem \ref{theorem 9.18,00}
 containing $p$).

\begin{lemma}  
                  \label{lemma 4.2,1}

For any $q>  p>1 $,
$\kappa\geq 2$, and $u\in C^{\infty}_{0}$   we have 
$$
( D^{2}u)^{\sharp}
\leq N(d,\delta,q,\kappa)\Big(\bM(|\cL_{ 0 } u|^{q})\big)^{1/q} 
$$
\begin{equation}
                         \label{4.2,1}
+\big(N(d,\delta,q,p,\kappa)a^{\shharp} 
+N(d,\delta, q,p  )\kappa^{-1}\big)\big(\bM(|D^{2} u|^{q})+\bM(|\partial_{t}u|^{q})\big)^{1/q}\Big).
\end{equation}
\end{lemma}  
 
Since $\partial_{t}u=\cL_{ 0 } u-a^{ij}D_{ij}u$,
similar estimate holds for $\partial_{t}u$.
Therefore,   on account of taking $q$ as close to $1$
as needed,
as in the case of \eqref{9.18,12} we get the following. 
\begin{lemma}   
                        \label{lemma 4.2,6}
For any $p>1$, $w\in A_{p}$, $\kappa\geq2 $,
and $u\in C^{\infty}_{0}$
$$
\|D^{2}u\|_{L_{p}(w)}^{p}+\|\partial_{t}u\|_{L_{p}(w)}^{p}
\leq N([w]_{A_p})\Big(N(d,\delta,p,\kappa)
\|\cL_{ 0 } u\|_{L_{p}(w)}^{p} 
$$
\begin{equation}
                      \label{4.2,7}
+\big(N(d,\delta,p,\kappa)a^{\shharp} +N(d,\delta,p )\kappa^{-1}\big)^{p}\Big( \|D^{2}u\|_{L_{p}(w)}^{p}+\|\partial_{t}u\|_{L_{p}(w)}^{p}\big)\Big). 
\end{equation}
\end{lemma}

\begin{theorem}
                \label{theorem 4.30,1}
Let $p \in(1,\infty)$, $ \beta>0$. Then
there is a constant $\hat a=\hat a(d,\delta,p,\beta) >0$ and a constant
$N=N(d,\delta,p,\beta)$ such that if
\begin{equation}
                       \label{4.30,2}
a^{\shharp}\leq \hat a,
\end{equation}
then for any $\lambda\geq0$ and
$u\in \dot E^{1,2}_{p,\beta}$
\begin{equation}
                       \label{4.30,3}
\|\lambda u,\sqrt\lambda Du,\partial_{t}u,D^{2}u\|_{\dot E_{p,\beta}
(\bR^{d+1})}\leq N\|\cL_{ 0 } u-\lambda u\|
_{\dot E_{p,\beta}(\bR^{d+1})}.
\end{equation}
Furthermore, for any $\lambda>0$
and $f\in \dot E_{p,\beta}(\bR^{d+1})$
there exists a unique solution 
$u\in \dot E^{1,2}_{p,\beta}$
of the equation $\cL_{ 0 } u-\lambda u=f$.
\end{theorem}

Proof. To prove \eqref{4.30,3}
observe that the case of arbitrary $\lambda\geq0$
is reduced to $\lambda=0$ by Agmon's method and interpolation (see
Lemma \ref{lemma 3.31.4}). Then it is not hard to see
that we may assume that $u$ has
compact support, so that $u\in W^{1,2}_{p}$.

Now take $\alpha\in[0,d+2)$ such that $\alpha+p\beta>d+2$
and let $w =w_{\alpha}\wedge 1$,
so that $w$ is an $A_{1}$-weight owing to Lemma \ref{lemma 1.29.3}.
Approximating $u$ by smooth functions
in $W^{1,2}_{p}$ we conclude that \eqref{4.2,7}
holds for our $u$.

  By using
Lemma \ref{lemma 1.30.2}, we get from
\eqref{4.2,7}  and H\"older's inequality  that
 for $w=w_{\alpha}\wedge 1$ we have 
$$
\Big(\dashint_{C_{1}}|D^{2}u|^{p}\,dxdt\Big)^{1/p}\leq
N\|D^{2}u\|_{L_{p}( w)}
 \leq N\|\cL_{ 0 } u \|
_{\dot E_{p,\beta}(\bR^{d+1})} 
$$
$$
+\big(N(d,\delta,p,\kappa,\beta)a^{\shharp} +N(d,\delta,
 p, \beta )\kappa^{-1}\big)\|\partial_{t}u,D^{2}u\|_{\dot E_{p,\beta}
(\bR^{d+1})}. 
$$

Shifting $w$, using the self-similarity and also
treating $\partial_{t}u$ in the same way give      
$$
\|\partial_{t}u,D^{2}u\|_{\dot E_{p,\beta}
(\bR^{d+1})}\leq N\|\cL_{ 0 } u \|
_{\dot E_{p,\beta}(\bR^{d+1})} 
$$
$$
+\big(N(d,\delta,p,\kappa,\beta)a^{\shharp}
 +N(d,\delta, p, \beta )\kappa^{-1}\big)\|\partial_{t}u,D^{2}u\|_{\dot E_{p,\beta}
(\bR^{d+1})}, 
$$
which easily yields \eqref{4.30,3}
with $\lambda=0$.

The solvability part of the theorem,
as usual, is proved by the method of continuity starting from Theorem
\ref{theorem 1.31.3}.
\qed
 
In light of Theorem  \ref{theorem 12.2.1}  and
Remark \ref{remark 4.2,1}, Lemma \ref{lemma 4.2,6}  leads
to the following.

\begin{lemma}
                    \label{lemma 4.2,2}
For any $p,q> 1$,
$\kappa\geq 2$, and $u\in C^{\infty}_{0}$  we have 
$$
\|\partial_{t}u,D^{2}u\|_{\L_{q,p}}\leq
  N_{0}(d,\delta,p,q,\kappa) \|\cL_{ 0 } u\|_{\L_{q,p}} 
$$
\begin{equation}
                         \label{4.2,2}
+\big(N_{1}(d,\delta,p,q,\kappa)a^{\shharp} +N_{2}(d,\delta, p,q)\kappa^{-1}\big)\|
\partial_{t}u,D^{2}u\|_{\L_{q,p}}.
\end{equation}
 
\end{lemma}

\begin{corollary}
                 \label{corollary 4.2,1}
Let $p,q>1$ and assume that
\begin{equation}
                         \label{4.2,3}
 a^{\shharp} \leq (1/4)
N_{1}^{-1}\big(d,\delta,p,q,4
N_{2} (d,\delta, p,q)\big).
\end{equation}
Then for any $u\in C^{\infty}_{0}$  we have
\begin{equation}
                                 \label{4.3,3}
\|\partial_{t}u,D^{2}u\|_{\L_{q,p}}\leq
  2N_{0}(d,\delta,p,q,4
N_{2} (d,\delta, p,q)) \|\cL_{ 0 } u\|_{\L_{q,p}}, 
\end{equation}
where $N=N(d,\delta, p,q )$.
\end{corollary}

We restate this corollary as follows.

\begin{theorem}
                    \label{theorem 4.3,1}
Let $p,q\in(1,\infty) $. Then
there are   constants $\hat a =\hat a(d,\delta,p,q)>0$ and $N=N(d,\delta, p,q)$ such that if 
\begin{equation}
                         \label{4.3,5}
 a^{\shharp} \leq \hat a,
\end{equation}
then for any $u\in C^{\infty}_{0}$
we have
\begin{equation}
                         \label{4.2,5}
\|\partial_{t}u,D^{2}u\|_{\L_{q,p}}\leq
  N  \|\cL_{ 0 } u\|_{\L_{q,p}}. 
\end{equation}

\end{theorem}

Next, as in the case of Theorem \ref{theorem 4.2,1} we introduce $\lambda\geq0$  to see that
for $\lambda\geq0$ and $u\in C^{\infty}_{0}$
$$
\|\partial_{t}u,D^{2}u,\sqrt\lambda Du,\lambda u\|_{\L_{q,p}}\leq
N\|\cL_{ 0 } u-\lambda u\|_{\L_{q,p}}. 
$$
In this way we arrive at the first assertion in the following theorem
if  $u\in C^{\infty}_{0}$
and $s=-\infty$.   

\begin{theorem}
                 \label{theorem 4.3,3}
Let $p,q\in( 1,\infty)$, $s\in[-\infty,\infty)$ and assume \eqref{4.3,5}.  
Then there exists $N =N (d,\delta,p,q)$   such that for any $u\in \W^{1,2}_{q,p}
(\bR^{d+1}_{s})$
and $\lambda\geq 0$ 
\begin{equation}
                            \label{4.3,4}
\|\partial_{t}u,D^{2}u,\sqrt\lambda Du,\lambda u\|_{\L_{q,p}(\bR^{d+1}_{s})}\leq
N\|\cL_{ 0 } u-\lambda u\|_{\L_{q,p}(\bR^{d+1}_{s})}. 
\end{equation}
Furthermore, for any $f\in \L_{q,p}(\bR^{d+1}_{s})$
and $\lambda>0$
there exists a unique $u\in \W^{1,2}_{q,p}(\bR^{d+1}_{s})$
satisfying in $\bR^{d+1}_{s}$
$$
\cL_{ 0 } u-\lambda u+f=0.
$$
\end{theorem}

 For $s=-\infty$ the a priori estimate \eqref{4.3,4} extends
from $C^{\infty}_{0}$ to $\W^{1,2}_{q,p}$
due to the denseness of the former in the latter and the solvability result is obtained by the method of continuity starting from
Theorem \ref{theorem 2.14.1} or
\ref{theorem 4.2,1}. General $s$ is treated as in the proof of 
Theorem  
\ref{theorem 4.2,1}.

In the future we will need the following
interior estimate.

\begin{theorem}
                 \label{theorem 4.29,3}
Let $p,q\in( 1,\infty)$ and assume \eqref{4.3,5}.
Let $0<r<\rho<\infty$,    
$u\in \W^{1,2}_{q,p}(C_{\rho})$.
Set
$f:=\cL_{ 0 } u $.
 Then   
 \begin{equation}
                                                 \label{06.5.25.20}
 \|\partial_{t}u, D^{2}u\|_{\L_{q,p}(C_{r})}
\leq N\big(
\|f\|_{\L_{q,p}(C_{\rho})} 
+(\rho-r)^{-2}\|u-c\|_{\L_{q,p}(C_{\rho})}\big), 
\end{equation} 
where $N=N(d,\delta,p,q)$ and $c$ is any number.

\end{theorem}

Proof. We follow the proof of Lemma 2.4.4
of \cite{Kr_08}.   
For obvious reasons we may assume that $u\in C^{1,2}(\bar{C}_{\rho})$ and $c=0$.
Then, let $\chi(s)$ be an infinitely differentiable function
on $\bR$ such that $\chi(s)=1$ for $s\leq0$ and $\chi(s)=0$
for $s\geq1$. For $m=0,1,2,...$ introduce  ($r_{0}=r$)
 $$
r_{m}=r+(\rho-r)\sum_{j=1}^{m}2^{-j},\quad
\xi_{m}(x)=\chi\big(2^{m+1}(\rho-r)^{-1}(|x|-r_{m}) \big),
$$
 $$
\eta_{m}(t)=\chi\big(2^{2m+2}(\rho-r)^{-2}(t-r^{2}_{m})\big),
\quad \zeta_{m}(t,x)=\xi_{m}(x)\eta_{m}(t).    
$$ 
As is easy to check, for
$$
C(m)=C_{r_{m}}=(0,r_{m}^{2})\times B_{r_{m}}
 $$
it holds that
 $$
\zeta_{m}=1\quad\text{on}\quad C(m),\quad\zeta_{m}=0\quad\text{on}
\quad C_{\rho}\setminus C(m+1).
 $$
Also (observe that $N2^{m+1}=N_{1}2^{m}$ with $N_{1}=2N$) 
 $$
|D\zeta_{m }|\leq N2^{m}(\rho-r)^{-1},\quad
|\partial_{t}\zeta_{m}|\leq N2^{2m}(\rho-r)^{-2},
  $$
\begin{equation}
                                                                \label{06.5.25.1}
|D^{2}\zeta_{m}|\leq N2^{2m}(\rho-r)^{-2}.
\end{equation} 

Next, the function $\zeta_{m}u$ is in $\W^{1,2}_{q,p}(\bR^{d+1}_{0})$  {}
and satisfies 
 $$
\cL_{ 0 }(\zeta_{m} u) =\zeta_{m} f
+u\cL\zeta_{m} +
2a^{ij}D_{i}\zeta_{m}D_{j}u .
 $$
By Theorem \ref{theorem 4.3,1}
and the above-mentioned properties of $\zeta_{m}$
$$
A_{m} :=\|\partial_{t}(\zeta_{m} u), D^{2}(\zeta_{m}u)\|_{\L_{q,p}(Q(m+1))}
             \label{06.5.25.2} 
 \leq N\|f\|_{\L_{q,p}(C_{\rho})} +N2^{2m}(\rho-r)^{-2}I+NJ_{m},
$$ 
where
 $$
I=\|u \|_{\L_{q,p}(C_{\rho})},
\quad J_{m}=\| 
2a^{ij}D_{i}\zeta_{m}D_{j}u\|_{\L_{q,p}(C_{\rho})}
 $$ 
 $$
 \leq  N2^{m}(\rho-r)^{-1}\|Du \|_{\L_{q,p}
(Q(m+1))}.
 $$
By interpolation inequalities  for  
any $\varepsilon>0$
$$
\|Du \|_{\L_{q,p}
(Q(m+1))}\leq \|D(\zeta_{m+1}u) \|_{\L_{q,p}
(\bR^{d+1}_{0})} 
$$
$$
\leq \varepsilon (\rho-r) 2^{- m}
A_{m+1}+N\varepsilon^{-1}(\rho-r)^{-1}
2^{ m}I. 
$$

It follows that for any $\varepsilon\in(0,1]$
$$
A_{m}\leq N\|f\|_{\L_{q,p}(C_{\rho})} +N
\varepsilon^{-1}2^{2m}(\rho-r)^{-2}I
+\varepsilon A_{m+1}. 
$$
 
Now we take $\varepsilon=1/8$ and get 
 \begin{gather}\nonumber
\varepsilon^{m}A_{m} 
\leq N\varepsilon^{m}\|f\|_{\L_{q,p}(C_{\rho})}+N\varepsilon^{m-1}
  2^{2m}(\rho-r)^{-2} I+
\varepsilon^{m+1} A_{m+1},
\\ 
                                                                \label{06.5.25.3}
A_{0}+\sum_{m=1}^{\infty}\varepsilon^{m}A_{m} 
\leq N \|f\|_{\L_{q,p}(C_{\rho})}+N (\rho-r)^{-2} I +
\sum_{m=1}^{\infty}\varepsilon^{m}  A_{m }.
\end{gather} 
Here the series converges because owing to \eqref{06.5.25.1},
 $$
A_{m}\leq N(1+4^{m}(\rho-r)^{-2})\|u\|_{\W^{1,2}_{q,p}(C_{\rho})}.
 $$
Therefore upon cancelling like terms in \eqref{06.5.25.3},
we see that $A_{0}$ is less than the right-hand side of \eqref{06.5.25.20}.
Since its left-hand side is obviously less than $A_{0}$, the theorem is proved.
\qed

Next, we refer to the arguments in Remark
\ref{remark 2.17.1}. Its conclusions  are based
on some facts which in no way used
the specific order of integration defining $\L_{q,p}$
and $\W^{1,2}_{q,p}$
and
are valid word for word in our setting
despite the order of integration here is reversed.

\begin{theorem}
                         \label{theorem 4.1,1}
Let $p,q\in(1,\infty)$  and
$$
\frac{d}{p}+\frac{2}{q}=\beta>1.
$$
Let $\bR^{d}$-valued $b\in \L_{q,p}$
and \eqref{4.3,5} be satisfied. Then
there exists a constant $\hat b>0$, depending
only on $d,\delta,p,q,\beta$ such that, if
$$
\|b\|_{\L_{\beta q,\beta p}}\leq \hat b,
$$
then there exists $N =N (d,\delta,p,q,\beta)$   such that for any $u\in \W^{1,2}_{q,p}$
and $\lambda\geq0$ 
\begin{equation}
                            \label{4.8,3}
\|\partial_{t}u,D^{2}u,\sqrt\lambda Du,\lambda u\|_{\L_{q,p}}\leq
N\|\cL u-\lambda u\|_{\L_{q,p}}, 
\end{equation}
where
$$
\cL u=\partial_{t}u+a^{ij}D_{ij}u+b^{i}D_{i}u .
$$
Furthermore, for any $f\in \L_{q,p}$
and $\lambda>0 $
there exists a unique $u\in \W^{1,2}_{q,p}$
satisfying 
$$
\cL u-\lambda u+f=0.
$$
\end{theorem}   
 
Of course, this theorem is also
true for equations in $\bR^{d+1}_{s}$  
as in Theorem~\ref{theorem 4.3,3}.
Observe that in this theorem the norm
$\|b\|_{L_{q\beta,p\beta}}$ is completely
analogous to the norm $\|b\|_{L_{d}}$
in the elliptic case and the Sobolev space
theory. In particular, because $d/(p\beta)+2/(q\beta)=1$, the function
\begin{equation}
                         \label{5.9,1}
\frac{1}{|x|+\sqrt{|t|}}I_{|x|\le	 1, |t| \leq1}
\end{equation}
does belong to $\L_{q'\beta,p'\beta}$
if $p'<p$, $q'<q$, but not
for $p'=p,q'=q$.

\subsection{Basic estimate in mixed-norm 
homogeneous Morrey-So\-bo\-lev spaces}

Theorem \ref{theorem 4.1,1} is about
the solvability of parabolic
equations in the mixed-norm Sobolev spaces.
Here we concentrate on the mixed-norm 
Morrey-Sobolev spaces.

Take $p,q\in(1,\infty)$ and
   $\beta> 0$ and introduce the
Morrey space $\dot \E_{q,p,\beta} $
as the set of $g\in  \L_{q,p,\loc}$ 
\index{$A$@Sets of functions!$\dot \E_{q,p,\beta}$}%
\index{$N$@Norms!$"|"|f"|"|_{\dot \E_{q,p,\beta}}$}%
(using the  definition of $\L_{q,p}$ 
in \eqref{4.29,5}) such that  
\begin{equation}
                             \label{2.16.50}
\|g\|_{\dot \E_{q,p,\beta} }:=
\sup_{\rho >0,C\in\bC_{\rho}}\rho^{\beta}
\dashnorm g  \|_{ \L_{q,p}(C)} <\infty .
\end{equation}  
Define
$$
\dot \E^{1,2}_{q,p,\beta} =\{u:u,Du,D^{2}u,
\partial_{t}u\in \dot \E_{q,p,\beta} \},
$$
where $Du,D^{2}u,
\partial_{t}u$ are Sobolev
\index{$A$@Sets of functions!$\dot E^{1,2}_{q,p,\beta}$}%
 derivatives,
and 
provide $\dot \E^{1,2}_{q,p,\beta} $ with an obvious norm.
Recall that the subsets of these spaces
consisting of functions independent of
$t$ are denoted by $\dot E_{p,\beta}$
and $\dot E^{2}_{p,\beta}$, respectively.
\index{$A$@Sets of functions!$\dot E_{p,\beta}$}%
\index{$A$@Sets of functions!$\dot E^{2}_{p,\beta}$}%

We take the operator $\cL_{ 0 }$
from the beginning of Subsection
\ref{section 4.29.1}.

\begin{theorem}
                    \label{theorem 4.29,4}
Let $p,q\in (1,\infty),\beta\in (1,2) $,
$u\in \dot \E^{1,2}_{q,p, \beta}$. 
Assume that the conditions \eqref{4.30,2}
 with $p\wedge q$ in place of $p$  and \eqref{4.3,5} are satisfied.
Then
\begin{equation}
                        \label{4.29,3}
\|\partial_{t}u,D^{2}u\|_{\dot \E_{q,p,\beta}}
\leq N(d,\delta,p,q,\beta)\|f\|_{\dot \E_{q,p,\beta}},
\end{equation}
where $f=\cL_{ 0 } u$.
\end{theorem}

Proof. First we observe that it suffices to
prove \eqref{4.29,3} for $u$ with compact support. For such $u$ and   $\lambda\in(0,1]$
define
$$
f_{\lambda}=\cL_{ 0 } u-\lambda u,\quad g_{\lambda}=f_{\lambda}I_{C_{2}},\quad h_{\lambda}=f_{\lambda}I_{C_{2}^{c}}
$$ 
and let $G_{\lambda},H_{\lambda}$ be solutions
of class $\W^{1,2}_{q,p}(\bR^{d+1})$
of the equations 
\begin{equation}
                        \label{4.29,30}
\cL_{ 0 } G_{\lambda}-\lambda G_{\lambda}=g_{\lambda},\quad 
\cL_{ 0 } H_{\lambda}-\lambda H_{\lambda}=h_{\lambda}. 
\end{equation}
Since $f_{\lambda}$ is compactly supported
element of $\dot \E_{q,p,\beta}$,
it belongs to $\L_{q,p}$ and $G_{\lambda},
H_{\lambda}$ exist by Theorem 
\ref{theorem 4.3,3}. Furthermore,
$f_{\lambda}\in \dot E_{p\wedge q,\beta}
(\bR^{d+1})$,
so that by Theorem \ref{theorem 4.30,1}
the equations in \eqref{4.29,30} have
solutions in $\dot E^{1,2}_{p\wedge q,
 \beta }(\bR^{d+1})$.
These are the same $G_{\lambda},H_{\lambda}$, since these solutions
in different spaces were obtained
by the method of continuity
starting from the heat equation
for which the explicit formulas
for solutions are available.

Next, note that $\hat H_{\lambda}
:=e^{t}H_{\lambda}$ satisfies
$$
\cL_{ 0 } \hat H_{\lambda}-(\lambda+1)
\hat H_{\lambda}=e^{t}h_{\lambda}:=\hat
h_{\lambda}
$$
and $\hat h_{\lambda}\in \L_{q,p}(\bR^{d+1}_{0})$ with its norm dominated by a constant independent of $\lambda$, since
$h_{\lambda}$ has compact support.
It follows from Theorem \ref{theorem 4.3,3}
that 
\begin{equation}
                  \label{1.5,2}
\sup_{\lambda\in(0,1]}
\|H_{\lambda}\|_{\L_{q,p}(\bR^{d+1}_{0})}
\leq \sup_{\lambda\in(0,1]}\|\hat H_{\lambda}\|_{\L_{q,p}(\bR^{d+1}_{0})}<\infty. 
\end{equation}

Concerning $G_{\lambda}$ note that
by Theorem \ref{theorem 4.3,3} 
$$
\| D^{2}G_{\lambda}\|_{\L_{q,p}(C_{1})}\leq N\|g_{\lambda}
\|_{\L_{q,p}}\leq N\|f_{\lambda}\|_{\L_{q,p}(C_{2})} \leq N\|f_{\lambda}\|_{\dot \E_{q,p,\beta}}. 
$$

Now observe that by applying Theorem \ref{theorem 4.29,3} we get 
$$
\|D^{2}H_{\lambda}\|_{\L_{q,p}(C_{1})}
\leq N\|\lambda H_{\lambda}\|_{\L_{q,p}
(C_{2})}+N\|H_{\lambda}- H_{\lambda}(0)
\|_{\L_{q,p}
(C_{2})}, 
$$
where owing to Theorem \ref{theorem 2.3.1}  
($1<\beta<2$) 
 the last term is dominated by
$$
N\|\partial_{t}H_{\lambda},D^{2}H_{\lambda}\|_{\dot E_{p\wedge q,\beta}},
$$
which, in turn, by Theorem
\ref{theorem 4.30,1} is dominated by 
$$
N\|h_{\lambda}\|_{\dot E_{p\wedge q,\beta}}\leq N\|f_{\lambda}\|_{\dot E_{p\wedge q,\beta}}\leq N\|f_{\lambda}\|_{\dot \E_{q,p,\beta}}.
$$

Since $u=G_{\lambda}+H_{\lambda}$
and \eqref{1.5,2} holds
it follows that 
$$
\|D^{2}u\|_{\L_{q,p}(C_{1})}
\leq \lim_{\lambda\downarrow 0}(
\|D^{2}G_{\lambda}\|_{\L_{q,p}(C_{1})}
+\|D^{2}H_{\lambda}\|_{\L_{q,p}(C_{1})}) 
$$
$$
\leq N\lim_{\lambda\downarrow 0}\|f_{\lambda}\|_{\dot \E_{q,p,\beta}}
\leq N\|f \|_{\dot \E_{q,p,\beta}}. 
$$
Parabolic scaling (preserving 
\eqref{4.30,2} and \eqref{4.3,5})
shows that for any
$\rho>0$ 
$$
\rho^{\beta}\|D^{2}u\|_{\L_{q,p}(C_{\rho})}
\leq N\|f \|_{\dot \E_{q,p,\beta}} 
$$
and then shifting the origin leads to
\eqref{4.29,3} in what concerns $D^{2}u$.
Since $\partial_{t}u=f-a^{ij}D_{ij}u$,
the theorem is proved. \qed

\begin{remark}
                                   \label{remark 8.16,1}
 One can show that the constraint: $\beta<2$, can be dropped.
\end{remark}   
\def\bim{\protect\eqref{4.30,2},
\protect\eqref{4.3,5}\,}
\subsection{Localizing  \bim and local Morrey-Sobolev spaces}

In many situations it is more convenient to work in 
local Morrey-Sobolev spaces rather
then in $\dot \E_{q,p,\beta}$, which
does not contain constants apart
from the case that $\beta=d/p+2/q$.

Take $p,q\in(1,\infty)$ and
   $\beta\geq 0$ and introduce the
Morrey space $\E_{q,p,\beta} $ 
as the set of $g\in  \L_{q,p,\loc}$ 
\index{$A$@Sets of functions!$\E_{q,p,\beta}$}%
\index{$N$@Norms!$"|"|g"|"|_{\E_{q,p,\beta}}$}%
such that  
\begin{equation}
                             \label{8.11.020}
\|g\|_{\E_{q,p,\beta} }:=
\sup_{\rho\leq 1,C\in\bC_{\rho}}\rho^{\beta}
\dashnorm g  \|_{ \L_{q,p}(C)} <\infty .
\end{equation}  
Define
$$
\E^{1,2}_{q,p,\beta} =\{u:u,Du,D^{2}u,
\partial_{t}u\in \E_{q,p,\beta} \},
$$
where $Du,D^{2}u,
\partial_{t}u$ are Sobolev derivatives,
and 
\index{$A$@Sets of functions!$\E^{1,2}_{q,p,\beta}$}%
provide $\E^{1,2}_{q,p,\beta} $ with an obvious norm.
The subsets of these spaces
consisting of functions independent of
$t$ we denoted by $E_{p,\beta}$
and $E^{2}_{p,\beta}$, respectively.
\index{$A$@Sets of functions!$E_{p,\beta}$}%
\index{$A$@Sets of functions!$E^{2}_{p,\beta}$}%

Observe that if    $\beta\leq
d/p+2/q$ and the support of $g$ is in some $C\in\bC_{1}$,
then the $\dot \E_{q,p,\beta}$-norm of $g$ coincides with its $\E_{q,p,\beta}$-norm.
One more useful observation is that
\begin{equation}
                               \label{3.30.3}
\|g\|_{\E_{q,p,\beta} }=\sup_{C\in \bC_{1}}
\|gI_{C}\|_{\E_{q,p,\beta} }.
\end{equation}

 \begin{remark} 
                \label{remark 7,2.10}
If $\beta\geq d/p+2/q$, then
$$
\|g\|_{\E_{q,p,\beta}}= \sup_{C\in \bC_{1}}\|g\|_{L_{ p,  q}(C)},
$$
H\"older's inequality easily shows that,
if $\beta\leq d/p+2/q$, then 
$$
\|g\|_{\E_{q,p,\beta}}\leq \sup_{C\in \bC_{1}}\|g\|_{L_{\bar p,\bar q}(C)},
$$
 where $(\bar p,\bar q) =\alpha (p,q)$ and $\beta\alpha=
d/p+2/q$.
\end{remark}

\begin{remark}
                        \label{remark 3.27.10}
If $u$ vanishes outside $C\in\bC_{1}$  and  if $\beta\leq d/p +2/q$, one easily sees that
for any $r\geq 1$ and $C'\in \bC_{r}$
$$
r^{\beta }\dashnorm u\|_{\L_{q,p}(C')}\leq N(d)
\|u\|_{\L_{q,p}(C)}. 
$$
\end{remark}

\begin{remark}
                \label{remark 4.14,1}
If $f(t,x)\geq 0$ is nonnegative, $r\in(0,\infty),\mu
\in(0,1)$, then
with 
$\zeta=|C_{r\mu}|^{-1}I_{C_{r\mu}}$ we have
$$
\dashnorm I_{C_{r}}f\|_{\L_{q,p}}
=N(d)r^{-d/p-2/q}\Big\|\int_{C_{(1+\mu)r}} (I_{C_{r}}f)  \zeta(s-\cdot,y-\cdot)\,dyds\Big\|_{\L_{q,p}} 
$$
$$
\leq Nr^{-d/p-2/q}\int_{C_{(1+\mu)r}}\|\zeta(s-\cdot,y-\cdot)f\|_{\L_{q,p}}\,dyds   
$$
$$
\leq N(d)(1+\mu)^{d+2}   \mu  ^{d/p+2/q-(d+2)}\sup_{C\in
\bC_{\mu r}}\dashnorm
f\|_{\L_{q,p}(C )}. 
$$

Since $\mu\leq 1$ and $p,q>1$
it follows that for $\rho\leq r$ 
$$
 \sup_{C\in\bC_{r}}\dashnorm f\|_{\L_{q,p}(C)}\leq N(d) (\rho/r)^{d/p+2/q-(d+2)} \sup_{C\in\bC_{\rho}}\dashnorm f\|_{\L_{q,p}(C)}. 
$$
\end{remark} 

For $\rho>0$, $p ,q  \in(1,\infty)$ introduce 
\index{$S$@Miscelenea!$a^{\shharp}$@$a^{\shharp}_{\rho}$}%   
\begin{equation}
                         \label{6.3.1}
 a^{\shharp}_{\rho} = \sup
_{r\leq\rho}
 \sup
_{ C\in\bC_{r} }\dashint_{C}|a(t,x)-\tilde a_{C}(t)|\,dxdt ,
\end{equation} 
\index{$S$@Miscelenea!$\bar  b_{\rho_{b}}$@$\bar b_{p ,q ,\rho}$}% 
\begin{equation}
                           \label{3.14.2}
\bar b_{q ,p ,\rho }=\sup_{r\leq\rho }r
\sup_{C\in \bC_{r}} 
\dashnorm b \|_{\L_{q ,p  }(C)}.
\end{equation}

Fix  $p ,q ,\beta    $ such that
\begin{equation}
                        \label{3.21.01}
p ,q  \in(1,\infty),\quad 1<\beta\leq\frac{d}{p }+\frac{2}{q } ,\quad \beta<2.
\end{equation}

\begin{lemma}
                          \label{lemma 5.7,1}
There exist a   constant $\kappa
=\kappa(d )\in(0,1)$ and a (large)
constant $K=K(d )$ such that for any $\rho>0$  there
exists an $\bS_{\delta}$-valued $\check a$
on $\bR^{d+1}$
such that $a=\check a$ in $C_{\kappa \rho }$
and $\check a^{ \shharp}\leq 
Ka_{\rho}^{ \shharp}$.
\end{lemma}
  
We are not proving this lemma here
and plan to do this in a separate
publication because the line of arguments
needed is somewhat far from the main
reasoning in this article.

 Recall that 
$$
\cL_{ 0 }  u=\partial_{t}u+a^{ij}D_{ij}u,
$$
and fix some
$$
\rho_{a},\rho_{b}\in(0,1].
$$

Let $\zeta_{0}\in C^{\infty}_{0}$
have suport in $C_{\kappa \rho_{a}} $ and
equal $1$ in $(1/2)C_{\kappa \rho_{a}}$,
where $\kappa$ is taken from Lemma \ref{lemma 5.7,1}. 
\begin{lemma}
                     \label{lemma 5.8,1}
There exists
$\hat a=\hat a(d,\delta,p,q,\beta)>0$
such that, if $a^{\shharp}_{\rho_{a}}\leq 
\hat  a$, 
\index{$S$@Miscelenea!$a^{\shharp}$@$\hat a$}%
then for any $u\in \E^{1,2}_{q,p,\beta}$ 
\begin{equation}
                        \label{5.8.3}
\|\zeta_{0}\partial_{t}u,\zeta_{0}D^{2}u\|_{  \E_{q,p,\beta}}
\leq N  \|f\|_{  \E_{q,p,\beta}}+ N_{1}\|Du,u\|_{  \E_{q,p,\beta}} , 
\end{equation}
where $f=\cL_{ 0 } u$, $N =N (d,\delta,p,q,\beta )$,
$N_{1} =N (d,\delta,p,q,\beta,\rho_{a} )$. 
\end{lemma}

Proof. Observe that
\begin{equation}
                        \label{5.8,2}
\|\zeta_{0}\partial_{t}u,\zeta_{0}D^{2}u\|_{ 
 \E_{q,p,\beta}}\leq
\|\partial_{t}(\zeta_{0}u),D^{2}(\zeta_{0}u)\|_{ 
{ \dot \E_{q,p,\beta}}}+N\|Du,u\|_{  \E_{q,p,\beta}}
\end{equation}
and, owing to Lemma \ref{lemma 5.7,1}
and Theorem \ref{theorem 4.29,4}, the first term is estimated through
$$
(\partial_{t}+\check a^{ij}D_{ij})(\zeta_{0}u)
=\cL_{ 0 }(\zeta_{0}u)=\zeta_{0}f + g,
$$
where, thanks to  $\beta\leq d/p+2/q$, the $\dot \E_{q,p,\beta}$-norm of $g$
is estimated by the last term in
\eqref{5.8,2}. \qed

Naturally, \eqref{5.8.3} holds for
any shift of $\zeta_{0}$ and this
along with Remark \ref{remark 4.14,1} proves
the following with    $ \|Du,u\|_{  \E_{q,p,\beta}}$
in place of $ 
\|u\|_{  \E_{q,p,\beta}}$. The term 
$\|Du \|_{  \E_{q,p,\beta}}$ is eliminated
by using interpolation inequalities.

\begin{theorem}
                     \label{theorem 5.8,1}
There exists
$\hat a=\hat a(d,\delta,p,q,\beta)>0$
such that, if $a^{\shharp}_{\rho_{a}}\leq 
\hat  a$, then for any $u\in \E^{1,2}_{q,p,\beta}$ 
\begin{equation}
                        \label{5.8.40}
\| \partial_{t}u, D^{2}u\|_{  \E_{q,p,\beta}}
\leq N  \|\cL_{ 0 } u\|_{  \E_{q,p,\beta}}+N\| u\|_{  \E_{q,p,\beta}} , 
\end{equation}
where $N =N (d,\delta,p,q,\beta,\rho_{a})$.
 
\end{theorem}

\begin{remark}
                      \label{remark 5.14,1}
Of course, it would be desirable to
exclude $\rho_{a}$ from the arguments
of the factor of $ \|\cL_{ 0 } u\|_{  \E_{q,p,\beta}}$ in \eqref{5.8.40} as in somewhat similar elliptic situation of Lemma
\ref{lemma 2.20.5}.
\end{remark}

Next, as a few times before, we use Agmon's
idea along
with the interpolation Lemma \ref{lemma 3.31.4} to treat 
$\cL_{ 0 }-\lambda$. 
After that to deal with the first-order
terms, note that for $\zeta\in C^{\infty}_{0}(
2C_{1})$ such that $\zeta=1$ on $C_{1}$
by Remark \ref{remark 10.8.2} we obtain
for $\rho\leq 1$ and $u\in C^{\infty}_{0}$ that   
$$
\rho^{\beta}\|b^{i} D_{i}u\|_{\L_{q,p}(C_{\rho})}
\leq  \|\zeta b^{i} D_{i}(\zeta u)\|_{\dot \E_{q,p,\beta}(C_{\rho})}
\leq N\|b  \|_{ \E_{q\beta,p\beta,1}}
\|u\|_{ \E^{1,2}_{q,p,\beta }},
$$
implying that
\begin{equation}
                             \label{5.8,5}
 \|b^{i} Du\|_{\E_{q,p,\beta}}
\leq N\|b  \|_{ \E_{q\beta,p\beta,1}}
\|u\|_{ \E^{1,2}_{q,p,\beta }}. 
\end{equation}
It is to be noted that Remark \ref{remark 10.8.2} treats mixed norm in the order
different from what is in this section.
However, just repeating the arguments
leading to Remark \ref{remark 10.8.2}
shows that it is valid in our current setting
as well. We have already made a few comments to that effect.

By considering the cut-off mollifiers of
$u\in \E^{1,2}_{q,p}$ we extend 
\eqref{5.8,5} to all $u\in \E^{1,2}_{q,p}$.
In this way we come to the a priori estimate
below.  Note 1 in the norm of $b$.

\begin{theorem}
                     \label{theorem 5.8,2}
There 
\index{$S$@Miscelenea!$a^{\shharp}$@$\hat a$}%
\index{$S$@Miscelenea!$a^{\shharp}$@$\hat b$}%
exist
$$
\hat a=\hat a(d,\delta,p,q,\beta)>0,\quad\hat b=\hat b(d,\delta,p,q,\beta,\rho_{a} )>0,
$$ 
$$\lambda_{0}= \hat \lambda(d,\delta,p,q,\beta,\rho_{a} ) >0,\quad
N=N(d,\delta,p,q,\beta,\rho_{a}  ),
$$
such that, if 
$$
a^{\shharp}_{\rho_{a}}\leq 
\hat  a,\quad \bar b_{q\beta,p\beta,1} \leq \hat b,
$$
 then for any $u\in \E^{1,2}_{q,p,\beta}$  
 and $\lambda\geq\lambda_{0}$ 
\begin{equation}
                        \label{5.8.4}
\|\lambda u,\sqrt\lambda Du, \partial_{t}u, D^{2}u\|_{  \E_{q,p,\beta}}
\leq N \|f\|_{  \E_{q,p,\beta}}, 
\end{equation}
where   $f=\cL_{ 0 } u+b^{i}D_{i}u-\lambda u$.
Furthermore, for any
$f\in \E_{q,p,\beta}$ and $\lambda\geq\lambda_{0}$ there exists a unique
$u\in \E^{1,2}_{q,p,\beta}$ such that in
$\bR^{d+1}$
$$
\cL_{ 0 } u+b^{i}D_{i}u-\lambda u=f.
$$
 
\end{theorem}

Proof. 
To prove the existence of solutions we use
the method of continuity. In light of the a priori estimate  \eqref{5.8.4},
for this method to work, we  need only prove
the solvability for the heat equation
$$
\partial_{t}u+\Delta u-\lambda u=f.
$$
If $f$ is bounded and with each derivatives bounded, the solution $u$
admits an explicit representation,
is bounded and has bounded
derivatives.
Therefore, $u\in \E^{1,2}_{q,p,\beta}$
and admits estimate \eqref{5.8.4}.
For general $f\in \E_{q,p,\beta}$
it suffices to use its mollifiers
and repeat with obvious changes the corresponding
argument in the proof
of Theorem \ref{theorem 1.20.3}. \qed

In the following theorem we further
localize conditions on $b$.
Recall that $\rho_{b}\in(0,1]$.

 \begin{theorem}
                     \label{theorem 5.8,20}
There 
\index{$S$@Miscelenea!$a^{\shharp}$@$\hat a$}%
\index{$S$@Miscelenea!$a^{\shharp}$@$\hat b$}%
exist
$$
\hat a=\hat a(d,\delta,p,q,\beta)>0,\quad\hat b=\hat b(d,\delta,p,q,\beta,\rho_{a} )>0,
$$ 
$$\lambda_{0}= \hat \lambda(d,\delta,p,q,\beta,\rho_{a} ) >0,\quad
N=N(d,\delta,p,q,\beta,\rho_{a}  ),
$$
such that, if $\lambda\geq \lambda_{0}\rho_{b}^{-2}$ and
$$
a^{\shharp}_{\rho_{a}}\leq 
\hat  a,\quad\bar b_{q\beta,p\beta,\rho_{b}}\leq \hat b,
$$
 then for any $u\in \E^{1,2}_{q,p,\beta}$ 
\begin{equation}
                        \label{5.10.2}
\|\lambda u,\sqrt\lambda Du, \partial_{t}u, D^{2}u\|_{  \E_{q,p,\beta}}
\leq N \rho_{b}^{-\alpha}\|f\|_{  \E_{q,p,\beta}}, 
\end{equation}
where   
$$
f=\cL_{ 0 } u+b^{i}D_{i}u-\lambda u,\quad \alpha=
d+2+\beta-\frac{d}{p}-\frac{2}{q}.
$$
Furthermore, for any
$f\in \E_{q,p,\beta}$ and $\lambda\geq\lambda_{0}\rho_{b}^{-2}$ there exists a unique
$u\in \E^{1,2}_{q,p,\beta}$ such that in
$\bR^{d+1}$
$$
\cL_{ 0 } u+b^{i}D_{i}u-\lambda u=f.
$$
 
\end{theorem}

Proof. Owing to Theorem \ref{theorem 5.8,2}
our assertions are true if $\rho_{b}=1$.
Then, having \eqref{5.10.2} for $\rho_{b}=1$,
for $\rho_{b}<1$ we make a parabolic dilation
(not affecting $a^{\shharp}_{\rho_{a}}\leq 
\hat  a$) and see that for $\lambda\geq\lambda_{0}\rho_{b}^{-2}$, where
$\lambda_{0}$ is from Theorem~\ref{theorem 5.8,2},
$$
\sup_{\rho\leq\rho_{b}}\rho^{\beta}\sup
_{C\in\bC_{\rho}}\dashnorm\lambda u,\sqrt\lambda Du, \partial_{t}u, D^{2}u\|_{ \L_{q,p}(C)} 
$$
\begin{equation}
                        \label{5.8.5}
\leq N(d,\delta,p,q,\beta,\rho_{a})\sup_{\rho\leq\rho_{b}}\rho^{\beta}\sup
_{C\in\bC_{\rho}}\dashnorm  f\|_{ \L_{q,p}(C)}
\leq N(d,\delta,p,q,\beta,\rho_{a})\|f\|_{  \E_{q,p,\beta}}. 
\end{equation}
By Remark \ref{remark 4.14,1} the left-hand
side of \eqref{5.8.5} times $N(d)\rho_{b}^{-\alpha}$ dominates the left-hand side
of \eqref{5.10.2}. This proves the theorem.
\qed

\begin{remark}
                      \label{remark 5.10,3}
This theorem allows $b$ to be bounded. Also
the function \eqref{5.9,1}, call it $g$,  belongs, for instance, to
$E_{p,p,\beta}$ as long as $ 1<p<d+2 $,
so that, if $|b|=cg$ with sufficiently small
$c$, $b$ is a valid choice for Theorem \ref{theorem 5.8,20}
 for $\rho_{b}=1$ and appropriate $\beta$. 

Indeed, if $|x|,\sqrt{|t|}\leq
2\rho$, then $C_{\rho}(t,x)\subset
(-5\rho^{2},5
\rho^{2})\times B_{3\rho}$ and 
$$
\rho^{d+2}\dashnorm g\|^{p}_{L_{p,p}( C_{\rho}(t,x))}=N\|g\|^{p}_{L_{p,p}( C_{\rho}(t,x))}
\leq N\int_{B_{3\rho}}
 \int_{0}^{5\rho^{2}}\frac{1}{(|x|+\sqrt t)^{p}}\,dt \,dx
$$
$$
=N\rho^{d+2-p}\int_{|y|<3}
|y|^{2-p}\int_{0}^{5/|y|^{2}}\frac{dt}{
(1+\sqrt t)^{p }}dy. 
$$
The last factor of $N\rho^{d+2-p}$ is finite
if $p>2$ because
$$
\int_{0}^{\infty}\frac{dt}{
(1+\sqrt t)^{p }}<\infty 
$$
and $2-p>-d$. It is also 
obviously finite if $p\in(1,2]$.

In case $\max(|x|,\sqrt{|t|})\geq
2\rho$ for $(s,y)\in C_{\rho}(t,x)$, we have 
$$
|y|+\sqrt{|s|}\geq \max(|y|,\sqrt{|s|})
\geq \max(|x|,\sqrt{|t|}) 
$$
$$
-\max(|x-y|,\sqrt{|t-s|})\geq \rho.
$$
We see that in both cases
$$
\dashnorm g\|_{L_{p,p}( C_{\rho}(t,x))}\leq
N\rho^{-1}.
$$

\end{remark}

\begin{remark}
                       \label{remark 5.14,5}
Everything in this section is also
true if we use the traditional way
to define the mixed-norm spaces integrating
first in $x$ and then in $t$. However, there are cases when \eqref{3.21.01} holds and
\begin{equation}
                              \label{5.14,4}
\sup_{\rho \leq1,C\in\bC_{\rho}}\rho 
\dashnorm b  \|_{ \L_{q\beta,p\beta}(C)} <\infty \quad
\text{and}\quad \sup_{\rho \leq  1 ,C\in\bC_{\rho}}\rho 
\dashnorm b \|_{ L_{p\beta,q\beta}(C)} =\infty .
\end{equation}
Therefore, developing the theory of
solvability in $\L_{q,p,\beta}$ was worth
the effort. 

For instance,  
take $p,q,\beta>1$,  such that 
$$
\beta\leq \frac{d}{p}+\frac{2}{q},\quad 2\leq
q\beta<p\beta<d+1.
$$
  Set $p_{0}=\beta p,q_{0}=\beta q$ and
let
  $b$ be such that
$$
|b(t,x)|=f(t,x):=I_{ t>0 }|x|^{-1}\Big| \frac{\sqrt t}{|x|}-1 \Big|^{-1/p_{0}}.
$$
Observe that for any constant $\lambda>0$
we have $f(\lambda^{2}t,\lambda x)=\lambda^{-1}f(t,x)$, implying that for any $C_{r}(t,x)$
\begin{equation}
                          \label{8.27.1}
\dashnorm f\|_{\L_{q_{0},p_{0}}(C_{r}(t,x))}=r^{-1}\dashnorm f\|_{\L_{q_{0},p_{0}}(C_{1}(t/r^{2},x/r))}.
\end{equation}

Next, 
$$
\int_{0}^{1}|f(t,x)|^{q_{0}}\,dt=|x|^{2-q_{0}}
\int_{0}^{1/|x|^{2}}|\sqrt{s}-1|^{-q/p}\,ds=:
|x|^{2-q_{0}}J(1/|x|^{2}).
$$
Here $q/p<1$, so that the singularity in $J$
is ingterable and $J(t)\sim t^{1-q/(2p)}$
as $t\to\infty$, so that
$J(t)\leq N t^{1-q/(2p)}$. It follows that 
$$
\Big(\int_{0}^{1}|f(t,x)|^{q_{0}}\,dt\Big)^{1/q_{0}}
\leq N|x|^{1/p_{0}-1},
$$
where $N$ is independent of $x$, and 
since $1-p_{0}  >-d$, the right-hand side
to the power $p_{0}$ is summable over $B_{1}$.
Together with \eqref{8.27.1} this yield
\begin{equation}
                                \label{1.23.1}
\dashnorm f\|_{\L_{q\beta ,p\beta }(C_{r})}\leq Nr^{-1}\quad r>0.
\end{equation}

If $t>0$ and $\rho(t,x)\leq 3r$, then
$$
\dashnorm f\|_{\L_{q \beta,p \beta}(C_{r}(t,x))}\leq N
\dashnorm f\|_{\L_{q\beta ,p\beta }(C_{4r} )}\leq Nr^{-1}.
$$
In case $t>0$ and $\rho(t,x)> 3r$
use \eqref{8.27.1} which shows that
\begin{equation}
                                \label{8.27.4}
\dashnorm f\|_{\L_{q \beta,p \beta}(C_{r}(t,x))}=r^{-1}\dashnorm f\|_{\L_{q \beta,p \beta}(C_{1}(\tilde t,\tilde x))}
=Nr^{-1}\|f\|_{\L_{q \beta,p \beta}(C_{1}(\tilde t,\tilde x))},
\end{equation}
where $\tilde t=t/r^{2},\tilde x=x/r$, so that
$\rho(\tilde t,\tilde x)>3$. This case we split
into two subcases: A) $\tilde t>1$, B) $|\tilde
x|>2$.

In case A)
$$
\int_{\tilde t}^{\tilde t+1}|f(t,x)|^{q_{0}}\,dt=|x|^{2-q_{0}}
\int_{\tilde t/|x|^{2}}^{(\tilde t+1)/|x|^{2}}|\sqrt{s}-1|^{-q/p}\,ds 
$$
$$
=I_{\tilde t>|x|^{2}}...+I_{\tilde t\leq|x|^{2}}...
=J_{1}+J_{2}\leq J_{1}+ J(2) .
$$
Futhermore, $J_{1}(\tilde t)$ has negative derivative for $\tilde t>|x|^{2}$. Therefore,
$$
J_{1}\leq |x|^{2-q_{0}}\int_{1}^{2/|x|^{2}}
|\sqrt{s}-1|^{-q/p}\,ds .
$$
It folows as above that
$$
\Big(\int_{\tilde t}^{\tilde t+1}|f(t,x)|^{q_{0}}\,dt\Big)^{1/q_{0}}
\leq N|x|^{1/p_{0}-1}+N,
$$
where the integrals over unit balls
of the right-hand side to the power $p_{0}$
are dominated by a finite constant.
This provides the desired estimate
of the left-hand side of \eqref{8.27.4}
in case A).

In case B) one can use the same arguments
augmented by the fact that on $B_{1}(\tilde x)$
we have $|x|\geq 1$.
 
 Finally, notice that for $t\in(0,1)$
$$
\int_{|x|\leq 1}f^{p_{0}}(t,x)\,dx
=N\int_{0}^{1}
\rho^{d-1-p_{0}}\Big| \frac{\sqrt t}{\rho}-1 \Big|^{-1 }\,d\rho=\infty.
$$
This proves the second relation in \eqref{5.14,4}. 
 
\end{remark}

 \scriptsize
 
\printindex

\end{document}